\documentclass[a4paper,review]{cas-sc}

\usepackage[section]{placeins} 
\usepackage{svg} 
\usepackage{amssymb}
\usepackage{amsthm}
\usepackage{amsmath}
\usepackage{upgreek}
\usepackage{pdflscape}
\usepackage{listings}
\usepackage{multirow}
\usepackage[percent]{overpic}
\usepackage{multicol}
\usepackage{color}
\usepackage{stmaryrd}
\usepackage{xfrac}
\usepackage[capitalise]{cleveref}
\usepackage{booktabs}
\usepackage[section]{placeins} 
\usepackage[section]{algorithm}
\usepackage{calc}
\usepackage[titletoc]{appendix}
\usepackage{siunitx}
\usepackage{changepage}
\usepackage[sort&compress,square,numbers]{natbib}

\usepackage{graphicx}
\usepackage{graphics}
\usepackage{wrapfig}
\usepackage{float}
\usepackage{subcaption}
\usepackage{pdfpages}
\usepackage{alphalph}

\usepackage[percent]{overpic}
\usepackage{varwidth}
\usepackage{tikz}
\usepackage{pgfplots}
\usepackage{adjustbox}
\usetikzlibrary{arrows,matrix,positioning,fit}
\usetikzlibrary{shapes,positioning}
\usetikzlibrary{backgrounds}
\usepackage{tikz-layers}
\usepgfplotslibrary{fillbetween}
\usetikzlibrary{intersections}
\pgfplotsset{compat=1.14}
\usepgfplotslibrary{colorbrewer}
\usepgfplotslibrary{patchplots}
\usepgfplotslibrary[colorbrewer]
\usetikzlibrary{pgfplots.colorbrewer}
\usetikzlibrary[pgfplots.colorbrewer]
\usepgfplotslibrary{units}
\usetikzlibrary{spy}
\usepackage{pgfplotstable}
\usepackage{arrayjobx}
\graphicspath{ {./figs/} }
\usetikzlibrary{external}
\usetikzlibrary{shapes.misc}

\tikzset{cross/.style={cross out, draw=red, minimum size=2*(#1-\pgflinewidth), inner sep=0pt, outer sep=0pt},
cross/.default={1pt}}

\pgfplotsset{
discard if not/.style 2 args={
    x filter/.append code={
        \edef\tempa{\thisrow{#1}}
        \edef\tempaa{#2}
        \ifx\tempa\tempaa
        \else
            
        \fi
    }
}}

\newlength\myheight
\newlength\mydepth
\settototalheight\myheight{Xygp}
\newcommand*\inlinegraphics[1]{%
  \settototalheight\myheight{Xygp}%
  \settodepth\mydepth{Xygp}%
  \raisebox{-\mydepth}{\includegraphics[height=\myheight]{#1}}%
}
\newcommand\orcid[1]{\href{https://orcid.org/#1}{\inlinegraphics{orcid_16x16.png}}}

\makeatletter
\def\BState{\State\hskip-\ALG@thistlm}
\makeatother

\newdefinition{definition}{Definition}[section]
\newtheorem*{remark}{Remark}

\newcommand\pxvar[2]{\partial_{#2} #1}

\usepackage{algorithmic}

\usepackage{mathtools}
\newcommand{\defeq}{\vcentcolon=}

\newcommand{\abf}{\mathbf{a}}

\newcommand{\ebf}{\mathbf{e}}

\newcommand{\obf}{\mathbf{o}}

\newcommand{\rbf}{\mathbf{r}}

\usepackage{bm}

\newcommand{\pibf}{\bm{\pi}}

\newcommand{\ebar}{\bar{e}}

\newcommand{\Ebb}{\mathbb{E}}

\newcommand{\Rbb}{\mathbb{R}}

\newcommand{\Acal}{\mathcal{A}}
\newcommand{\Bcal}{\mathcal{B}}

\newcommand{\Ocal}{\mathcal{O}}

\newcommand{\Scal}{{\mathcal{S}}}

\newcommand\RCOMMENT[1]{\hfill\(\triangleright\) #1}

\DeclareMathOperator*{\argmin}{\arg\!\min}

\begin{document}

\title[mode=title]{DynAMO: Multi-agent reinforcement learning for dynamic anticipatory mesh optimization with applications to hyperbolic conservation laws}
\shorttitle{DynAMO: Multi-agent reinforcement learning for dynamic anticipatory mesh optimization}
\shortauthors{T.~Dzanic \textit{et al.}}

\author[1]{T.~Dzanic}[orcid=0000-0003-3791-1134]
\cormark[1]
\cortext[cor1]{Corresponding author}
\ead{dzanic1@llnl.gov}
\author[1]{K.~Mittal}[orcid=0000-0002-2062-852X]
\author[2]{D.~Kim}[orcid=0000-0001-9892-0611]
\author[1]{J.~Yang}[orcid=0000-0003-2234-6224]
\author[1]{S.~Petrides}[orcid=0000-0002-1284-5495]
\author[2]{B.~Keith}[orcid=0000-0002-6969-6857]
\author[1]{R.~Anderson}[orcid=0000-0002-3508-9944]

\address[1]{Lawrence Livermore National Laboratory, Livermore, CA 94550, United States of America}
\address[2]{Division of Applied Mathematics, Brown University, Providence, RI 02912, United States of America}

\begin{abstract}
We introduce DynAMO, a reinforcement learning paradigm for Dynamic Anticipatory Mesh Optimization. Adaptive mesh refinement is an effective tool for optimizing computational cost and solution accuracy in numerical methods for partial differential equations. However, traditional adaptive mesh refinement approaches for time-dependent problems typically rely only on instantaneous error indicators to guide adaptivity. As a result, standard strategies often require frequent remeshing to maintain accuracy. In the DynAMO approach, multi-agent reinforcement learning is used to discover new local refinement policies that can anticipate and respond to future solution states by producing meshes that deliver more accurate solutions for longer time intervals. By applying DynAMO to discontinuous Galerkin methods for the linear advection and compressible Euler equations in two dimensions, we demonstrate that this new mesh refinement paradigm can outperform conventional threshold-based strategies while also generalizing to different mesh sizes, remeshing and simulation times, and initial conditions.
\end{abstract}

\begin{keywords}
Adaptive mesh refinement \sep 
Finite element methods \sep
Scientific machine learning \sep 
Reinforcement learning \sep 
Hyperbolic conservation laws \sep 
\end{keywords}


\maketitle

\section{Introduction}
\label{sec:intro}
The numerical approximation of partial differential equations (PDEs) has broadly relied on discretization techniques such as finite volume, finite difference, and finite element methods. These techniques present a general framework for targeting a wide variety of governing equations arising in engineering and scientific applications. While their formulations and numerical properties vary, the methods above share the common task of effectively discretizing the problem domain onto a mesh upon which the governing equations can be solved numerically. In many applications, mesh discretization has a considerable impact on accuracy and efficiency. In particular, for systems with a large spatio-temporal variation of scales --- e.g., systems with regions of steep gradients and small-scale features that require high resolution --- the use of a uniform mesh discretization can result in an exceedingly inefficient numerical approximation. These multi-scale systems are common in scientific and engineering applications, with a wide variety of examples ranging from within fluid dynamics and astrophysics to solid mechanics and bioengineering \cite{TAKAKI2005263,Berger2017,kumar2020,kitzmann2016discontinuous,van2003finite}. An efficient treatment of these problems often requires some form of adaptive mesh refinement (AMR).

AMR strategies aim to balance computational cost and solution accuracy by selectively increasing the resolution of the mesh in regions where the largest discretization errors are generated while keeping the resolution low in regions where the errors are sufficiently small. Ideally, the end result achieves maximal error reduction with minimal computational cost. AMR is widely used to solve PDEs, with the typical approach involving some form of an instantaneous \emph{error indicator} to dictate which mesh elements are to be \emph{refined} and/or coarsened (i.e., \emph{de-refined}). However, a long-standing limitation of AMR strategies for time-dependent PDEs lies in the general observation that optimizing a mesh for computational errors at a fixed point in time will deliver a mesh that is suboptimal at future times.
Indeed, adapted meshes can quickly become inadequate as a system evolves and introduce large errors to the discrete solution that persist indefinitely. To mitigate this effect, it is often necessary to compute with a dynamic mesh that is re-adapted frequently. Unfortunately, the computational cost associated with mesh adaptation is not negligible and can quickly negate the practical benefits of AMR unless the refinement frequency is low and the solution does not vary much over the time period between re-adaptation steps.

The core belief motivating this work is that a more effective dynamic AMR paradigm can be achieved with \emph{anticipatory} refinement. Namely, we seek a refinement strategy to predict future error propagation and respond by preemptively refining the mesh. Compared to standard AMR approaches based only on instantaneous indicators, an anticipatory refinement strategy has the potential to achieve far superior accuracy and performance owing to its ability to refine the mesh \emph{before} future discretization errors appear. Most importantly, this would allow for longer time intervals between mesh adaptation, which is particularly beneficial for compute architectures where bandwidth is the bottleneck. 

Developing an anticipatory refinement strategy for time-dependent PDEs presents significant challenges. Except for trivial governing equations for which analytic solutions exist, the future spatio-temporal behavior of a PDE solution is not known \textit{a priori}, and directly predicting the evolution of the error is computationally infeasible. Instead of trying to directly predict the behavior of the error by techniques such as solving an adjoint problem~\citep{Becker2001}, we find that using surrogate models to approximate the error evolution offers a more tractable approach. One such approach is through error transport equations of \citet{Tyson2016}, which present an evolution equation for the error functional~\citep{Tyson2019, Tyson2019b, Wang2023}, or through even simple retrospective error estimates based on re-propagating the solution and measuring the error evolution. Further within the realm of surrogate modeling, machine learning (ML) techniques have shown promise as an effective tool for error estimation~\citep{Moosavi2017, Drohmann2015}. This success has driven efforts towards combing ML and AMR techniques, with the earliest works dating as far back as the early 1990s~\citep{dyck1992determining}. To date, ML has been used to design marking strategies~\cite{paszynski2021deep, Sualec2023, gillette2022learning, bohn2021recurrent}, goal-oriented AMR techniques~\citep{chen2020output, chen2021output, chakraborty2021multigoaloriented, roth2022neural}, and mesh density function-based mesh adaptation~\citep{dyck1992determining, chedid1996automatic, pfaff2020learning,huang2021machine, zhang2020meshingnet, song2022m2n, chan2022locally}.
However, an effective anticipatory refinement strategy for time-dependent PDEs requires optimizing long-term objectives, i.e., solution accuracy/cost over time. This in itself is a sequential decision-making problem as this objective at any given time is dependent on all of the refinement decisions leading up to that point. For this class of problems, reinforcement learning (RL) methods have shown promise as they are inherently tailored towards optimizing long-term objectives. The use of RL for AMR was first proposed in \citet{yang2021reinforcement}, where RL methods with variable-size global state and action spaces were applied to account for the changing mesh topology. 
This approach was further investigated in \citet{yang2023reinforcement}, where the authors used a multi-agent graph neural network with a team reward to generalize to arbitrary unstructured meshes, in \citet{freymuth2023swarm}, where local individual rewards were used to produce higher refinement levels, and in \citet{foucart2023deep}, where AMR was formulated as a partially observable Markov decision process in which a single agent applied to each element can observe only local surrounding information. For all of these works, AMR via RL was applied to scalar linear time-dependent or stationary PDEs. 

Given the potential of RL for optimizing long-term objectives in complex environments, the goal of this work is to introduce DynAMO, a new reinforcement learning paradigm for \textit{Dynamic Anticipatory Mesh Optimization}, in order to establish and explore the feasibility of using RL to guide anticipatory mesh refinement strategies for complex nonlinear systems of PDEs. In particular, we look at learning refinement strategies for discontinuous Galerkin finite element methods (FEMs), including $h$-refinement and $p$-refinement, for general hyperbolic conservation laws, with applications to the linear advection equation and the compressible Euler equations. Due to the local domain of influence of these governing equations, we consider a decentralized partially observable Markov decision process model of dynamic mesh optimization through a multi-agent reinforcement learning viewpoint, with independent agents corresponding to elements in the mesh that observe local surrounding information. Furthermore, we introduce novel and highly generalizable observations and reward functions that  1) enable anticipatory mesh refinement for arbitrary nonlinear hyperbolic conservation laws without analytic solutions; 2) show invariance to problem scale and mesh resolution; and 3) allow for user-controlled error/cost targets at evaluation time. We demonstrate the efficacy of DynAMO in multi-dimensional numerical experiments ranging from simple linear transport to complex nonlinear shock interactions, thus establishing that anticipatory refinement can unlock levels of efficiency that have remained out of reach with previous, conventional AMR techniques.

The remainder of this manuscript is organized as follows. In \cref{sec:preliminaries}, some background material on the overarching topic of this work is presented, including an overview of hyperbolic conservation laws, finite element methods, adaptive mesh refinement, and reinforcement learning. An in-depth description of the proposed DynAMO approach is given in \cref{sec:methodology}, followed by numerical implementation details in \cref{sec:implementation}. The results of the numerical experiments for $h$- and $p$-refinement on the advection and Euler equations are then shown in \cref{sec:results}, after which we use \cref{sec:conclusion} to draw our final conclusions.

\section{Preliminaries}\label{sec:preliminaries}

In this section, we present overviews of hyperbolic conservation laws, finite element methods, adaptive
mesh refinement, and reinforcement learning that are necessary for the subsequent sections.

\subsection{Hyperbolic conservation laws}
The primary applications of this work pertain to approximations of hyperbolic conservation laws of the form 
\begin{equation}\label{eq:hyp}
     \begin{cases}
        \pxvar{\mathbf{u}}{t} + \boldsymbol{\nabla}{\cdot}\mathbf{F}(\mathbf{u}) = \bm{0}, \quad \mathrm{for}\ (\mathbf{x}, t)\in\Omega\times\mathbb{R}_+, \\
        \mathbf{u}(\mathbf{x}, 0)= \mathbf{u}_0(\mathbf{x}), \quad \quad  \mathrm{for}\ \mathbf{x}\in \Omega,
    \end{cases}
\end{equation}
where $\Omega\subset\mathbb{R}^d$ is a $d$-dimensional spatial domain, $\mathbf{u}\in\mathbb{R}^m$ is a solution of $m$ variables, and $\mathbf{F}(\mathbf{u})\in \mathbb{R}^{m\times d}$ is the associated flux. We consider both linear scalar conservation laws as well as nonlinear systems of equations. For the former, the quintessential example is the scalar advection equation, given as
\begin{equation}
    \pxvar{u}{t} + \boldsymbol{\nabla}{\cdot}\left (\mathbf{c} u \right) = 0,
\end{equation}
where $\mathbf{c} \in \mathbb{R}^d$ is some constant velocity. For the latter, the compressible Euler equations of gas dynamics can be employed, given in the form of \cref{eq:hyp} as
\begin{equation}\label{eq:euler}
    \mathbf{u} = \begin{bmatrix}
            \rho \\ \mathbf{m} \\ E
        \end{bmatrix} \quad  \mathrm{and} \quad \mathbf{F}(\mathbf{u}) = \begin{bmatrix}
            \mathbf{m}^T\\
            \mathbf{m}\otimes\mathbf{v} + P\mathbf{I}\\
        (E+P)\mathbf{v}^T
    \end{bmatrix},
\end{equation}
where $\rho$ is the density, $\mathbf{m} \in \mathbb{R}^d$ is the momentum vector, and $E$ is the total energy. Furthermore, $\mathbf{I}$ denotes the identity matrix in $\mathbb{R}^{d\times d}$ and $\mathbf{v} = \mathbf{m}/\rho$ denotes the velocity. The pressure $P$ can be computed as
\begin{equation}
    P = (\gamma-1)\left(E - \frac{1}{2}\rho\mathbf{v}{\cdot}\mathbf{v}\right),
\end{equation}
where the specific heat ratio is taken as $\gamma = 1.4$ in this work.

The conservation law in \cref{eq:hyp} can also be represented in a quasi-linear formulation, given as
\begin{equation}
    \pxvar{\mathbf{u}}{t} + \mathbf{A}{\cdot}\nabla\mathbf{u} = \bm{0},
\end{equation}
where $\mathbf{A} \in \mathbb{R}^{m\times m \times d}$ is the flux Jacobian tensor defined
\begin{equation}
    \mathbf{A} = \frac{\partial \mathbf{F} (\mathbf{u})}{\partial \mathbf{u}},
\end{equation}
and the product $\mathbf{A}\cdot\nabla\mathbf{u}\in\mathbb{R}^m$ is defined $(\mathbf{A}\cdot\nabla\mathbf{u})_i=\sum_{j,k}A_{ijk}(\partial u_j/\partial x_k)$. For the advection equation, this yields an identical formulation as the conservative form due to the linearity of the flux, with the flux Jacobian simply reducing to the velocity, i.e., $\mathbf{A} = \mathbf{c}$. For the Euler equations, the flux Jacobian becomes much more complex as the flux is nonlinear with respect to a vector-valued solution, and its derivation is left to \cref{app:fluxjacobian}. However, this nonlinearity can cause discrepancies between the conservative formulation and the quasi-linear form, particularly so in the vicinity of discontinuities. Nevertheless, even in the nonlinear case, the structure of the flux Jacobian still possesses some useful information about the underlying system as it gives an approximation of the characteristic propagation velocity of the solution. 

An important property of hyperbolic conservation laws is that this characteristic propagation velocity is finite, such that the domain of influence of any point is locally supported in space and time. An upper bound for the extent of the domain of influence $\mathcal D$ for some arbitrary point $\mathbf{x}_0$ over some time period $[t, t + \Delta \tau]$ can readily be estimated as
\begin{equation}
    \mathcal D (\mathbf{x}_0, \Delta \tau) = \mathbb B^d (\mathbf{x}_0, \Delta \tau \lambda_{\max}),
\end{equation}
where $\mathbb B^d (\mathbf x', r) $ is a $d$-ball centered at $\mathbf x'$ with radius $r$ and $\lambda_{\max}$ is an estimate of the maximum wavespeed of the system. This domain of influence effectively ensures that the solution at $\mathbf{x}_0$ cannot be affected by the solution at any point farther than $r = \Delta \tau \lambda_{\max}$ away, which can greatly reduce the level of information necessary to adequately predict the local evolution of the system over a sufficiently small time interval. 

\subsection{Finite element methods}
To numerically approximate hyperbolic conservation laws of the form of \cref{eq:hyp}, we utilize finite element methods to discretize the governing equations, particularly the nodal discontinuous Galerkin (DG) method \citep{Hesthaven2008}. In this approach, the domain $\Omega$ is partitioned into $N$ elements $\Omega_k$ such that $\Omega = \bigcup_N\Omega_k$ and $\Omega_i\cap\Omega_j=\emptyset$ for $i\neq j$. Within each element $\Omega_k$, the discrete solution $\mathbf{u}_h(\mathbf{x})$ is approximated via a set of $n_s$ nodal interpolating polynomials as 
\begin{equation}
    \mathbf{u}_h(\mathbf{x}) = \sum_{i=1}^{n_s} \mathbf{u}_{i}\phi_i(\mathbf{x})\subset V_h,
\end{equation}
where $\mathbf{x}_i$ is a set of solution nodes, $\phi_i(\mathbf{x})$ are their associated nodal basis functions which possess the property $\phi_i(\mathbf{x}_j) = \delta_{ij}$, and $V_h$ is the piece-wise polynomial space spanned by the nodal basis functions. We use the short-hand notation $\mathbf{u}_i = \mathbf{u}_h(\mathbf{x}_i)$ for brevity and denote the order of the approximation as $\mathbb P_p$, where $p$ is the maximal order of $\mathbf{u}_h(\mathbf{x})$.

The problem is solved through the weak formulation by applying the DG approximation and integrating \cref{eq:hyp} with respect to a test function $\mathbf{w}_h(\mathbf{x})$, which resides in the same finite element space as $\mathbf{u}_h(\mathbf{x})$, yielding
\begin{equation}\label{eq:semi-disc}
    \sum_{k=1}^{N} \left \{ \int_{\Omega_k} \partial_t \mathbf{u}_h{\cdot} \mathbf{w}_h\ \mathrm{d}V  + \int_{\partial \Omega_k} \hat{\mathbf{F}}(\mathbf{u}_h^-, \mathbf{u}_h^+, \mathbf{n}) {\cdot} \mathbf{w}_h\ \mathrm{d}S  -  \int_{\Omega_k} \mathbf{F}(\mathbf{u}_h){\cdot} \nabla \mathbf{w}_h\ \mathrm{d}V\right \} = 0.
\end{equation}
As the approximate solution is discontinuous across element boundaries, the flux function must be replaced by a numerical flux $\hat{\mathbf{F}}(\mathbf{u}_h^-, \mathbf{u}_h^+, \mathbf{n})$ which is dependent on the solution within the element of interest, denoted by the $\mathbf{u}_h^-$, and its face-adjacent neighbor, denoted by $\mathbf{u}_h^+$, along with the outward facing normal vector $\mathbf{n}$. For the hyperbolic conservation laws considered in this work, the Rusanov approximate Riemann solver \citep{Rusanov1962} is used, such that the numerical flux can be calculated as
\begin{equation}
    \hat{\mathbf{F}}(\mathbf{u}_h^-, \mathbf{u}_h^+, \mathbf{n}) = \frac{1}{2}\left (\mathbf{F}(\mathbf{u}_h^-){\cdot}\mathbf{n} + \mathbf{F}(\mathbf{u}_h^+){\cdot}\mathbf{n} \right) - \frac{1}{2}\lambda_{\max} \left (\mathbf{u}_h^+ - \mathbf{u}_h^-\right).
\end{equation} 
For the advection equation, the maximum wavespeed is simply taken as $\lambda_{\max} = |\mathbf{c}{\cdot}\mathbf{n}|$, which reduces to the upwind flux, whereas for the Euler equations, the Davis wavespeed estimate \citep{Davis1988} is used.

\subsection{Adaptive mesh refinement}\label{ssec:amr}

While AMR techniques for various classes of governing equations and numerical methods differ in their specific implementations, they typically fall within the widely-used paradigm of the \texttt{SOLVE}\textrightarrow\texttt{ESTIMATE}\textrightarrow\texttt{MARK}\textrightarrow\texttt{REFINE} loop. A complete pass through of the above sequence for time-dependent PDEs consists of solving the governing equations over some time interval, estimating some error indicator functional (or functionals) for each element in the mesh, marking which elements are to be refined and de-refined, and applying the resulting refinement actions onto the mesh. A more in-depth overview of the steps in the AMR sequence is presented below.

\texttt{SOLVE}. In this step, the semi-discrete form of the governing time-dependent PDE is evolved using an appropriate temporal integration method, the specific choice of which is typically dependent on the characteristics of the physical system being solved. For each \texttt{SOLVE} step in the loop, the simulation is advanced over the time interval $[t, t + T]$, 
where $T$ is referred to as the remesh time as this is the time period over which the mesh is fixed. We remark here that $T$ is not the temporal integration time step $\Delta t$, the latter of which is dependent on the numerical solver and is largely irrelevant in the AMR procedure. In most applications, the remesh time interval is many times larger than the simulation time step as frequent remeshing (i.e., low $T$) incurs a large computational cost overhead. 

\texttt{ESTIMATE}. Once the simulation has been advanced to the next remesh time, the \texttt{ESTIMATE} step consists of computing a local error indicator for each element in the mesh. The indicator functional maps the solution to a set of element-wise indicator values $e_i$, i.e.,
\begin{equation}
    \mathbf{u}_h(\mathbf{x}) \mapsto e_i ~~ \text{for all~} i \in \{1,\dotsc, N\}.
\end{equation}
The indicator is assumed to be deterministic, such that an identical input state would result in an identical indicator distribution throughout the domain. Furthermore, the \texttt{ESTIMATE} step may result in multiple indicators for each element.

\texttt{MARK}. Given the distribution of the error indicator (or indicators) for the elements in the mesh, the \texttt{MARK} step maps this distribution to a set of refinement decisions for the elements through a policy $\pi (e_i) \to \{0,1\}\ \forall \ i \in \{1,\dotsc, N\}$, where the resulting decisions $0$ and $1$ correspond to de-refinement and refinement actions, respectively. We remark here that the output of this policy may not necessarily be binary as multiple refinement methods may be allowed, a notion that will be further explained in the subsequent step of the refinement loop, and that an element's refinement decision may depend on other elements' decisions. This step of the refinement process is largely reliant on heuristics, with a wide variety of approaches, including absolute threshold-based methods, i.e., 
\begin{equation}
    \pi(e_i) = \begin{cases}
        1, \quad \mathrm{if}\ e_i > \theta,\\
        0, \quad \mathrm{else},
    \end{cases}
\end{equation}
where $\theta$ here is some arbitrary global error threshold, and relative threshold-based methods, i.e.,
\begin{equation}
    \pi(e_i) = \begin{cases}
        1, \quad \mathrm{if}\ \frac{e_{\max} - e_{i}}{e_{\max} - e_{\min}} > \theta,\\
        0, \quad \mathrm{else},
    \end{cases}
\end{equation}
where $\theta \in [0, 1]$ here is some arbitrary relative error threshold and $e_{\max}$ and $e_{\min}$ are the global maximum and minimum error indicator values, respectively. 

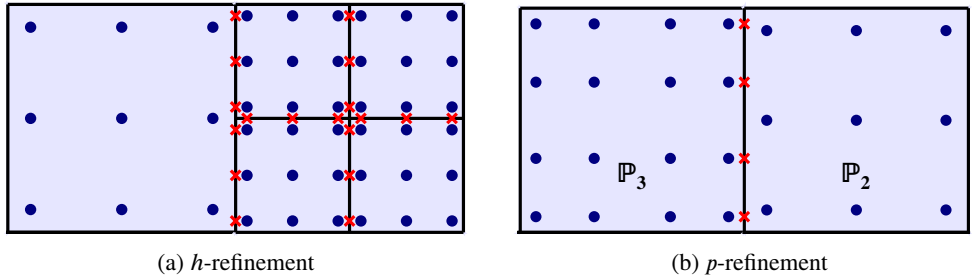
\begin{figure}[htbp!]
    \centering
    \subfloat[$h$-refinement]{
    \adjustbox{width=0.4\linewidth,valign=b}{\begin{tikzpicture}[scale=1]
\coordinate (0) at (0,0);
\coordinate (1) at (0.5,0);
\coordinate (2) at (1,0);
\coordinate (3) at (0,0.5);
\coordinate (4) at (0.5,0.5);
\coordinate (5) at (1,0.5);
\coordinate (6) at (0,1);
\coordinate (7) at (0.5,1);
\coordinate (8) at (1,1);
\coordinate (100) at (0.1,0.1);
\coordinate (101) at (0.5,0.1);
\coordinate (102) at (0.9,0.1);
\coordinate (103) at (0.1,0.5);
\coordinate (104) at (0.5,0.5);
\coordinate (105) at (0.9,0.5);
\coordinate (106) at (0.1,0.9);
\coordinate (107) at (0.5,0.9);
\coordinate (108) at (0.9, 0.9);

\draw [blue!10!white, fill=blue!10!white] plot [smooth] coordinates {(0) (1) (2)} -- (4) -- cycle;
\draw [blue!10!white, fill=blue!10!white] plot [smooth] coordinates {(2) (5) (8)} -- (4) -- cycle;
\draw [blue!10!white, fill=blue!10!white] plot [smooth] coordinates {(8) (7) (6)} -- (4) -- cycle;
\draw [blue!10!white, fill=blue!10!white] plot [smooth] coordinates {(6) (3) (0)} -- (4) -- cycle;

\draw[black] plot [smooth cycle] coordinates {(0) (1) (2)} ;
\draw[black] plot [smooth] coordinates {(8) (7) (6)};
\draw[black] plot [smooth] coordinates {(6) (3) (0)};

\foreach \i in {100,...,108}
{
    \draw (\i) node[circle, fill=blue!50!black, inner sep=0.5pt] {} ;
}

\coordinate (9) at (1,0);
\coordinate (10) at (1.25,0);
\coordinate (11) at (1.5,0);
\coordinate (12) at (1.75,0);
\coordinate (13) at (2.0,0);

\coordinate (14) at (1,0.25);
\coordinate (15) at (1.25,0.25);
\coordinate (16) at (1.5,0.25);
\coordinate (17) at (1.75,0.25);
\coordinate (18) at (2.0,0.25);

\coordinate (19) at (1,0.5);
\coordinate (20) at (1.25,0.5);
\coordinate (21) at (1.5,0.5);
\coordinate (22) at (1.75,0.5);
\coordinate (23) at (2.0,0.5);

\coordinate (24) at (1,0.75);
\coordinate (25) at (1.25,0.75);
\coordinate (26) at (1.5,0.75);
\coordinate (27) at (1.75,0.75);
\coordinate (28) at (2.0,0.75);

\coordinate (29) at (1,1);
\coordinate (30) at (1.25,1);
\coordinate (31) at (1.5,1);
\coordinate (32) at (1.75,1);
\coordinate (33) at (2.0,1);

\coordinate (209) at (1.05,0.05);
\coordinate (210) at (1.25,0.05);
\coordinate (211) at (1.45,0.05);
\coordinate (212) at (1.05,0.25);
\coordinate (213) at (1.25,0.25);
\coordinate (214) at (1.45,0.25);
\coordinate (215) at (1.05,0.45);
\coordinate (216) at (1.25,0.45);
\coordinate (217) at (1.45,0.45);

\coordinate (309) at (1.55,0.05);
\coordinate (310) at (1.75,0.05);
\coordinate (311) at (1.95,0.05);
\coordinate (312) at (1.55,0.25);
\coordinate (313) at (1.75,0.25);
\coordinate (314) at (1.95,0.25);
\coordinate (315) at (1.55,0.45);
\coordinate (316) at (1.75,0.45);
\coordinate (317) at (1.95,0.45);

\coordinate (409) at (1.05,0.55);
\coordinate (410) at (1.25,0.55);
\coordinate (411) at (1.45,0.55);
\coordinate (412) at (1.05,0.75);
\coordinate (413) at (1.25,0.75);
\coordinate (414) at (1.45,0.75);
\coordinate (415) at (1.05,0.95);
\coordinate (416) at (1.25,0.95);
\coordinate (417) at (1.45,0.95);

\coordinate (509) at (1.55,0.55);
\coordinate (510) at (1.75,0.55);
\coordinate (511) at (1.95,0.55);
\coordinate (512) at (1.55,0.75);
\coordinate (513) at (1.75,0.75);
\coordinate (514) at (1.95,0.75);
\coordinate (515) at (1.55,0.95);
\coordinate (516) at (1.75,0.95);
\coordinate (517) at (1.95,0.95);

\draw [blue!10!white, fill=blue!10!white] plot [smooth] coordinates {(9) (13)} -- (21) -- cycle;
\draw [blue!10!white, fill=blue!10!white] plot [smooth] coordinates {(13) (33)} -- (21) -- cycle;
\draw [blue!10!white, fill=blue!10!white] plot [smooth] coordinates {(29) (33)} -- (21) -- cycle;
\draw [blue!10!white, fill=blue!10!white] plot [smooth] coordinates {(9) (29)} -- (21) -- cycle;

\draw[black] plot [smooth] coordinates {(9) (13)} ;
\draw[black] plot [smooth] coordinates {(13) (33)};
\draw[black] plot [smooth] coordinates {(33) (29)};

\draw[black] plot [smooth] coordinates {(11) (21)} ;
\draw[black] plot [smooth] coordinates {(19) (21)} ;
\draw[black] plot [smooth] coordinates {(23) (21)} ;
\draw[black] plot [smooth] coordinates {(31) (21)} ;

\draw[black] plot [smooth] coordinates {(29) (9)};

\foreach \i in {209,...,217}
{
    \draw (\i) node[circle, fill=blue!50!black, inner sep=0.5pt] {} ;
}
\foreach \i in {309,...,317}
{
    \draw (\i) node[circle, fill=blue!50!black, inner sep=0.5pt] {} ;
}
\foreach \i in {409,...,417}
{
    \draw (\i) node[circle, fill=blue!50!black, inner sep=0.5pt] {} ;
}
\foreach \i in {509,...,517}
{
    \draw (\i) node[circle, fill=blue!50!black, inner sep=0.5pt] {} ;
}

\coordinate (600) at (1.0,0.05);
\coordinate (601) at (1.0,0.25);
\coordinate (602) at (1.0,0.45);
\coordinate (603) at (1.0,0.55);
\coordinate (604) at (1.0,0.75);
\coordinate (605) at (1.0,0.95);
\coordinate (606) at (1.5,0.05);
\coordinate (607) at (1.5,0.25);
\coordinate (608) at (1.5,0.45);
\coordinate (609) at (1.5,0.55);
\coordinate (610) at (1.5,0.75);
\coordinate (611) at (1.5,0.95);
\coordinate (612) at (1.05,0.5);
\coordinate (613) at (1.25,0.5);
\coordinate (614) at (1.45,0.5);
\coordinate (615) at (1.55,0.5);
\coordinate (616) at (1.75,0.5);
\coordinate (617) at (1.95,0.5);

\foreach \i in {600,...,617}
{
    \draw (\i) node[cross, fill=red, inner sep=0.5pt] {} ;
}

\end{tikzpicture}}}
    \subfloat[$p$-refinement]{
    \adjustbox{width=0.41\linewidth,valign=b}{\begin{tikzpicture}[scale=1]
\coordinate (0) at (0,0);
\coordinate (1) at (0.33333333333, 0);
\coordinate (2) at (0.66666666667,0);
\coordinate (3) at (1.0,0);

\coordinate (4) at (0.0,0.33333333333);
\coordinate (5) at (0.33, 0.33333333333);
\coordinate (6) at (0.67,0.33333333333);
\coordinate (7) at (1.0,0.33333333333);

\coordinate (8) at (0.0,0.66666666667);
\coordinate (9) at (0.33, 0.66666666667);
\coordinate (10) at (0.67,0.66666666667);
\coordinate (11) at (1.0,0.66666666667);

\coordinate (12) at (0,1);
\coordinate (13) at (0.33333333333,1);
\coordinate (14) at (0.66666666667,1);
\coordinate (15) at (1.0,1);

\coordinate (100) at (0.07,0.07);
\coordinate (101) at (0.33, 0.07);
\coordinate (102) at (0.67,0.07);
\coordinate (103) at (0.93,0.07);

\coordinate (104) at (0.07,0.33);
\coordinate (105) at (0.33, 0.33);
\coordinate (106) at (0.67,0.33);
\coordinate (107) at (0.93,0.33);

\coordinate (108) at (0.07,0.67);
\coordinate (109) at (0.33, 0.67);
\coordinate (110) at (0.67,0.67);
\coordinate (111) at (0.93,0.67);

\coordinate (112) at (0.07,0.93);
\coordinate (113) at (0.33, 0.93);
\coordinate (114) at (0.67,0.93);
\coordinate (115) at (0.93,0.93);

\coordinate (116) at (1.1,0.1);
\coordinate (117) at (1.5,0.1);
\coordinate (118) at (1.9,0.1);
\coordinate (119) at (1.1,0.5);
\coordinate (120) at (1.5,0.5);
\coordinate (121) at (1.9,0.5);
\coordinate (122) at (1.1,0.9);
\coordinate (123) at (1.5,0.9);
\coordinate (124) at (1.9, 0.9);

\draw [blue!10!white, fill=blue!10!white] plot [smooth] coordinates {(0) (1) (2) (3)} -- (6) -- cycle;
\draw [blue!10!white, fill=blue!10!white] plot [smooth] coordinates {(3) (7) (11) (15)} -- (6) -- cycle;
\draw [blue!10!white, fill=blue!10!white] plot [smooth] coordinates {(15) (14) (13) (12)} -- (6) -- cycle;
\draw [blue!10!white, fill=blue!10!white] plot [smooth] coordinates {(12) (8) (4) (0)} -- (6) -- cycle;
\draw [blue!10!white, fill=blue!10!white] plot [smooth] coordinates {(6) (10) (9)} -- (5) -- cycle;

\draw[black] plot [smooth cycle] coordinates {(0) (1) (2) (3)} ;
\draw[black] plot [smooth] coordinates {(15) (14) (13) (12)} ;
\draw[black] plot [smooth] coordinates {(12) (8) (4) (0)} ;

\coordinate (16) at (1.0,0.5);
\coordinate (17) at (1.5, 0);
\coordinate (18) at (2.0, 0);

\coordinate (19) at (1.5,0.5);
\coordinate (20) at (2, 0.5);

\coordinate (21) at (1.5, 1);
\coordinate (22) at (2, 1);

\draw [blue!10!white, fill=blue!10!white] plot [smooth] coordinates {(3) (17) (18)} -- (19) -- cycle;
\draw [blue!10!white, fill=blue!10!white] plot [smooth] coordinates {(18) (20) (22)} -- (19) -- cycle;
\draw [blue!10!white, fill=blue!10!white] plot [smooth] coordinates {(22) (21) (15)} -- (19) -- cycle;
\draw [blue!10!white, fill=blue!10!white] plot [smooth] coordinates {(15) (16) (3)} -- (19) -- cycle;

\draw[black] plot [smooth cycle] coordinates {(3) (17) (18)} ;
\draw[black] plot [smooth] coordinates {(18) (20) (22)} ;
\draw[black] plot [smooth] coordinates {(22) (21) (15)} ;

\draw[black] plot [smooth] coordinates {(3) (7) (11) (15)} ;

\foreach \i in {100,...,115}
{
    \draw (\i) node[circle, fill=blue!50!black, inner sep=0.5pt] {} ;
}

\coordinate (200) at (1,0.07);
\coordinate (201) at (1,0.33);
\coordinate (202) at (1,0.67);
\coordinate (203) at (1,0.93);

\draw (200) node[cross, fill=red, inner sep=0.5pt] {} ;
\draw (201) node[cross, fill=red, inner sep=0.5pt] {} ;
\draw (202) node[cross, fill=red, inner sep=0.5pt] {} ;
\draw (203) node[cross, fill=red, inner sep=0.5pt] {} ;

\foreach \i in {116,...,124}
{
    \draw (\i) node[circle, fill=blue!50!black, inner sep=0.5pt] {} ;
}

\draw (0.5, 0.25) node[scale=0.4, text=black] {$\boldsymbol{\mathbb P_3}$};
\draw (1.5, 0.25) node[scale=0.4, text=black] {$\boldsymbol{\mathbb P_2}$};

\end{tikzpicture}}}
    
    \caption{\label{fig:hp_ref_schematic} Schematic of one-level $h$-refinement (left) and $p$-refinement on a two-dimensional quadrilateral discontinuous Galerkin finite element mesh with equal number of quadrature nodes and basis functions. Blue circles represent volume quadrature nodes, red crosses represent surface quadrature nodes.}
\end{figure}

\texttt{REFINE}. With the output of the \texttt{MARK} step, the appropriate refinement and de-refinement actions can be then be applied to the mesh. There exist a variety of mesh refinement methods for AMR, with the most ubiquitous approach being $h$-refinement. For this method, show on the left-hand side of \cref{fig:hp_ref_schematic}, refining an element of choice subdivides it into a set of sub-elements, typically referred to as the \emph{child elements}, whereas de-refining a group of sub-elements coalesces them into a singular element, typically referred to as the \emph{parent element}. The number of successive refinement actions performed between a child element and its native parent element is referred to as the \emph{refinement depth}. The use of $h$-refinement typically results in a non-conforming mesh for which the numerical method of choice must appropriately handle (e.g., the use of constrained degrees of freedom in FEM). 

Due to the element-local polynomial approximation used in finite element methods, another, albeit not as common, AMR approach can be achieved with $p$-refinement, shown on the right-hand side of \cref{fig:hp_ref_schematic}. For this method, the order of the polynomial basis within an element is increased or decreased appropriately in accordance with the output of the marking policy, such that effective order of accuracy and number of degrees of freedom within elements can vary. Similarly to $h$-refinement, this adaptation technique results in a degree of non-conformity, although in the basis instead of in the mesh. 

Additional refinement methods exist separately from $h$-refinement and $p$-refinement, although they are not considered in this work. For FEM, the joint use of $h$- and $p$-refinement, dubbed $hp$-refinement, may be yield superior performance for certain applications as it is typically assumed that $h$-refinement is better suited for regions with discontinuous features whereas $p$-refinement is better suited for regions with smooth features. Furthermore, mesh deformation, dubbed $r$-refinement, may be used in applications where dynamically varying features and geometries are encountered. 

\subsection{Reinforcement learning}
Reinforcement learning enables an agent or a team of agents to solve sequential decision-making problems in an unknown environment by interacting with the environment and receiving rewards over many episodes \citep{sutton2018reinforcement}. 
In a single-agent scenario, the environment is formalized as a stationary Markov decision process (MDP) \citep{puterman2014markov}, a discrete-time stochastic process $\lbrace \Scal, \Acal, P, R, \gamma \rbrace$ with the following components: 
\begin{itemize}
    \item A state space $\Scal$.
    \item An action space $\Acal$.
    \item A transition function $P \colon \Scal \times \Acal \times \Scal \rightarrow [0,1]$ that defines the distribution $P(s_{\tau+1}|s_{\tau}, a_{\tau})$ of next states $s_{\tau+1} \in \Scal$ given an action $a_{\tau} \in \Acal$ taken by the agent at state $s_{\tau} \in \Scal$, and an initial state distribution $P_0 \colon \Scal \rightarrow [0,1]$, from which the initial state of each episode is sampled.
    \item A reward function $R \colon \Scal \times \Acal \rightarrow \Rbb$ which includes a temporal discount factor $\gamma \in (0,1)$ that determines the value of future rewards.
\end{itemize}
We use $\tau = 0,1,\dotsc$ to denote the discrete time index. An agent's goal is to find a policy $\pi \colon \Scal \times \Acal \rightarrow [0,1]$, which determines the probability distribution of actions $\pi(a_{\tau}|s_{\tau})$ to take at each state, to maximize the expected cumulative discounted reward $J(\pi) \defeq \Ebb [ \sum_{\tau=0}^{\infty} \gamma^{\tau} R(s_{\tau},a_{\tau}) ]$, where the expectation is taken over the environment transition and policy.
In practice, the transition function contains a termination criterion, formalized by an absorbing state with zero reward, that triggers a reset of the environment for a new episode, so that each episode has a finite time horizon.

An alternative RL approach is the multi-agent viewpoint, which is the method of choice in this work. The environment is formalized as a decentralized partially observable Markov decision process (Dec-POMDP) \citep{oliehoek2016concise}, denoted by $\lbrace \Scal, \Ocal, \Acal, P, R, O, \gamma \rbrace$, which extends an MDP in the following ways: 1) multiple agents, indexed by $i = 1,\dotsc,N$, interact with the shared environment; 2) each agent $i$ has a limited partial observation $o^i_{\tau} = O(s_{\tau}) \in \Ocal$, produced by an observation function $O$ from the true state; 3) actions are decentralized, meaning that each agent $i$ acts according to its own policy $\pi^i(a^i_{\tau}|o^i_{\tau})$; and 4) the environment transition $P(s_{\tau+1} | s_{\tau}, a^1_{\tau},\dotsc, a^N_{\tau})$ depends on the actions of all agents. For the purposes of this work, the Dec-POMDP definition is specialized in the following two ways. First, agents are homogeneous, meaning that they have the same observation and action spaces. Second, all agents share the same definition of reward function $R$, but each agent receives its own individual reward value $R(s,a^i,i)$ which depends on its own action and situation within the state. Each agent $i$ learns a policy $\pi^i(a^i|o^i)$ to maximize the  objective
\begin{align}\label{eq:objective}
    J(\pi^i) \defeq \Ebb_{s_0 \sim P_0, \abf_{\tau} \sim \pibf(\cdot | s_{\tau}), s_{\tau+1} \sim P(\cdot|s_{\tau},\abf_{\tau})} \left[ \sum_{\tau=0}^{\infty} \gamma^{\tau} R(s_{\tau},a^i_{\tau},i) \right] \, ,
\end{align}
where $\abf = (a^1,\dotsc,a^N)$ denotes the collection of all agents' actions and $\pibf(\abf|s) = \prod_{i=1}^N \pi(a^i|o^i)$ denotes the joint policy.
\section{Methodology}\label{sec:methodology}
Given the drawbacks of standard AMR techniques for the approximation of complex hyperbolic conservation laws, the objective of this work is to 
explore a reinforcement learning approach for anticipatory mesh refinement strategies. As reinforcement learning is quite general, a wide variety of approaches are feasible for integrating AMR within the RL framework. As such, to have a robust and reliable policy for AMR, certain desirable attributes and features would ideally be satisfied:
\begin{enumerate}
    \item The policy should be able to generalize to arbitrary hyperbolic conservation laws without requiring substantial domain specific knowledge, including complex nonlinear systems of equations which do not possess closed-form expressions for the error (i.e., analytic solutions). 
    \item The user should be able to control the relative degree of refinement (i.e., error and cost targets) of the policy \emph{at evaluation time.}
    \item The resulting policy should be insensitive to or, ideally, completely invariant to mesh/problem scale, choice of error estimator, mesh resolution, simulation time, and remeshing time interval. 
\end{enumerate}
With these characteristics in mind, we present the proposed DynAMO approach, a \textit{Dynamic Anticipatory Mesh Optimization} method for hyperbolic conservation laws using multi-agent reinforcement learning. The goal of DynAMO is to learn more optimal long-term refinement strategies which consider the underlying physics of the problem, allowing for more efficient refinement decisions which account for the spatio-temporal evolution of the solution, longer remeshing time intervals, and the preemptive refinement of mesh regions prior to the introduction of discretization error. Therefore, novel formulations for the observation and reward functions are introduced to best achieve these target goals and features. We remark here that the primary purpose of this work is to evaluate the ability of RL and the proposed observation and reward formulations as methods of achieving more efficient anticipatory mesh refinement policies for complex nonlinear systems of equations. As such, the focus of this work will be on discovering policies for one-level $h$- and $p$-refinement in finite element methods on structured (but non-conforming), two-dimensional, periodic meshes. While the extension to unstructured meshes and curved mesh elements, boundary conditions and complex geometries, and multi-level refinement poses a necessary task for the proposed approach, these difficulties are orthogonal to the study in this work and have been addressed by the formulations presented by the authors in \citet{yang2023reinforcement} and by \citet{foucart2023deep} and \citet{freymuth2023swarm} in the context of simpler linear problems. Furthermore, extensions to alternate refinement actions (such as conforming splitting in unstructured meshes) and mesh movement in $r$-refinement may be achieved through increasing the size of the action space or changing the action space to a continuous output, respectively.

\subsection{Dynamic mesh optimization as a Dec-POMDP}
DynAMO is instantiated as a Dec-POMDP model of dynamic mesh optimization
with the following components: 1) a set of agents corresponding to each element in the mesh; 2) an observation space for each agent consisting of a local spatial window centered on the agent's element that observes derived quantities of the solution in the surrounding elements at the given time; 3) an action space corresponding to refinement and de-refinement operations; 4) a transition function which applies these actions and advances the PDE until the next remesh time; and 5) a reward function which encodes the efficacy of the refinement/de-refinment actions on the local error. 

\subsubsection{Agent}
We consider a multi-agent reinforcement learning approach to AMR, where each element $\Omega_i$ is represented by an agent, denoted by an identifier $i \in \lbrace 1,\dotsc, N \rbrace$. Crucially, this viewpoint enables each individual element to make its own refinement or coarsening decision based on its own local observation and independently from other agents. Furthermore, all agents are allowed to act simultaneously at each discrete remeshing step in contrast to single-agent approaches \citep{yang2021reinforcement,foucart2023deep}, which is critical for maintaining exactly the same mode of execution between training and deployment. For $p$-refinement, $N$ is constant throughout the simulation, such that the number of agents remains constant. However, for $h$-refinement, the number of elements varies according to mesh refinement and de-refinement.
To avoid the problem of agent creation and deletion, posthumous credit assignment in multi-agent reinforcement learning, 
and agent co-operation (for coarsening actions in $h$-refinement), we instead take an alternate approach where the agents are assigned to the initial coarse elements, and any additional elements stemming from subsequent refinement actions on these initial coarse elements also remain under the scope of the original agent. This allows for a constant number of agents throughout the simulation even in the $h$-refinement case, avoiding the difficulties of agent creation and deletion.
Note that while the number of agents is the same during a single episode, it does not have to be the same in different episodes, such that different mesh resolution may be used for training and evaluation.

\subsubsection{Observation space}
The choice of the agent observation plays arguably the most important role in the design of an anticipatory mesh refinement approach. For time-dependent PDEs such as hyperbolic conservation laws, the observation must have two components: an effective and robust approach for approximating the error distribution within the domain and a general method for predicting the spatio-temporal evolution of the solution and, by extension, the error. These components, which will be described momentarily, form the quantities that are observed by each individual agent. 

\begin{figure}
    \centering
    \adjustbox{width=0.9\linewidth,valign=b}{\includegraphics[]{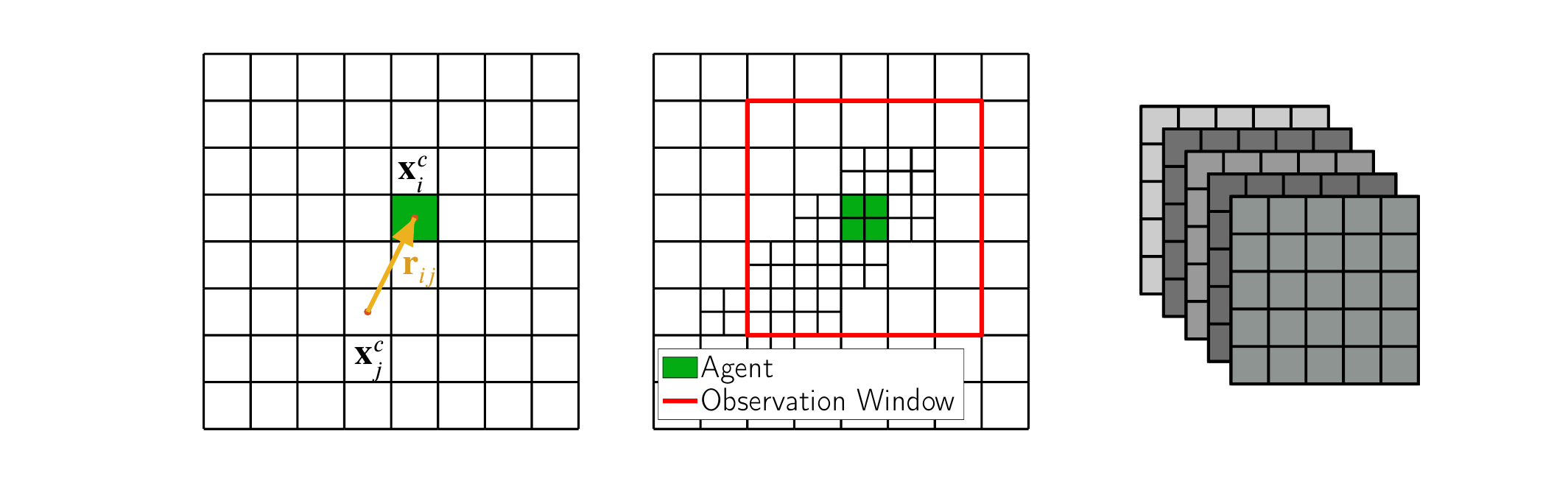}}
    \caption{\label{fig:obs_schematic} Schematic for an example $h$-refinement observation setup with an initial mesh of $N = 8^2$ agents/elements, showing the mesh and example displacement vector (left), the observation window of an agent on a refined mesh (center), and the observation channels of an agent with 5 observable quantities (right).
    }
\end{figure}

Given the local domain of influence of hyperbolic conservation laws, it is natural to assume that a local observation window is sufficient to adequately predict the evolution of the underlying physics. As such, we utilize a formulation where each agent $i$ has an individual observation $o^{ij} \in \Rbb^{k_x \times k_y \times n}$ which consists of a $k_x \times k_y$ spatial window with one channel for each of the $n$ observable quantities. The spatial window is centered on the element of agent $i$ and contains the values of $k_x = 2n_x + 1$ and $k_y = 2n_y + 1$ elements, where $n_x$ and $n_y$ specify the number of neighbor elements along each direction, such that $1 \leq j \leq k_x k_y$. An example schematic for this approach is shown in \cref{fig:obs_schematic}. For $h$-refinement, these numbers correspond to the initial coarse elements (i.e., agents). This spatial extent is fixed at training time, which constrains the extent of the local domain of influence. As will be later presented, this gives an effective bound maximum remesh time $T$ for the approach. 

The first observable quantity is an error indicator metric for each element in the observation window. For brevity, we present the observation for an arbitrary agent $i$ with an arbitrary neighbor index $j$, where $j$ corresponds to any one of the $k_x k_y$ agent/element indices which reside in the spatial window of the observation $o^i$. As most governing equations do not possess an analytic solution, it is unreasonable to expect that an analytic expression exists for the error in practical applications. As such, we utilize an \emph{error estimator} to approximate the distribution of error within the domain. While the details of the particular error estimators used in this work are left to \cref{sec:implementation}, we assume that the error estimator returns a set of quantities $e^i_\tau\ \forall \ i \in \lbrace 1,\dotsc, N \rbrace$, where $e^i_\tau$ represents an element-wise constant approximation of the error in element $\Omega_i$ at some time index $\tau$. However, observing the raw output of an error estimator poses some drawbacks as this quantity is not scale-invariant and, in a sense, not invariant to the mesh resolution. Instead, we observe a normalized error estimate as
\begin{align}\label{eq:obs}
    o^{ij}_{\tau, 0} &\defeq - \frac{\log_{10}(e^j_\tau)}{\log_{10}(e_{\tau,\max})},
\end{align}
where a scaled maximum error is defined as
\begin{align}\label{eq:error_threshold_max}
    e_{\tau, \max} &\defeq \alpha \lvert \ebf_\tau \rvert_{ \infty}.
\end{align}
This normalized error observation includes a user-defined parameter $\alpha$, which can be considered as a form of a relative error threshold that will be used in the reward function to indicate which elements are to be refined or coarsened. A key feature of the proposed approach here is that $\alpha$ \emph{can be varied at evaluation time}, such that the normalized error distribution in the observation can be adjusted by the user. It will be later shown that this allows for the user to control the relative degree of refinement at evaluation time, one of the previously presented critical features of an effective AMR policy. We note here that, in theory, one may observe multiple error components for systems of conservation laws, but we only consider one component in this work which will be discussed in \cref{sec:implementation}.

While the normalized error observation is critical to the AMR policy, it does not contain any of the required information to predict its spatio-temporal evolution. Therefore, it is necessary to include some information about the local state in the observation. In \citet{yang2021reinforcement}, the policy observed the high-order FEM solution on an equispaced grid and the physics of the problem were ostensibly learned through the RL approach. However, this method has some notable disadvantages, with the main drawback being that the RL policy has to learn the physics of the underlying PDE essentially from scratch. Furthermore, this approach is not invariant to mesh/problem scale or remesh time. We instead propose an alternate approach which relies on feature engineering based on some useful properties of hyperbolic systems to form a non-dimensional observable quantity which possesses significantly better generalization properties.

We introduce this approach in the form of an example for the linear advection equation, where the advection velocity is given by a constant vector $\mathbf{c}$.  We define a displacement vector $\mathbf{r}_{ij}$ as
\begin{equation}
    \mathbf{r}_{ij} = \mathbf{x}^c_i - \mathbf{x}^c_j,
\end{equation}
where $\mathbf{x}^c_i$ is the centroid of the element $\Omega_i$. The propagation velocity along this displacement vector can be simply computed as $\mathbf{c} \cdot \mathbf{r}_{ij}/\|\mathbf{r}_{ij}\|_2$. If this quantity is again normalized by $\|\mathbf{r}_{ij}\|_2$, one yields an estimate for the (inverse) time required for a feature to propagate from $\mathbf{x}^c_j$ to $ \mathbf{x}^c_i$. By then introducing the remesh time $T$ as 
\begin{equation*}
    \frac{\mathbf{c} \cdot \mathbf{r}_{ij}}{\|\mathbf{r}_{ij}\|^2_2}T,
\end{equation*}
we recover a \emph{non-dimensional} quantity which estimates the likelihood, although not in a probability sense, that element $\Omega_i$ will encounter a feature presently in $\Omega_j$ in the span of the remesh time $T$, i.e., in the range $[t_\tau, t_\tau +T]$. This novel observation formulation possesses some very useful properties, namely that it is invariant to problem/mesh scale and remesh time and that the quantity is effectively re-normalized around a critical value of unity. 

With this example in mind, we introduce a general form of this novel observable quantity for arbitrary hyperbolic conservation laws. Given $m$ solution variables in the PDE, the next $m$ observable quantities are defined as
\begin{equation}\label{eq:fluxjac_obs}
    o^{ij}_{\tau, k} = \left \langle \mathbf{A}_{kl} \frac{\mathbf{u}_k}{\mathbf{u}_l} \right \rangle_j \cdot \frac{\mathbf{r}_{ij}}{\|\mathbf{r}_{ij}\|^2_2}T \quad \forall \ k \in \lbrace 1, \dotsc, m \rbrace,
\end{equation}
where $\mathbf{A}_{kl}$ is the $k$-th component of the flux Jacobian, defined in \cref{sec:preliminaries}, with respect to the chosen solution variable index $l$ at time $t_\tau$. The notation $\langle \cdot \rangle_j$ denotes the average of the quantity $\cdot$ across the element $\Omega_j$. Similarly to the error observation, one may observe multiple components, but we only consider one component in this work which coincides with the solution component used for the error observation. Furthermore, this component is assumed to be non-zero, such that the quantity in \cref{eq:fluxjac_obs} is always well-defined (e.g., density and total energy in the Euler equations). It must be noted that for nonlinear conservation laws, this non-dimensional likelihood quantity is no longer exact (i.e., a value of unity does not necessarily mean that a feature will propagate across the given distance), but it instead approximates a linearization of the exact quantity around $(\mathbf x, t) = (\mathbf{x}^c_i, t_\tau)$.

Finally, while this non-dimensional likelihood quantity gives an essential representation of the spatio-temporal evolution of the error, the linearized formulation may not present enough information to adequately predict strongly nonlinear interactions (e.g., shocks in the Euler equations). For such problems for the Euler equations that do exhibit strong discontinuities, we also append the element-wise averaged conserved ($[\rho, \mathbf{m}, E]$) and primitive variables ($[\rho, \mathbf{v}, P]$) to the observation as this was found to yield moderately better results. 

\begin{remark}[Aggregation for $h$-refinement]
For $h$-refinement, as the agents are defined on the coarse elements, it is necessary to aggregate the observables from any child elements up to the coarse parent level (i.e., the agent level). For the error field, this is performed as an $L^2$ norm across the child elements, whereas for the other derived quantities such as the flux Jacobian and solution, a simple area-weighted average was performed. Any error normalization and element-wise averaging was performed prior to the aggregation. 
\end{remark}


\subsubsection{Action space}
Each agent $i$ then chooses its individual action $a^i$ from the same discrete action space $\mathcal A$. In the proposed approach, we utilize an absolute action space for refinement/de-refinement, where the action corresponds to an output state of the mesh independently of the input mesh. As such, for the one-level refinement considered in this work, the action space is simply taken as a binary output,
\begin{equation}
    \mathcal{A} \defeq \lbrace 0, 1 \rbrace = \lbrace \texttt{coarse}, \texttt{fine} \rbrace,
\end{equation}
which dictates the refinement level for the chosen agent. For the \texttt{coarse} action, the agent will set that given element's refinement level to the coarse state. Likewise for the \texttt{fine} action, the agent will set that given element's refinement level to the fine state. For $p$-refinement, these actions simply correspond to setting the approximation order for the element, where for some base approximation order $p$, the \texttt{coarse} action sets the element to a $\mathbb P_p$ approximation and the \texttt{fine} action elevates it to a $\mathbb P_{p+1}$ approximation. If the current state of the agent is such that the element approximation is of order $p+1$ and the \texttt{coarse} action is chosen, the polynomial reduction is performed with an $L^2$ projection. For $h$-refinement, since the agents are defined over the coarse mesh, the \texttt{coarse} action sets the refinement level of that element to the base level while the \texttt{fine} action sets the refinement level to one higher than the base level, subdividing the initial coarse element into four congruent child elements. Similarly as with $p$-refinement, the reduction from the fine child elements to the coarse parent element is performed using an $L^2$ projection. An example of these refinement actions for both $h$- and $p$-refinement is shown in \cref{fig:hp_ref_schematic}. 

\subsubsection{Transition function}

The initial state of each episode is defined by the initial coarse mesh and the initial conditions that are sampled from a class of parameterized functions, which are defined in \Cref{sec:results}.
Given the agents' actions $\abf_\tau$, the environment transitions from current state $s_\tau$ to next state $s_{\tau+1}$ according to the conditional distribution $P(s_{\tau+1} | s_\tau, \abf_\tau)$ via the following steps:
\begin{enumerate}
    \item Apply refinement to each element whose agent's action is \texttt{refine} and current refinement level is coarse and de-refinement to each element whose agent's action is \texttt{coarse} and current refinement level is fine.
    \item Advance the PDE in simulation time from $t_\tau$ to $t_{\tau + 1} = t_\tau + T$.
    \item If final time is reached ($t = t_f$), reset the environment. 
\end{enumerate}
The transition is deterministic for the PDEs in this work, but the overall methodology applies equally to stochastic transitions in the case of stochastic PDEs.

\subsubsection{Reward}

After an environment transitions with actions $\abf_\tau$, each agent $i$ receives an individual reward $r^i_{\tau+1} \in \mathbb{R}$. Similarly to the observation, the reward function must be carefully crafted to allow for generalizability between different problems, meshes, and target error/cost. We utilize a formulation in which the agents are penalized for doing the ``wrong'' actions, such that the optimal reward obtainable by the policy is zero. In brief, we penalize the agents if they choose the \texttt{coarse} action when refinement is more desirable (i.e., when errors are high) and vice versa. The reward function is explicitly defined as 
\begin{align}\label{eq:reward}
r^i_{\tau+1}(s_\tau,\abf_\tau) =
    \begin{cases}
    - p_{\text{ur}} \bigl\lvert \log_{10}\left(\hat{e}^i_{\tau+1} / e_{\tau, \max}\right) \bigr\rvert, \quad &\text{if $\hat{e}^i_{\tau+1} > e_{\tau, \max}$ and $a^i_\tau = \texttt{coarse}$,} \\
    - p_{\text{or}} \bigl\lvert \log_{10}(\hat{e}^i_{\tau+1} / e_{\tau, \min}) \bigr\rvert, \quad &\text{if $\hat{e}^i_{\tau+1} < e_{\tau, \min}$ and $a^i_\tau = \texttt{fine}$,}\\
    0, &\text{else},
    \end{cases}
\end{align}
where $p_{ur}$ and $p_{or}$ are strictly positive hyperparameters representing the penalty factors for under-refinement and over-refinement, respectively, and $e_{\tau, \max}$ and $e_{\tau, \min}$ are maximum and minimum error thresholds, respectively, computed at time $t_\tau$. The use of threshold values at $t_\tau$ rather than at $t_{\tau+1}$ ensures that the reward depends only on the normalization applied to the observations that the policy used to produce action $a^i_\tau$. Crucially, this reward formulation is time-independent and does not depend on any global spatial information, which means policies trained with this reward can in principle generalize to any simulation duration and mesh size (retaining the same observation normalization) at evaluation time. 

For the thresholds in \cref{eq:reward}, the maximum error threshold $e_{\tau, \max}$ is computed identically to \cref{eq:error_threshold_max} while the minimum error threshold is computed as 
\begin{equation}\label{eq:error_threshold_min}
    e_{\tau, \min} = e_{\tau, \max}^\beta, 
\end{equation}
where $\beta > 1$ is a hyperparameter representing the error threshold hysteresis. We remark that this formulation is under the assumption that $e_{\tau, \max} < 1$, such that $e_{\tau, \min} < e_{\tau, \max}$, which is typically fair to assume if the problem setup is appropriately normalized. These hyperparameters are chosen purely for the purpose of training the policy and, unlike the user-defined parameter $\alpha$ in the observation, cannot affect the policy at evaluation time. In this work, the parameter $\beta$ is fixed. 

\begin{figure}[htbp!]
    \centering
    \subfloat[Raw error distribution]{
    \adjustbox{width=0.4\linewidth,valign=b}{\begin{tikzpicture} [spy using outlines={rectangle, height=3cm,width=2.3cm, magnification=3, connect spies}]
    \begin{axis}
    [   axis lines = none,
        ybar interval, 
        ymax=1,
        ymin=0,
        xmax=1.05,
        xmin=0,
        ticks=none,
        x post scale=1.2,]
        
    \addplot [color=black!70, fill=blue!10] coordinates { 
        (0.0, 0.2189)
        (0.05, 0.5104)
        (0.1, 0.6607)
        (0.15, 0.7177)
        (0.2, 0.7153)
        (0.25, 0.6768)
        (0.3, 0.6181)
        (0.35, 0.5498)
        (0.4, 0.4789)
        (0.45, 0.4098)
        (0.5, 0.3449)
        (0.55, 0.2856)
        (0.6, 0.2325)
        (0.65, 0.1855)
        (0.7, 0.1443)
        (0.75, 0.1086)
        (0.8, 0.0777)
        (0.85, 0.0512)
        (0.9, 0.0283)
        (0.95, 0.0087)};
    
    \end{axis}

    \begin{axis}
    [   axis line style={latex-latex},
        axis y line=left,
        axis x line=left,
        xmode=linear,
        ymode=linear,
        xlabel = {$\mathbf{e}_{\tau}$},
        ylabel = {$n$},
        xmin = 0, xmax = 1.05,
        ymin = 0, ymax = 1,
        xticklabels={,,},
        yticklabels={,,},
        x tick style={draw=none},
        y tick style={draw=none},
        x post scale=1.2,
        extra x ticks={1},
        extra x tick style={major x tick style={draw,black}, xticklabel=$|\mathbf{e}_{\tau}|_\infty$},
        clip mode=individual,
        label style={font=\large},
    ]
    \draw[domain=0:1, smooth, variable=\x, black, style={very thick}]  plot ({\x}, {\x*10*exp(-\x*5) - \x*0.0674});

    \end{axis}

\end{tikzpicture}}}
    \subfloat[Normalized error distribution]{
    \adjustbox{width=0.4\linewidth,valign=b}{\begin{tikzpicture} [spy using outlines={rectangle, height=3cm,width=2.3cm, magnification=3, connect spies}]
    \begin{axis}
    [   axis lines = none,
        xmin = 0, xmax = 1.05,
        ymin = 0, ymax = 1,
        x post scale=1.2
    ]
    \addplot[domain=0:0.667, name path=f1, smooth, variable=\x, black, style={very thick}]  ({\x}, {\x*10*1.5*exp(-1.5*\x*5) - 1.5*\x*0.0674});
    \addlegendentry{$\alpha_1$}
    
    \addplot[domain=0:0.8333, name path=f2, smooth, variable=\x, black, style={very thick, dashed}]  ({\x}, {\x*10*1.2*exp(-\x*5*1.2) - \x*0.0674*1.2});
    \addlegendentry{$\alpha_2$}
    
    \addplot[domain=0:1, name path=f3, smooth, variable=\x, black, style={very thick, dotted}]  ({\x}, {\x*10*exp(-\x*5) - \x*0.0674});
    \addlegendentry{$\alpha_3$}
    
    \path[name path=axis] (axis cs:0,0) -- (axis cs:1,0);
    
    \addplot [thick, color=black, fill=red, fill opacity=0.2]
    fill between[of=f1 and axis,
        soft clip={domain=0.5:0.667}
    ];
    
    \addplot [thick, color=black, fill=red, fill opacity=0.2]
    fill between[of=f2 and axis,
        soft clip={domain=0.5:0.8333},
    ];
    \addplot [thick, color=black, fill=red, fill opacity=0.2]
    fill between[of=f3 and axis,
        soft clip={domain=0.5:1.0},
    ];
    
    \addplot[black!50, style={semithick}]  coordinates { 
        (0.5, 0.0)
        (0.5, 0.3767)};

    \draw[->] (0.45, 0.5) -- (0.525, 0.65);
    \draw (0.55, 0.7) node[scale=1] {Decreasing $\alpha$};

    \addplot [thick, color=black, fill=blue, fill opacity=0.2]
    fill between[of=f1 and axis,
        soft clip={domain=0.0:0.25}
    ];
    \addplot[black!50, style={semithick}]  coordinates { 
        (0.25, 0.0)
        (0.25, 0.55)};
        
    \addplot [thick, color=black, fill=blue, fill opacity=0.2]
    fill between[of=f2 and axis,
        soft clip={domain=0.0:0.225},
    ];
    \addplot[black!50, style={semithick}]  coordinates { 
        (0.225, 0.0)
        (0.225, 0.682)};

    \addplot [thick, color=black, fill=blue, fill opacity=0.2]
    fill between[of=f3 and axis,
        soft clip={domain=0.0:0.2},
    ];
    \addplot[black!50, style={semithick}]  coordinates { 
        (0.2, 0.0)
        (0.2, 0.722)};
    
    \draw[->] (0.25, 0.8) -- (0.1, 0.8);
    \draw (0.175, 0.85) node[scale=1] {Increasing $\beta$};

    \end{axis}

    \begin{axis}
    [   axis line style={latex-latex},
        axis y line=left,
        axis x line=left,
        xmode=linear,
        ymode=linear,
        xlabel = {$\bar{\mathbf{e}}_{\tau}$},
        ylabel = {$n$},
        xmin = 0, xmax = 1.05,
        ymin = 0, ymax = 1,
        xticklabels={,,},
        yticklabels={,,},
        x tick style={draw=none},
        y tick style={draw=none},
        x post scale=1.2,
        extra x ticks={0.5},
        extra x tick style={major x tick style={draw,black}, xticklabel=$1.0$},
        clip mode=individual,
        label style={font=\large},
    ]    
    \end{axis}
\end{tikzpicture}}}
    \caption{\label{fig:normalization} Schematic for an example raw error distribution (left) transformed to a normalized error distribution (right) as a function of the hyperparameters $\alpha$ and $\beta$. Red region denotes normalized error values for which the policy penalizes de-refinement, blue region denotes normalized error values for which the policy penalizes refinement.}
\end{figure}
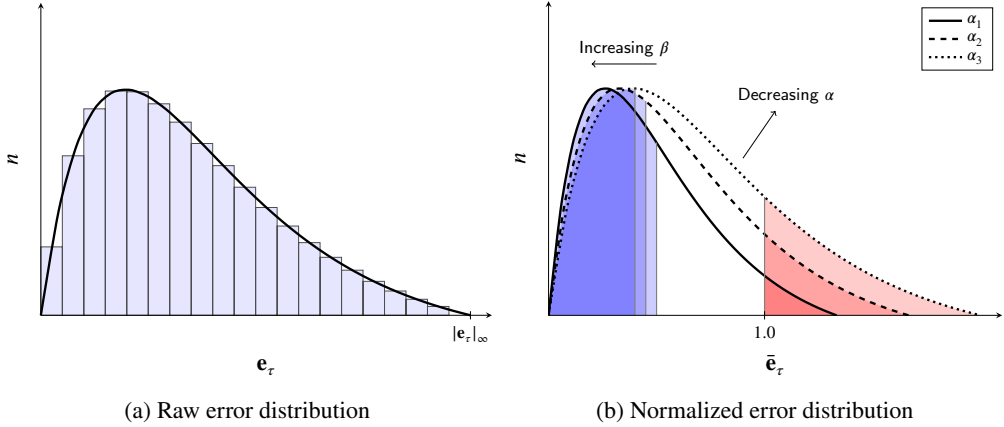

An example of how this error normalization and reward functional operate is shown in \cref{fig:normalization}. Given a raw error distribution, the error is re-normalized around a unit value per \cref{eq:error_threshold_max} based on the threshold parameter $\alpha$, much like in the observation. From the reward functional, the policy is penalized for applying any de-refinement actions on elements with normalized errors above this critical value, shown by the red region in the figure. Furthermore, the policy is also penalized for applying refinement actions to elements with normalized error values below the minimum error threshold computed by \cref{eq:error_threshold_min}, shown by the blue region in the figure. 
Therefore, the policy is effectively incentivized to maintain the element error distributions within the middle (white) region. The key feature of the observation and reward functional is that the raw error represented by the white region \emph{is controlled by the threshold parameter $\alpha$}. As a result, by varying the parameter $\alpha$, the user can control the raw target error at evaluation time even though the reward is not used then since the policy is optimized with respect to a normalized error in the observation which itself can be varied at evaluation time. 

Note that the reward is computed with respect to a modified error estimate $\hat{e}^i_{\tau + 1}$ instead of the instantaneous error estimate $e^i_{\tau + 1}$. We claim that the instantaneous error at time $t_{\tau+1}$ is not necessarily a good indicator of the appropriateness of the action at time $t_{\tau}$, nor is any combination of the instantaneous error estimates $e^i_{\tau}$ and $e^i_{\tau + 1}$. A simple demonstration of the failure mode introduced by using the instantaneous error in the reward can be exemplified through the advection of some compactly-supported characteristic feature of interest, where the feature itself on the order of the element scale. Over the remesh time $T$, this feature may move through multiple elements (recall that longer remesh times are a desirable feature of the AMR policy), such that it is possible that elements which have encountered the feature in that time interval will still have near-zero error estimates since the feature is no longer currently in that element. As such, a reward based on the instantaneous error may penalize agents which correctly refine elements that experience complex features in the time interval $[t_\tau, t_{\tau +1}]$. To mitigate this failure mode, we utilize a modified error estimate $\hat{e}^i_{\tau + 1}$ in the reward, where the modified error estimate is computed as
\begin{equation}
    \hat{e}^i_{\tau + 1} = \underset{t\in [t_\tau, t_{\tau +1}]}{\max} \hat{e}^i(t),
\end{equation}
i.e., the maximum error estimate for the element across the time range $[t_\tau, t_{\tau +1}]$. This maxima is computed discretely across the discrete time steps of the temporal integration method, which are sufficiently small such that $T >> \Delta t$. Much like with the observation, the modified error field for the $h$-refinement case is aggregated to the coarsest (agent) level. 

An algorithmic overview of the DynAMO approach which combines these individual components is shown in \cref{alg:dynamo}.

\begin{algorithm}[htb!]
\caption{\label{alg:dynamo}AMR environment step subroutine in DynAMO.}
\label{alg:step}
\begin{algorithmic}[1]
\REQUIRE $\abf_\tau$: actions of all agents
\STATE Update mesh with actions $\abf_\tau$.
\STATE Advance PDE in simulation time on new mesh from $t_\tau$ to $t_{\tau + 1} = t_\tau + T$, computing the maximum element-wise error $\hat{\ebf}$ over the time interval. 
\STATE Compute new element errors $\ebf_{\tau+1}$.
\STATE Compute reward $\rbf_{\tau+1}$ from $\ebf_{\tau+1}$ and thresholds $\ebar_{\tau,\max}, \ebar_{\tau,\min}$ using \cref{eq:error_threshold_max}, \cref{eq:error_threshold_min}, and \cref{eq:reward}.
\STATE Compute new error thresholds $\ebar_{\tau+1,\max}$ and $\ebar_{\tau+1,\min}$ from $\ebf_\tau$ using \cref{eq:error_threshold_min} and \cref{eq:error_threshold_max}.
\STATE Compute next observations $\obf_{\tau+1}$ from $\ebf_{\tau+1}$ and new error thresholds using \cref{eq:obs} and \cref{eq:fluxjac_obs}.
\ENSURE $\rbf_{\tau+1}, \obf_{\tau+1}, \texttt{done}$
\end{algorithmic}
\end{algorithm}

\subsection{Independent Proximal Policy Optimization algorithm}

We employ Independent PPO, a simple yet effective method for fully-decentralized multi-agent reinforcement learning based on Proximal Policy Optimization (PPO) \citep{schulman2017proximal}. Each agent $i$ learns a policy $\pi^i_{\mu}(a^i|o^i)$, parameterized by a neural network with trainable weights $\mu$, to maximize the objective in \cref{eq:objective}. A value function $V_{\phi}(o^i_{\tau})$, parameterized by another neural network with trainable weights $\phi$, is used to estimate the long-term value of each observation to update the policy at each training step. To improve efficiency, we employ \emph{parameter-sharing} and \emph{experience-sharing} for all agents, which means individual policies $\pi^i_{\mu}, \forall\ i \in \lbrace 1, \dotsc, N \rbrace$ are represented by one and the same policy network $\pi_{\mu}$ that is trained using experiences from all agents. This enables faster learning and convergence than completely independent learning using individual experiences \citep{tan1993multi}.
Crucially, note that each agent still acts differently based on its own local information, since each action sample $a^i_\tau \sim \pi_{\mu}(\cdot|o^i_\tau)$ for agent $i$ is conditioned on its own individual observation.

Although independent learning faces challenges such as the lack of centralized global information and environment non-stationary from each agent's perspective due to the changing policies of other agents, it has significant benefits in the context of dynamic mesh adaptation. As previously mentioned, the spatial domain of dependency for each element is bounded over the time interval $T$, which means it is sufficient for each element's refinement decision to be based on an adequately large local window instead of the full global state. Furthermore, the use of individual rewards means each agent receives unambiguous feedback on the optimality of its own action, which circumvents the multi-agent credit assignment problem that arises in the case of a team reward \citep{chang2004all}. A pseudo-code overview of the DynAMO training routine using Independent PPO is shown in \cref{alg:training}.




\begin{algorithm}[htb!]
\caption{\label{alg:training}DynAMO training routine.}
\begin{algorithmic}[1]
\FOR{each training iteration}
    \STATE $\obf_0 \leftarrow$ env.reset() \RCOMMENT{Reset environment}
    \FOR{time step $\tau = 1,\dotsc,n_t$}
        \STATE $\abf_\tau \leftarrow \pibf(\abf_\tau|\obf_\tau)$ \RCOMMENT{Get actions for all agents}
        \STATE $\rbf_{\tau+1}, \obf_{\tau+1}, \texttt{done} \leftarrow$ env.step($\abf_\tau$) \RCOMMENT{\Cref{alg:step}}
        \STATE Store transition $(\obf_\tau, \abf_\tau, \rbf_{\tau+1}, \texttt{done})$ into trajectory buffer $\Bcal$
        \IF{\texttt{done = True}}
            \STATE $\obf_{\tau+1} \leftarrow \texttt{env.reset()}$ \RCOMMENT{Reset environment}
        \ENDIF
    \ENDFOR
    \STATE Train on minibatches of transitions from $\Bcal$ by PPO
\ENDFOR
\end{algorithmic}
\end{algorithm}

\section{Implementation}\label{sec:implementation}
\begin{figure}[!htbp]
    \centering
    \adjustbox{width=\linewidth,valign=b}{\includegraphics[]{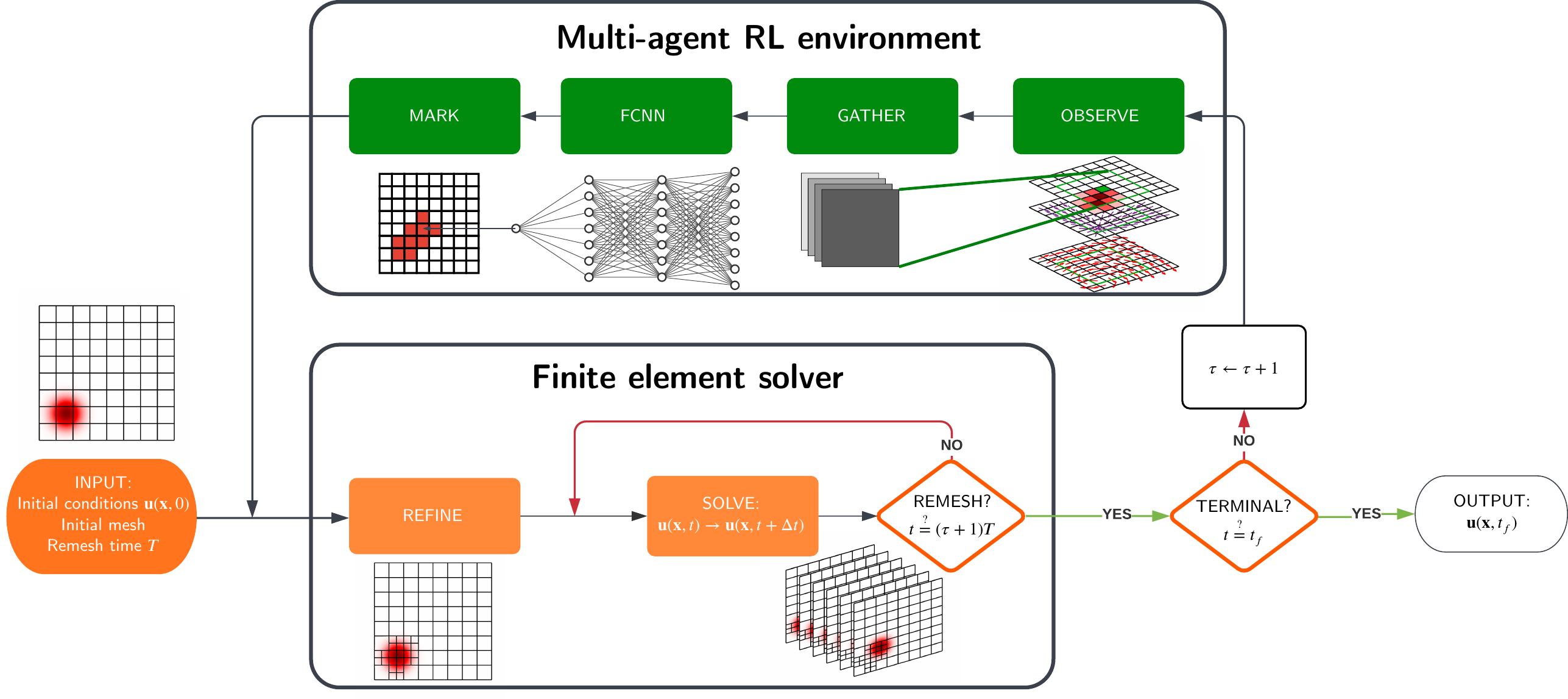}}
    \caption{\label{fig:dynamo_schematic}Schematic of the DynAMO approach.}
\end{figure}

With the DynAMO approach set forth, we present here an in-depth overview of the implementation details of the approach, including numerical frameworks and experimental setups. A schematic of the proposed approach is shown in \cref{fig:dynamo_schematic}.

\subsection{Numerical solver}
The numerical solver component of the proposed approach was implemented using the modular open-source C++ FEM library MFEM \citep{mfem} and its Python interface PyMFEM \citep{pymfem-web}. The FEM approach of choice was a nodal DG \citep{Hesthaven2008} approximation with Gauss--Lobatto solution nodes and the Rusanov approximate Riemann solver~\citep{Rusanov1962}. The Rusanov approach was chosen primarily due to its robustness, although better results may be obtained using less dissipative Riemann solvers. Temporal integration was performed using an explicit fourth-order, four-stage Runge--Kutta scheme with a time step based on the Courant--Friedrichs--Lewy (CFL) condition. The CFL number was generally set to $\text{CFL} = 0.5$ using the standard RKDG estimate of \citet{Cockburn2001}. The base approximation order was set to $\mathbb P_2$ except for problems in the Euler equations which exhibit discontinuities where a $\mathbb P_1$ approximation with the Barth--Jespersen limiter \citep{Barth1989} was used, with the limiting performed with respect to the density field. 

For a given analytic initial condition, the solution was initialized by interpolating the initial conditions to the solution nodes. To compute the estimated error for the observation and reward, we utilized a separate error estimator for $h$- and $p$-refinement. We remark here that these estimators may not satisfy the efficiency and the reliability condition of traditional \textit{a posteriori} error estimators \citep{carstensen2014axioms}, but they were found to be sufficient to yield highly-performant policies. These estimators primarily target discretization error, although more sophisticated estimators, such as ones based on the solution residual, may achieve more accurate results for hyperbolic conservation laws. However, the particular choice of error estimator is separate from the main aim of this work, and a fair comparison can only be made against baseline approaches using the same error estimators.
For $p$-refinement, the error was estimated through a measure of the difference between the polynomial approximation and itself projected to polynomial space of one degree lower, i.e.,
\begin{equation*}
    e^i =\|\mathbf{u}_h^i-\Pi_{p-1}\mathbf{u}_h^i\|_{L^2(\Omega_i)},
\end{equation*}
where $\Pi_{p-1}$ is an $L^2$ projection operator onto the polynomial space of degree $p-1$. Note that $p$ is the current order of the element, which in the $p$-refinement case may not be the base order. For $h$-refinement, the error was estimated by an interface solution jump-based error estimator using a bulk $L^2$-projection similar to the Zienkiewicz--Zhu error estimator \cite{zienkiewicz1992superconvergent1,zienkiewicz1992superconvergent2}.
For an arbitrary element $\Omega_i$, let $\mathcal{N}_e=\{\Omega_+,\Omega_-\}$ be the neighboring elements of an interior edge $e$ and $R$ be the smallest rectangle that contains $\mathcal{N}_e$. We define a polynomial reconstruction operator $\mathcal{R}_e:\mathbb{P}_{p_+}(\Omega_+)\times\mathbb{P}_{p_-}(\Omega_-)\rightarrow \mathbb{P}_r(R)$, where $r=\max\{p_+,p_-\}$ and
\begin{equation*}
    \mathcal{R}_e(\mathbf{u}_h):=\argmin_{\mathbf{v}\in\mathbb{P}_r(R)}\|\mathbf{v} - \mathbf{u}_h\|_{L^2(\mathcal{N}_e)}.
\end{equation*}
Then, the error estimate $e^i$ is defined as
\begin{equation*}
    e^i=\left(\frac{1}{N_{e,i}}\sum_{e\subset \partial \Omega_i\backslash\partial\Omega}\|\mathbf{u}_h-\mathcal{R}_e(\mathbf{u}_h)\|_{L^2(\Omega_i)}^2\right)^{1/2},
\end{equation*}
where $N_{e,i}$ is the number of interior edges contained in the boundary of $\Omega_i$.

While these error estimates were used for the observation and reward, the efficacy of the policies (and the comparison to baseline approaches) was evaluated with respect to the ``true'' error. For the advection equation, this was simply the analytic solution, whereas for the Euler equations, the true error was approximated using a reference simulation with resolution equivalent to one refinement level higher than the maximum refinement level allowable (e.g., for $p$-refinement, the reference simulation was performed at $\mathbb P_{p+2}$).

\subsection{Reinforcement learning framework}
The reinforcement learning framework was implemented using RLLib \citep{rllib}. The PPO approach was implemented using the default hyperparameter values in RLLib except the following: a learning rate of $10^{-4}$, a rollout fragment length of 20, and a batch size of 1000 for each epoch of stochastic gradient descent with minibatches of size 50. The network architecture was a fully-connected neural network with two hidden layers of 256 neurons using a hyperbolic tangent activation function. Purely fully-connected layers were used instead of convolutional neural network (CNN) layers (or a mixture thereof) to maintain the directional structure in the input  and to allow for straightforward extension to unstructured observations through approaches such as graph neural networks~\citep{yang2023reinforcement}. For each policy, training was performed across up to 16 CPUs over the span of 24 hours, which covers approximately $10^4 - 10^5$ training episodes. The most performant policy (i.e., the policy with the highest batch-averaged mean reward) over the training period was used for evaluation. 

The DynAMO-specific hyperparameters were fixed across the various policies trained. The error threshold value at training time was set to $\alpha_{\text{train}} = 0.1$ and the error threshold hysteresis was set to $\beta = 1.2$. Furthermore, the under-refinement and over-refinement penalty factors were set as $p_{ur} = 10$ and $p_{or} = 5$, respectively. Unless otherwise stated, the error threshold value $\alpha$ at evaluation time was set to its training value. For the Euler equations, the solution component for the observed quantities in \cref{eq:obs} and \cref{eq:fluxjac_obs} was the total energy, which was also the component used for evaluating the error against the reference solution. The positivity of total energy in the Euler equations ensures that the quantity in \cref{eq:fluxjac_obs} remains well-defined. For the observation, the element neighborhood size was set to $n_x = n_y = 8$, such that the size of the observation window was $k_x = k_y = 17$. 

\begin{figure}[!htbp]
    \centering
    \adjustbox{width=0.4\linewidth,valign=b}{\includegraphics[]{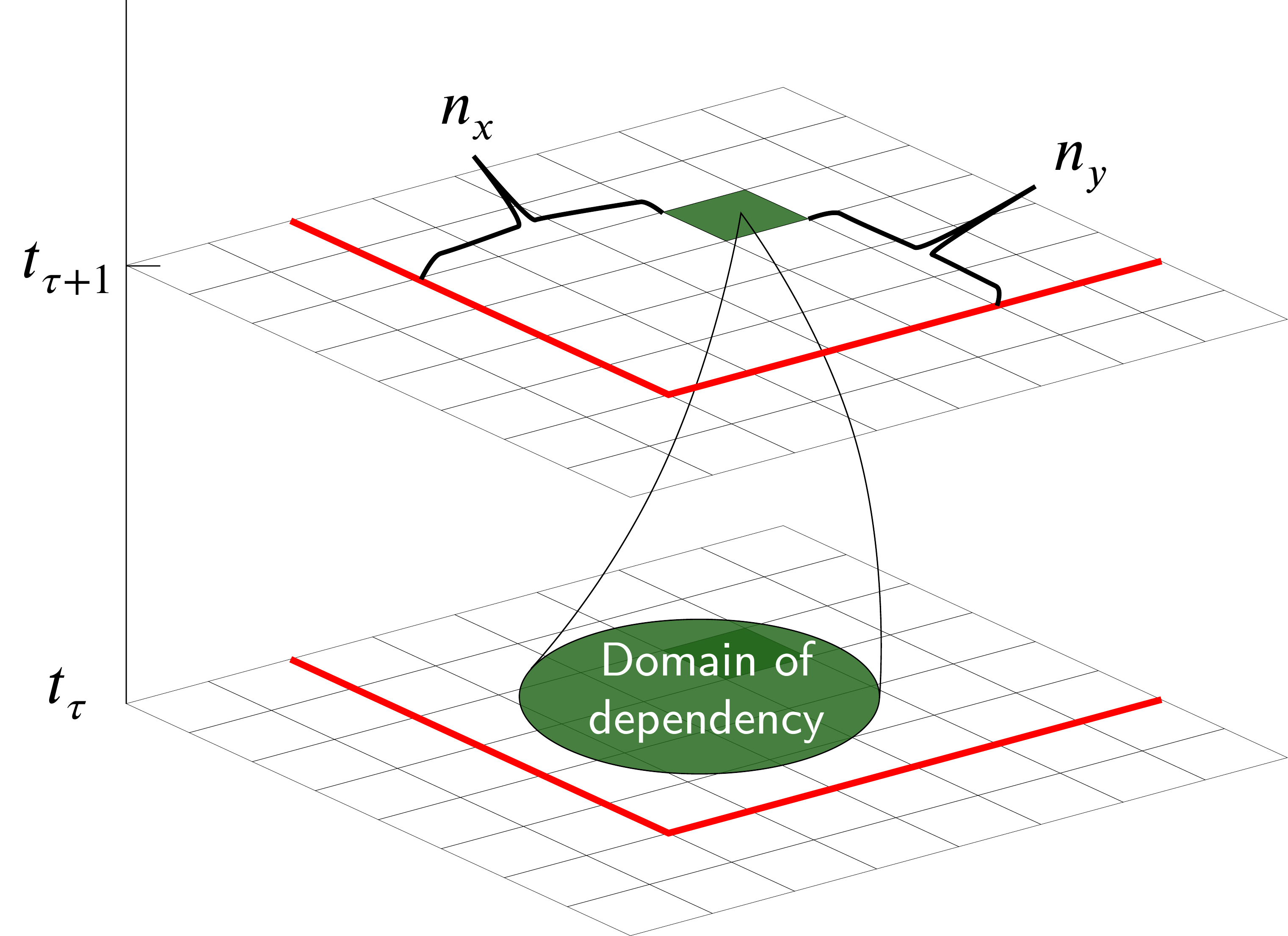}}
    \caption{\label{fig:dod_diagram}The domain of dependency and observation window for an element over a remeshing time interval $T=t_{\tau+1}-t_\tau$.}
\end{figure}

To ensure that domain of influence of a given element does not exceed the observable region of the policy, as exemplified in \cref{fig:dod_diagram}, this observation window size effectively constrained the maximum remeshing time interval $T$ via a CFL-like condition as 
\begin{equation*}
    T \leq \min \left [ \frac{n_x \Delta x}{\lambda_{\max}},  \frac{n_y \Delta y}{\lambda_{\max}} \right ],
\end{equation*}
where $\Delta x$ and $\Delta y$ are the corresponding characteristic mesh spacings along the $x$ and $y$ directions, respectively. The specific remesh time for the experiments was chosen close to this limit during training -- note that this inequality does not necessarily always have to be satisfied, but the performance of the approach would likely suffer in scenarios where it is not. As such, since the policy relies on a fixed observation size (due to the architecture of the network), the remesh time can be modified accordingly if needed to ensure that this condition is met at evaluation time. For the purposes of training, we utilize a fixed episode length of 4 RL steps, corresponding to 4 remesh time intervals, such that the final simulation time was set as $t_f = 4T$. However, this episode length and simulation time was freely varied during evaluation. We remark here that the proposed approach is primarily tailored towards solvers utilizing explicit time stepping, where the observation window can encompass the domain of influence over many solver time steps. Extensions to implicit time schemes, which allow for significantly larger time steps, can be readily performed with the caveat that the observation window must be sufficiently large to account for time steps which themselves may be larger than a standard remesh time interval used for explicit time schemes. Although the assumption of a local domain of influence in hyperbolic conservation laws does not necessarily hold for an implicit time scheme, it  may still offer a suitable approximation for the local solution behavior.

As the purpose of the proposed approach is to increase the efficiency of AMR, it is important to understand the computational overhead of the RL framework at evaluation time to ensure that its cost do not outweigh its benefits. We remark that this work primarily pertains to evaluating the capabilities of RL for AMR, and due to the use of the RLLib framework on top of the FEM solver, the current implementation of the proposed approach is not necessarily tailored towards peak computational efficiency. However, even in its current form, the computational cost overhead of the RL framework was a fraction of the cost of the FEM solve (e.g., for the Euler equations on a $N = 64^2$ mesh, the cost of inference was less than $20\%$ of the cost of the FEM solve). Furthermore, due to the use of an independent multi-agent approach which allows for a highly-parallelizable implementation, the cost of the proposed method can still be drastically decreased with parallelization, a more computationally efficient implementation of the network, and the use of GPU-accelerated computing. 
\section{Results}\label{sec:results}
The proposed DynAMO approach was evaluated on a series of numerical experiments showcasing applications to linear and nonlinear hyperbolic conservation laws with $h$- and $p$-refinement. While the linear advection equation may not be a rigorous test for the approach, it possesses some useful features such as anisotropic solution and error evolution and an analytic solution that allows for the calculation of an exact error. As such, it serves as a suitable preliminary example to present certain features of DynAMO and an initial efficacy comparison to baseline methods. Afterward, the approach is more rigorously evaluated on the compressible Euler equations for both smooth and discontinuous problems, showcasing the viability of DynAMO for complex nonlinear systems. For each set of numerical experiments, a separate policy is trained, such that the results of four policies are shown (one each for advection with $p$-refinement, advection with $h$-refinement, Euler with $p$-refinement, and Euler with $h$-refinement). The goal is to train the policy on smaller representative problems and show that this approach can effectively generalize to larger and unseen problems. 

The results are presented for a series of problems, including \textit{in-distribution} examples, where the problem conditions are similar (but not necessarily exactly identical) to episodes encountered during training, and \textit{out-of-distribution} examples, where the approach is applied to unseen problems, including different initial conditions (which are not representative of the initial conditions used during training), longer simulation time/episode lengths, longer remesh time intervals, and different initial mesh resolutions. To compare the proposed approach to baseline methods, we consider a normalized metric for the computational cost $\bar{c}$ and error $\bar{e}$, defined as 
\begin{equation}
    \bar{c} =  \frac{c - c_{\text{coarse}}}{c_{\text{fine}} - c_{\text{coarse}}} \quad \text{and} \quad 
    \bar{e} = \frac{e - e_{\text{fine}}}{e_{\text{coarse}} - e_{\text{fine}}},
\end{equation}
where the subscripts \textit{coarse} and \textit{fine} refer to the results of simulations performed on completely unrefined and completely refined meshes, respectively. The computational cost metric $c$ is computed as the cumulative total degrees of freedom in the mesh summed between remesh steps. 
From these quantities, we define an \emph{efficiency} metric $0\leq\epsilon \leq 1$ as 
\begin{equation}
    \epsilon = 1 - \sqrt{\bar{c}^2 + \bar{e}^2},
\end{equation}
where a higher efficiency represents an AMR policy that can achieve a higher accuracy for a given computational cost or, conversely, a given accuracy for a lower computational cost. This efficiency metric is the primary point of comparison between DynAMO and baseline methods. For this comparison, we utilize the absolute threshold-based method described in \cref{ssec:amr} as the baseline, which for the Euler equations was applied to the total energy error estimate for consistency with the DynAMO approach. 

\subsection{Advection with $p$-refinement}
\begin{figure}[htbp!]
        \centering
        \adjustbox{width=0.45\linewidth,valign=b}{    \begin{tikzpicture}[spy using outlines={rectangle, height=3cm,width=2.3cm, magnification=3, connect spies}]
		\begin{axis}[name=plot1,
		    axis x line=left,
            axis y line=left,
    		xlabel={Env. steps},
		    ylabel={$-\bar{r}$},
    		xmin=0,
    		ymin=0,
    		ymax=10,
    		style={font=\normalsize}]
    		
			\addplot[color=black, style={thick}]
				table[x=num_env_steps_sampled,
                      y expr=-\thisrow{episode_reward_mean},
                      col sep=comma,unbounded coords=jump]{./figs/data/adv_p_training.csv};
    		    		
		\end{axis} 		
	\end{tikzpicture}}
        \caption{\label{fig:adv_pref_reward} Batch-averaged (negative) reward with respect to number of environment steps during the training process for $p$-refinement on the advection equations.}
    \end{figure}
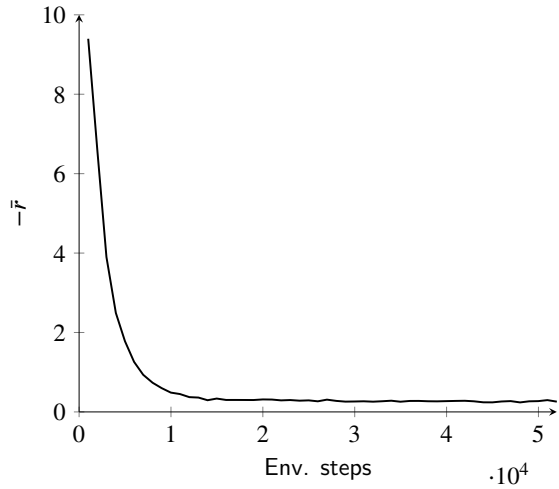

The DynAMO approach was first applied to the task of $p$-refinement for the advection equation, which allows for an evaluation of the approach for AMR on simple linear PDEs with a constant number of mesh elements. We consider the advection of smooth features for which $p$-refinement is particularly well-suited for. For this task, the initial conditions for the training episodes were sampled from a set of ``ring''-like shapes, defined as
\begin{equation}\label{eq:adv_rings}
    u(x,y,0) = 1 + \exp \left(-w\left (\sqrt{(x - x_0)^2 + (y - y_0)^2} - r_0\right)^2 \right),
\end{equation}
where the parameters $x_0$, $y_0$, $r_0$, and $w$ were uniformly distributed across the ranges $x_0 \in [0.3, 0.7]$, $y_0 \in [0.3, 0.7]$, $r_0 \in [0.1, 0.3]$, and $w \in [50, 150]$. This initial condition was chosen to give an extra degree of complexity to the error distribution, adding a region of relatively low error surrounded by a region of high error. Furthermore, the domain was set to $\Omega = [0,1]^2$, the remesh time was set as $T = 0.3$ with a total number of 4 RL steps, the initial mesh resolution was set as $N=24^2$, and the propagation velocity magnitude was sampled across the range $\|\mathbf{c}\|_2 \in [0.7, 1]$, with the direction of propagation uniformly distributed across the azimuthal direction. The efficacy of the agent over the training period is shown through the training reward curve in \cref{fig:adv_pref_reward}.

    \begin{figure}[htbp!]
        \centering
        \subfloat[Remesh at $t = 0$]{
        \adjustbox{width=0.24\linewidth,valign=b}{\includegraphics{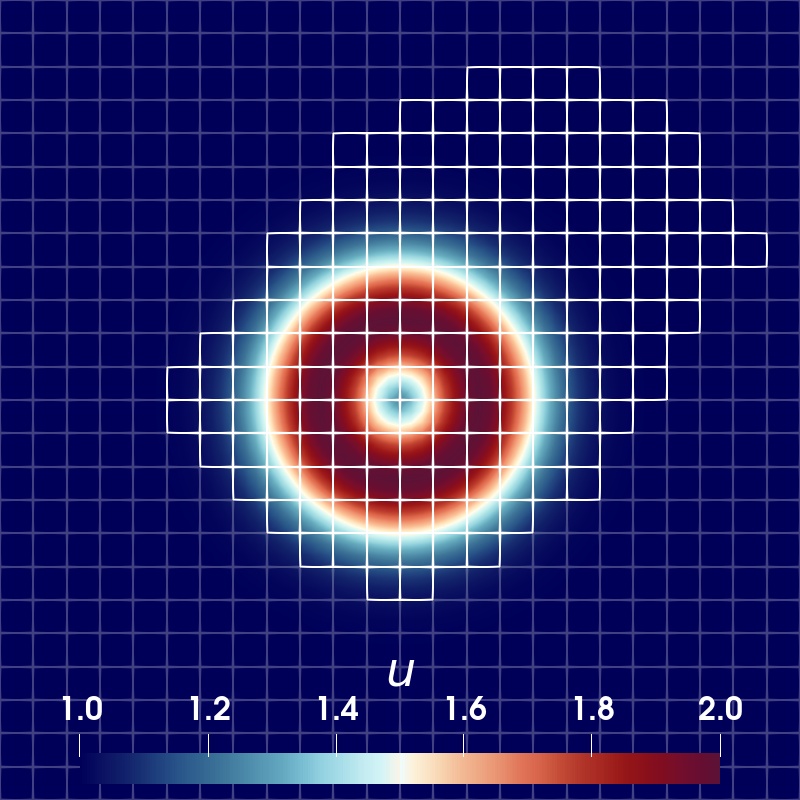}}}
        \subfloat[Solution at $t = T$]{
        \adjustbox{width=0.24\linewidth,valign=b}{\includegraphics{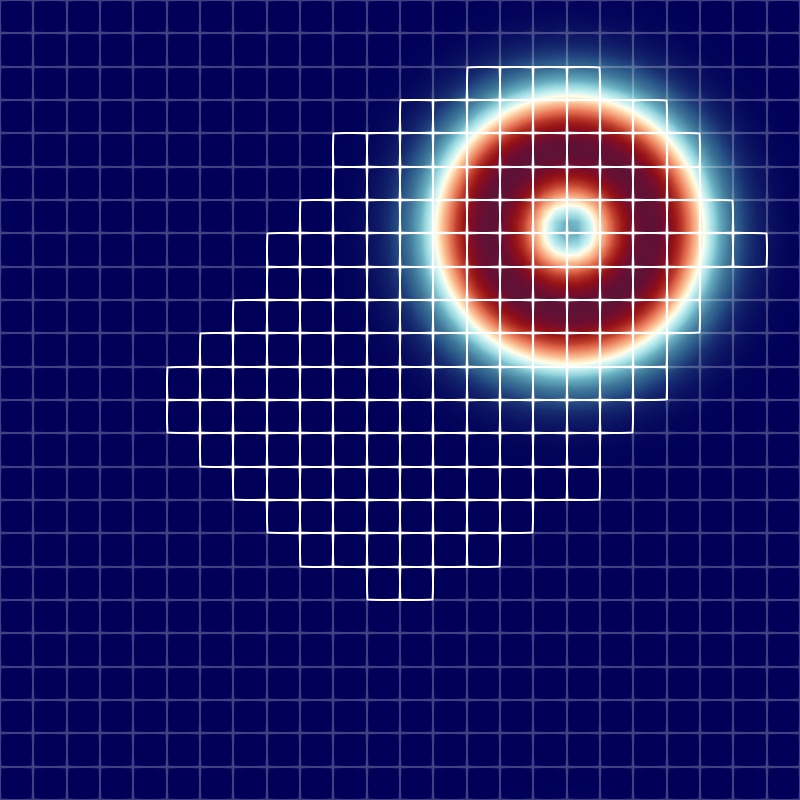}}}
        \subfloat[Remesh at $t = T$]{
        \adjustbox{width=0.24\linewidth,valign=b}{\includegraphics{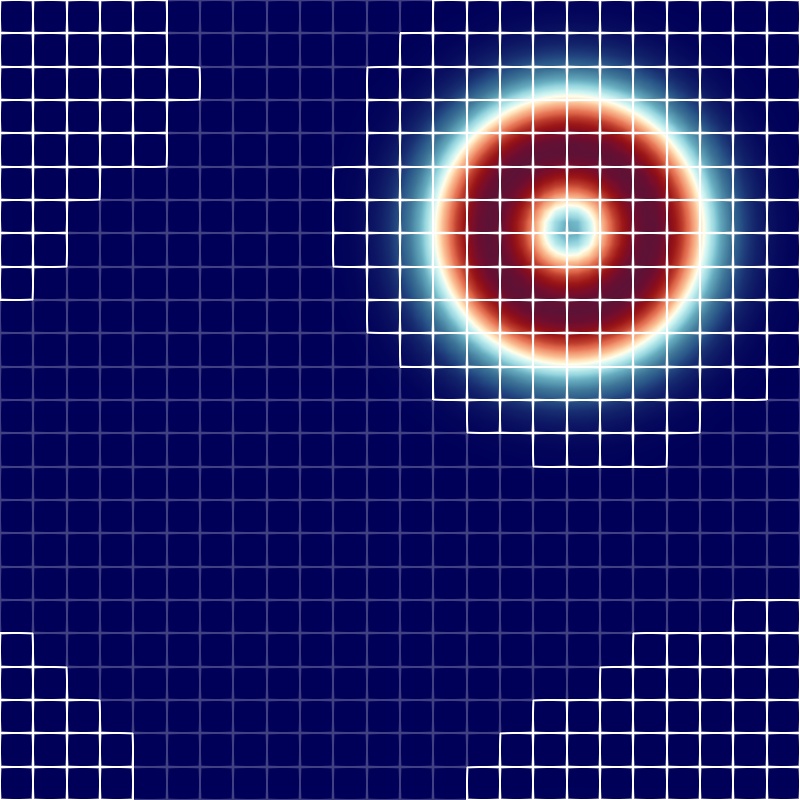}}}
        \subfloat[Solution at $t = 2T$]{
        \adjustbox{width=0.24\linewidth,valign=b}{\includegraphics{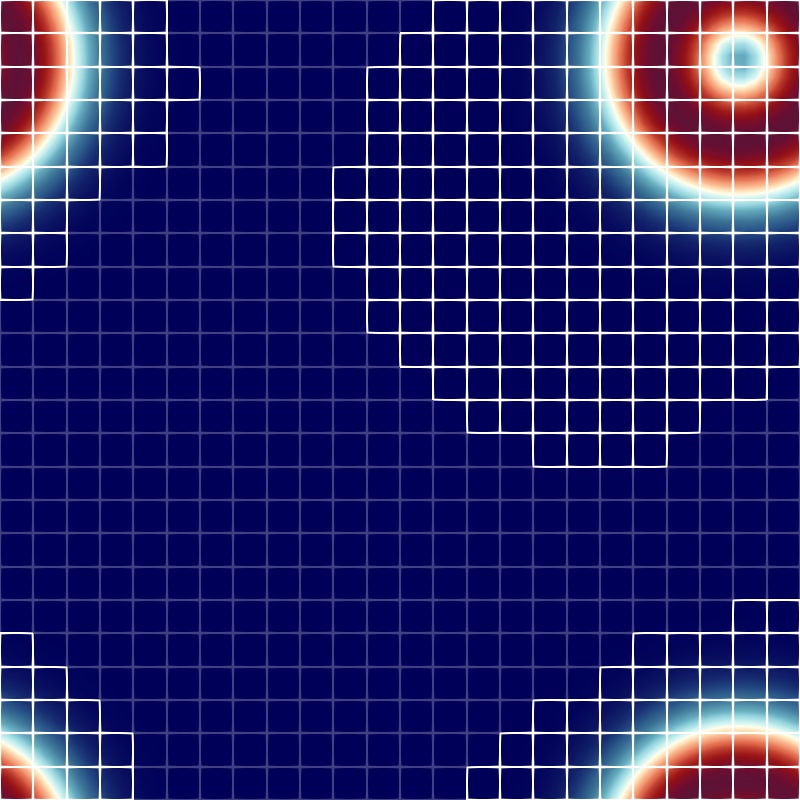}}}
        \newline
        \subfloat[Remesh at $t = 2T$]{
        \adjustbox{width=0.24\linewidth,valign=b}{\includegraphics{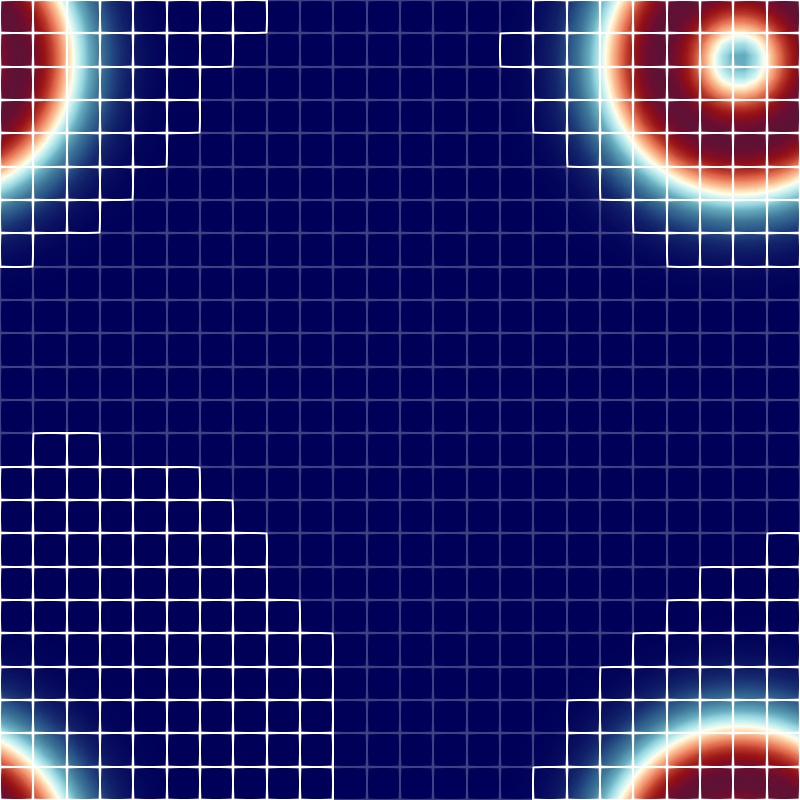}}}
        \subfloat[Solution at $t = 3T$]{
        \adjustbox{width=0.24\linewidth,valign=b}{\includegraphics{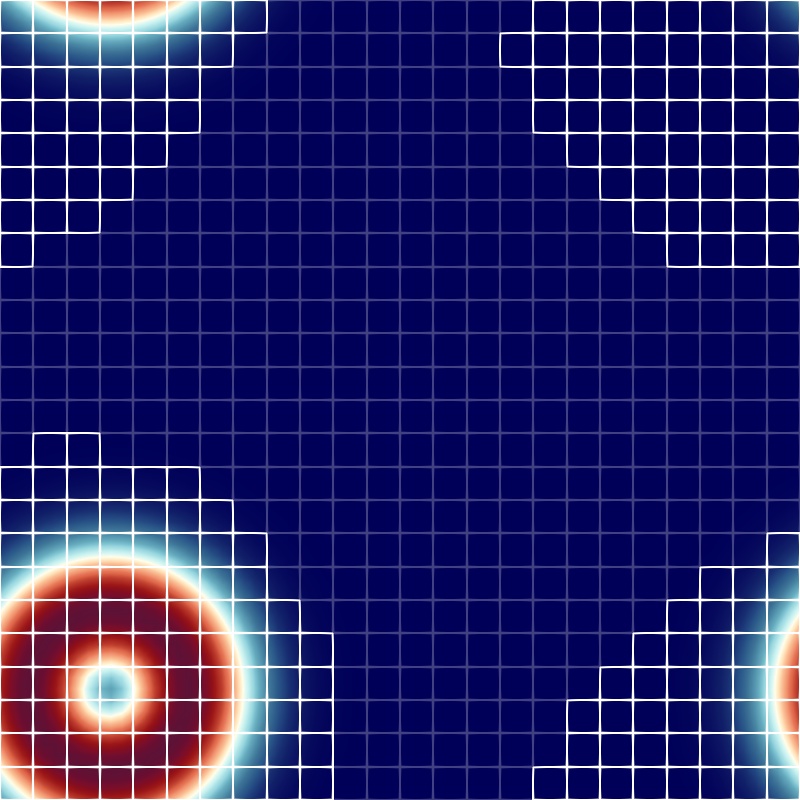}}}
        \subfloat[Remesh at $t = 3T$]{
        \adjustbox{width=0.24\linewidth,valign=b}{\includegraphics{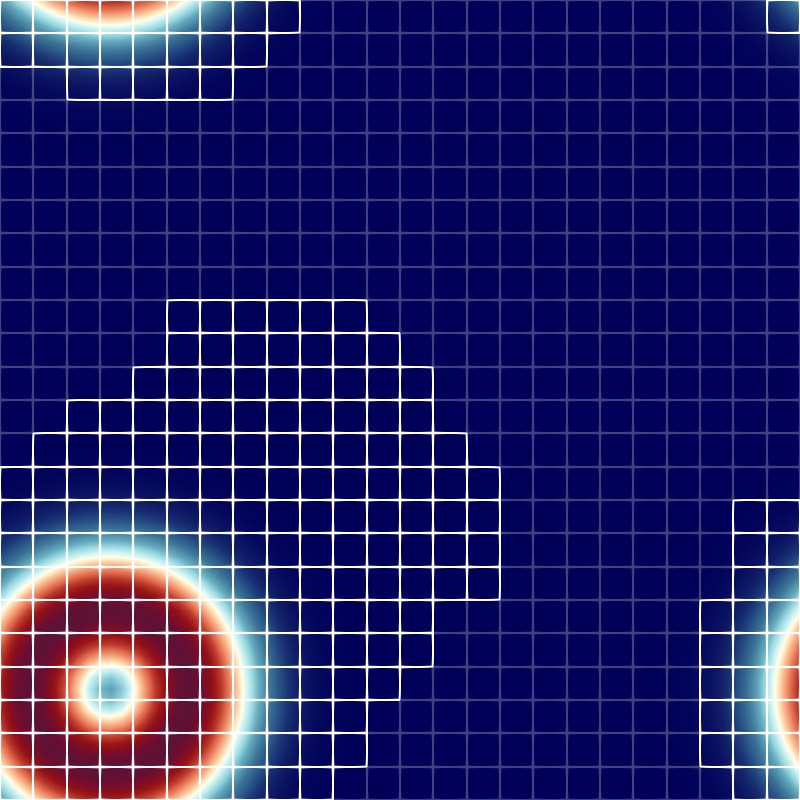}}}
        \subfloat[Solution at $t = 4T$]{
        \adjustbox{width=0.24\linewidth,valign=b}{\includegraphics{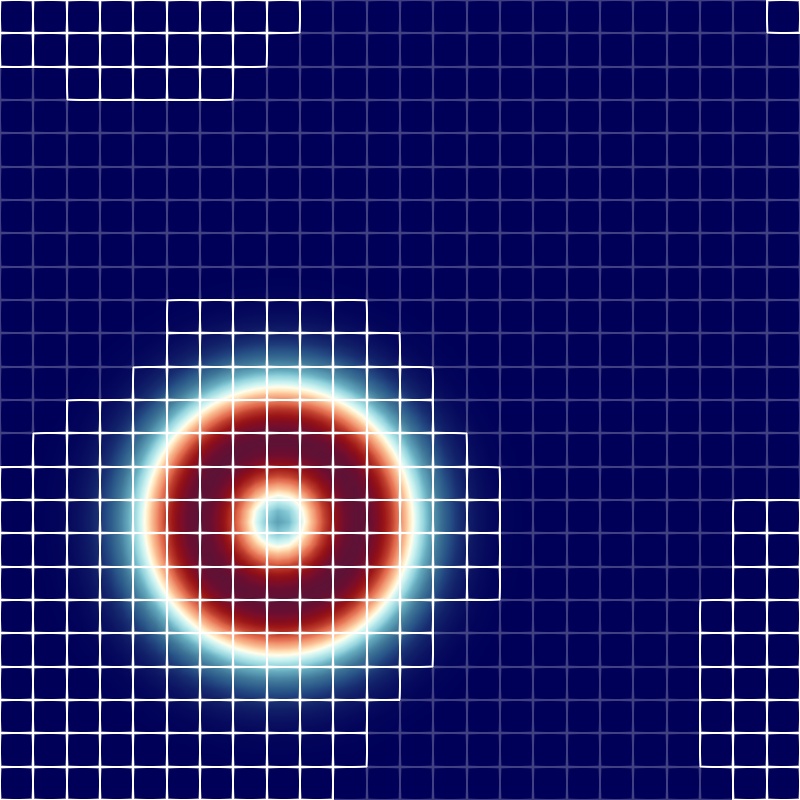}}}
        \newline
        \caption{\label{fig:advpex1_amr_drl} Solution contours overlaid with $p$-adapted mesh at varying remesh intervals using DynAMO for the \textit{in-distribution} advecting ring case. Highlighted elements represent $p$-refinement.} 
    \end{figure}
    
\subsubsection{Sample problem}
\begin{figure}[htbp!]
        \centering
        \subfloat[Remesh at $t = 0$]{
        \adjustbox{width=0.24\linewidth,valign=b}{\includegraphics{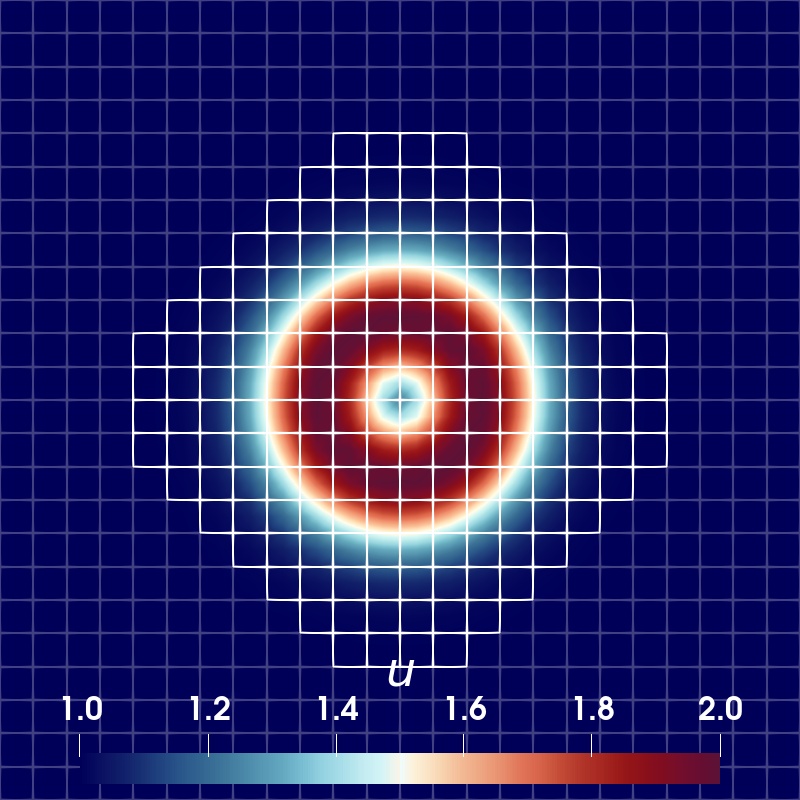}}}
        \subfloat[Solution at $t = T$]{
        \adjustbox{width=0.24\linewidth,valign=b}{\includegraphics{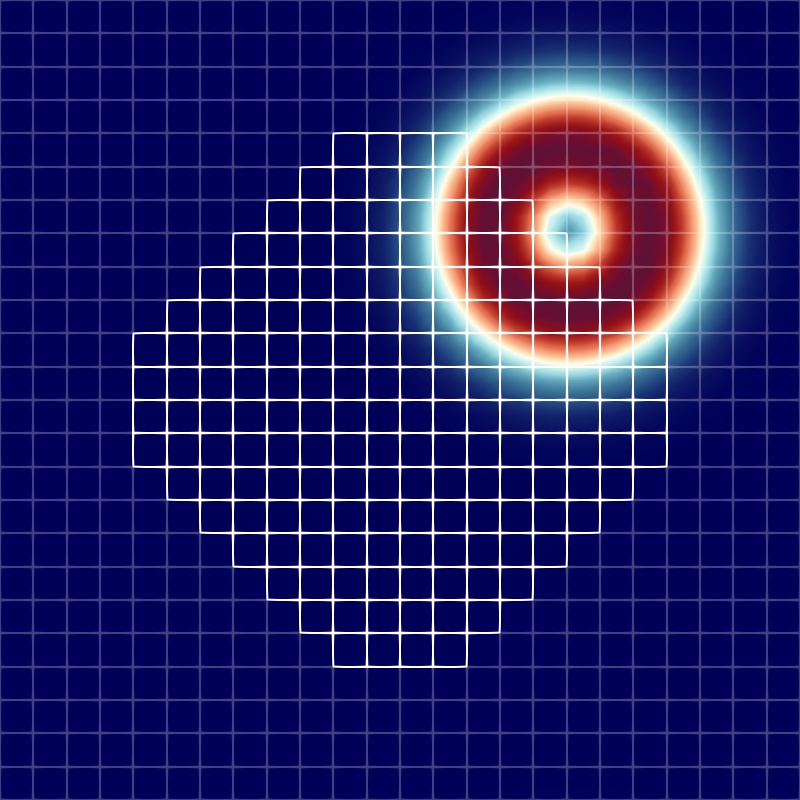}}}
        \subfloat[Remesh at $t = T$]{
        \adjustbox{width=0.24\linewidth,valign=b}{\includegraphics{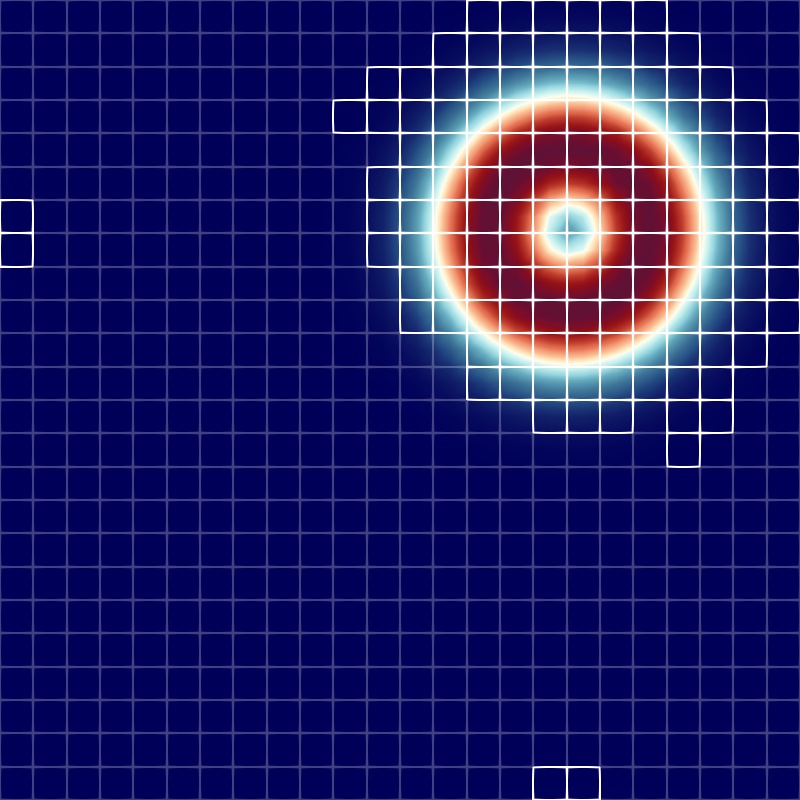}}}
        \subfloat[Solution at $t = 2T$]{
        \adjustbox{width=0.24\linewidth,valign=b}{\includegraphics{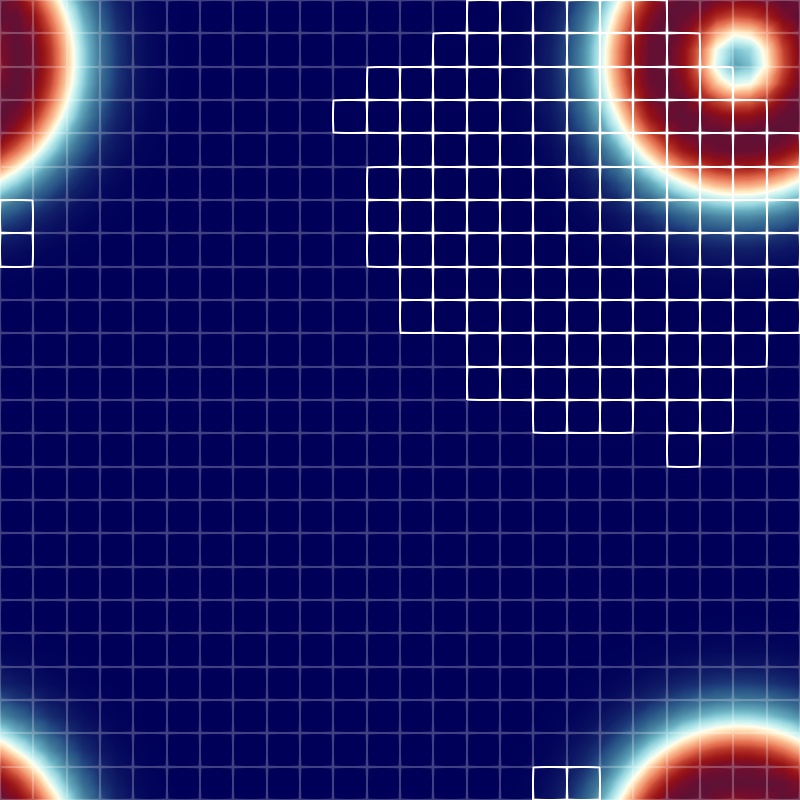}}}
        \newline
        \subfloat[Remesh at $t = 2T$]{
        \adjustbox{width=0.24\linewidth,valign=b}{\includegraphics{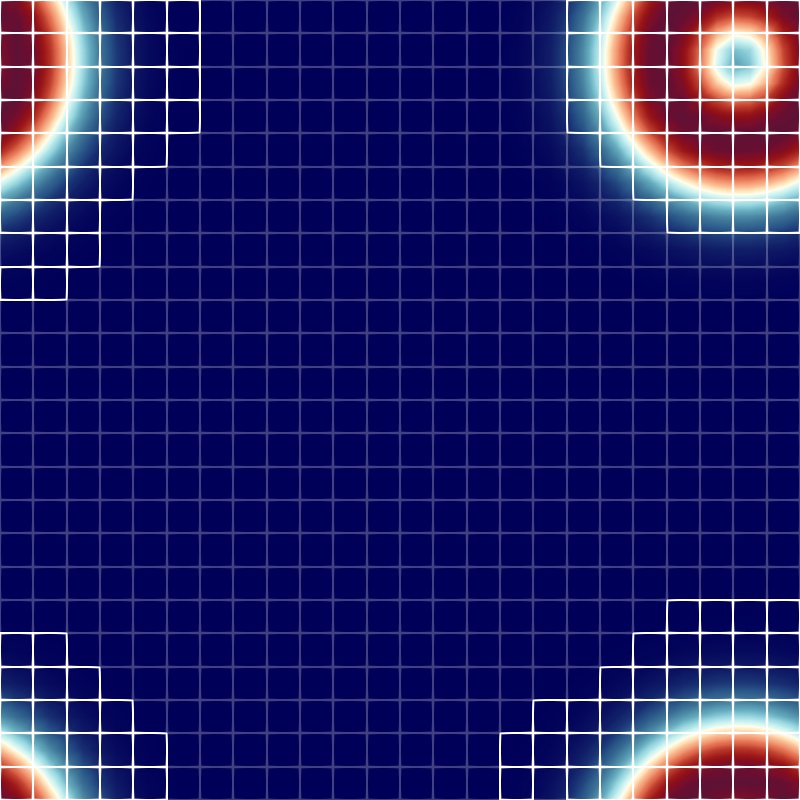}}}
        \subfloat[Solution at $t = 3T$]{
        \adjustbox{width=0.24\linewidth,valign=b}{\includegraphics{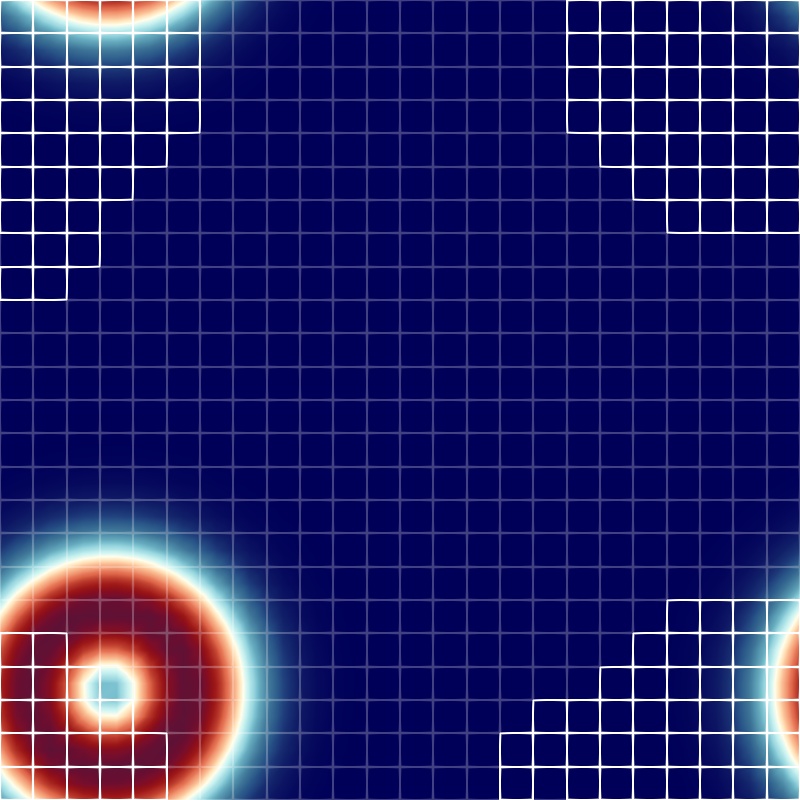}}}
        \subfloat[Remesh at $t = 3T$]{
        \adjustbox{width=0.24\linewidth,valign=b}{\includegraphics{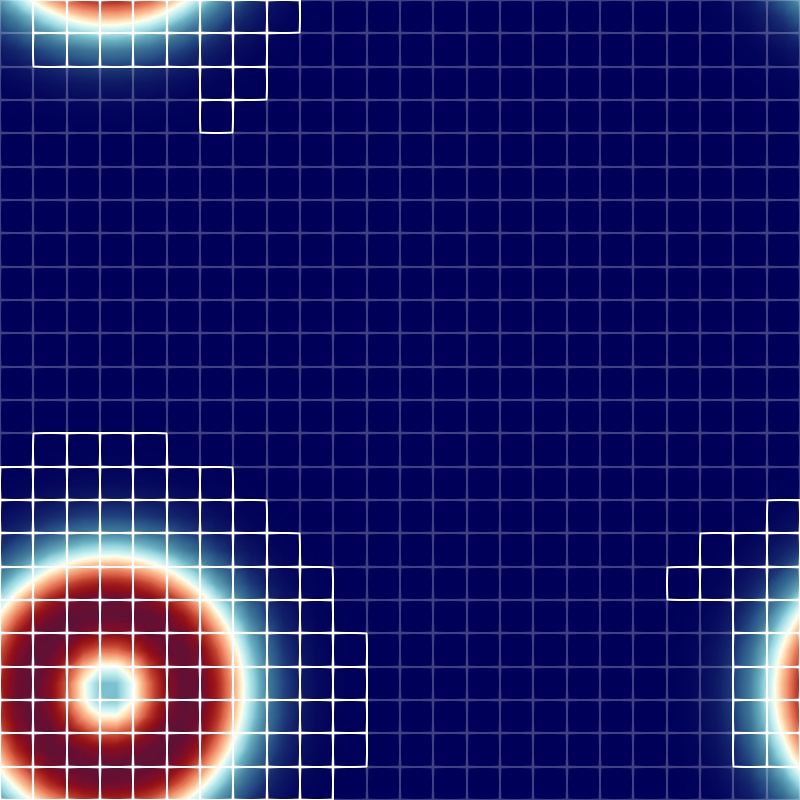}}}
        \subfloat[Solution at $t = 4T$]{
        \adjustbox{width=0.24\linewidth,valign=b}{\includegraphics{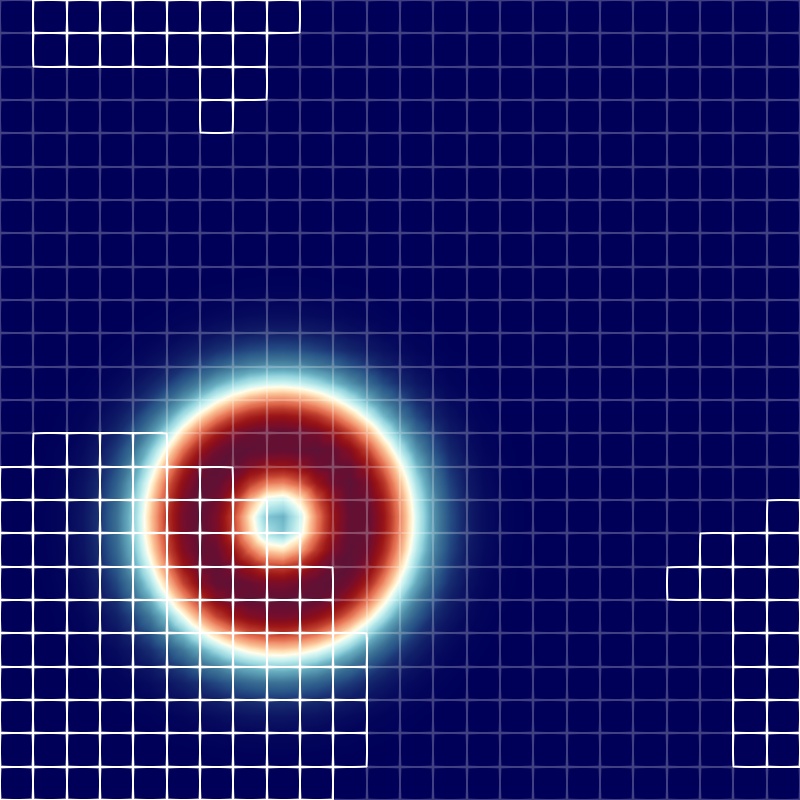}}}
        \newline
        \caption{\label{fig:advpex1_amr_threshold} Solution contours overlaid with $p$-adapted mesh at varying remesh intervals using the threshold policy ($\theta = 10^{-3}$) for the \textit{in-distribution} advecting ring case. Highlighted elements represent $p$-refinement.} 
    \end{figure}

To showcase the features of the DynAMO approach and present an initial comparison to baseline methods, we first present results for a single \textit{in-distribution} test case, where the parameters were set as $x_0 = y_0 = 0.5$, $r_0 = 0.1$, and $\|\mathbf{c}\|_2 = 1$ with a propagation angle of $\pi/4$. The evolution of the solution and $p$-adapted mesh across 4 remesh times as computed by the DynAMO approach and the threshold policy is shown in \cref{fig:advpex1_amr_drl} and \cref{fig:advpex1_amr_threshold}, respectively. The error threshold for the DynAMO policy was set identical to the training value, i.e., $\alpha = \alpha_{\text{train}}$, while the error threshold for the threshold policy was set to $\theta = 10^{-3}$ as this will later be shown to be the value at which this policy achieves the highest efficiency. It can immediately be seen that at $t = 0$, the threshold policy based on the instantaneous error estimate results in an isotropic refinement pattern centered about the ring, which cannot account for the propagation of the ring until the next remesh time. Consequently, the majority of the ring exits the refinement region by $t = T$, causing an increase in error. In contrast, the DynAMO policy shows a highly anisotropic refinement pattern, with the refinement region stretched in the direction of the propagation velocity. As such, it is able to preemptively refine the region through which the ring will propagate between the remesh intervals. In fact, the refinement decisions by the DynAMO policy qualitatively seem to cover \emph{exactly} the region necessary to account for the propagation distance between the remesh intervals without excessively refining ahead of the ring. This behavior persisted throughout the entire simulation time, with the DynAMO policy showing excellent preemptive refinement while the threshold policy was unable to adequately account for the spatio-temporal evolution of the solution and subsequent error. As a result, the DynAMO policy was able to use much longer remesh intervals without sacrificing accuracy in comparison to standard AMR approaches based on instantaneous error estimates.
    
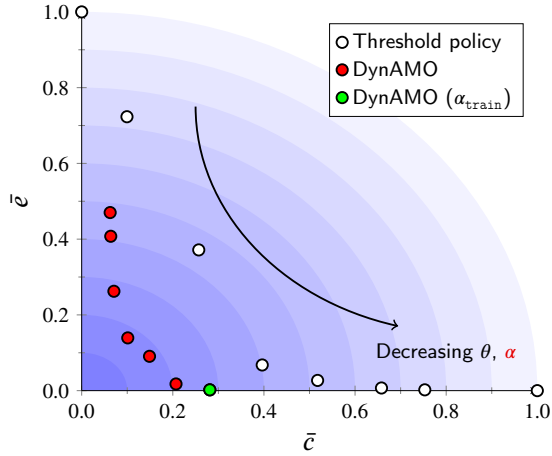
\begin{figure}[htbp!]
        \centering
        \adjustbox{width=0.45\linewidth,valign=b}{\begin{tikzpicture}[spy using outlines={rectangle, height=3cm,width=2.3cm, magnification=3, connect spies}]

    \begin{axis}
    [   axis lines = none,
        xmin = 0, xmax = 1,
        ymin = 0, ymax = 1
    ]
        \fill[fill=blue!5] (1.0, 0) arc[start angle=0, end angle=90, radius=1.0] -- (0,0) -- (1.0, 0);
        \fill[fill=blue!10] (0.9, 0) arc[start angle=0, end angle=90, radius=0.9] -- (0,0) -- (0.9, 0);
        \fill[fill=blue!15] (0.8, 0) arc[start angle=0, end angle=90, radius=0.8] -- (0,0) -- (0.8, 0);
        \fill[fill=blue!20] (0.7, 0) arc[start angle=0, end angle=90, radius=0.7] -- (0,0) -- (0.7, 0);
        \fill[fill=blue!25] (0.6, 0) arc[start angle=0, end angle=90, radius=0.6] -- (0,0) -- (0.6, 0);
        \fill[fill=blue!30] (0.5, 0) arc[start angle=0, end angle=90, radius=0.5] -- (0,0) -- (0.5, 0);
        \fill[fill=blue!35] (0.4, 0) arc[start angle=0, end angle=90, radius=0.4] -- (0,0) -- (0.4, 0);
        \fill[fill=blue!40] (0.3, 0) arc[start angle=0, end angle=90, radius=0.3] -- (0,0) -- (0.3, 0);
        \fill[fill=blue!45] (0.2, 0) arc[start angle=0, end angle=90, radius=0.2] -- (0,0) -- (0.2, 0);
        \fill[fill=blue!50] (0.1, 0) arc[start angle=0, end angle=90, radius=0.1] -- (0,0) -- (0.1, 0);

    \end{axis}
        
    \begin{axis}
    [   axis line style={latex-latex},
        axis y line=left,
        axis x line=left,
        xmode=linear,
        ymode=linear,
        xlabel = {$\bar{c}$},
        ylabel = {$\bar{e}$},
        xmin = 0, xmax = 1,
        ymin = 0, ymax = 1,
        xtick = {0,0.2,0.4,0.6,0.8,1.0},
        ytick = {0,0.2,0.4,0.6,0.8,1.0},
        minor x tick num=1,
        minor y tick num=1,
        legend cell align={left},
        legend style={at={(0.97, 0.97)},anchor=north east},
        clip mode=individual,
        x tick label style={/pgf/number format/.cd, fixed, fixed zerofill, precision=1, /tikz/.cd},
        y tick label style={/pgf/number format/.cd, fixed, fixed zerofill, precision=1, /tikz/.cd},
        label style={font=\large},
    ]
        \addplot[color=black, fill=white, style={thick}, only marks, mark=*, mark options={scale=1.2}]
        table[x = cost_mean, y = err_mean, col sep=comma]{./figs/data/advp_ex1_id_threshold.csv};
        \addlegendentry{Threshold policy};
        
        \addplot[color=black, fill=red, style={thick}, only marks, mark=*, mark options={scale=1.2}]
        table[x = cost, y = err, col sep=comma]{./figs/data/advp_ex1_id_drl_alphasweep.csv};
        \addlegendentry{DynAMO};
        
        \addplot[color=black, fill=green, style={thick}, only marks, mark=*, mark options={scale=1.2}]
        coordinates {(0.28125,0.0021230)};
        \addlegendentry{DynAMO ($\alpha_{{\tt train}}$)};

        \draw[->, black, thick] (0.25,.75) arc (180:255:0.6);
        \draw (0.8, 0.1) node[scale=1] {Decreasing $\theta$, \textcolor{red!90!black}{$\alpha$}};
        
    \end{axis}
\end{tikzpicture}}
        \caption{\label{fig:adv_pref_pareto} Pareto plot of normalized cost vs. error for $p$-refinement on the advection equations with DynAMO and the threshold policy for the \textit{in-distribution} advecting ring case. Contours of efficiency shown on the background. Red markers represent a sweep of $\alpha \in [0.1, 0.8]$ at evaluation time, green marker represents $\alpha$ set to the training value of $\alpha_{\texttt{train}} = 0.1$. Peak efficiency achieved at $\alpha=0.5$.}
    \end{figure}

A more quantitative comparison between the two approaches was performed by analyzing the error, cost, and efficiency for the given test case. A Pareto plot showing the normalized cost vs. error distributions at varying threshold values for the DynAMO and threshold policies is shown in \cref{fig:adv_pref_pareto}. The threshold policy shows a clear Pareto front, with the optimal efficiency achieved at $\theta = 10^{-3}$. As expected, the results of the DynAMO policy sit well within the Pareto front of the threshold policy, showing higher efficiency through a much lower error as a result of the policy's ability to achieve preemptive refinement. This much lower error, represented by the green marker in the figure, was achieved using an identical error threshold at evaluation time as with training time. Furthermore, a key feature of the proposed DynAMO approach is that this error threshold $\alpha$ can be varied at evaluation time to achieve different error and cost targets. The efficacy of this feature was evaluated by changing the error threshold at evaluation time, with the values of $\alpha$ varied between $0.1$ and $0.8$. It can be seen that this results in a similar Pareto front to the threshold policy, where the user can dictate the trade-off between error and cost. However, the key difference is that the DynAMO policy \emph{unlocks regions of efficiency that are unattainable by conventional AMR approaches}. The result is that the DynAMO policy can achieve significantly lower error for a similar computational cost in comparison to these standard AMR approaches. We remark here that the value of the error threshold $\alpha$ used at training time may not necessarily be the value that obtains the highest efficiency at evaluation time, which is exemplified in this case where the optimal efficiency was achieved at $\alpha=0.5$ whereas the training value was set as $\alpha_{\text{train}}=0.1$. This is simply a result of the trade-off between error and cost, where policies that prioritize one over the other will achieve lower efficiency than ones that balance them equivalently and is the expected behavior for a policy that allows users to control error or cost targets. The results from the policy using the training value $\alpha_{\text{train}}$ will just be one point along that Pareto front.

\subsubsection{In-distribution experiments}
While this single example showcases the features and benefits of the DynAMO approach, it is of more interest to observe how its efficacy holds up over a large sample of operating conditions. To this end, a comparison was performed between the DynAMO policy and the threshold policy over 100 \textit{in-distribution} samples of the advecting rings problem, where the operating conditions were uniformly sampled over the training range. The efficiency, normalized error, and normalized cost for the threshold policy at varying values of the absolute threshold $\theta$ is compared to the DynAMO policy (evaluated at $\alpha_{\text{train}}$) in \cref{tab:advection_pref_indistribution}. At its peak, the threshold policy only achieves a mean efficiency of roughly $0.4$, with normalized mean error and cost values of $0.204$ and $0.465$, respectively. In comparison, the DynAMO policy at $\alpha = \alpha_{\text{train}}$, which is not necessarily even the value at which the policy achieves peak efficiency, showed a noticeable increase in mean efficiency from $0.409$ to $0.539$ ($31.8\%$). The extent of the computational benefits provided by the DynAMO approach in comparison to the most efficient threshold policy is best seen in the relative error and cost reduction. For effectively the same mean computational cost, the use of the DynAMO policy resulted in $93.1\%$ less mean error than the threshold policy, such that the DynAMO policy could achieve essentially the accuracy of a fully-refined mesh with less than half of the cost. Achieving this level of accuracy with the threshold policy would require over double the computational cost.

    \begin{figure}[htbp!] 
        \centering
        \begin{tabular}{cccc}
        \toprule
        Method &  Efficiency & Normalized error & Normalized cost\\ 
        \midrule
        Threshold ($\theta = 10^{-2}$) & 0.076 (0.946) & 0.850 (1.008) & 0.087 (0.060) \\
        Threshold ($\theta = 10^{-3}$) &\textbf{ 0.409 (0.313)} & \textbf{0.204 (0.416)} & \textbf{0.465 (0.130) }\\
        Threshold ($\theta = 10^{-4}$) & 0.233 (0.302) & 0.099 (0.370) & 0.715 (0.151) \\
        Threshold ($\theta = 10^{-5}$) & 0.141 (0.185) & 0.057 (0.234) & 0.838 (0.114) \\
        Threshold ($\theta = 10^{-6}$) & 0.083 (0.115) & 0.039 (0.168) & 0.904 (0.081) \\
        Threshold ($\theta = 10^{-7}$) & 0.053 (0.063) & 0.019 (0.090) & 0.943 (0.058) \\
        \midrule
        DynAMO & \textbf{0.539 (0.124)} & \textbf{0.014 (0.052)}& \textbf{0.457 (0.126)} \\
        DynAMO/Optimal $\theta$ & \textcolor{green!70!black}{+31.8\%} & \textcolor{green!70!black}{-93.1\%}& \textcolor{green!70!black}{-1.7\%} \\
        \bottomrule
        \end{tabular}
        \captionof{table}{\label{tab:advection_pref_indistribution} Comparison of the mean efficiency, normalized error, and normalized cost for $p$-refinement on the advection equation with DynAMO and the threshold policy for the advecting rings over 100 \textit{in-distribution} runs using uniform random initial conditions. Standard deviation is shown in parentheses.}
    \end{figure}

\subsubsection{Generalization experiments}
The generalization capabilities of the proposed DynAMO approach were then evaluated by applying the policy to \textit{out-of-distribution} problems, ones with characteristics that were not encountered during training. Three experiments were conducted, testing the ability of the approach to generalize to finer meshes, different shapes for the initial conditions, and longer simulation times. The mean efficiency for the threshold policy at varying values of the absolute threshold $\theta$ and the DynAMO policy evaluated at $\alpha_{\text{train}}$ is shown in \cref{tab:advection_pref_ood} for the three experiments over 100 samples. For the generalization to finer meshes, the initial mesh resolution was increased from $N=24^2$ to $N=96^2$, such that the base finite element space consisted of nearly $10^5$ degrees of freedom. It can be seen that while the increase in mesh resolution increased the peak efficiency of the threshold policy, the efficiency of the DynAMO policy also increased appropriately, such that the relative efficiency benefits of DynAMO stayed approximately the same. This positive result can be attributed to the highly generalizable nature of the novel observation proposed in this work, which is insensitive to mesh resolution and problem scale.

    \begin{figure}[htbp!] 
        \centering
        \begin{tabular}{ccccc}
        \toprule
        Method &  In-distribution & Finer mesh & Different shapes & Longer sim. time \\ 
        \midrule
        Threshold ($\theta = 10^{-2}$) & 0.076 (0.946) & 0.113 (0.186) & 0.000 (0.002) & 0.121 (0.185)\\
        Threshold ($\theta = 10^{-3}$) & \textbf{0.409 (0.313)} & 0.337 (0.255) & 0.121 (0.120) & \textbf{0.429 (0.396)}\\
        Threshold ($\theta = 10^{-4}$) & 0.233 (0.302) & \textbf{0.567 (0.244)} & \textbf{0.483 (0.260)} & 0.281 (0.171)\\
        Threshold ($\theta = 10^{-5}$) & 0.141 (0.185) & 0.487 (0.192) & 0.435 (0.282) & 0.160 (0.120)\\
        Threshold ($\theta = 10^{-6}$) & 0.083 (0.115) & 0.291 (0.190) & 0.307 (0.149) & 0.090 (0.080)\\
        Threshold ($\theta = 10^{-7}$) & 0.053 (0.063) & 0.184 (0.172) & 0.193 (0.112) & 0.050 (0.054)\\
        \midrule
        DynAMO & \textbf{0.539 (0.124)} & \textbf{0.774 (0.135)} & \textbf{0.568 (0.107)} & \textbf{0.840 (0.116)} \\
        DynAMO/Optimal $\theta$ & \textcolor{green!70!black}{+31.8\%} & \textcolor{green!70!black}{+36.5\%} & \textcolor{green!70!black}{+17.6\%}  & \textcolor{green!70!black}{+95.8\%} \\
        \bottomrule
        \end{tabular}
        \captionof{table}{\label{tab:advection_pref_ood} Comparison of the mean efficiency for $p$-refinement on the advection equation with DynAMO and the threshold policy for the advecting rings over 100 \textit{out-of-distribution} runs using uniform random initial conditions with finer mesh resolution, different advecting shapes, and longer simulation time. Standard deviation shown in parentheses. In-distribution results from \cref{tab:advection_pref_indistribution} shown for comparison.}
    \end{figure}

To show that the trained policy is not overfitting to the solution features observed during training, a second generalization experiment was conducted by changing the class of initial conditions from ``ring''-type shapes to ``bump''-type shapes of varying size and width, characterized by the Gaussian function
\begin{equation}
    u(x,y,0) = h \exp \left(-w \left((x - x_0)^2 + (y - y_0)^2 \right) \right),
\end{equation}
where the parameters $x_0$, $y_0$, $h$, and $w$ were uniformly distributed across the ranges $x_0 \in [0.3, 0.7]$, $y_0 \in [0.3, 0.7]$, $h \in [0.2, 1]$, and $w \in [25, 100]$. The efficiency comparison in \cref{tab:advection_pref_ood} shows that the efficacy of the DynAMO policy does not deteriorate with changing initial conditions, but in fact marginally increases. However, due to the simple nature of the bump shape, the efficiency of the threshold policy also increased, such that the efficiency benefits of DynAMO were not as pronounced.

Finally, the ability to generalize to longer simulation times was tested by increasing the simulation length (and total number of remeshing cycles) by a factor of eight. We note here that for simple linear problems such as advection, these experiments may not necessarily show generalization to out-of-distribution problems as the observed features do not vary much over time, but they act as a verification that the fixed episode length used for training does not affect evaluation for longer episode lengths and showcase the type of problems where anticipatory mesh refinement is expected to yield the most benefits. For this experiment, it is expected that an anticipatory refinement approach would yield far superior performance as it would be able to mitigate the introduction of discretization error much more effectively, error which would continuously accumulate over the course of the simulation. As can be seen from \cref{tab:advection_pref_ood}, this was indeed the case, with the DynAMO policy achieving nearly double the efficiency of the standard threshold policy. 

\subsection{Advection with $h$-refinement}
\begin{figure}[htbp!]
        \centering
        \adjustbox{width=0.45\linewidth,valign=b}{    \begin{tikzpicture}[spy using outlines={rectangle, height=3cm,width=2.3cm, magnification=3, connect spies}]
		\begin{axis}[name=plot1,
		    axis x line=left,
            axis y line=left,
    		xlabel={Env. steps},
		    ylabel={$-\bar{r}$},
    		xmin=0,
    		ymin=0,
    		ymax=10,
    		style={font=\normalsize}]
    		
			\addplot[color=black, style={thick}]
				table[x=num_env_steps_sampled,
                      y expr=-\thisrow{episode_reward_mean},
                      col sep=comma,unbounded coords=jump]{./figs/data/adv_h_training.csv};
    		    		
		\end{axis} 		
	\end{tikzpicture}}
        \caption{\label{fig:adv_href_reward} Batch-averaged (negative) reward with respect to number of environment steps during the training process for $h$-refinement on the advection equations.}
    \end{figure}
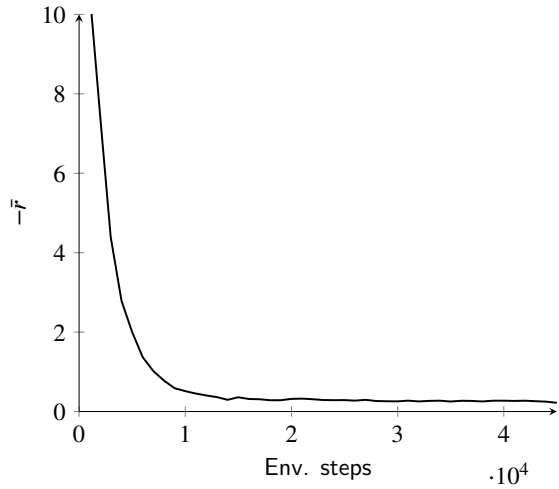

With the successful application of the DynAMO approach to the advection equation with $p$-refinement, the more complex case of $h$-refinement, where the number of elements varies throughout the simulation, was attempted. For consistency with the results of the $p$-refinement experiments, we utilize effectively the same experimental setup, although the results are expected to vary due to differences in the choice of error estimators and the efficacy of $h$-refinement versus $p$-refinement for smooth solutions. The purpose of these experiments is primarily to validate the proposed methodology for $h$-refinement where mesh inhomogeneity is addressed from the coarse agent level. The training reward curve is presented in \cref{fig:adv_href_reward}, showing nearly identical behavior to the $p$-refinement case.

\subsubsection{Sample problem}
\begin{figure}[htbp!]
        \centering
        \subfloat[Remesh at $t = 0$]{
        \adjustbox{width=0.24\linewidth,valign=b}{\includegraphics{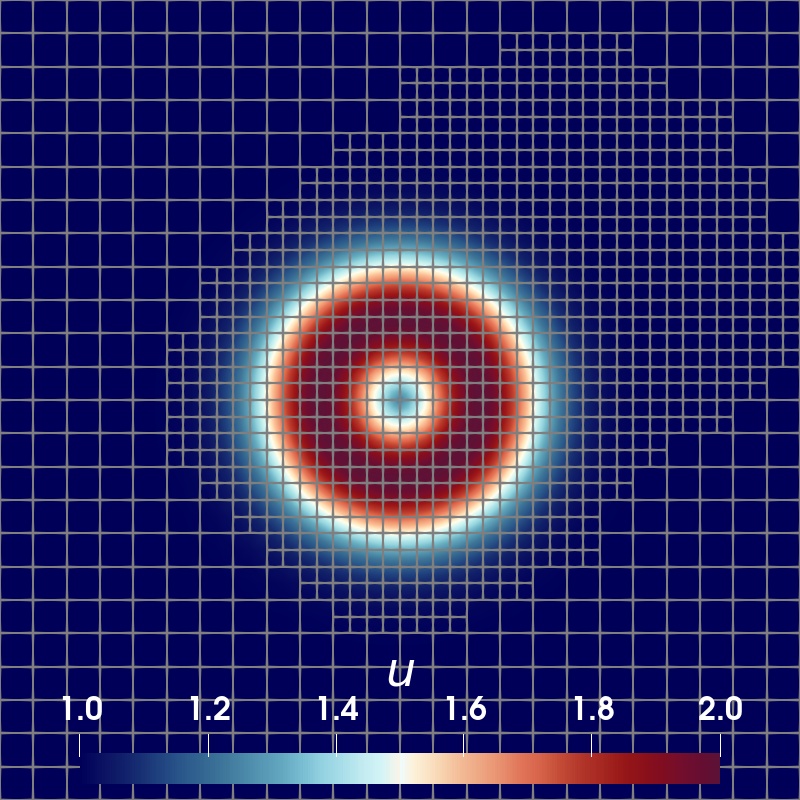}}}
        \subfloat[Solution at $t = T$]{
        \adjustbox{width=0.24\linewidth,valign=b}{\includegraphics{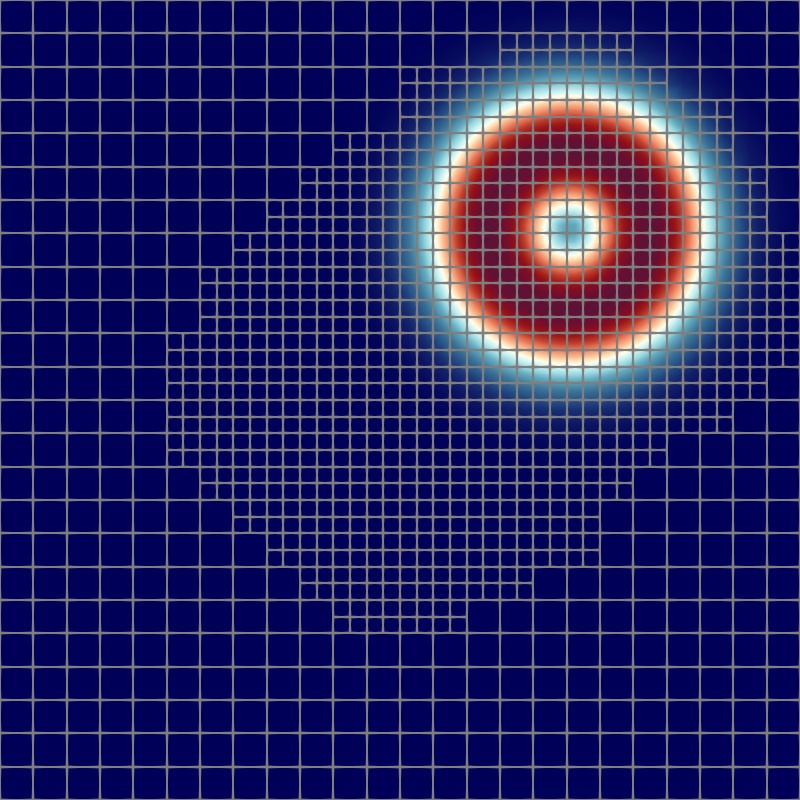}}}
        \subfloat[Remesh at $t = T$]{
        \adjustbox{width=0.24\linewidth,valign=b}{\includegraphics{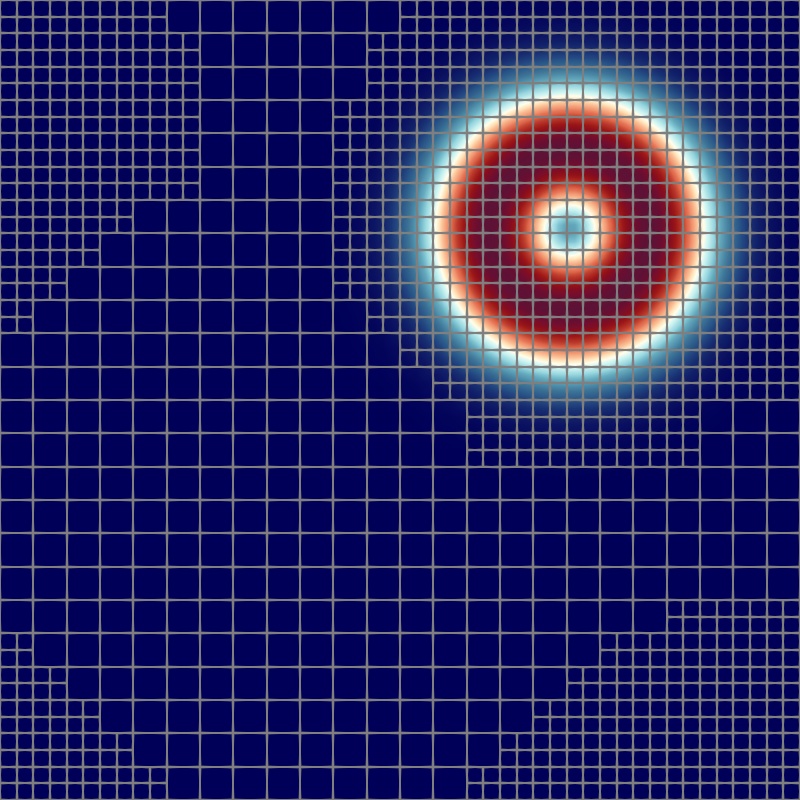}}}
        \subfloat[Solution at $t = 2T$]{
        \adjustbox{width=0.24\linewidth,valign=b}{\includegraphics{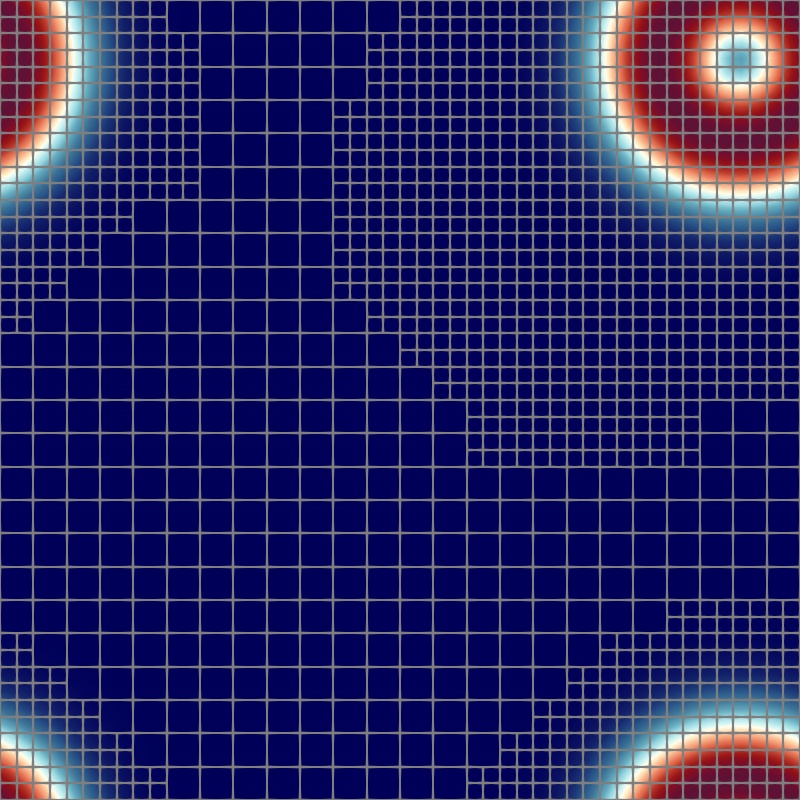}}}
        \newline
        \subfloat[Remesh at $t = 2T$]{
        \adjustbox{width=0.24\linewidth,valign=b}{\includegraphics{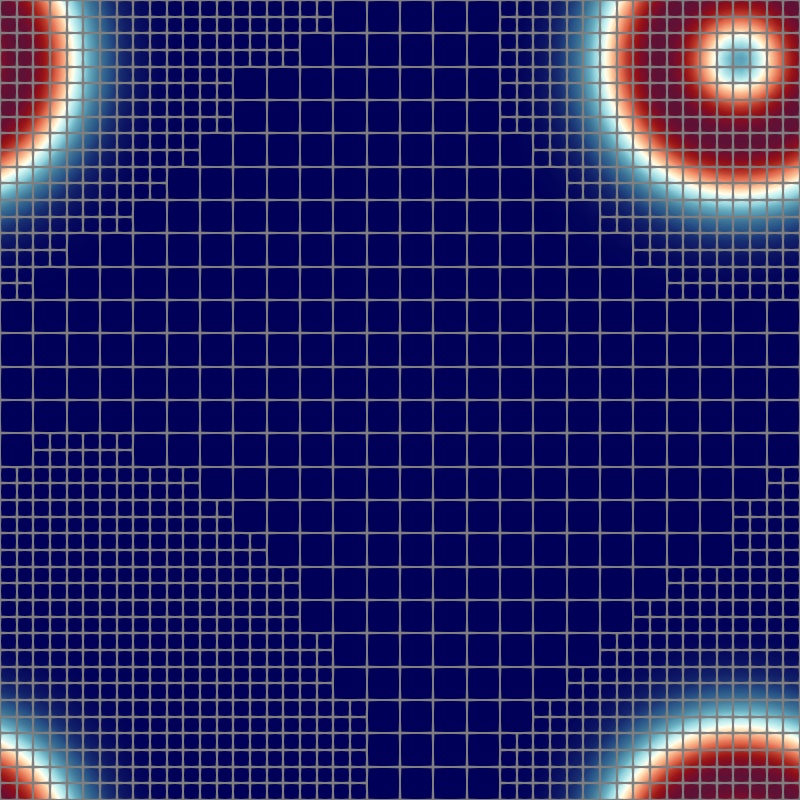}}}
        \subfloat[Solution at $t = 3T$]{
        \adjustbox{width=0.24\linewidth,valign=b}{\includegraphics{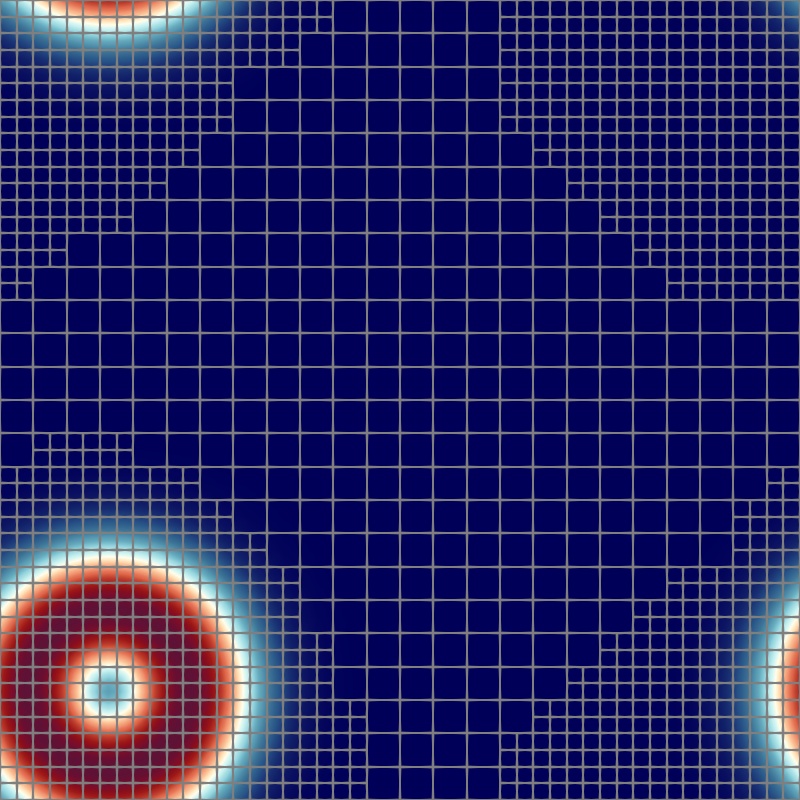}}}
        \subfloat[Remesh at $t = 3T$]{
        \adjustbox{width=0.24\linewidth,valign=b}{\includegraphics{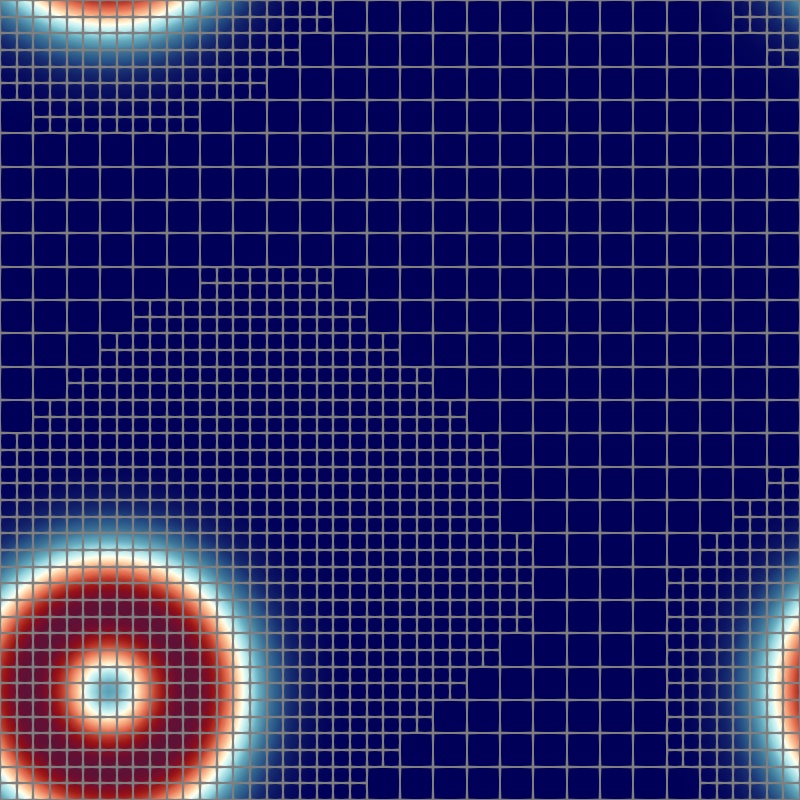}}}
        \subfloat[Solution at $t = 4T$]{
        \adjustbox{width=0.24\linewidth,valign=b}{\includegraphics{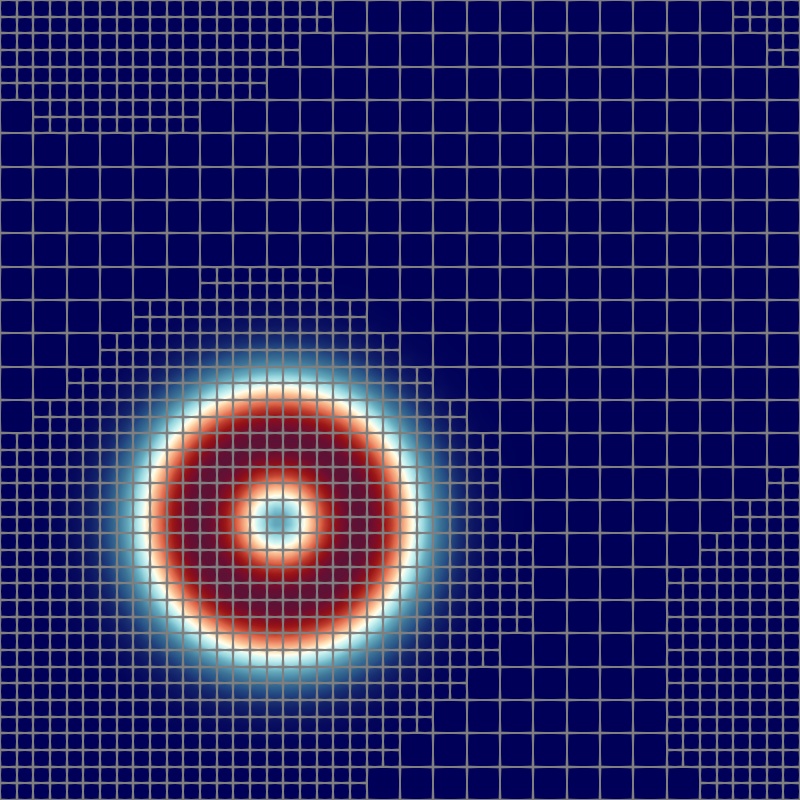}}}
        \newline
        \caption{\label{fig:advhex1_amr_drl} Solution contours overlaid with $h$-adapted mesh at varying remesh intervals using DynAMO for the \textit{in-distribution} advecting ring case.} 
    \end{figure}

\begin{figure}[htbp!]
        \centering
        \subfloat[Remesh at $t = 0$]{
        \adjustbox{width=0.24\linewidth,valign=b}{\includegraphics{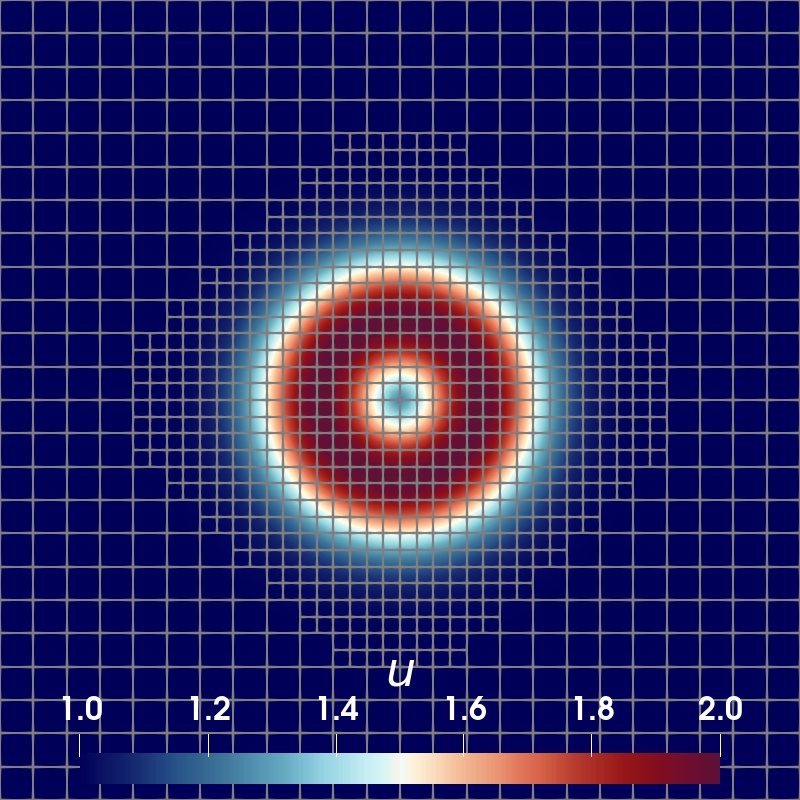}}}
        \subfloat[Solution at $t = T$]{
        \adjustbox{width=0.24\linewidth,valign=b}{\includegraphics{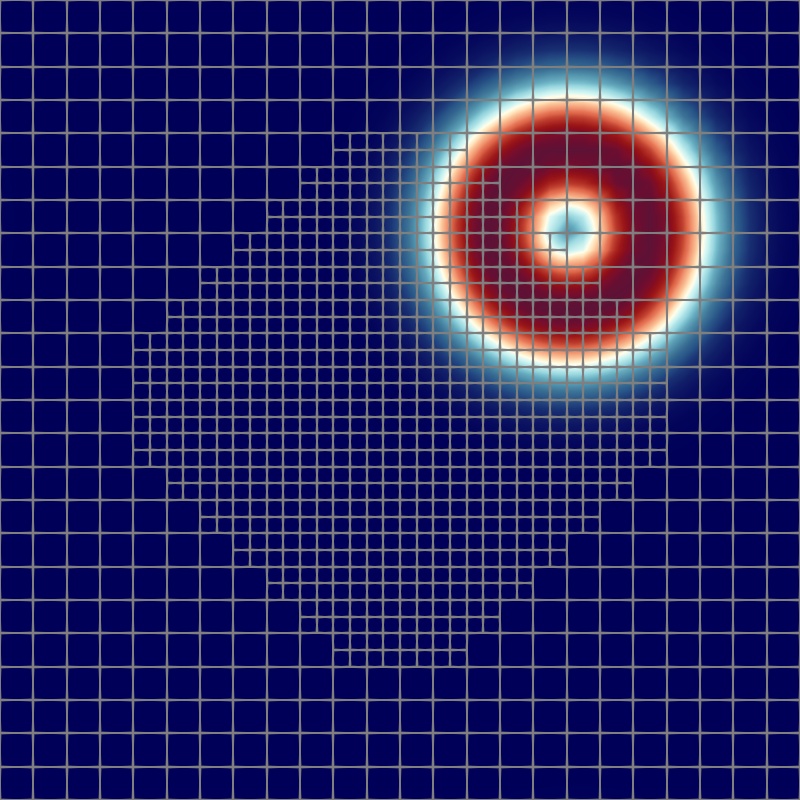}}}
        \subfloat[Remesh at $t = T$]{
        \adjustbox{width=0.24\linewidth,valign=b}{\includegraphics{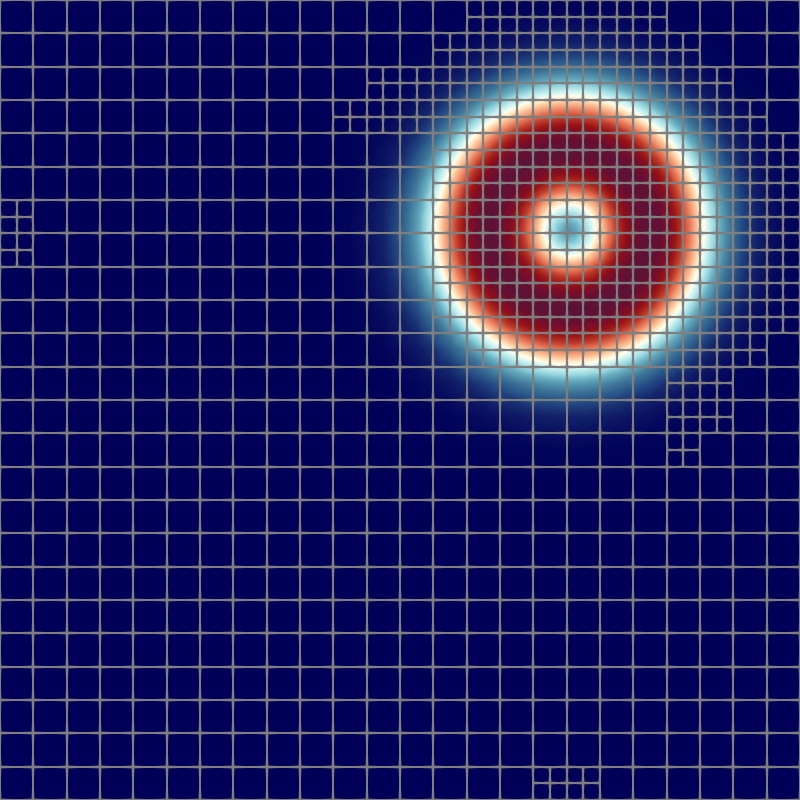}}}
        \subfloat[Solution at $t = 2T$]{
        \adjustbox{width=0.24\linewidth,valign=b}{\includegraphics{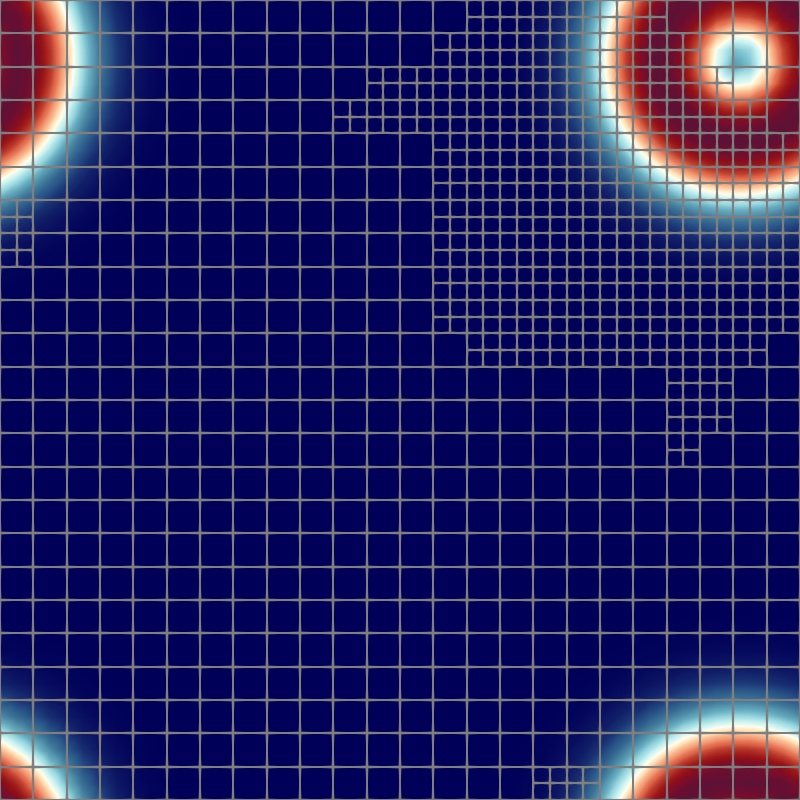}}}
        \newline
        \subfloat[Remesh at $t = 2T$]{
        \adjustbox{width=0.24\linewidth,valign=b}{\includegraphics{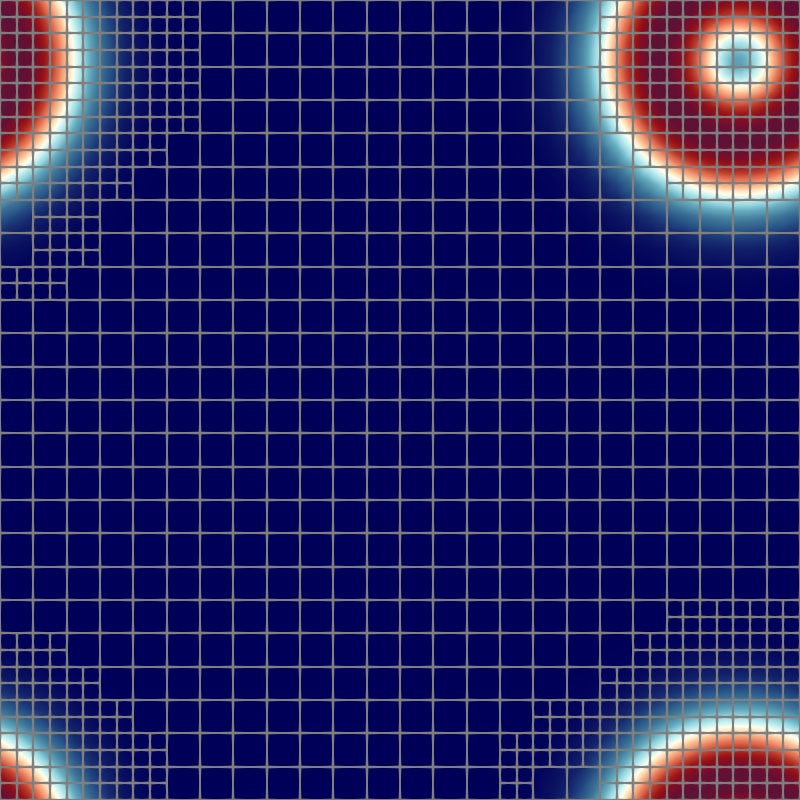}}}
        \subfloat[Solution at $t = 3T$]{
        \adjustbox{width=0.24\linewidth,valign=b}{\includegraphics{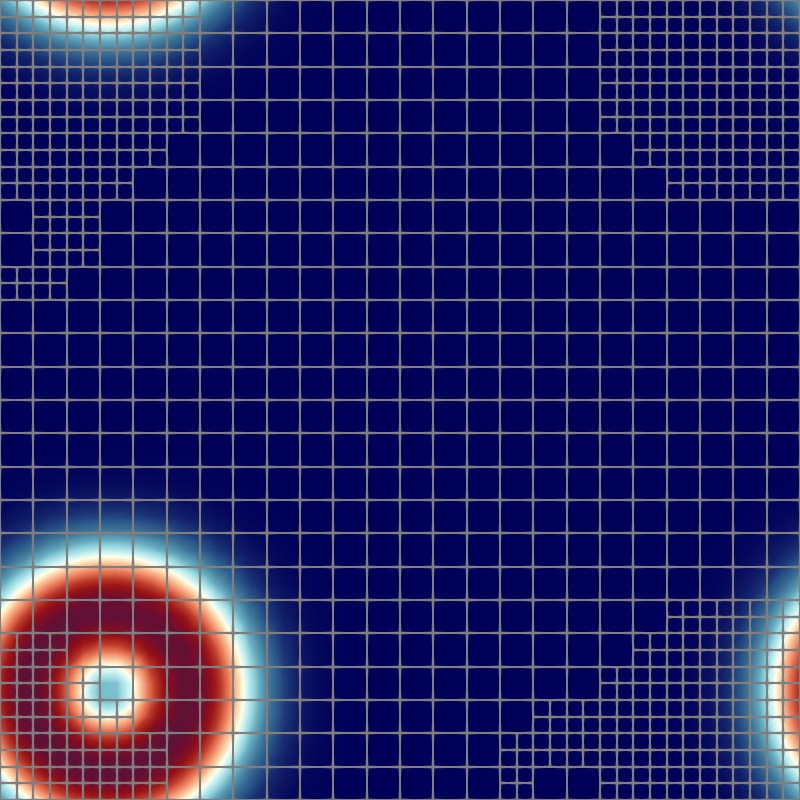}}}
        \subfloat[Remesh at $t = 3T$]{
        \adjustbox{width=0.24\linewidth,valign=b}{\includegraphics{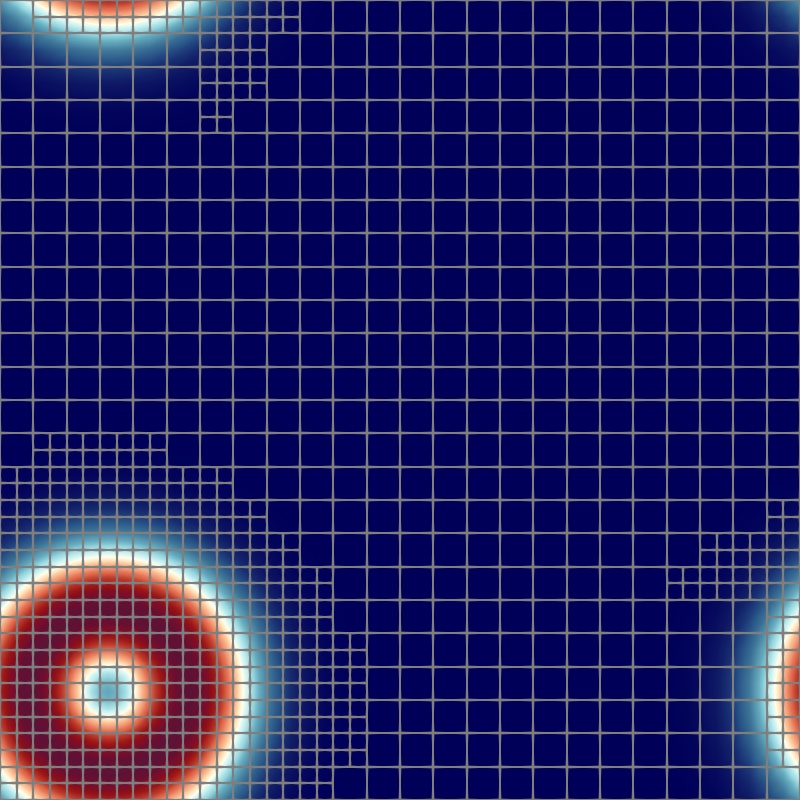}}}
        \subfloat[Solution at $t = 4T$]{
        \adjustbox{width=0.24\linewidth,valign=b}{\includegraphics{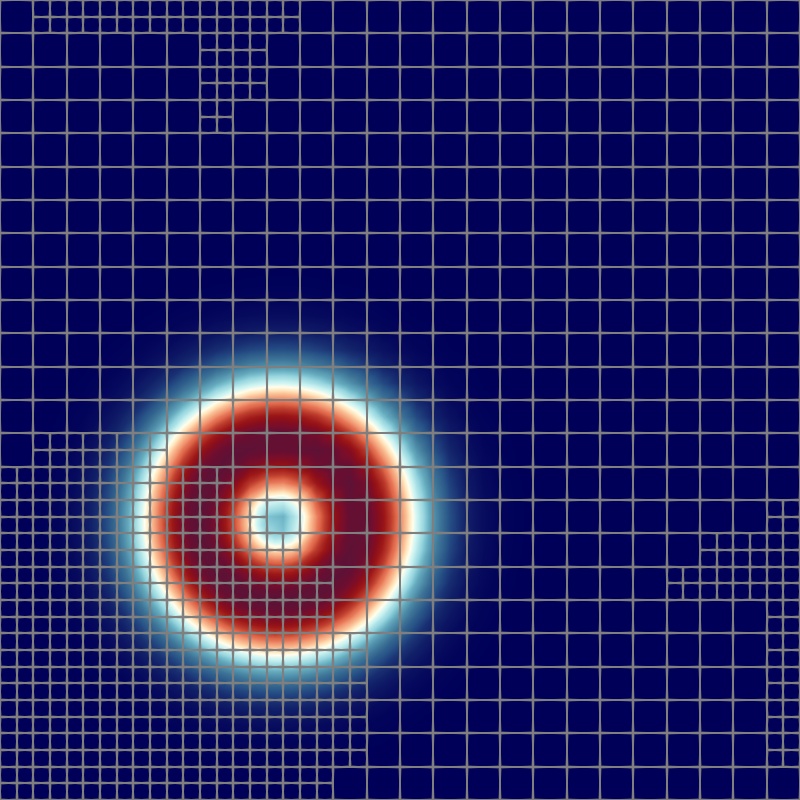}}}
        \newline
        \caption{\label{fig:advhex1_amr_threshold} Solution contours overlaid with $h$-adapted mesh at varying remesh intervals using the threshold policy ($\theta = 10^{-3}$) for the \textit{in-distribution} advecting ring case.} 
    \end{figure}

The training procedure was repeated identically to the $p$-refinement experiments and the resulting approach was compared to the threshold policy for the same advecting ring case, with the initial conditions given by \cref{eq:adv_rings}. The $h$-adapted mesh and solution evolution across 4 remesh times as computed by the DynAMO policy ($\alpha = \alpha_{\text{train}}$) and the threshold policy ($\theta = 10^{-3}$) is shown in \cref{fig:advhex1_amr_drl} and \cref{fig:advhex1_amr_threshold}, respectively. Much like with the $p$-refinement case, it can be seen that the threshold policy based on the instantaneous error estimator results in an isotropic refinement pattern centered on the ring, which does not adequately account for the movement of the ring between remesh intervals. This resulting introduction of discretization error cause irregular refinement patterns later in the simulation. In contrast, the $h$-refinement DynAMO policy showed similar behavior to the $p$-refinement DynAMO policy, with an anistropic refinement pattern stretched in the direction of propagation. This predictive refinement persisted across the entire simulation time, such that the feature of interest remained within the refined region throughout. Furthermore, it was observed that the DynAMO policy tended to produce contiguous regions of refined elements, which would be the proper refinement pattern for linearly advecting compactly-supported features.
    

    \begin{figure}[htbp!]
        \centering
        \adjustbox{width=0.45\linewidth,valign=b}{\begin{tikzpicture}[spy using outlines={rectangle, height=3cm,width=2.3cm, magnification=3, connect spies}]

    \begin{axis}
    [   axis lines = none,
        xmin = 0, xmax = 1,
        ymin = 0, ymax = 1
    ]
        \fill[fill=blue!5] (1.0, 0) arc[start angle=0, end angle=90, radius=1.0] -- (0,0) -- (1.0, 0);
        \fill[fill=blue!10] (0.9, 0) arc[start angle=0, end angle=90, radius=0.9] -- (0,0) -- (0.9, 0);
        \fill[fill=blue!15] (0.8, 0) arc[start angle=0, end angle=90, radius=0.8] -- (0,0) -- (0.8, 0);
        \fill[fill=blue!20] (0.7, 0) arc[start angle=0, end angle=90, radius=0.7] -- (0,0) -- (0.7, 0);
        \fill[fill=blue!25] (0.6, 0) arc[start angle=0, end angle=90, radius=0.6] -- (0,0) -- (0.6, 0);
        \fill[fill=blue!30] (0.5, 0) arc[start angle=0, end angle=90, radius=0.5] -- (0,0) -- (0.5, 0);
        \fill[fill=blue!35] (0.4, 0) arc[start angle=0, end angle=90, radius=0.4] -- (0,0) -- (0.4, 0);
        \fill[fill=blue!40] (0.3, 0) arc[start angle=0, end angle=90, radius=0.3] -- (0,0) -- (0.3, 0);
        \fill[fill=blue!45] (0.2, 0) arc[start angle=0, end angle=90, radius=0.2] -- (0,0) -- (0.2, 0);
        \fill[fill=blue!50] (0.1, 0) arc[start angle=0, end angle=90, radius=0.1] -- (0,0) -- (0.1, 0);

    \end{axis}
        
    \begin{axis}
    [   axis line style={latex-latex},
        axis y line=left,
        axis x line=left,
        xmode=linear,
        ymode=linear,
        xlabel = {$\bar{c}$},
        ylabel = {$\bar{e}$},
        xmin = 0, xmax = 1,
        ymin = 0, ymax = 1,
        xtick = {0,0.2,0.4,0.6,0.8,1.0},
        ytick = {0,0.2,0.4,0.6,0.8,1.0},
        minor x tick num=1,
        minor y tick num=1,
        legend cell align={left},
        legend style={at={(0.97, 0.97)},anchor=north east},
        clip mode=individual,
        x tick label style={/pgf/number format/.cd, fixed, fixed zerofill, precision=1, /tikz/.cd},
        y tick label style={/pgf/number format/.cd, fixed, fixed zerofill, precision=1, /tikz/.cd},
        label style={font=\large},
    ]
        \addplot[color=black, fill=white, style={thick}, only marks, mark=*, mark options={scale=1.2}]
        table[x = cost_mean, y = err_mean, col sep=comma]{./figs/data/advh_ex1_id_threshold.csv};
        \addlegendentry{Threshold policy};
        
        \addplot[color=black, fill=red, style={thick}, only marks, mark=*, mark options={scale=1.2}]
        table[x = cost, y = err, col sep=comma]{./figs/data/advh_ex1_id_drl_alphasweep.csv};
        \addlegendentry{DynAMO};

        \addplot[color=black, fill=green, style={thick}, only marks, mark=*, mark options={scale=1.2}]
        coordinates {(0.33463,0.00038329)} ;
        \addlegendentry{DynAMO ($\alpha_{{\tt train}}$)};

        \draw[->, black, thick] (0.25,.75) arc (180:255:0.6);
        \draw (0.8, 0.1) node[scale=1] {Decreasing $\theta$, \textcolor{red!90!black}{$\alpha$}} ;
        
    \end{axis}
\end{tikzpicture}}
        \caption{\label{fig:adv_href_pareto} Pareto plot of normalized cost vs. error for $p$-refinement on the advection equations with DynAMO and the threshold policy for the \textit{in-distribution} advecting ring case. Contours of efficiency shown on background. Red markers represent a sweep of $\alpha \in [0.1, 2.0]$ at evaluation time, green marker represents $\alpha$ set to the training value of $\alpha_{\texttt{train}} = 0.1$. Peak efficiency achieved at $\alpha=0.8$. }
    \end{figure}
    
A quantitative comparison of the error, cost, and efficiency between the DynAMO policy and the threshold policy for the advecting ring case is shown in the Pareto plot in \cref{fig:adv_href_pareto}. By varying the error threshold parameter, both the threshold policy and the DynAMO policy showed the typical Pareto front. However, much like with the $p$-refinement case, the DynAMO policy achieved noticeably higher efficiency than the threshold policy. For a given computational cost, the DynAMO policy was able to achieve a significantly lower error. Furthermore, by varying the error threshold parameter $\alpha$ at evaluation time, various points on that Pareto front could be achieved, such that error and cost targets could be specified by the user. We remark here that again the DynAMO policy did not necessarily achieve peak efficiency for the error threshold parameter used at training time, which resulted in refinement decisions that were biased more towards error reduction than cost reduction. Nevertheless, the efficiency of the DynAMO policy evaluated over a wide range of values of the error threshold parameter was noticeably higher than the threshold policy at its peak efficiency. These results, which show a marked consistency with the results of the $p$-refinement policy, indicate that the proposed multi-agent reinforcement learning approach for $h$-refinement using coarse level agents is highly effective and adequately sidesteps the agent creation/deletion problem without sacrificing performance.

\subsubsection{In-distribution experiments}
To verify that these observations extend over many sampled initial conditions, the mean efficiency, normalized error, and normalized cost over 100 randomly sampled \textit{in-distribution} conditions was calculated, the results of which are presented in \cref{tab:advection_href_indistribution} for the threshold policy at varying values of the error threshold $\theta$ and the DynAMO policy at the training error threshold value $\alpha = \alpha_{\text{train}}$. It can be seen that the DynAMO approach still retains efficiency benefits over the baseline threshold policy, with nearly a $20\%$ increase in the mean efficiency. Furthermore, the error-to-cost ratio of the DynAMO approach was significantly lower, showing a nearly $99\%$ decrease in the error for approximately $30\%$ higher computational cost. This accuracy was on par with the accuracy of a fully-refined mesh, which would require nearly double the computational cost. We remark here that the efficiency increase of DynAMO in the $h$-refinement case was lower than the $p$-refinement case, but this was primarily due to the choice of error threshold $\alpha$. At the training value $\alpha = \alpha_{\text{train}}$, the results of the DynAMO policy for $h$-refinement were farther along the Pareto front, which is evidenced by \cref{fig:adv_href_pareto} and the relative error in \cref{tab:advection_href_indistribution}. Much like with the $p$-refinement policy, even higher efficiencies could be obtained by modifying the value of $\alpha$ as the taken training value biases towards higher error reduction instead of peak efficiency. 

    \begin{figure}[htbp!] 
        \centering
        \begin{tabular}{cccc}
        \toprule
        Method &  Efficiency & Normalized error & Normalized cost \\ 
        \midrule
        Threshold ($\theta = 10^{-2}$) & 0.174 (0.135) & 0.815 (0.157) & 0.084 (0.056) \\
        Threshold ($\theta = 10^{-3}$) & \textbf{0.407 (0.193)} & \textbf{0.417 (0.212)} & \textbf{0.401 (0.098)} \\
        Threshold ($\theta = 10^{-4}$) & 0.299 (0.159) & 0.121 (0.166) & 0.675 (0.140) \\
        Threshold ($\theta = 10^{-5}$) & 0.191 (0.147) & 0.056 (0.144) & 0.797 (0.128) \\
        Threshold ($\theta = 10^{-6}$) & 0.121 (0.116) & 0.029 (0.126) & 0.872 (0.098) \\
        Threshold ($\theta = 10^{-7}$) & 0.079 (0.072) & 0.012 (0.048) & 0.920 (0.071) \\
        \midrule
        DynAMO & \textbf{0.480 (0.106)} & \textbf{0.006 (0.044)} & \textbf{0.519 (0.108)} \\
        DynAMO/Optimal $\theta$ & \textcolor{green!70!black}{+17.9\%} & \textcolor{green!70!black}{-98.6\%} & \textcolor{red!70!black}{+29.4\%} \\
        \bottomrule
        \end{tabular}
        \captionof{table}{\label{tab:advection_href_indistribution} Comparison of the mean efficiency, normalized error, and normalized cost for $h$-refinement on the advection equation with DynAMO and the threshold policy for the advecting rings over 100 \textit{in-distribution} runs using uniform random initial conditions. Standard deviation shown in parentheses.}
    \end{figure}
    
\subsubsection{Generalization experiments}
    \begin{figure}[htbp!] 
        \centering
        \begin{tabular}{ccccc}
        \toprule
        Method &  In-distribution & Finer mesh & Different shapes & Longer sim. time \\ 
        \midrule
        Threshold ($\theta = 10^{-2}$) & 0.174 (0.135) & 0.092 (0.153) & 0.000 (0.002) & 0.098 (0.137)\\
        Threshold ($\theta = 10^{-3}$) & \textbf{0.407 (0.193)} & 0.283 (0.199) & 0.040 (0.145) & \textbf{0.399 (0.197)}\\
        Threshold ($\theta = 10^{-4}$) & 0.299 (0.159) & \textbf{0.483 (0.210)} & \textbf{0.484 (0.163)} & 0.328 (0.142)\\
        Threshold ($\theta = 10^{-5}$) & 0.191 (0.147) & 0.475 (0.154) & 0.479 (0.140) & 0.206 (0.132)\\
        Threshold ($\theta = 10^{-6}$) & 0.121 (0.116) & 0.291 (0.160) & 0.351 (0.143) & 0.128 (0.099)\\
        Threshold ($\theta = 10^{-7}$) & 0.079 (0.072) & 0.196 (0.145) & 0.241 (0.128) & 0.076 (0.068)\\
        \midrule
        DynAMO & \textbf{0.480 (0.106)}  & \textbf{0.807 (0.129)} & \textbf{0.540 (0.110)} & \textbf{0.719 (0.081)}\\
        DynAMO/Optimal $\theta$ & \textcolor{green!70!black}{+17.9\%} & \textcolor{green!70!black}{+67.1\%}& \textcolor{green!70!black}{+11.5\%}& \textcolor{green!70!black}{+80.2\%} \\
        \bottomrule
        \end{tabular}
        \captionof{table}{\label{tab:advection_href_ood} Comparison of the mean efficiency for $h$-refinement on the advection equation with DynAMO and the threshold policy for the advecting rings over 100 \textit{out-of-distribution} runs using uniform random initial conditions with finer mesh resolution, different advecting shapes, and longer simulation time. Standard deviation shown in parentheses. In-distribution results from \cref{tab:advection_href_indistribution} shown for comparison.}
    \end{figure}

Finally, the generalization experiments were repeated for \textit{out-of-distribution} problems, such as finer meshes, different classes of initial conditions, and longer simulation times. A comparison of the threshold policy at varying values of the error threshold $\theta$ and the DynAMO policy at the training error threshold value $\alpha = \alpha_{\text{train}}$ performed over 100 samples is shown in \cref{tab:advection_href_ood} for the three experiments. When the initial mesh resolution was increased from $N = 24^2$ to $N = 96^2$, the DynAMO approach showed a notably higher efficiency ($67.1\%$) than the threshold policy. This efficiency gain was significantly higher than with the similar experiment on $p$-refinement, but this was primarily as a result of the threshold policy showing poorer performance on the fine mesh with $h$-refinement than with $p$-refinement. For the generalization to different initial condition shapes, the DynAMO approach retained the efficiency benefits, resulting in a similar increase in efficiency between DynAMO and the optimal threshold policy. Extensions to longer simulation times showed the most drastic results, with approximately $80\%$ higher efficiency with DynAMO compared to the optimal threshold policy due to superior error mitigation from preemptive refinement. These results are notably similar to the results of the $p$-refinement experiments, indicating that the proposed method for coarse-level observations in $h$-refinement does not degrade the accuracy and generalizability of the DynAMO approach. 

\subsection{Euler equations with $p$-refinement}
\begin{figure}[htbp!]
        \centering
        \adjustbox{width=0.45\linewidth,valign=b}{    \begin{tikzpicture}[spy using outlines={rectangle, height=3cm,width=2.3cm, magnification=3, connect spies}]
		\begin{axis}[name=plot1,
		    axis x line=left,
            axis y line=left,
    		xlabel={Env. steps},
		    ylabel={$-\bar{r}$},
    		xmin=0,
    		ymin=0,
    		ymax=10,
    		style={font=\normalsize}]
    		
			\addplot[color=black, style={thick}]
				table[x=num_env_steps_sampled,
                      y expr=-\thisrow{episode_reward_mean},
                      col sep=comma,unbounded coords=jump]{./figs/data/euler_p_training.csv};
    		    		
		\end{axis} 		
	\end{tikzpicture}}
        \caption{\label{fig:euler_pref_reward} Batch-averaged (negative) reward with respect to number of environment steps during the training process for $p$-refinement on the Euler equations.}
    \end{figure}
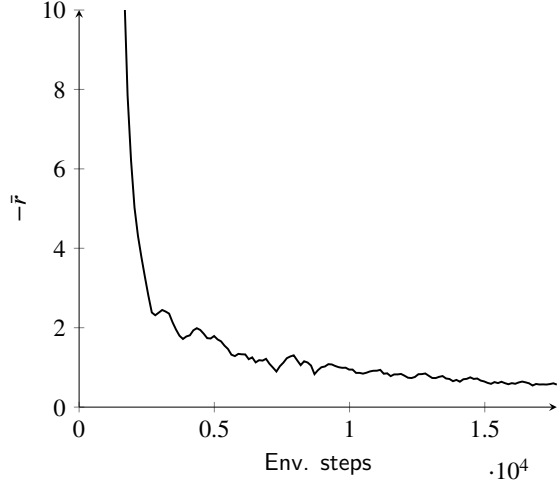

With the preliminary experiments for the linear advection equation showing promising results for both $p$- and $h$-refinement, the proposed DynAMO approach was extended to the much more complex task of AMR for the compressible Euler equations. Due to the inefficacy of $p$-refinement for discontinuous features, the $p$-refinement policy and experiments for the Euler equations focus strictly on smooth solutions while discontinuous solutions are later considered in the $h$-refinement experiments. For this class of solutions, common representative cases are convection-dominated and acoustic-type problems, where nonlinear flow physics are introduced through bulk propagation and pressure-wave interactions. To this end, the DynAMO policy for $p$-refinement on the Euler equations was trained on a class of solutions containing multiple convecting pressure pulses/waves, where the initial conditions were
\begin{subequations}
\begin{align}\label{eq:pressure_pulse}
    \rho (x, y, 0) &= 1, \\
    u(x, y, 0) &= u_0, \\
    v(x, y, 0) &= v_0, \\
    P(x, y, 0) &= \sum_{i = 1}^{n_p} h_i \exp \left(-w_i \left((x - x_{i,0})^2 + (y - y_{i,0})^2 \right) \right),
\end{align}
\end{subequations}
where the parameters $u_0$, $v_0$, $h$, $w$, $x_0$, and $y_0$ were uniformly distributed across the ranges $u_0 \in [0.0, 3.0]$, $v_0 \in [0.0, 3.0]$, $n_p \in \lbrace 1,2,3 \rbrace$, $h \in [0.05, 0.2]$, $w \in [200, 700]$, $x_0 \in [0.25, 1.25]$, and $y_0 \in [0.25, 1.25]$. The number of pressure pulses $n_p$ was chosen randomly, and the characteristic parameters of each pulse were independently sampled from their respective distributions. The domain size was set to $\Omega = [0, 1.5]^2$, the remesh time was set as $T = 0.05$ with a total number of 4 RL steps, and the initial mesh distribution was set as $N = 48^2$. The training reward curve over the training period is shown in \cref{fig:euler_pref_reward}.

\subsubsection{Sample problem}
\begin{figure}[htbp!]
        \centering
        \subfloat[Remesh at $t = 0$]{
        \adjustbox{width=0.16\linewidth,valign=b}{\includegraphics{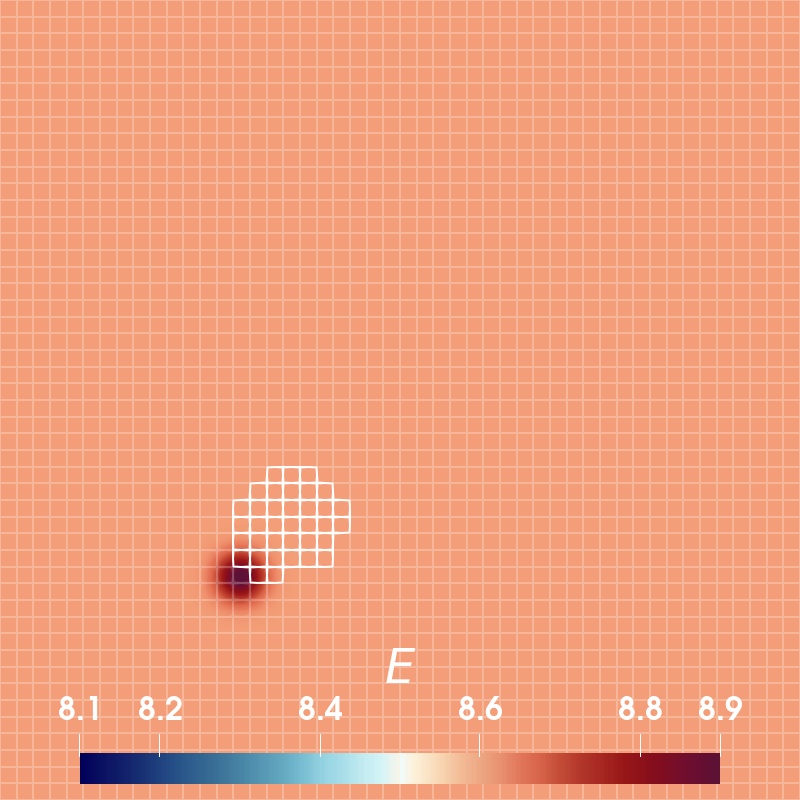}}}
        \subfloat[Solution at $t = T$]{
        \adjustbox{width=0.16\linewidth,valign=b}{\includegraphics{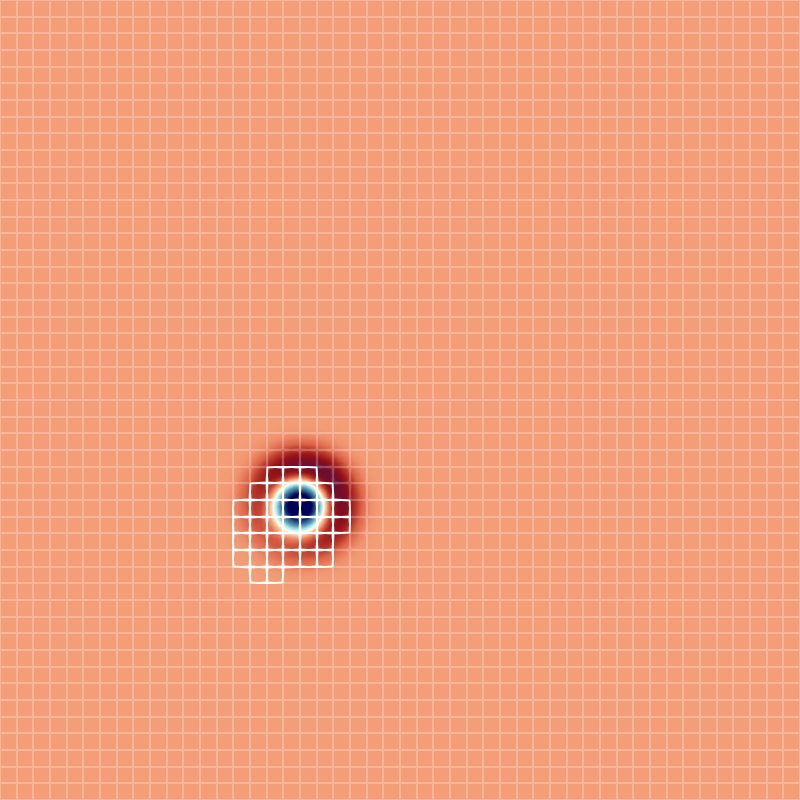}}}
        \subfloat[Remesh at $t = T$]{
        \adjustbox{width=0.16\linewidth,valign=b}{\includegraphics{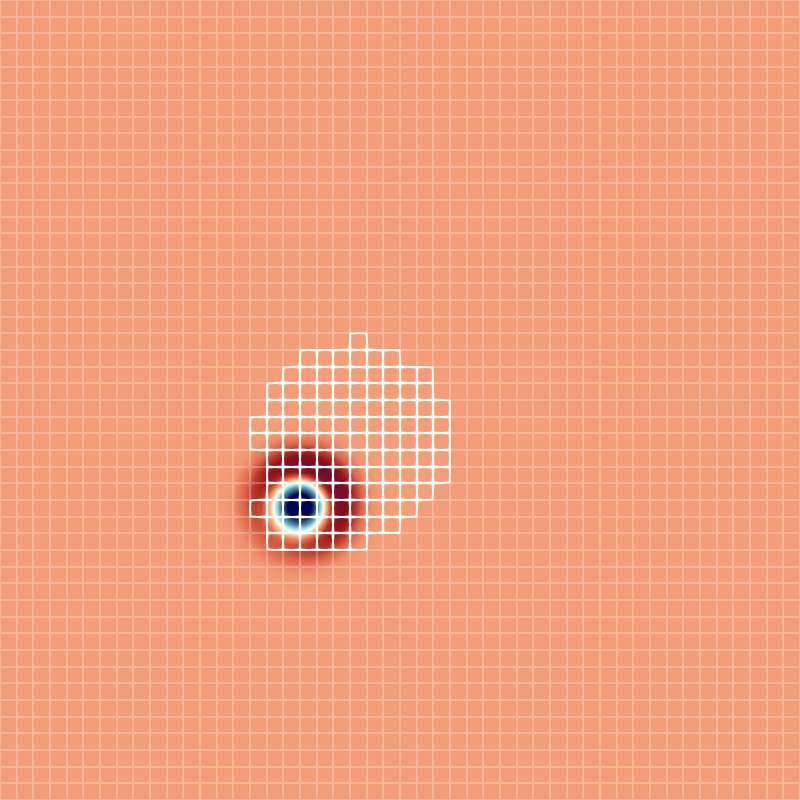}}}
        \subfloat[Solution at $t = 2T$]{
        \adjustbox{width=0.16\linewidth,valign=b}{\includegraphics{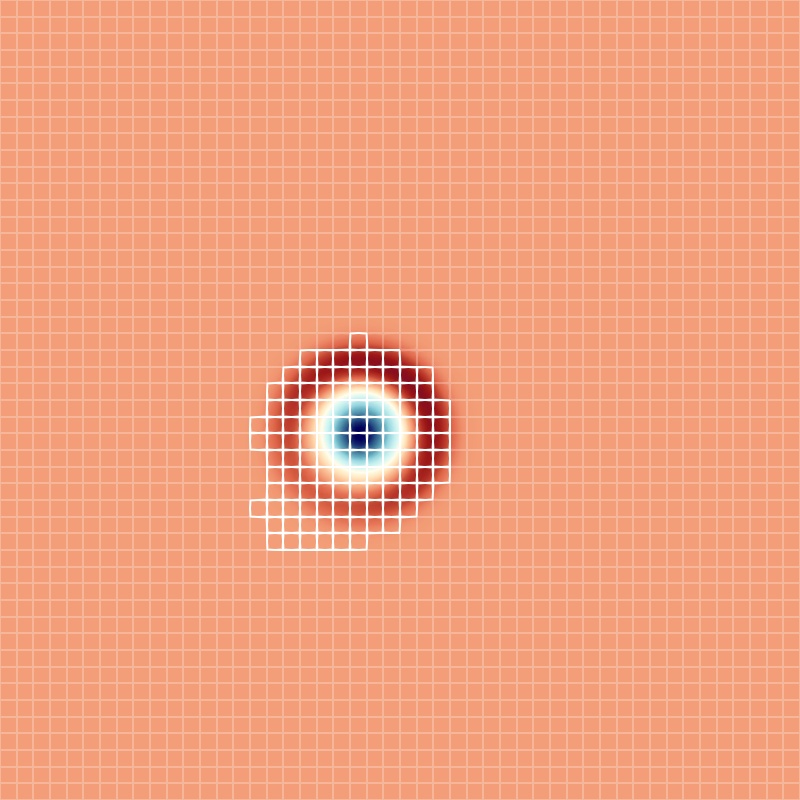}}}
        \subfloat[Remesh at $t = 2T$]{
        \adjustbox{width=0.16\linewidth,valign=b}{\includegraphics{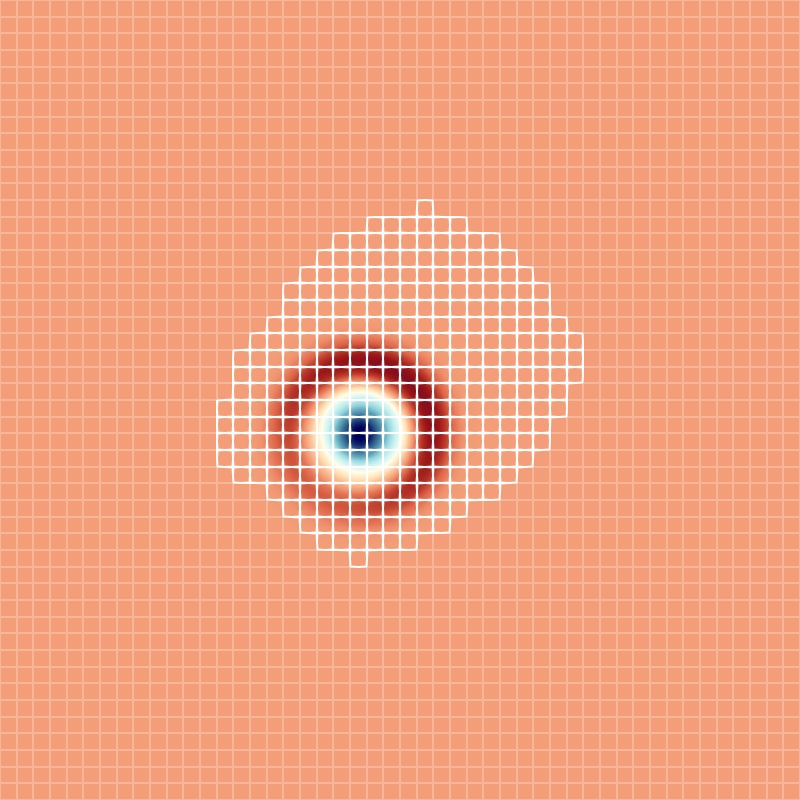}}}
        \subfloat[Solution at $t = 3T$]{
        \adjustbox{width=0.16\linewidth,valign=b}{\includegraphics{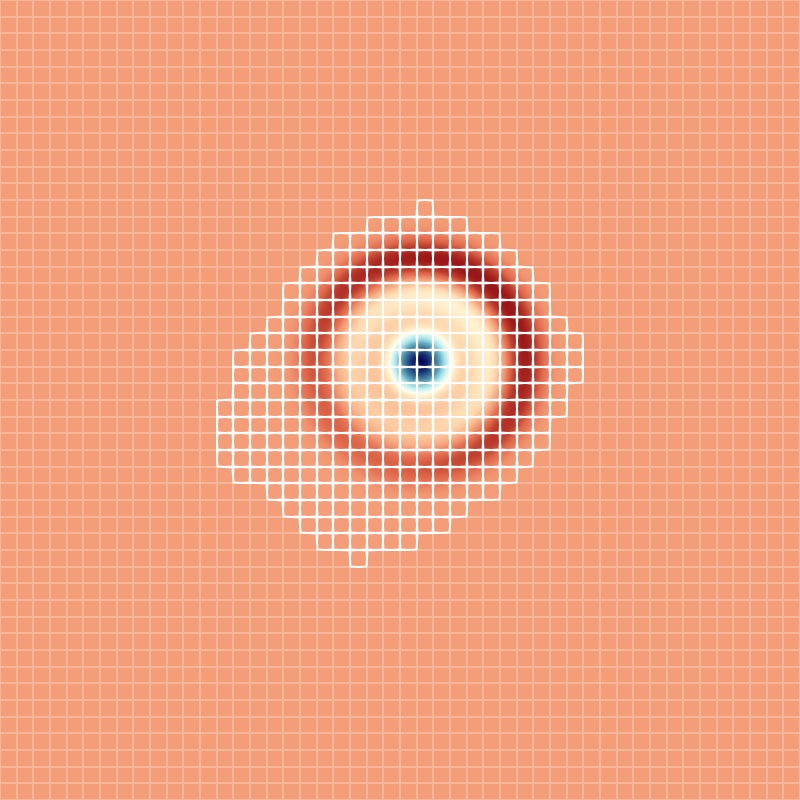}}}
        \newline
        \subfloat[Remesh, $3T$]{
        \adjustbox{width=0.16\linewidth,valign=b}{\includegraphics{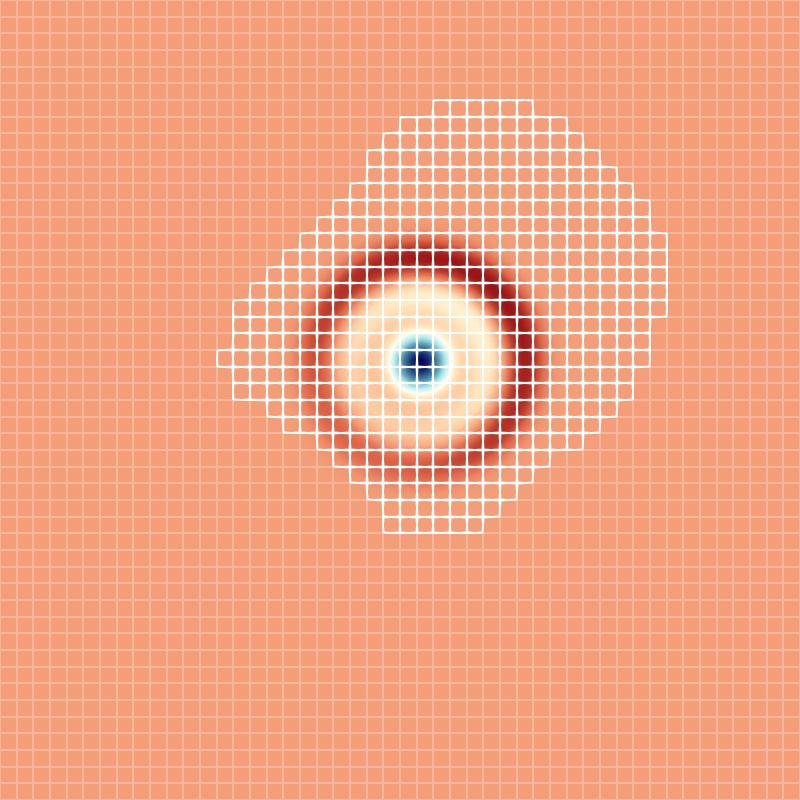}}}
        \subfloat[Solution, $4T$]{
        \adjustbox{width=0.16\linewidth,valign=b}{\includegraphics{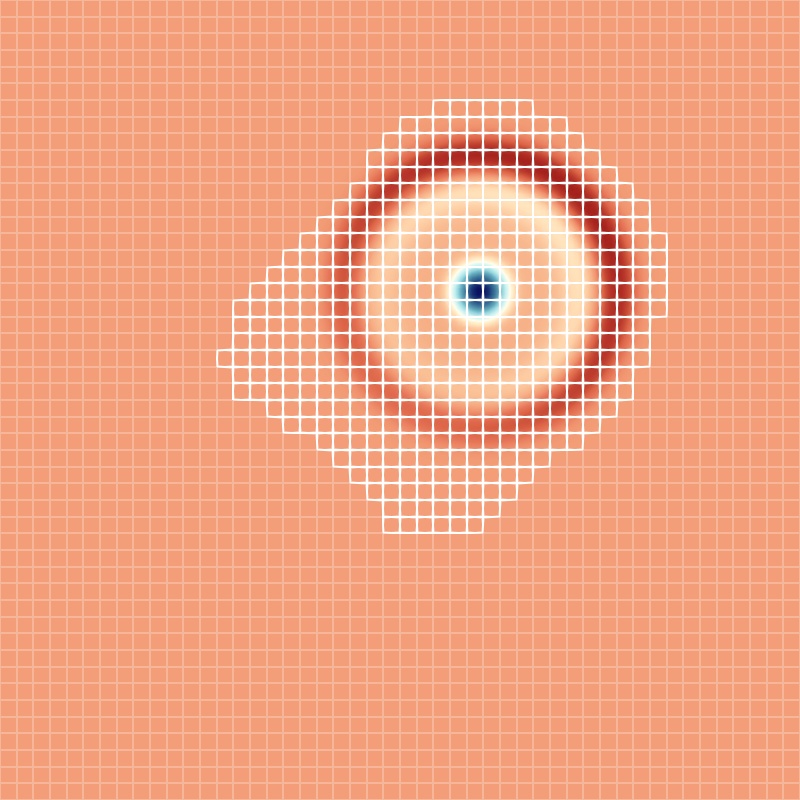}}}
        \subfloat[Remesh, $4T$]{
        \adjustbox{width=0.16\linewidth,valign=b}{\includegraphics{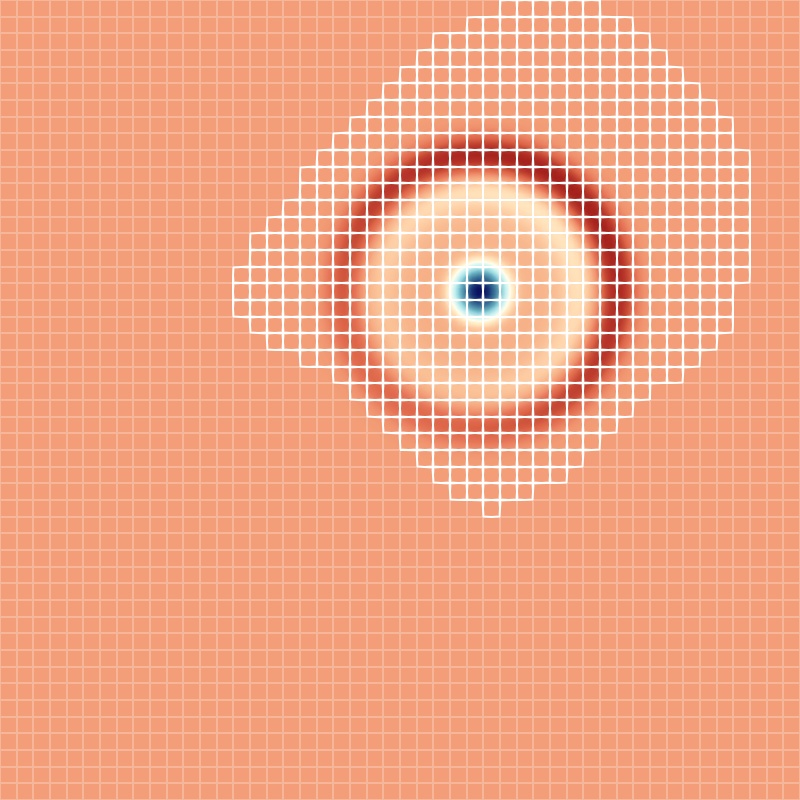}}}
        \subfloat[Solution, $5T$]{
        \adjustbox{width=0.16\linewidth,valign=b}{\includegraphics{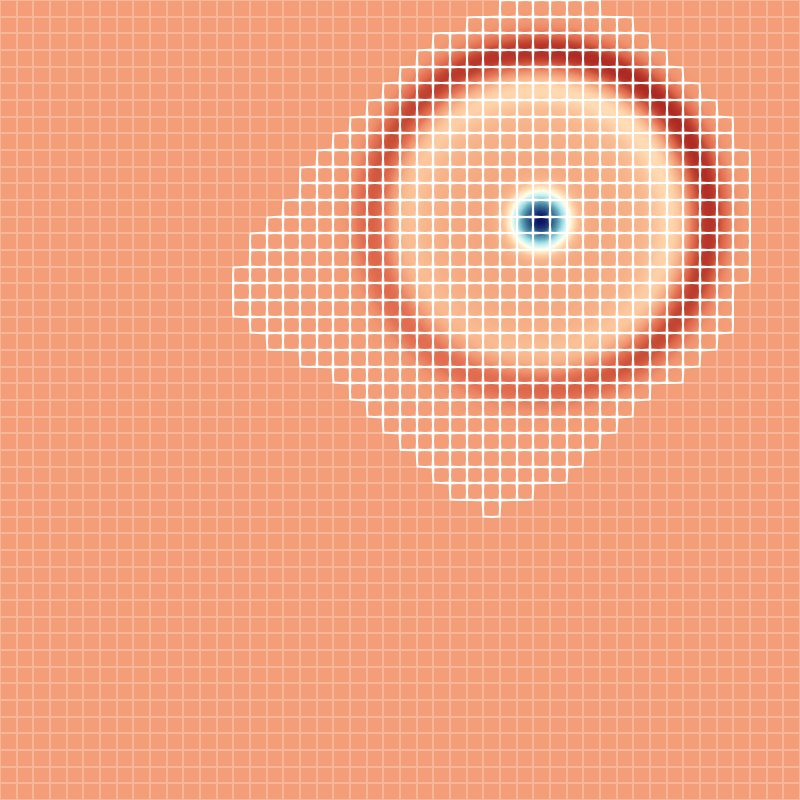}}}
        \subfloat[Remesh, $5T$]{
        \adjustbox{width=0.16\linewidth,valign=b}{\includegraphics{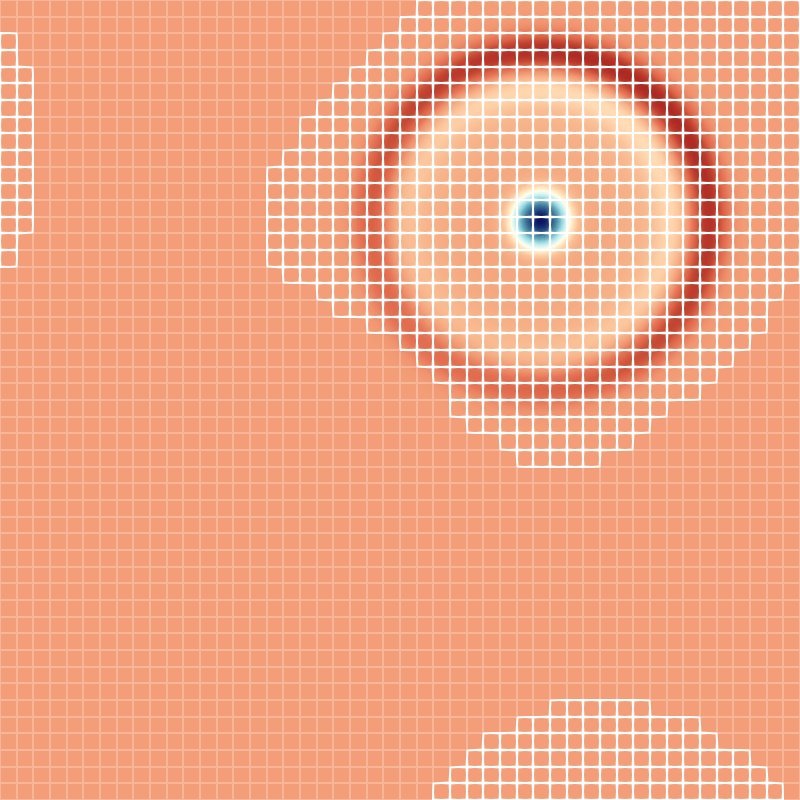}}}
        \subfloat[Solution, $6T$]{
        \adjustbox{width=0.16\linewidth,valign=b}{\includegraphics{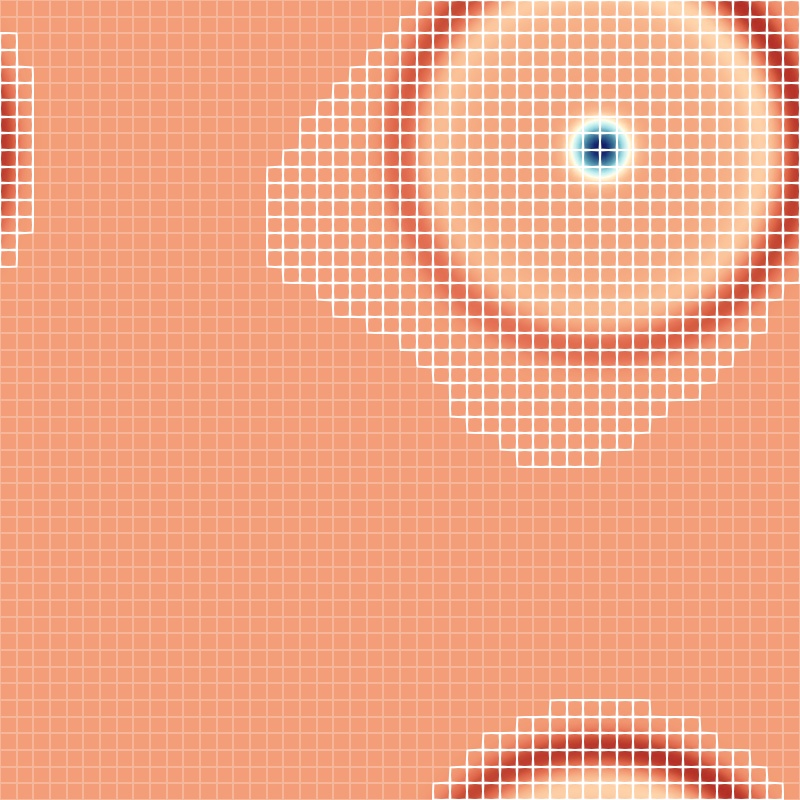}}}
        \newline
        \subfloat[Remesh, $6T$]{
        \adjustbox{width=0.16\linewidth,valign=b}{\includegraphics{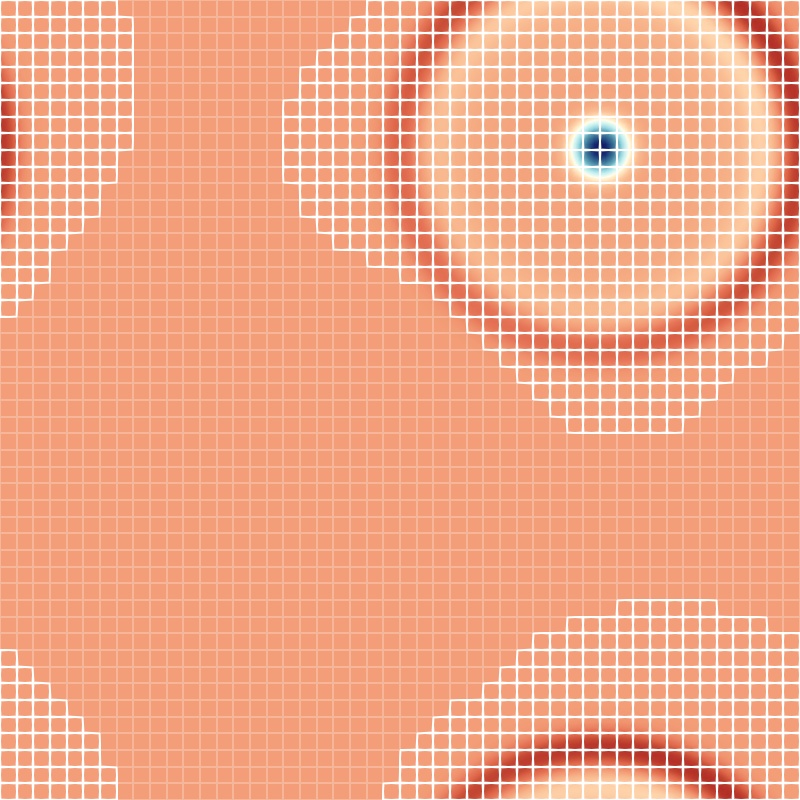}}}
        \subfloat[Solution, $7T$]{
        \adjustbox{width=0.16\linewidth,valign=b}{\includegraphics{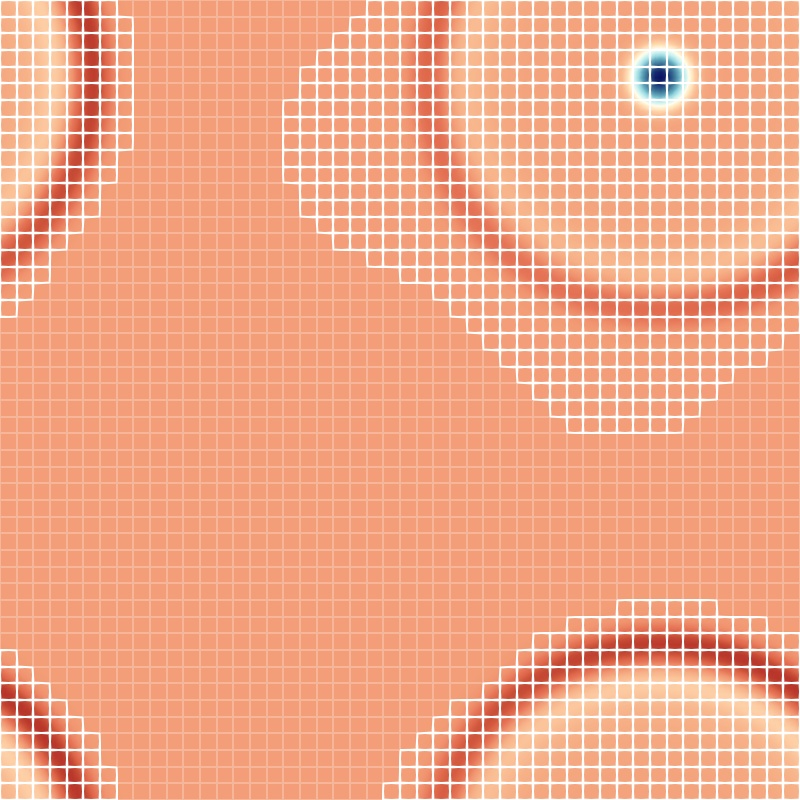}}}
        \subfloat[Remesh, $7T$]{
        \adjustbox{width=0.16\linewidth,valign=b}{\includegraphics{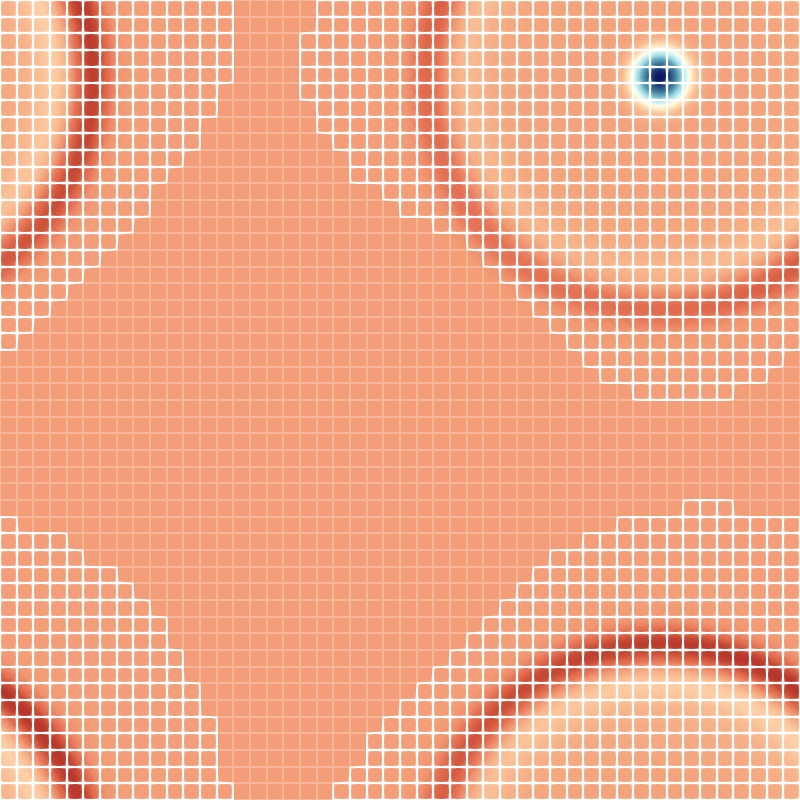}}}
        \subfloat[Solution, $8T$]{
        \adjustbox{width=0.16\linewidth,valign=b}{\includegraphics{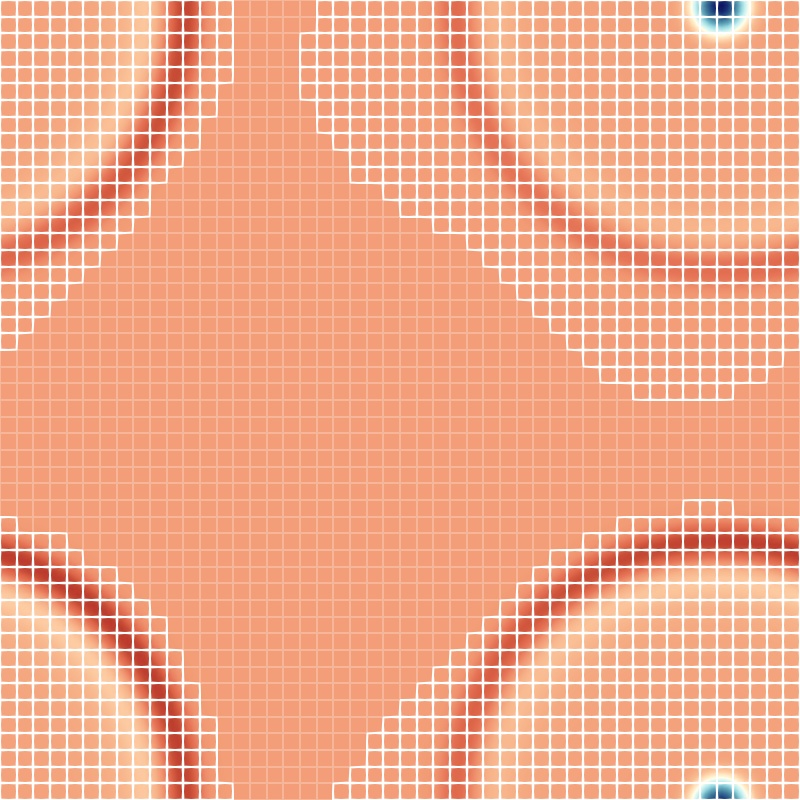}}}
        \subfloat[Remesh, $8T$]{
        \adjustbox{width=0.16\linewidth,valign=b}{\includegraphics{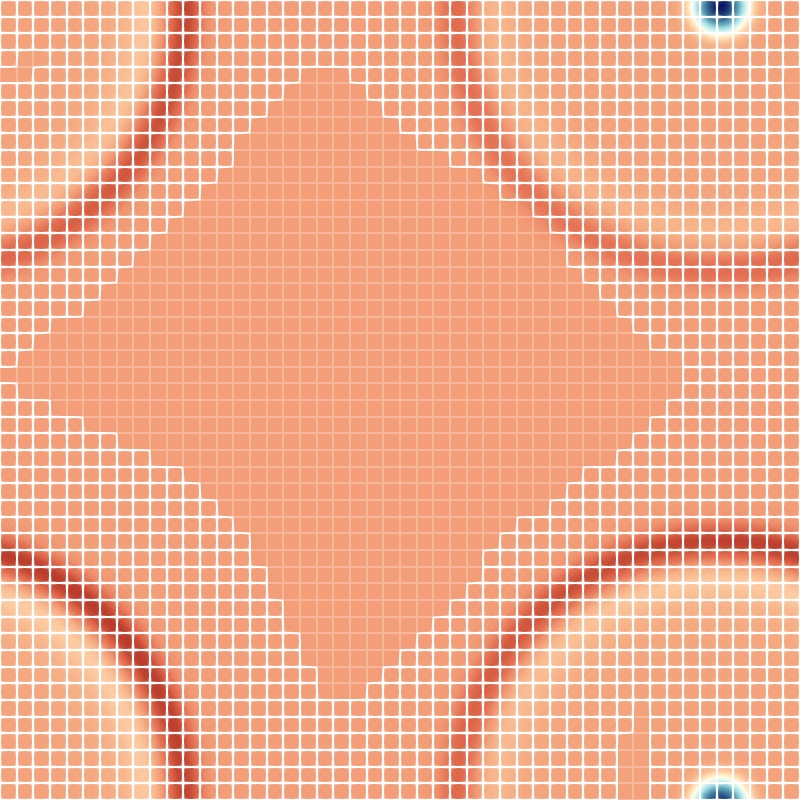}}}
        \subfloat[Solution, $9T$]{
        \adjustbox{width=0.16\linewidth,valign=b}{\includegraphics{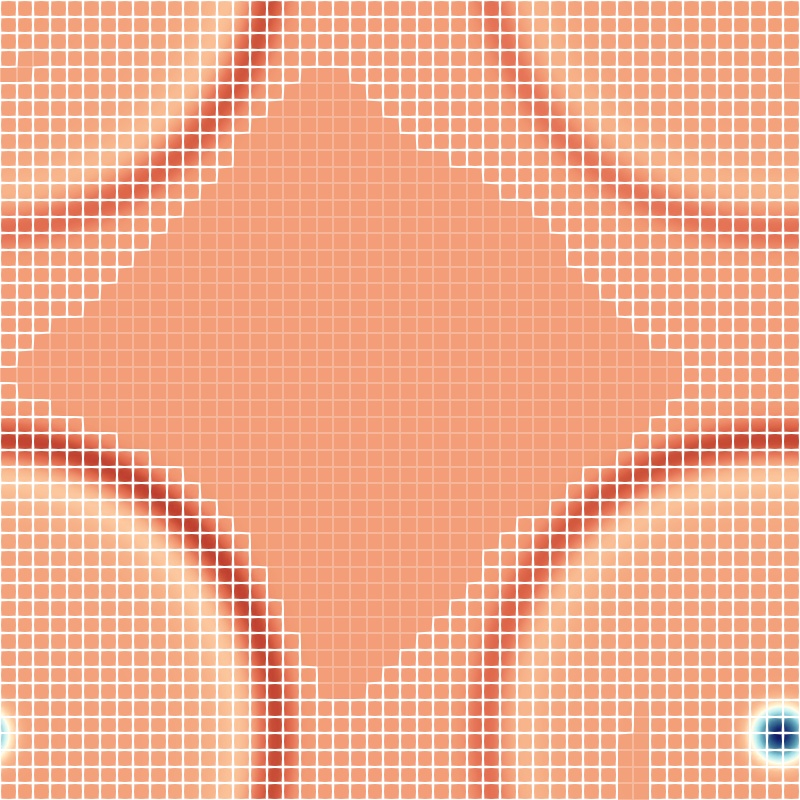}}}
        \newline
        \subfloat[Remesh, $9T$]{
        \adjustbox{width=0.16\linewidth,valign=b}{\includegraphics{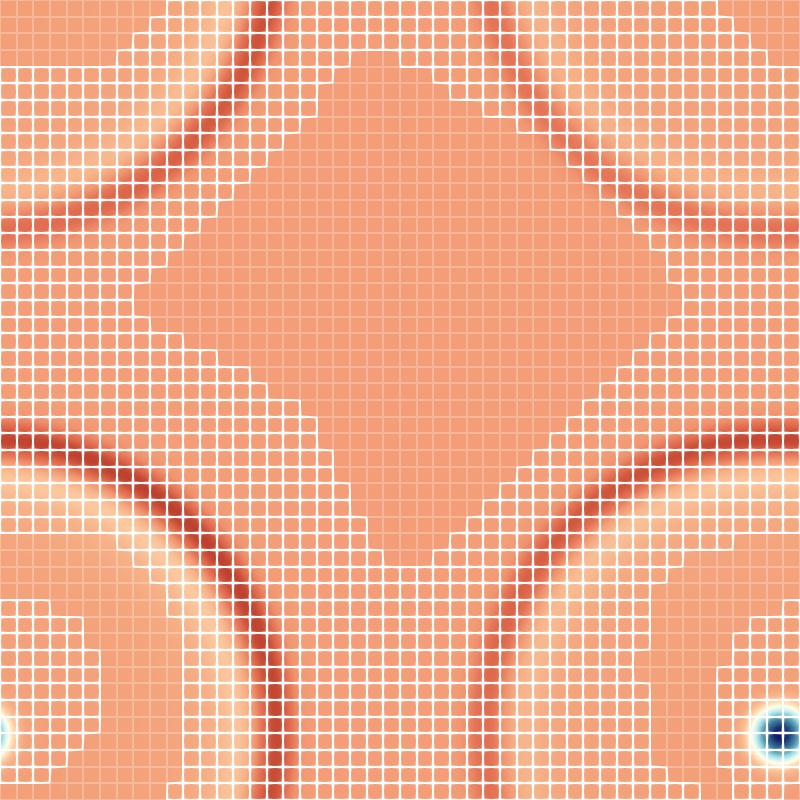}}}
        \subfloat[Solution, $10T$]{
        \adjustbox{width=0.16\linewidth,valign=b}{\includegraphics{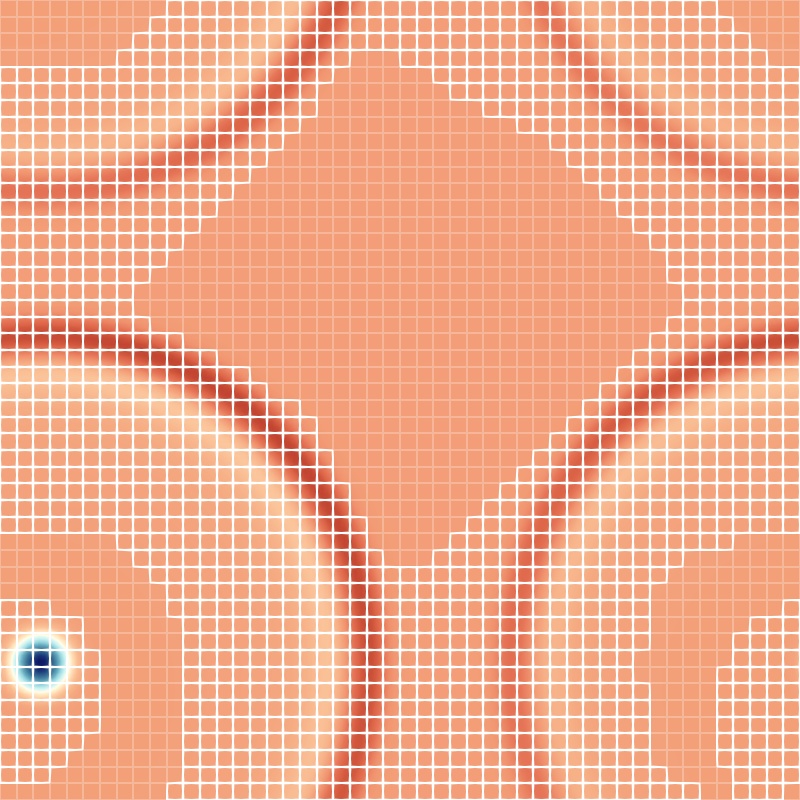}}}
        \subfloat[Remesh, $10T$]{
        \adjustbox{width=0.16\linewidth,valign=b}{\includegraphics{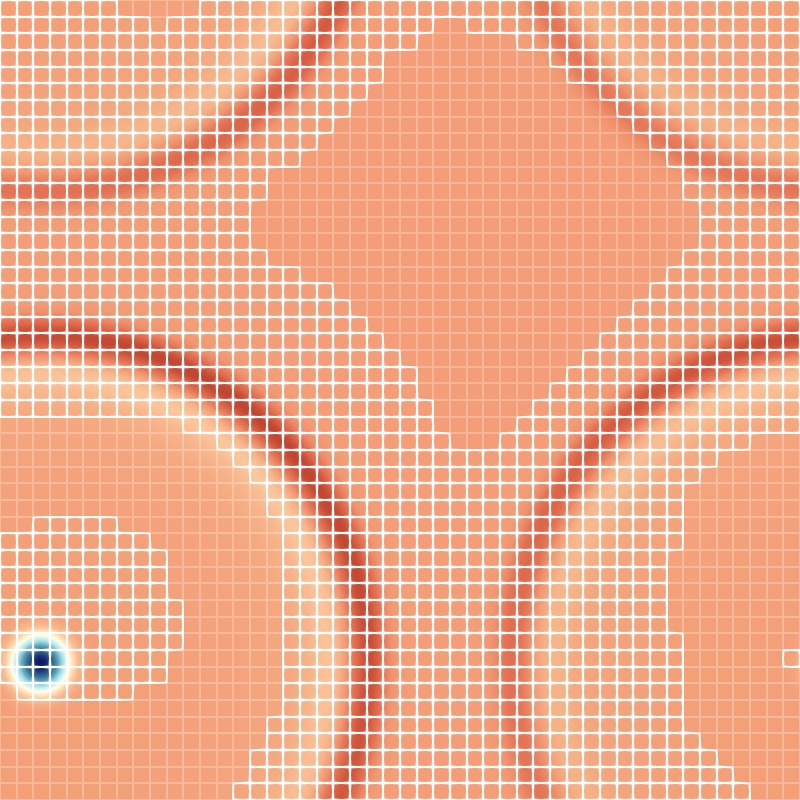}}}
        \subfloat[Solution, $11T$]{
        \adjustbox{width=0.16\linewidth,valign=b}{\includegraphics{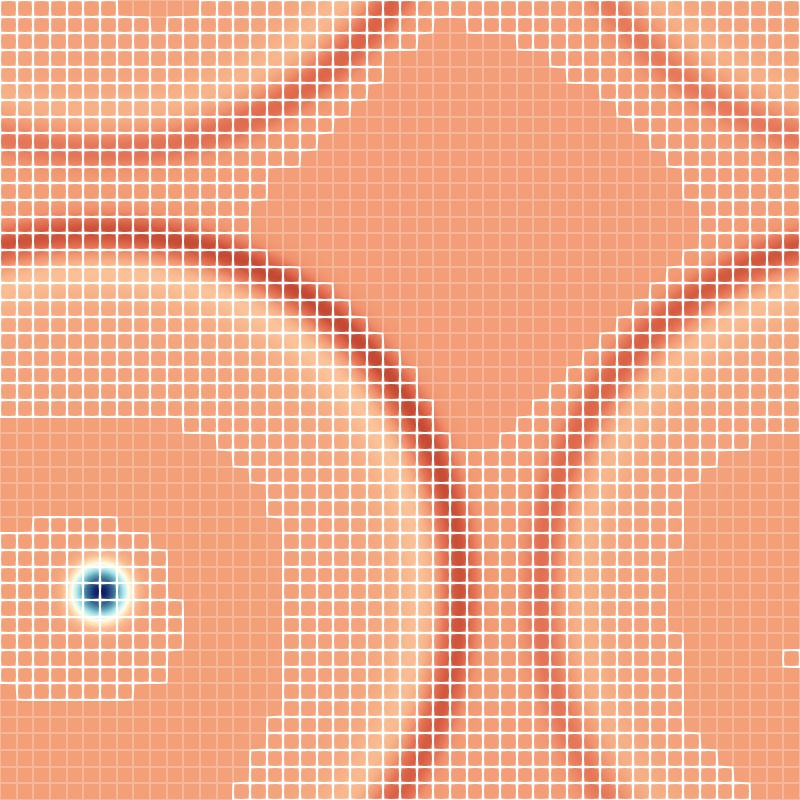}}}
        \subfloat[Remesh, $11T$]{
        \adjustbox{width=0.16\linewidth,valign=b}{\includegraphics{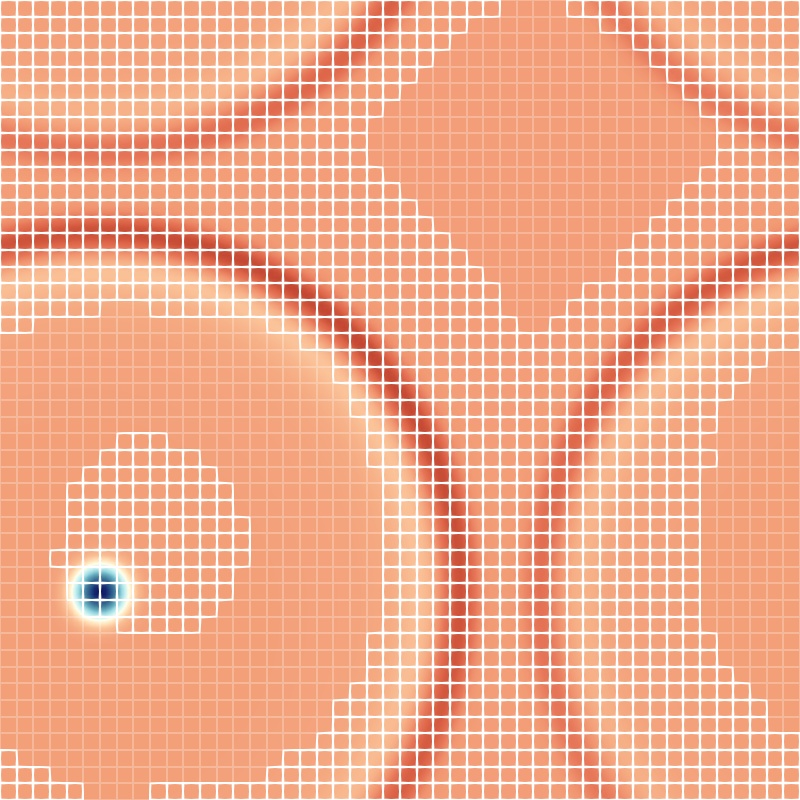}}}
        \subfloat[Solution, $12T$]{
        \adjustbox{width=0.16\linewidth,valign=b}{\includegraphics{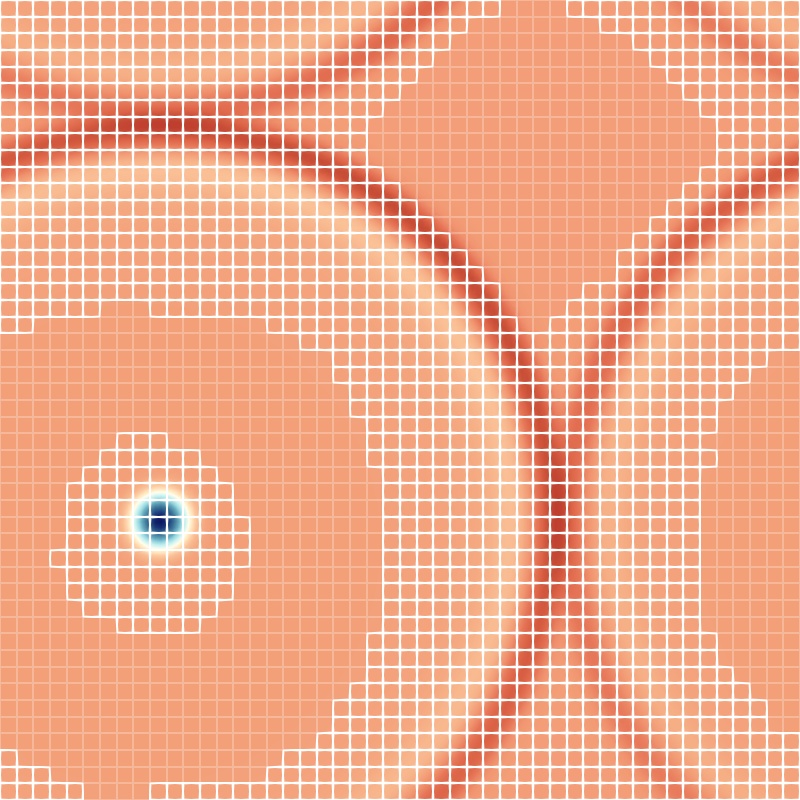}}}
        \newline
        \subfloat[Remesh, $12T$]{
        \adjustbox{width=0.16\linewidth,valign=b}{\includegraphics{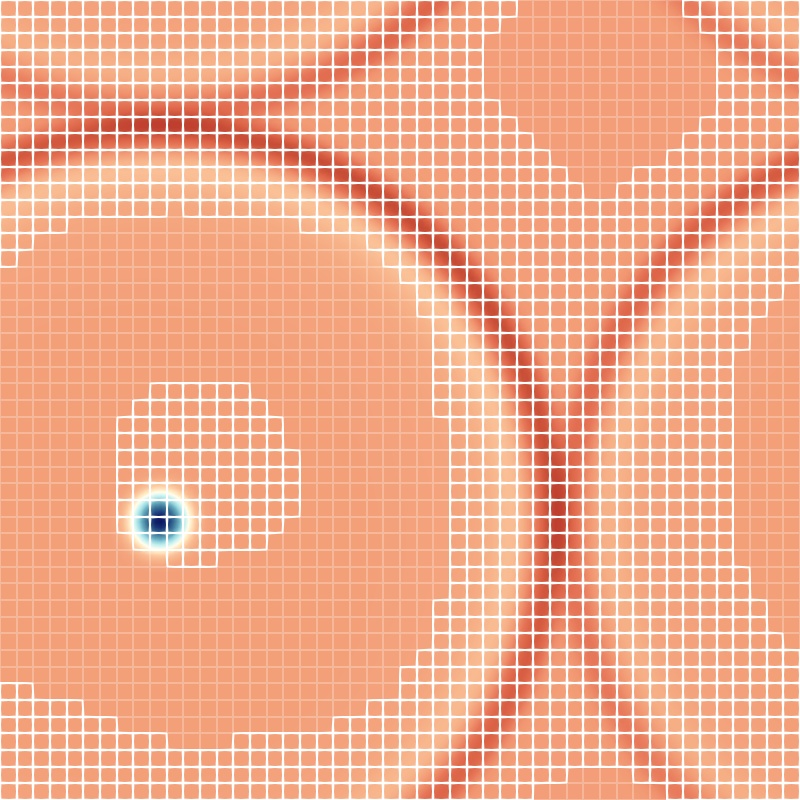}}}
        \subfloat[Solution, $13T$]{
        \adjustbox{width=0.16\linewidth,valign=b}{\includegraphics{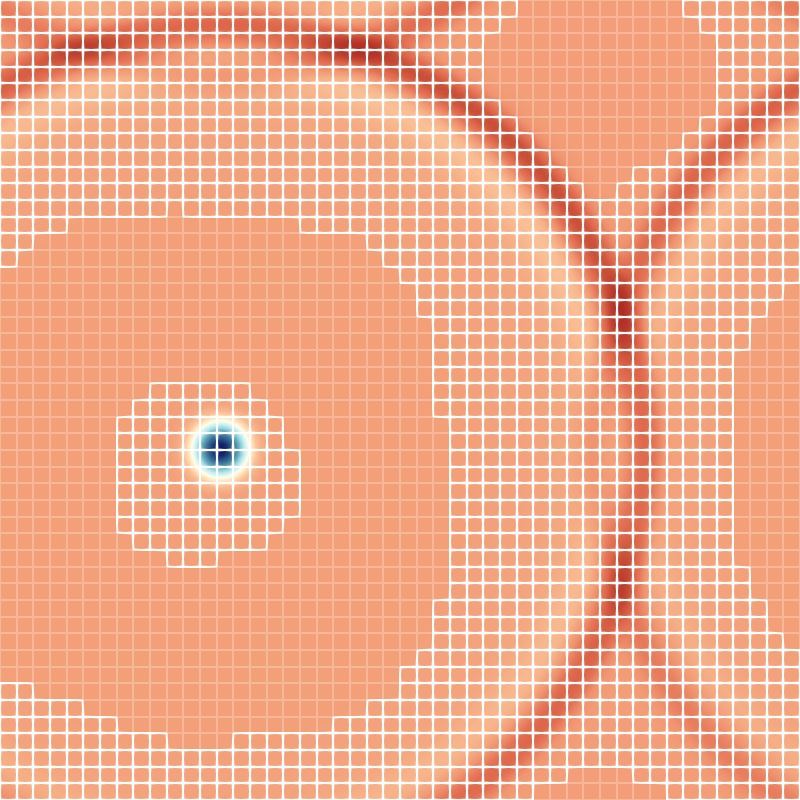}}}
        \subfloat[Remesh, $13T$]{
        \adjustbox{width=0.16\linewidth,valign=b}{\includegraphics{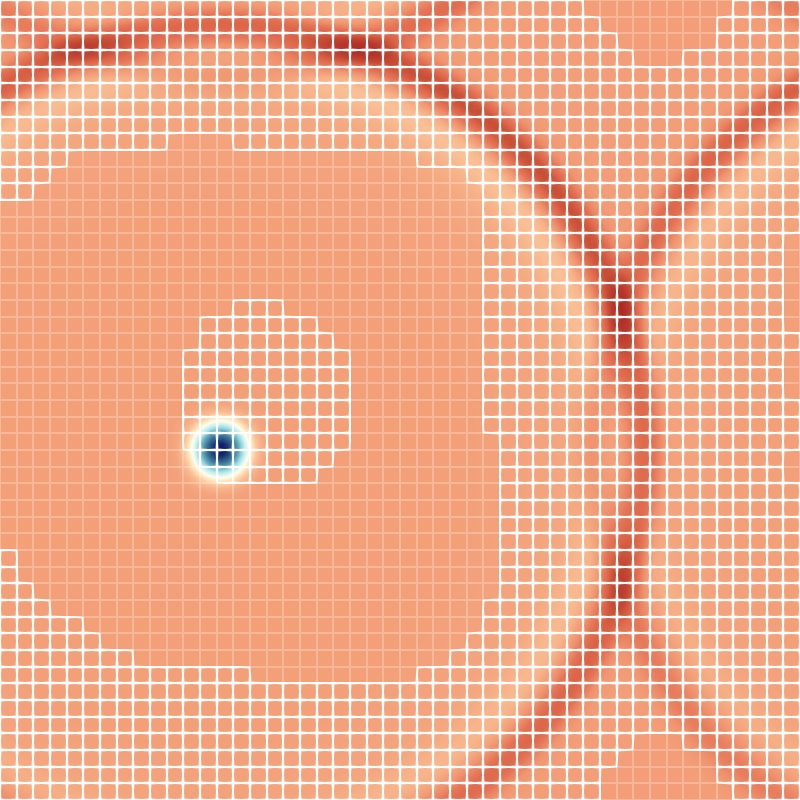}}}
        \subfloat[Solution, $14T$]{
        \adjustbox{width=0.16\linewidth,valign=b}{\includegraphics{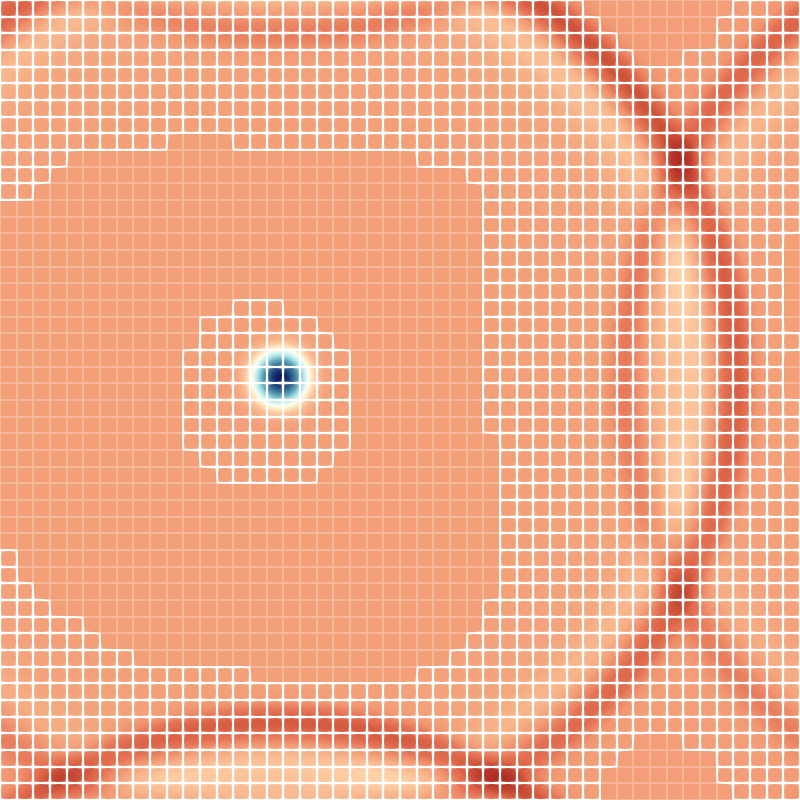}}}
        \subfloat[Remesh, $14T$]{
        \adjustbox{width=0.16\linewidth,valign=b}{\includegraphics{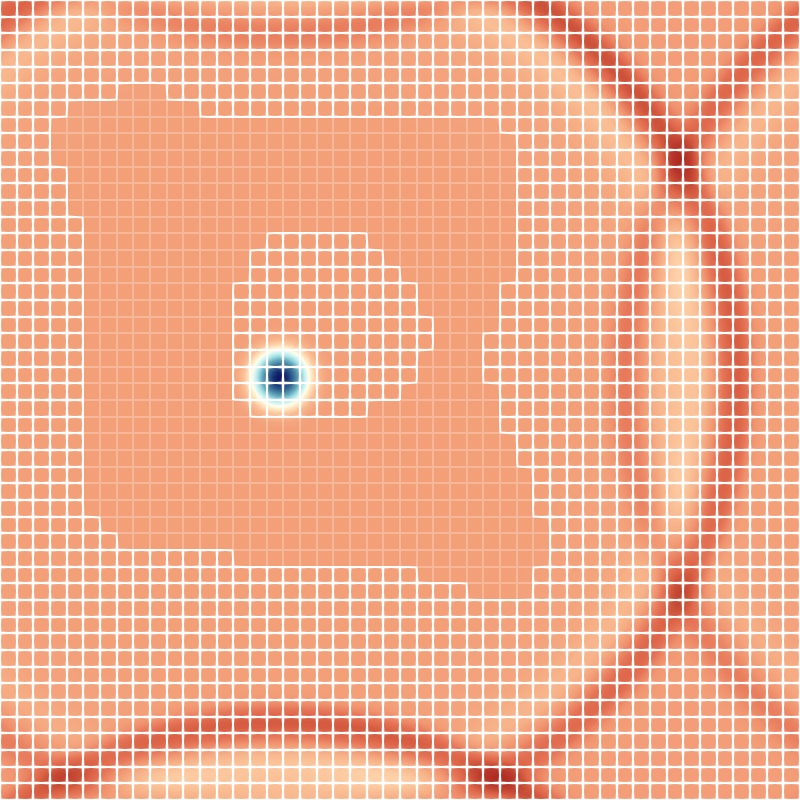}}}
        \subfloat[Solution, $15T$]{
        \adjustbox{width=0.16\linewidth,valign=b}{\includegraphics{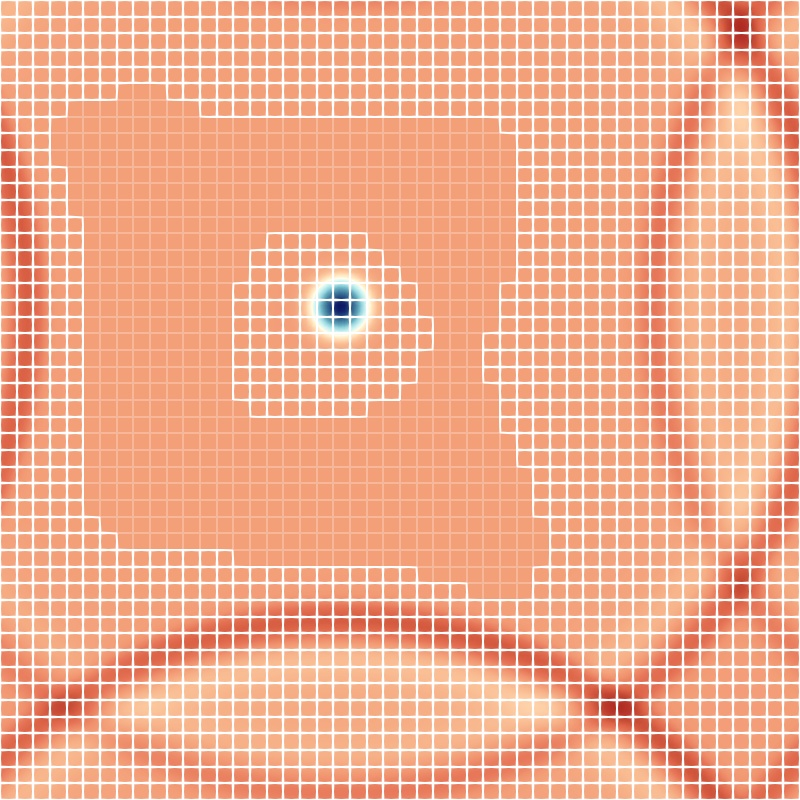}}}
        \newline
        \subfloat[Remesh, $15T$]{
        \adjustbox{width=0.16\linewidth,valign=b}{\includegraphics{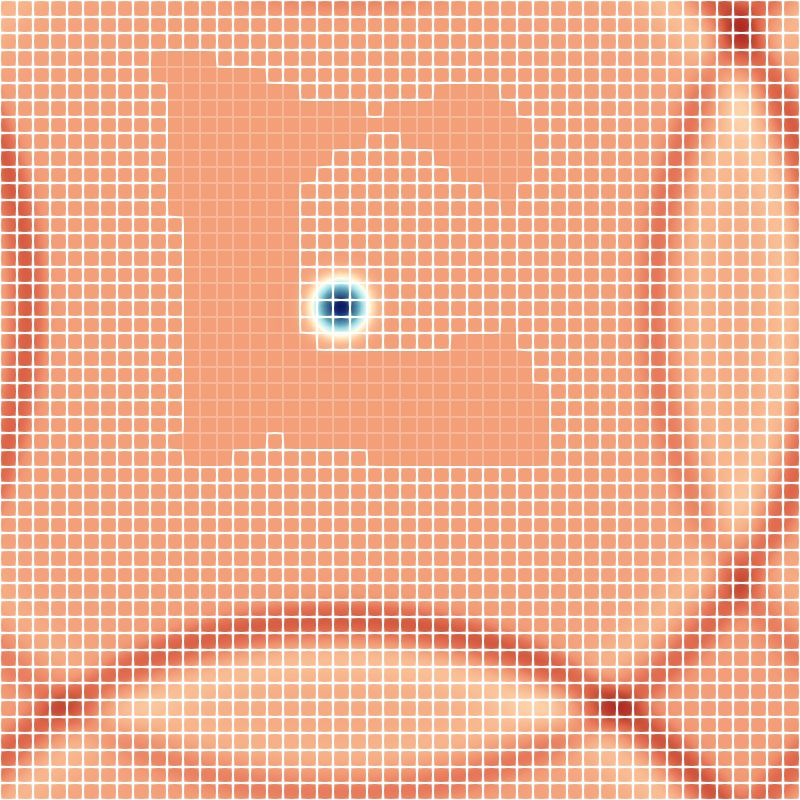}}}
        \subfloat[Solution, $16T$]{
        \adjustbox{width=0.16\linewidth,valign=b}{\includegraphics{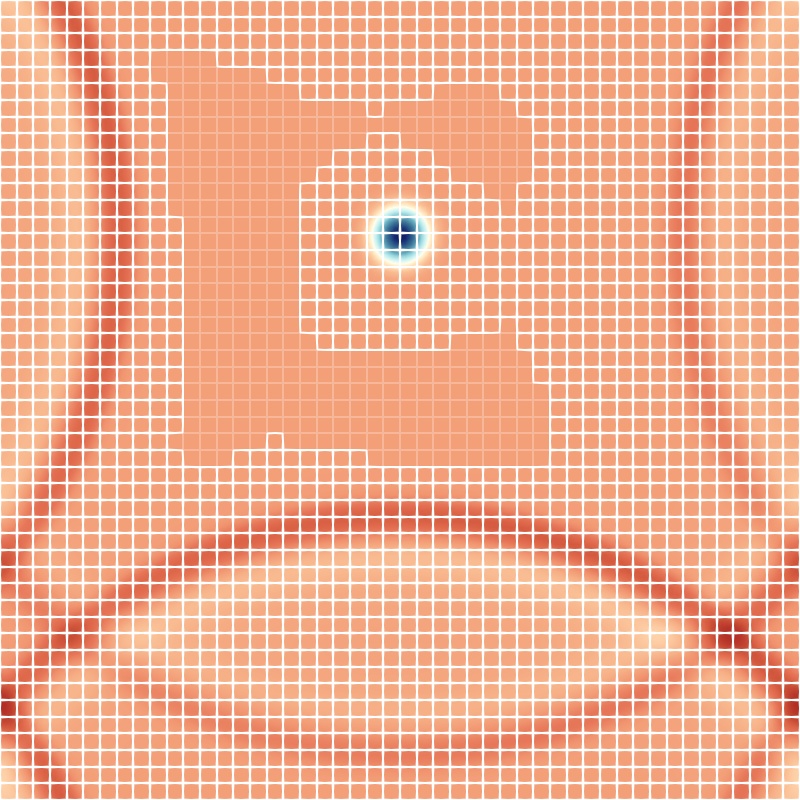}}}
        \subfloat[Remesh, $16T$]{
        \adjustbox{width=0.16\linewidth,valign=b}{\includegraphics{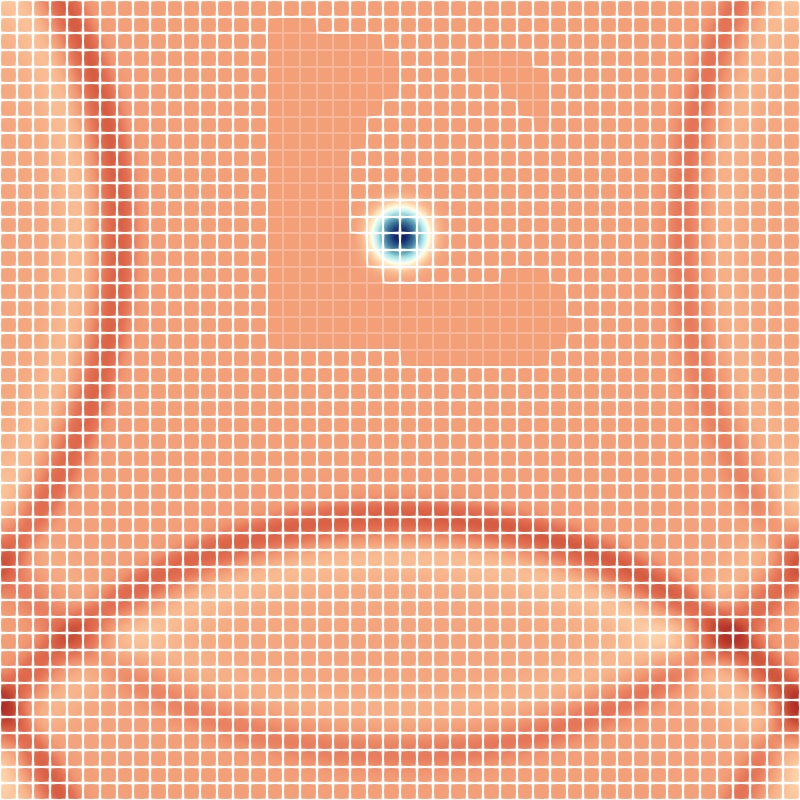}}}
        \subfloat[Solution, $17T$]{
        \adjustbox{width=0.16\linewidth,valign=b}{\includegraphics{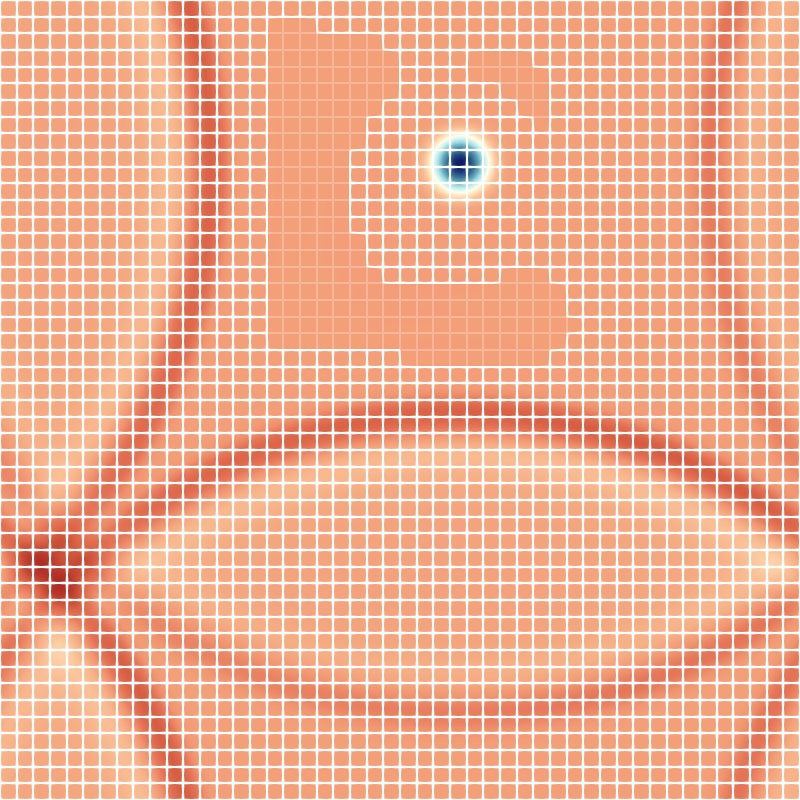}}}
        \subfloat[Remesh, $17T$]{
        \adjustbox{width=0.16\linewidth,valign=b}{\includegraphics{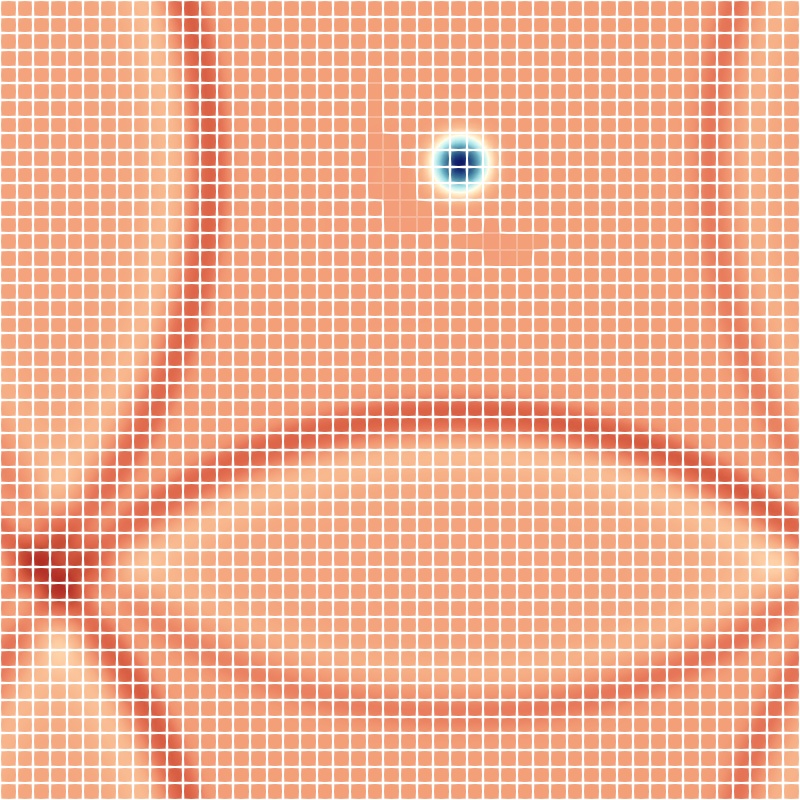}}}
        \subfloat[Solution, $18T$]{
        \adjustbox{width=0.16\linewidth,valign=b}{\includegraphics{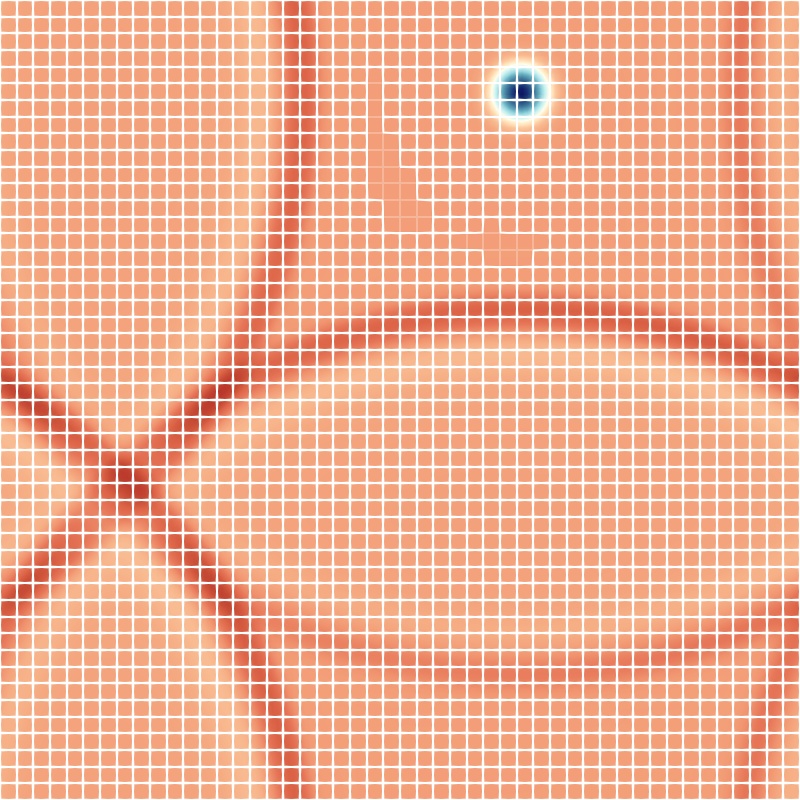}}}
        \newline
        \caption{\label{fig:pressurepulse_amr_drl} Contours of total energy overlaid with $p$-adapted mesh at varying remesh intervals using DynAMO for the convecting pressure pulse problem. Highlighted elements represent $p$-refinement.}
    \end{figure}

For these flow conditions, pressure sources, initialized as smooth Gaussian bumps, cause outward-running acoustic waves centered over a compact low-density/low-energy region. The combination of this acoustic propagation with a background convecting velocity yields interesting flow behavior, particular in the energy field where features with steep flow gradients are separated by relatively flat regions, presenting an ideal case for AMR. We first present an example of the DynAMO approach applied to this class of problems for a randomly selected \textit{in-distribution} set of initial conditions with the parameters equal to $u_0 = 2.25$, $v_0 = 2.67$, $h = 0.12$, $w = 580$, $x_0 = 0.3$, and $y_0 = 0.53$. The evolution of the total energy and $p$-adapted mesh as computed by the DynAMO approach (with $\alpha = \alpha_{\text{train}}$) across 18 remesh time intervals, much longer than the 4 remesh time intervals used for training, is shown in \cref{fig:pressurepulse_amr_drl}. The spatio-temporal evolution of the flow for this class of problems is notably more complex than the convecting density pulse case as a result of the nonlinear acoustic wave interactions. However, even with the significantly longer simulation time, the DynAMO approach was able to preemptively refine regions of the prior to features of interest appearing. At earlier times, the evolution of the acoustic wave was dominated by the background convective velocity, such that the refinement pattern was reminiscent of advection-type problems, with an elongated refinement region stretched in the direction of the propagation velocity. However, as the pressure gradients drove the flow outwards from the pulse center, the preemptive refinement capability of DynAMO was able to accurately account for both the wave propagation and convection. At later times, the acoustic wave propagated far enough from its origin such that a region of low variation in the solution existed between the low energy center and the wave front. The DynAMO policy appropriately de-refined this region without erroneously de-refining the regions across which the wave front propagated. Furthermore, the collisions of the acoustic wave fronts were preemptively predicted in the refinement patterns, indicating that the proposed approach can effectively account for nonlinear flow physics in which the characteristic direction of propagation can strongly vary across elements (e.g., colliding wave fronts). 

\subsubsection{In-distribution experiments}
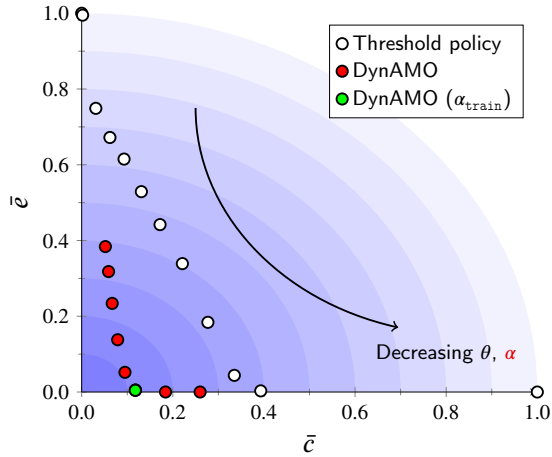
\begin{figure}[htbp!]
        \centering
        \adjustbox{width=0.45\linewidth,valign=b}{\begin{tikzpicture}[spy using outlines={rectangle, height=3cm,width=2.3cm, magnification=3, connect spies}]

    \begin{axis}
    [   axis lines = none,
        xmin = 0, xmax = 1,
        ymin = 0, ymax = 1
    ]
        \fill[fill=blue!5] (1.0, 0) arc[start angle=0, end angle=90, radius=1.0] -- (0,0) -- (1.0, 0);
        \fill[fill=blue!10] (0.9, 0) arc[start angle=0, end angle=90, radius=0.9] -- (0,0) -- (0.9, 0);
        \fill[fill=blue!15] (0.8, 0) arc[start angle=0, end angle=90, radius=0.8] -- (0,0) -- (0.8, 0);
        \fill[fill=blue!20] (0.7, 0) arc[start angle=0, end angle=90, radius=0.7] -- (0,0) -- (0.7, 0);
        \fill[fill=blue!25] (0.6, 0) arc[start angle=0, end angle=90, radius=0.6] -- (0,0) -- (0.6, 0);
        \fill[fill=blue!30] (0.5, 0) arc[start angle=0, end angle=90, radius=0.5] -- (0,0) -- (0.5, 0);
        \fill[fill=blue!35] (0.4, 0) arc[start angle=0, end angle=90, radius=0.4] -- (0,0) -- (0.4, 0);
        \fill[fill=blue!40] (0.3, 0) arc[start angle=0, end angle=90, radius=0.3] -- (0,0) -- (0.3, 0);
        \fill[fill=blue!45] (0.2, 0) arc[start angle=0, end angle=90, radius=0.2] -- (0,0) -- (0.2, 0);
        \fill[fill=blue!50] (0.1, 0) arc[start angle=0, end angle=90, radius=0.1] -- (0,0) -- (0.1, 0);

    \end{axis}
        
    \begin{axis}
    [   axis line style={latex-latex},
        axis y line=left,
        axis x line=left,
        xmode=linear,
        ymode=linear,
        xlabel = {$\bar{c}$},
        ylabel = {$\bar{e}$},
        xmin = 0, xmax = 1,
        ymin = 0, ymax = 1,
        xtick = {0,0.2,0.4,0.6,0.8,1.0},
        ytick = {0,0.2,0.4,0.6,0.8,1.0},
        minor x tick num=1,
        minor y tick num=1,
        legend cell align={left},
        legend style={at={(0.97, 0.97)},anchor=north east},
        clip mode=individual,
        x tick label style={/pgf/number format/.cd, fixed, fixed zerofill, precision=1, /tikz/.cd},
        y tick label style={/pgf/number format/.cd, fixed, fixed zerofill, precision=1, /tikz/.cd},
        label style={font=\large},
    ]
        \addplot[color=black, fill=white, style={thick}, only marks, mark=*, mark options={scale=1.2}]
        table[x = cm, y = em, col sep=comma]{./figs/data/euler_pref_2d_pa_base.csv};
        \addlegendentry{Threshold policy};
        
        \addplot[color=black, fill=red, style={thick}, only marks, mark=*, mark options={scale=1.2}]
        table[x = cm, y = em, col sep=comma]{./figs/data/euler_pref_2d_pa_base_rleps.csv};
        \addlegendentry{DynAMO};

        \addplot[color=black, fill=green, style={thick}, only marks, mark=*, mark options={scale=1.2}]
        coordinates {(1.177e-1, 4.526e-3)};
        \addlegendentry{DynAMO ($\alpha_{{\tt train}}$)};

        \draw[->, black, thick] (0.25,.75) arc (180:255:0.6);
        \draw (0.8, 0.1) node[scale=1] {Decreasing $\theta$, \textcolor{red!90!black}{$\alpha$}};
        
        \addplot[color=black, fill=white, style={thick}, only marks, mark=*, mark options={scale=1.2}]
        coordinates {(1.0, 0.0)};
        
    \end{axis}
\end{tikzpicture}}
        \caption{\label{fig:euler_pref_pareto_pa} Pareto plot of normalized cost vs. error for $p$-refinement on the Euler equations with DynAMO and the threshold policy for the in-distribution convecting pressure pulse case. Contours of efficiency shown on background.  Red markers represent a sweep of $\alpha \in [10^{-6}, 0.6]$ at evaluation time, green marker represents $\alpha$ set to the training value of $\alpha_{\texttt{train}} = 0.1$. Peak efficiency achieved at $\alpha=0.2$. 
        }
    \end{figure}

To quantify the capability of the DynAMO approach for this case and present a comparison to baseline methods, the normalized cost vs. error distributions are shown in the Pareto plot in \cref{fig:euler_pref_pareto_pa} for the DynAMO policy and the standard threshold policy. Both policies were evaluated over a range of threshold parameters to observe the effects on error, cost, and efficiency. For the threshold policy, where the threshold parameter was varied $\theta = 10^{-2}$ to $10^{-11}$, a linear relation was observed in the cost to error profile, with the most optimal value yielding a peak efficiency of roughly $0.65$. In contrast, the DynAMO policy was able to produce a far more efficient AMR approach. At the training threshold value $\alpha_{\text{train}} = 0.1$, the policy was able to achieve an efficiency of approximately $0.9$ with almost zero normalized error, yielding results of similar accuracy to a fully refined mesh. In contrast, the threshold policy could only achieve a normalized error of approximately $0.5$ for the same computational cost. Similarly, to achieve the same normalized error, the threshold policy required roughly quadruple the computational cost. When the error threshold for the DynAMO policy was varied at evaluation time from $\alpha = 10^{-6}$ to $0.6$, the quintessential convex Pareto curve was recovered with the sampled results much more closely distributed near the origin. Most importantly, all of the sampled points along that curve were well within the bounds of the threshold policy results, showcasing that the proposed approach can unlock efficiency levels that are out of reach with conventional AMR techniques. 

To verify that the increased efficiency of the proposed DynAMO approach persists across a wide distribution of test cases, a comparison of the DynAMO approach ($\alpha = \alpha_{\text{train}}$) and the threshold policy was performed across 100 randomly sampled \textit{in-distribution} initial conditions. The results of the comparison are shown in \cref{tab:euler_pref_indistribution} for varying threshold values for the threshold policy. At its peak, the threshold policy achieves a mean efficiency of $0.664$, with normalized mean error and cost values of $0.184$ and $0.277$, respectively. In comparison, the DynAMO policy achieved a noticeably higher mean efficiency, $0.882$, a relative increase of $32.8\%$. The benefits of the DynAMO approach were more pronounced when observing the relative decrease in mean error and computational cost. The DynAMO policy achieved a relative mean error of $0.005$, which effectively results in a level of accuracy identical to a fully-refined mesh. This error was $97.3\%$ less than the mean error provided by the most efficient threshold policy. Furthermore, the DynAMO policy achieved this level of error with a relative mean computational cost of $0.118$, which is $57.4\%$ less than the computational cost of the most efficient threshold policy that results in a mean error that is nearly $40$ times higher. To achieve the level of error provided by the DynAMO policy, the threshold policy would essentially have to utilize a fully-refined mesh, such that the mean computational cost would be almost $10$ times higher. Conversely, for the same mean computational cost as the DynAMO policy, the use of the threshold policy would result in more than two orders of magnitude increase in error.

    \begin{figure}[htbp!] 
        \centering
        \begin{tabular}{cccc}
        \toprule
        Method &  Efficiency & Normalized error & Normalized cost\\ 
        \midrule
        Threshold ($\theta = 10^{-2}$) &  0.000 (0.000) & 1.000 (0.000) &  0.000 (0.000) \\
        Threshold ($\theta = 10^{-3}$) &  0.005 (0.020) & 0.995 (0.020) & 0.002 (0.001)  \\
        Threshold ($\theta = 10^{-4}$) &  0.250 (0.048) & 0.749 (0.049) & 0.031 (0.006) \\
        Threshold ($\theta = 10^{-5}$) &  0.325 (0.023) & 0.672 (0.023) & 0.062 (0.006) \\
        Threshold ($\theta = 10^{-6}$) &  0.377 (0.033) & 0.615 (0.034) & 0.093 (0.008) \\
        Threshold ($\theta = 10^{-7}$) &  0.455 (0.032) & 0.529 (0.034) & 0.131 (0.009) \\
        Threshold ($\theta = 10^{-8}$) &  0.525 (0.038) & 0.442 (0.045) & 0.172 (0.010) \\
        Threshold ($\theta = 10^{-9}$) &  0.593 (0.040) & 0.339 (0.055) & 0.221 (0.011) \\
        Threshold ($\theta = 10^{-10}$) & \textbf{0.664 (0.029)} & \textbf{0.184 (0.057)} & \textbf{0.277 (0.011)} \\
        Threshold ($\theta = 10^{-11}$) &  0.661 (0.011) & 0.044 (0.026) & 0.335 (0.011) \\
        \midrule
        DynAMO & \textbf{0.882 (0.022)} &  \textbf{0.005 (0.008)} & \textbf{0.118 (0.023)} \\
        DynAMO/Optimal $\theta$ & \textcolor{green!70!black}{+32.8\%} & \textcolor{green!70!black}{-97.3\%} & \textcolor{green!70!black}{-57.4\%} \\
        \bottomrule
        \end{tabular} 
        \captionof{table}{\label{tab:euler_pref_indistribution} Comparison of the mean efficiency, normalized error, and normalized cost for $p$-refinement on the Euler equations with DynAMO and the threshold policy for the convecting pressure pulse problem over 100 \textit{in-distribution} runs using uniform random initial conditions. Standard deviation shown in parentheses.}
    \end{figure}

\subsubsection{Generalization experiments}
The efficacy and generalization capabilities of the DynAMO approach were further evaluated through experiments for \textit{out-of-distribution} problems, where the threshold parameter of the threshold policy was varied and the threshold parameter for the DynAMO policy was set to its training value. For brevity, certain results for the threshold policy are omitted as the recovered efficiency levels were not significant. 

The first generalization experiment was performed with respect to different flow physics through the convecting density pulse problem. This problem is given by the initial conditions
\begin{subequations}
\begin{align}\label{eq:density_pulse}
    \rho(x, y, 0) &= 1 + h \exp \left(-w \left((x - x_0)^2 + (y - y_0)^2 \right) \right), \\
    u(x, y, 0) &= u_0, \\
    v(x, y, 0) &= v_0, \\
    P(x, y, 0) &= 1,
\end{align}
\end{subequations}
where the parameters $u_0$, $v_0$, $h$, $w$, $x_0$, and $y_0$ were uniformly distributed across the same ranges as \cref{eq:pressure_pulse}. For this setup, a perturbation in the density field convects at a constant velocity and direction, much like in the advection problems. Besides showing the generalizability of the proposed observation and reward functions, the purpose of this example is to show a direct extension of the DynAMO approach from the advection equation to the Euler equations for a problem which effectively exhibits the same physics but achieves it from a more complex set of governing equations. Recall that while the evolution of the density field is linear with respect to the solution, the evolution of the total energy field, the solution component which is actually being observed by the policy, is not.

While the flow physics may seem relatively simple, this particular class of problems presents a very interesting test for the proposed policy as the initial distribution of the total energy, the observable quantity, is effectively identical between these cases and the pressure pulses that the policy was trained on, but the resulting flow fields are radically different. As such, it presents a validation of the ability of the proposed approach in learning the physics of the underlying governing equations as opposed to simply extrapolating from previously encountered flow problems. A comparison between the threshold policy and DynAMO policy over 100 randomly sampled initial conditions for the problem is shown in \cref{tab:euler_pref_ood}. Due to the compactly supported nature of the convecting density problem, the threshold policy was able to achieve a high mean efficiency, $0.799$. However, the DynAMO policy was able to achieve even higher, $0.955$, almost reaching the theoretical limit of unity. The relative efficiency benefits of the DynAMO approach over the optimal threshold policy remained similar to the \textit{in-distribution} case, such that it can be surmised that the approach may have the ability to effectively generalize to flow physics not encountered during the training process.

    \begin{figure}[htbp!] 
        \centering
        \begin{tabular}{ccccccc}
        \toprule
        Method & In-distribution & Convecting density &  Finer mesh  & Longer remesh time & Longer sim. time \\ 
        \midrule    
        Threshold ($\theta = 10^{-2}$)     &  0.000 (0.000) & 0.000 (0.000) &		0.000 (0.000) & 0.000 (0.000) & 0.000 (0.000) \\
        Threshold ($\theta = 10^{-4}$)     &  0.250 (0.048) & 0.015 (0.003) &		0.007 (0.006) & 0.011 (0.005) & 0.020 (0.011) \\
        Threshold ($\theta = 10^{-6}$)     &  0.377 (0.033) & 0.036 (0.007) &		0.200 (0.062) & 0.025 (0.022) & 0.114 (0.046) \\
        Threshold ($\theta = 10^{-8}$)     &  0.525 (0.038) & 0.073 (0.008) &		0.637 (0.092) & 0.075 (0.053) & 0.495 (0.025) \\
        Threshold ($\theta = 10^{-9}$)     &  0.593 (0.040) & 0.099 (0.009) &		0.738 (0.076) & 0.124 (0.076) & \textbf{0.516 (0.022)} \\
        Threshold ($\theta = 10^{-10}$)    &  \textbf{0.664 (0.029)} & 0.128 (0.011) &		0.824 (0.032) & 0.224 (0.084) & 0.507 (0.014) \\
        Threshold ($\theta = 10^{-11}$)    &  0.661 (0.011) & 0.163 (0.012) &		\textbf{0.843 (0.007)} & 0.308 (0.038) & 0.456 (0.015) \\
        Threshold ($\theta = 10^{-12}$)    &  0.607 (0.010) & \textbf{0.799 (0.012)} &		0.820 (0.008) & 0.348 (0.020) & 0.394 (0.014) \\
        Threshold ($\theta = 10^{-14}$)    &  0.493 (0.010) & 0.707 (0.013) &		0.757 (0.010) & 0.429 (0.040) & 0.284 (0.011) \\
        Threshold ($\theta = 10^{-15}$)    &  0.342 (0.020) & 0.459 (0.034) &		0.478 (0.047) & \textbf{0.470 (0.024)} & 0.183 (0.010) \\
        \midrule
        DynAMO & \textbf{0.882 (0.022)} & \textbf{0.955 (0.014)}  &  \textbf{0.910 (0.012)} & \textbf{0.568 (0.074)} & \textbf{0.726 (0.026)}  \\
        DynAMO/Optimal $\theta$ & \textcolor{green!70!black}{+32.8\%} & \textcolor{green!70!black}{+19.5\%} & \textcolor{green!70!black}{+7.9\%} & \textcolor{green!70!black}{+20.8\%} & \textcolor{green!70!black}{+40.1\%}  \\
        \bottomrule
        \end{tabular}
        \captionof{table}{\label{tab:euler_pref_ood} 
        Comparison of the mean efficiency for $p$-refinement on the Euler equations with DynAMO and the threshold policy for the two-dimensional convecting pressure-pulse over 100 \textit{out-of-distribution} runs using uniform random initial conditions with finer mesh resolution, longer remesh time, and longer simulation time. Standard deviation shown in parentheses. In-distribution results from \cref{tab:euler_pref_indistribution} shown for comparison. All sampled thresholds are not shown for brevity.
        }
    \end{figure}




    \begin{figure}[htbp!]
        \centering
        \subfloat[Remesh at $t = 0$]{
        \adjustbox{width=0.24\linewidth,valign=b}{\includegraphics{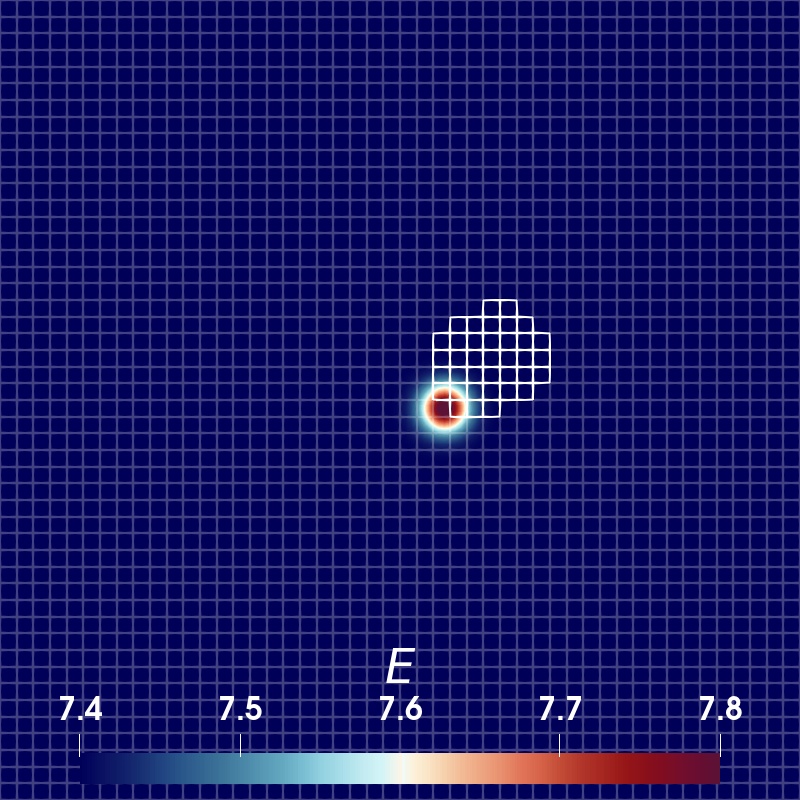}}}
        \subfloat[Solution at $t = T$]{
        \adjustbox{width=0.24\linewidth,valign=b}{\includegraphics{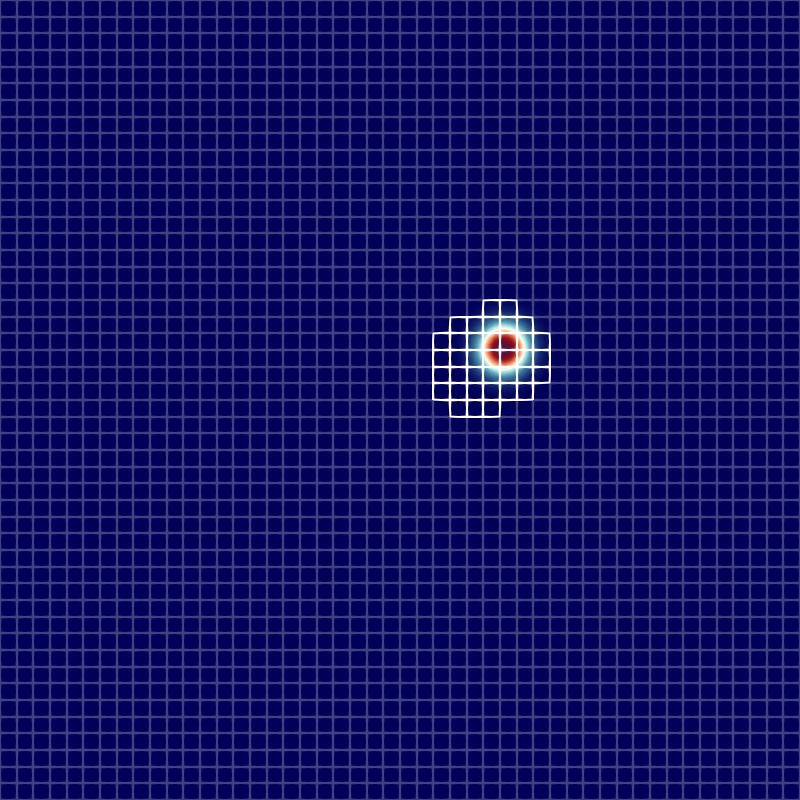}}}
        \subfloat[Remesh at $t = T$]{
        \adjustbox{width=0.24\linewidth,valign=b}{\includegraphics{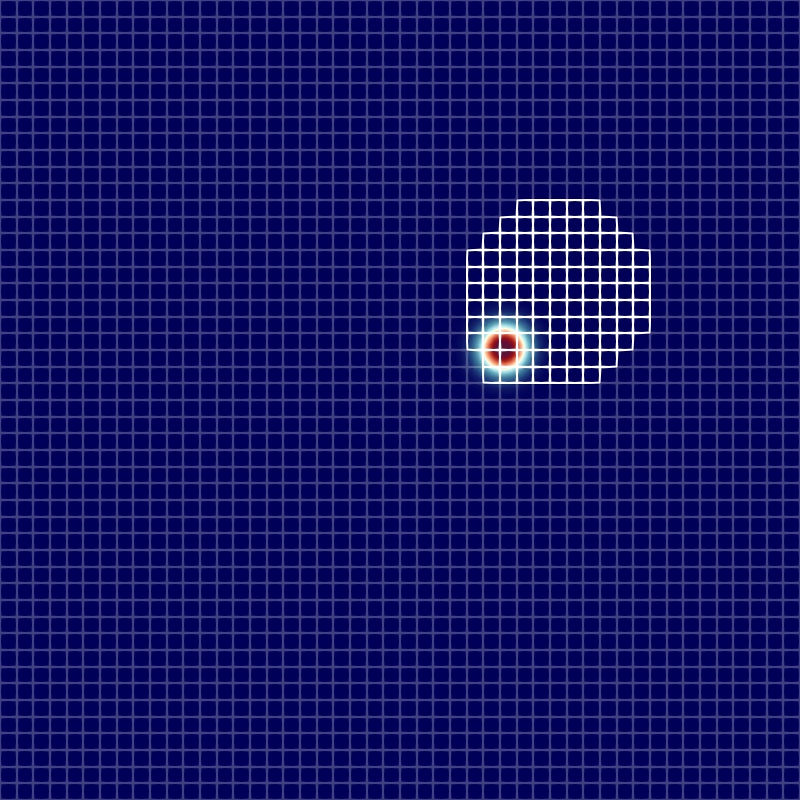}}}
        \subfloat[Solution at $t = 2T$]{
        \adjustbox{width=0.24\linewidth,valign=b}{\includegraphics{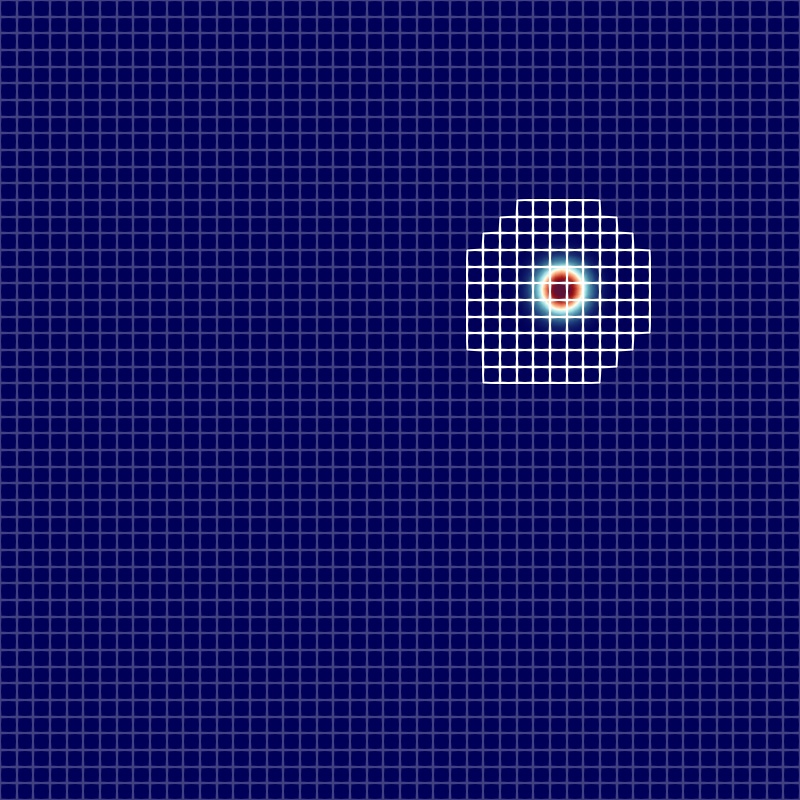}}}
        \newline
        \subfloat[Remesh at $t = 2T$]{
        \adjustbox{width=0.24\linewidth,valign=b}{\includegraphics{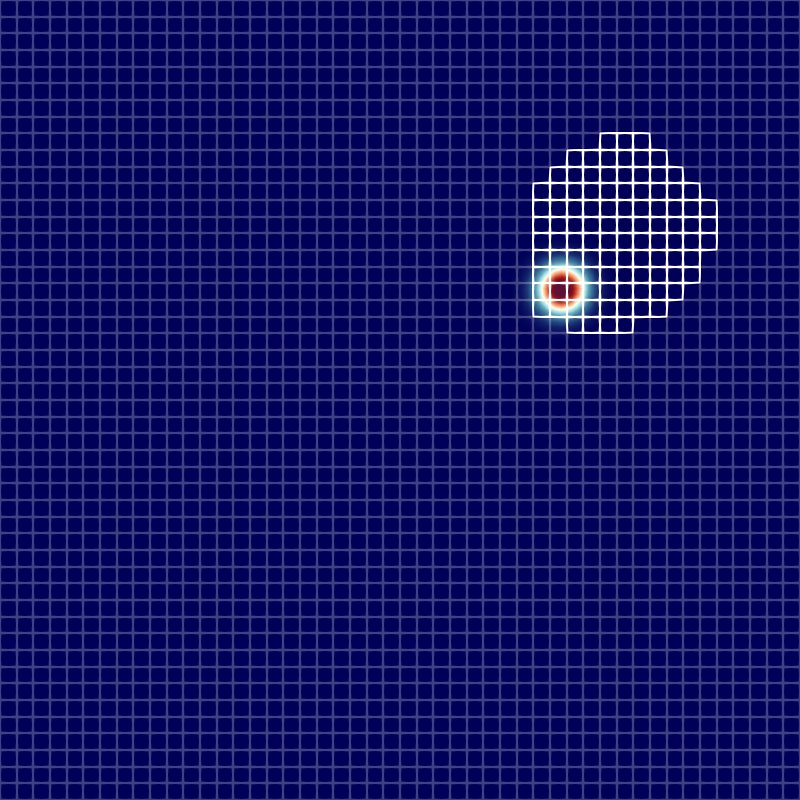}}}
        \subfloat[Solution at $t = 3T$]{
        \adjustbox{width=0.24\linewidth,valign=b}{\includegraphics{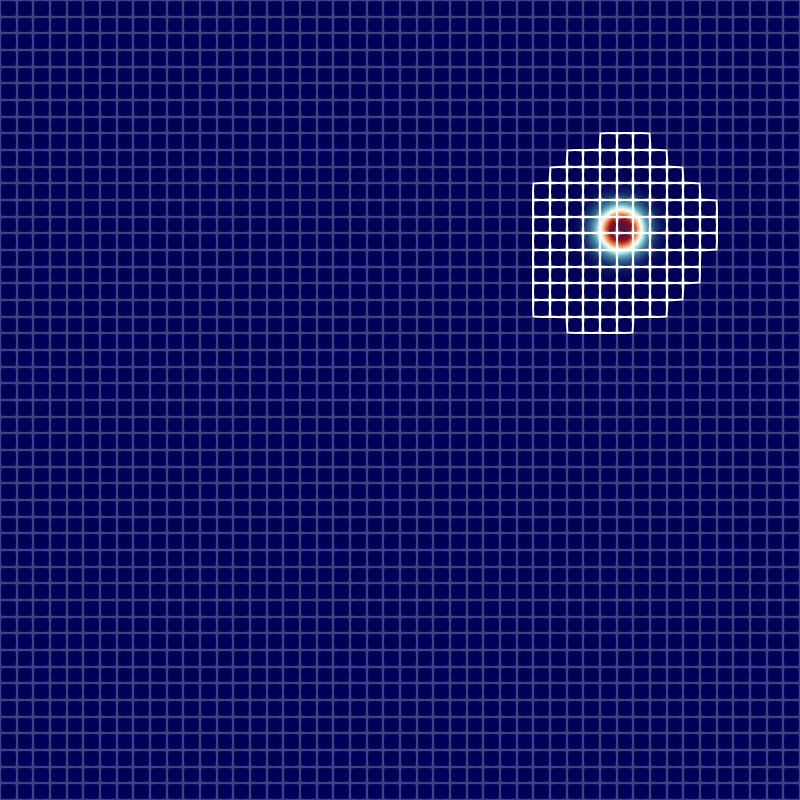}}}
        \subfloat[Remesh at $t = 3T$]{
        \adjustbox{width=0.24\linewidth,valign=b}{\includegraphics{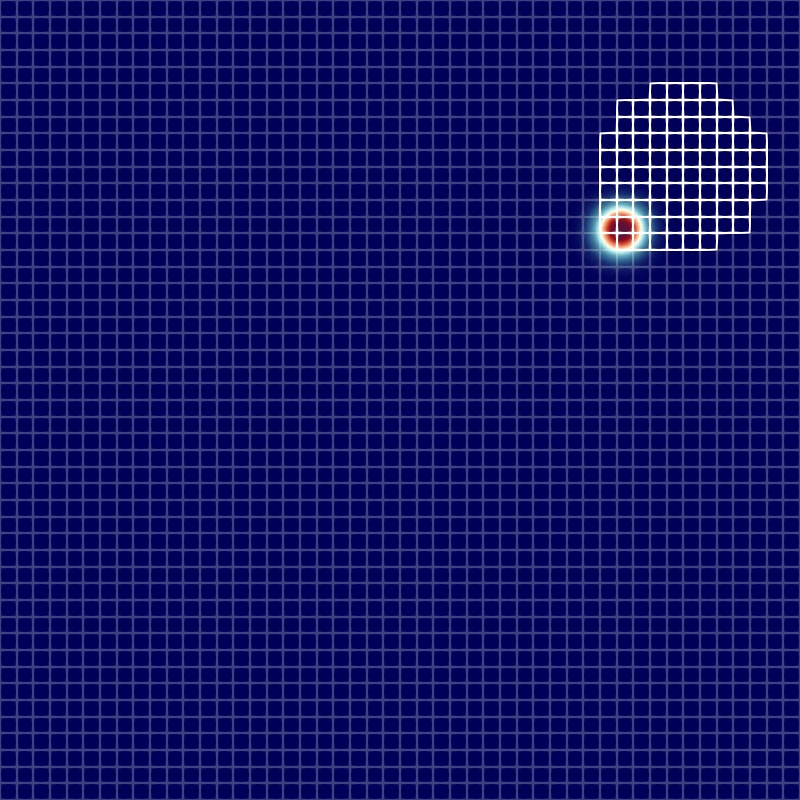}}}
        \subfloat[Solution at $t = 4T$]{
        \adjustbox{width=0.24\linewidth,valign=b}{\includegraphics{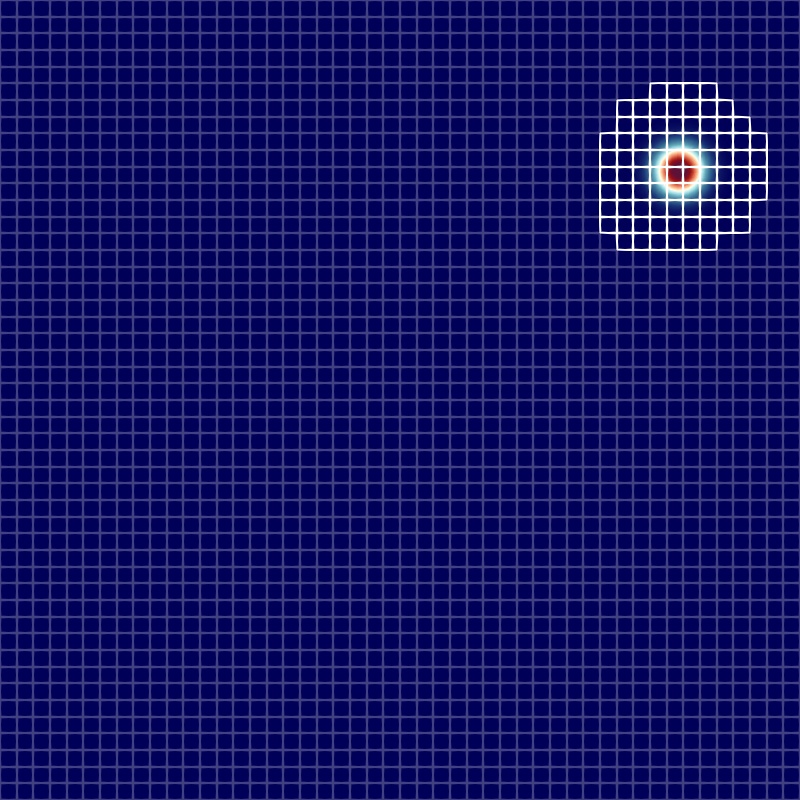}}}
        \newline
        \caption{\label{fig:densitywave_amr_drl} Contours of total energy overlaid with $p$-adapted mesh at varying remesh intervals using DynAMO for the convecting density pulse problem. Highlighted elements represent $p$-refinement.} 
    \end{figure}

To present an visualization of the approach for the convecting density pulse problem, a random initial condition was chosen from the distribution with parameters equal to $u_0 = 2.46$, $v_0 = 2.99$, $h = 0.08$, $w = 267$, $x_0 = 0.3$, and $y_0 = 0.53$. The evolution of the total energy and $p$-adapted mesh as obtained by the DynAMO policy trained on the pressure pulse cases is shown in \cref{fig:densitywave_amr_drl}. For brevity, we present only the first four remesh intervals, although this behavior was essentially identical at later simulation times. Furthermore, the results are shown with respect to the total energy, which is the observable field, instead of the density, but these distributions are nearly identical for a constant background velocity and pressure. It can be seen that similar behavior as with the policies trained on the advection equation is recovered, even with the more complex underlying governing equations and observations. The DynAMO policy predicts an elongated refinement region at the initial time, which adequately covers the propagation path of the density pulse. This preemptive refinement capability is maintained throughout the entire simulation time, such that the density pulse always remains within the refined region of the mesh. While the physics of the problem are simple, these results showcase a direct connection between linear and nonlinear equations for the proposed formulation.

The next generalization experiment was performed with respect to the mesh resolution. The base mesh resolution was increased from $N=48^2$ to $N=96^2$, such that the total number of degrees of freedom was over $3 \cdot 10^5$ at the coarsest level. The mean efficiency across 100 sampled initial conditions on the finer mesh level is shown in \cref{tab:euler_pref_ood} for both the DynAMO policy and the threshold policy. It can be seen that both policies recover a higher efficiency on the finer mesh which is consistent with previous results. For the threshold policy, this efficiency increased from $0.664$ to $0.843$, whereas for the DynAMO policy, this efficiency increased from $0.882$ to $0.910$. As a result, the relative benefit of DynAMO was lower, only $7.9\%$, but this is largely attributed to the marginal gains possible at the highest efficiency levels. However, the results indicate that the DynAMO policy can very effectively generalize to different mesh resolutions for the Euler equations.

The third generalization experiment was performed with respect to the remesh time $T$, which was doubled from its training value. With longer remesh time intervals, it is expected that the efficiency of policies based on instantaneous error estimators would degrade. This behavior was indeed observed for the threshold policy, also shown in \cref{tab:euler_pref_ood}, where the highest mean efficiency was $0.470$ over 100 randomly sampled initial conditions, noticeably less than the in-distribution result of $0.664$. Furthermore, this efficiency was only recovered at machine zero values of the threshold parameter, meaning that the policy refines any element with a non-constant total energy distribution. A notably more efficient approach was recovered with the DynAMO policy, with a mean efficiency of $0.568$. These relative efficiency benefits were on par with the increase observed for the \textit{in-distribution} experiments, indicating that the approach can also effectively generalize to longer remesh times. 

The final generalization experiment was performed with respect to the total simulation time, which was also doubled from its training value. The comparison of the threshold policy at varying values of the threshold parameter to the DynAMO approach over 100 randomly sampled initial conditions is shown in \cref{tab:euler_pref_ood}. It can be seen that the peak efficiency of the threshold approach decreased from $0.664$ to $0.516$, but the mean efficiency of the DynAMO approach decreased as well, from $0.882$ to $0.726$. This decrease in efficiency is somewhat expected at longer simulation times for flows whose complexity increases over time as the most appropriate refinement decision tends toward uniform refinement, which is seen in \cref{fig:pressurepulse_amr_drl}, such that the increase in normalized cost decreases the efficiency. However, we observe that the efficiency benefits of the DynAMO approach do not deteriorate for longer simulation times but, in fact, noticeably appreciate. This further supports the generalization abilities of the proposed approach.

\subsection{Euler equations with $h$-refinement}
Following the initial evaluation of the proposed DynAMO approach for the Euler equations on smooth problems, the approach was then extended to significantly more complex and highly nonlinear flow physics featuring shocks, rarefaction waves, and discontinuities. The discontinuous, compact nature of many of these features paired with their nonlinear evolution presents an ideal example for problems where anticipatory mesh refinement would excel while methods based on instantaneous estimators would be unsuccessful. Due to the lack of smoothness of these features, $h$-refinement was used for the policies as it is typically better suited to resolve such features.  

\begin{figure}[htbp!]
        \centering
        \adjustbox{width=0.35\linewidth,valign=b}{     \begin{tikzpicture}[spy using outlines={rectangle, height=3cm,width=2.3cm, magnification=3, connect spies}]
		\begin{axis}[name=plot1,
		    tick style={draw=none},
		    axis x line=left,
            axis y line=left,
            axis equal image,
            clip mode=individual,
    		xmin=0,
    		xmax=1,
    		xtick={0.0, 0.5, 1.0},
		    xlabel={$x$},
    		ymin=0,
    		ymax=1,
    		ytick={0.0, 0.5, 1.0},
		    ylabel={$y$},
    		ylabel style={rotate=-90},
    		style={font=\small},
    		scale = 1]

		    \draw[-] (axis cs:0.0, 0.0) -- (axis cs:1.0, 0.0);
		    \draw[-] (axis cs:1.0, 0.0) -- (axis cs:1.0, 1.0);
		    \draw[-] (axis cs:1.0, 1.0) -- (axis cs:0.0, 1.0);
		    \draw[-] (axis cs:0.0, 1.0) -- (axis cs:0.0, 0.0);
		    \draw[-, dashed] (axis cs:0.5, 0.0) -- (axis cs:0.5, 1.0);
		    \draw[-, dashed] (axis cs:0.0, 0.5) -- (axis cs:1.0, 0.5);
		    
		    ;
		  
		    \node at (axis cs:0.25,0.25) 
		    {\begin{tabular}{c} $\rho_{bl}$\\ $u_{bl}$ \\ $v_{bl}$ \\ $ P_{bl}$ \end{tabular}};
		    \node at (axis cs:0.25,0.75) 
		    {\begin{tabular}{c} $\rho_{tl}$\\ $u_{tl}$ \\ $v_{tl}$ \\ $ P_{tl}$ \end{tabular}};
		    \node at (axis cs:0.75,0.25) 
		    {\begin{tabular}{c} $\rho_{br}$\\ $u_{br}$ \\ $v_{br}$ \\ $ P_{br}$ \end{tabular}};
		    \node at (axis cs:0.75,0.75) 
		    {\begin{tabular}{c} $\rho_{tr}$\\ $u_{tr}$ \\ $v_{tr}$ \\ $ P_{tr}$ \end{tabular}};
			    
		\end{axis}

	\end{tikzpicture}}
        \caption{\label{fig:2dRP_IC} Schematic of the initial conditions and domain for the 2D Riemann problem.}
    \end{figure}
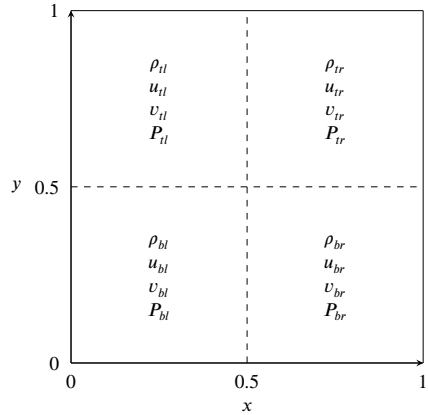

When relaxing the constraint of smoothness in the flow, a wide variety of problems with different flow physics become available for the Euler equations. To evaluate the efficacy of the approach for a sufficiently diverse set of flow physics, we consider the quintessential example of a two-dimensional Riemann problem \citep{SchulzRinne1993}, which possesses complex nonlinear flow behavior including shock waves, contact discontinuities, rarefaction waves, vortical flow, Kelvin--Helmholtz and Richtmyer--Meshkov instabilities, and shock-vortex interactions. This problem is initialized by subdividing the domain into four respective quadrants, shown in \cref{fig:2dRP_IC}, each with their own unique initial state. The partition between these states is representative of a physical diaphragm separating individual gases which is instantaneously broken, allowing for the interaction of these gases. 
Unless otherwise stated, we utilize a unit domain $\Omega = [0,1]^2$ with the diaphragms placed at $x_0 = y_0 = 0.5$, respectively. Furthermore, the domain is set as periodic, such that the numerical setup is representative of an infinite array of two-dimensional Riemann problems. While the enforcement of periodic boundary conditions is not standard for two-dimensional Riemann problems, this helps to increase the efficiency of the training process by increasing the overall complexity and variation of features encountered per episode (i.e., a single two-dimensional Riemann problem turns into four unique Riemann problems).
For brevity, we use the notation $tl$, $tr$, $bl$, and $br$ to denote the top-left, top-right, bottom-left, and bottom-right quadrants, respectively. 
    
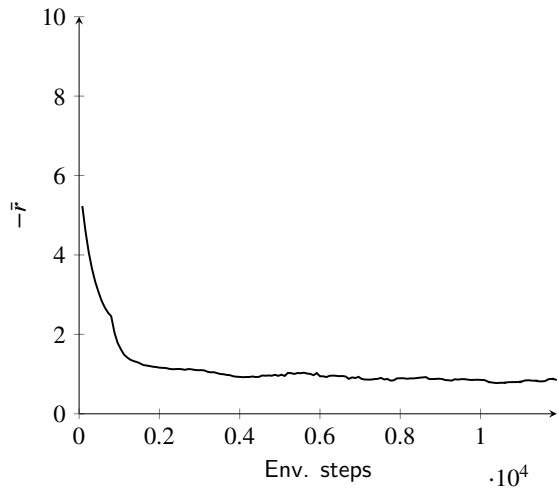
\begin{figure}[htbp!]
        \centering
        \adjustbox{width=0.45\linewidth,valign=b}{    \begin{tikzpicture}[spy using outlines={rectangle, height=3cm,width=2.3cm, magnification=3, connect spies}]
		\begin{axis}[name=plot1,
		    axis x line=left,
            axis y line=left,
    		xlabel={Env. steps},
		    ylabel={$-\bar{r}$},
    		xmin=0,
    		ymin=0,
    		ymax=10,
    		style={font=\normalsize}]
    		
			\addplot[color=black, style={thick}]
				table[x=num_env_steps_sampled,
                      y expr=-\thisrow{episode_reward_mean},
                      col sep=comma,unbounded coords=jump]{./figs/data/euler_h_training.csv};
    		    		
		\end{axis} 		
	\end{tikzpicture}}
        \caption{\label{fig:euler_href_reward} Batch-averaged (negative) reward with respect to number of environment steps during the training process for $h$-refinement on the Euler equations.}
    \end{figure}

The DynAMO policy was trained on a variety of two-dimensional Riemann problems generated by uniformly sampling the individual quadrant states from the distributions $\rho \in [0.2, 2]$, $u \in [-0.5,0.5]$, $v\in [-0.5,0.5]$, and $P \in [0.2,2]$. The remesh time was set as $T = 0.05$ with a total number of 4 RL steps, and the initial mesh distribution was set as $N = 32^2$. Furthermore, to add more variation to the flow problems during training, we vary the placement of the diaphragms within the ranges $x_0 \in [0.3, 0.7]$ and $y_0 \in [0.3, 0.7]$, but at evaluation time, these values were fixed at $x_0 = y_0 = 0.5$. The reward curve for the agent during training is shown in \cref{fig:euler_href_reward}.

\subsubsection{In-distribution experiments}
    \begin{figure}[htbp!] 
        \centering
        \begin{tabular}{cccc}
        \toprule
        Method &  Efficiency & Normalized error & Normalized cost\\ 
        \midrule
        Threshold ($\theta = 10^{-1}$) & 0.215 (0.083) & 0.763 (0.088) & 0.177 (0.049)\\
        Threshold ($\theta = 10^{-2}$) & \textbf{0.251 (0.056)} & \textbf{0.439 (0.111)} & \textbf{0.596 (0.061)}\\
        Threshold ($\theta = 10^{-3}$) & 0.180 (0.031) & 0.291 (0.099) & 0.760 (0.024)\\
        Threshold ($\theta = 10^{-4}$) & 0.162 (0.027) & 0.272 (0.093) & 0.787 (0.009)\\
        Threshold ($\theta = 10^{-5}$) & 0.157 (0.027) & 0.266 (0.090) & 0.795 (0.002)\\
        Threshold ($\theta = 10^{-6}$) & 0.157 (0.027) & 0.265 (0.090) & 0.796 (0.000)\\
        \midrule
        DynAMO & \textbf{0.486 (0.073)} & \textbf{0.142 (0.092)} & \textbf{0.482 (0.090)} \\
        DynAMO/Optimal $\theta$ & \textcolor{green!70!black}{+93.6\%} & \textcolor{green!70!black}{-67.7\%}  & \textcolor{green!70!black}{-19.1\%} \\
        \bottomrule
        \end{tabular}
        \captionof{table}{\label{tab:euler_href_indistribution} Comparison of the mean efficiency, normalized error, and normalized cost for $h$-refinement on the Euler equations with DynAMO and the threshold policy for the two-dimensional Riemann problem over 100 \textit{in-distribution} runs using uniform random initial conditions. Standard deviation shown in parentheses.}
    \end{figure}
    
An evaluation of the DynAMO policy with $\alpha = \alpha_{\text{train}}$ and a comparison to the baseline threshold policy at varying values of the threshold parameter was first performed across 100 uniformly sampled \textit{in-distribution} initial conditions. The normalized error, cost, and efficiency for the two approaches is shown in \cref{tab:euler_href_indistribution}. Due to the increased complexity of the flow physics, the optimal threshold policy only achieved a mean efficiency of $0.251$ with relatively large normalized error and cost values of $0.439$ and $0.596$, respectively. Furthermore, it can be seen that the threshold policy effectively stagnates around a minimum mean error of $0.265$ even with a decreasing threshold value. This behavior is as a result of the inability of instantaneous error estimators to account for the spatio-temporal evolution of the error. For this particular set of problems, this largely manifests at the first time step, where the estimator evaluated on the initially discontinuous state can, at best, predict a non-zero error only along (or one element adjacent to) the discontinuities. As such, the adapted mesh based on the instantaneous error estimator quickly becomes suboptimal, introducing a large degree of discretization error as the discontinuities evolve away from their initial positions. This effect is also encountered at later times as the discontinuities propagate across the mesh, although to a lesser extent. In contrast, the DynAMO approach showed a significantly better mean efficiency of $0.486$, a $93.6\%$ increase over the optimal threshold policy. This increased efficiency also came with a significantly lower mean error, $0.142$, a $67.7\%$ decrease in comparison to the optimal threshold policy. Most importantly, this error was much lower than the minimum error achievable with the threshold policy, showing that preemptive refinement decisions provided by the DynAMO policy can achieve levels of accuracy that standard AMR approaches may not be able to. Additionally, this increased accuracy came with a $19.1\%$ decrease in the computational cost as a result of the more efficient refinement decisions of the DynAMO policy.

\subsubsection{Sample problems}
To further validate the DynAMO approach, we consider canonical instances of two-dimensional Riemann problems that are used as example cases for evaluating numerical methods. In particular, we consider the six Riemann problems used in the work of \citet{Liska2003} which are derived from the work of \citet{Lax1998}. The initial conditions as well as the final simulation time for these six Riemann problems, denoted by case numbers 3, 4, 6, 12, 15, and 17, are shown in \cref{tab:euler_href_rpcases}. For these problems, the standard implementation requires Neumann boundary conditions to yield the expected flow patterns. As the current implementation of the DynAMO approach utilizes periodic meshes, we replicate the effects of the required boundary conditions by extending the domain by a factor of two to $\Omega = [0,2]^2$ and considering only the interior subdomain $\Omega' = [0.5,1.5]^2$ with the diaphragm placed at $x_0 = y_0 = 1$ . Due to the finite rate of propagation of flow perturbations as well as the self-similarity of Riemann problems with respect to the similarity parameter $x/t$, this approach yields identical results to applying Neumann boundary conditions on the subdomain $\Omega'$ over the time period of the problem. The characteristic mesh scale was kept identical by doubling the resolution to $N = 64^2$, and analysis of the cost and error was performed with respect to only the subdomain $\Omega'$. We remark here that these conditions are effectively \textit{out-of-distribution}, not only due to the increased domain size and mesh resolution but also due to the many of the quadrant states residing outside of the training distribution.

    \begin{figure}[htbp!] 
        \centering
        \begin{tabular}{c | llll | llll | l}
        \toprule
        Case & $\rho_l$ & $u_l$ & $v_l$ & $P_l$ & $\rho_r$ & $u_r$ & $v_r$ & $P_r$ & $t_f$\\ 
        \midrule
          & 0.5323 & 1.206 & 0.0   & 0.3   & 1.5    & 0.0 & 0.0   & 1.5 \\
        3 & 0.138  & 1.206 & 1.206 & 0.029 & 0.5323 & 0.0 & 1.206 & 0.3 & 0.3 \\
        \midrule
          & 0.5065 & 0.8939 & 0.0    & 0.35   & 1.1    & 0.0 & 0.0    & 1.1 \\
        4 & 1.1    & 0.8939 & 0.8939 & 1.1    & 0.5065 & 0.0 & 0.8939 & 0.35 & 0.25 \\
        \midrule
          & 2.0 &  0.75 & 0.5 & 1.0 & 1.0 &  0.75 & -0.5 & 1.0 \\
        6 & 1.0 & -0.75 & 0.5 & 1.0 & 3.0 & -0.75 & -0.5 & 1.0 & 0.3 \\
        \midrule
           & 1.0 & 0.7276 & 0.0 & 1.0 & 0.5313 & 0.0 & 0.0    & 0.4 \\
        12 & 0.8 & 0.0    & 0.0 & 1.0 & 1.0    & 0.0 & 0.7276 & 1.0 & 0.25 \\
        \midrule
           & 0.5197 & -0.6259 & -0.3 & 0.4 & 1.0    & 0.1 & -0.3   & 1.0 \\
        15 & 0.8    &  0.1    & -0.3 & 0.4 & 0.5313 & 0.1 & 0.4276 & 0.4 & 0.2 \\
        \midrule
           & 2.0    & 0.0 & -0.3   & 1.0 & 1.0    & 0.0 & -0.4    & 1.0 \\
        17 & 1.0625 & 0.0 & 0.2145 & 0.4 & 0.5197 & 0.0 & -1.1259 & 0.4 & 0.3 \\
        \bottomrule
        \end{tabular}
        \captionof{table}{\label{tab:euler_href_rpcases}Initial conditions for the two-dimensional Riemann problem tests from \citet{Liska2003}. Top and bottom quadrant states shown in upper and bottom rows, respectively.  }
    \end{figure}

\begin{figure}[htbp!]
        \centering
        \subfloat[Case 3]{
        \adjustbox{width=0.32\linewidth,valign=b}{\begin{tikzpicture}[spy using outlines={rectangle, height=3cm,width=2.3cm, magnification=3, connect spies}]

    \begin{axis}
    [   axis lines = none,
        xmin = 0, xmax = 1,
        ymin = 0, ymax = 1
    ]
        \fill[fill=blue!5] (1.0, 0) arc[start angle=0, end angle=90, radius=1.0] -- (0,0) -- (1.0, 0);
        \fill[fill=blue!10] (0.9, 0) arc[start angle=0, end angle=90, radius=0.9] -- (0,0) -- (0.9, 0);
        \fill[fill=blue!15] (0.8, 0) arc[start angle=0, end angle=90, radius=0.8] -- (0,0) -- (0.8, 0);
        \fill[fill=blue!20] (0.7, 0) arc[start angle=0, end angle=90, radius=0.7] -- (0,0) -- (0.7, 0);
        \fill[fill=blue!25] (0.6, 0) arc[start angle=0, end angle=90, radius=0.6] -- (0,0) -- (0.6, 0);
        \fill[fill=blue!30] (0.5, 0) arc[start angle=0, end angle=90, radius=0.5] -- (0,0) -- (0.5, 0);
        \fill[fill=blue!35] (0.4, 0) arc[start angle=0, end angle=90, radius=0.4] -- (0,0) -- (0.4, 0);
        \fill[fill=blue!40] (0.3, 0) arc[start angle=0, end angle=90, radius=0.3] -- (0,0) -- (0.3, 0);
        \fill[fill=blue!45] (0.2, 0) arc[start angle=0, end angle=90, radius=0.2] -- (0,0) -- (0.2, 0);
        \fill[fill=blue!50] (0.1, 0) arc[start angle=0, end angle=90, radius=0.1] -- (0,0) -- (0.1, 0);

    \end{axis}
        
    \begin{axis}
    [   axis line style={latex-latex},
        axis y line=left,
        axis x line=left,
        xmode=linear,
        ymode=linear,
        xlabel = {$\bar{c}$},
        ylabel = {$\bar{e}$},
        xmin = 0, xmax = 1,
        ymin = 0, ymax = 1,
        xtick = {0,0.2,0.4,0.6,0.8,1.0},
        ytick = {0,0.2,0.4,0.6,0.8,1.0},
        minor x tick num=1,
        minor y tick num=1,
        legend cell align={left},
        legend style={at={(0.97, 0.97)},anchor=north east},
        clip mode=individual,
        x tick label style={/pgf/number format/.cd, fixed, fixed zerofill, precision=1, /tikz/.cd},
        y tick label style={/pgf/number format/.cd, fixed, fixed zerofill, precision=1, /tikz/.cd},
        label style={font=\large},
    ]
        \addplot[color=black, fill=white, style={thick}, only marks, mark=*, mark options={scale=1.2}]
        table[x = cost_mean, y = err_mean, col sep=comma]{./figs/data/euler_href_RP3_threshold.csv};
        \addlegendentry{Threshold policy};
        
        \addplot[color=black, fill=red, style={thick}, only marks, mark=*, mark options={scale=1.2}]
        table[x = cost_mean, y = err_mean, col sep=comma]{./figs/data/euler_href_RP3_drl.csv};
        \addlegendentry{DynAMO};
        
    \end{axis}
\end{tikzpicture}}}
        \subfloat[Case 4]{
        \adjustbox{width=0.32\linewidth,valign=b}{\begin{tikzpicture}[spy using outlines={rectangle, height=3cm,width=2.3cm, magnification=3, connect spies}]

    \begin{axis}
    [   axis lines = none,
        xmin = 0, xmax = 1,
        ymin = 0, ymax = 1
    ]
        \fill[fill=blue!5] (1.0, 0) arc[start angle=0, end angle=90, radius=1.0] -- (0,0) -- (1.0, 0);
        \fill[fill=blue!10] (0.9, 0) arc[start angle=0, end angle=90, radius=0.9] -- (0,0) -- (0.9, 0);
        \fill[fill=blue!15] (0.8, 0) arc[start angle=0, end angle=90, radius=0.8] -- (0,0) -- (0.8, 0);
        \fill[fill=blue!20] (0.7, 0) arc[start angle=0, end angle=90, radius=0.7] -- (0,0) -- (0.7, 0);
        \fill[fill=blue!25] (0.6, 0) arc[start angle=0, end angle=90, radius=0.6] -- (0,0) -- (0.6, 0);
        \fill[fill=blue!30] (0.5, 0) arc[start angle=0, end angle=90, radius=0.5] -- (0,0) -- (0.5, 0);
        \fill[fill=blue!35] (0.4, 0) arc[start angle=0, end angle=90, radius=0.4] -- (0,0) -- (0.4, 0);
        \fill[fill=blue!40] (0.3, 0) arc[start angle=0, end angle=90, radius=0.3] -- (0,0) -- (0.3, 0);
        \fill[fill=blue!45] (0.2, 0) arc[start angle=0, end angle=90, radius=0.2] -- (0,0) -- (0.2, 0);
        \fill[fill=blue!50] (0.1, 0) arc[start angle=0, end angle=90, radius=0.1] -- (0,0) -- (0.1, 0);

    \end{axis}
        
    \begin{axis}
    [   axis line style={latex-latex},
        axis y line=left,
        axis x line=left,
        xmode=linear,
        ymode=linear,
        xlabel = {$\bar{c}$},
        ylabel = {$\bar{e}$},
        xmin = 0, xmax = 1,
        ymin = 0, ymax = 1,
        xtick = {0,0.2,0.4,0.6,0.8,1.0},
        ytick = {0,0.2,0.4,0.6,0.8,1.0},
        minor x tick num=1,
        minor y tick num=1,
        legend cell align={left},
        legend style={at={(0.97, 0.97)},anchor=north east},
        clip mode=individual,
        x tick label style={/pgf/number format/.cd, fixed, fixed zerofill, precision=1, /tikz/.cd},
        y tick label style={/pgf/number format/.cd, fixed, fixed zerofill, precision=1, /tikz/.cd},
        label style={font=\large},
    ]
        \addplot[color=black, fill=white, style={thick}, only marks, mark=*, mark options={scale=1.2}]
        table[x = cost_mean, y = err_mean, col sep=comma]{./figs/data/euler_href_RP4_threshold.csv};
        
        \addplot[color=black, fill=red, style={thick}, only marks, mark=*, mark options={scale=1.2}]
        table[x = cost_mean, y = err_mean, col sep=comma]{./figs/data/euler_href_RP4_drl.csv};
        
    \end{axis}
\end{tikzpicture}}}
        \subfloat[Case 6]{
        \adjustbox{width=0.32\linewidth,valign=b}{\begin{tikzpicture}[spy using outlines={rectangle, height=3cm,width=2.3cm, magnification=3, connect spies}]

    \begin{axis}
    [   axis lines = none,
        xmin = 0, xmax = 1,
        ymin = 0, ymax = 1
    ]
        \fill[fill=blue!5] (1.0, 0) arc[start angle=0, end angle=90, radius=1.0] -- (0,0) -- (1.0, 0);
        \fill[fill=blue!10] (0.9, 0) arc[start angle=0, end angle=90, radius=0.9] -- (0,0) -- (0.9, 0);
        \fill[fill=blue!15] (0.8, 0) arc[start angle=0, end angle=90, radius=0.8] -- (0,0) -- (0.8, 0);
        \fill[fill=blue!20] (0.7, 0) arc[start angle=0, end angle=90, radius=0.7] -- (0,0) -- (0.7, 0);
        \fill[fill=blue!25] (0.6, 0) arc[start angle=0, end angle=90, radius=0.6] -- (0,0) -- (0.6, 0);
        \fill[fill=blue!30] (0.5, 0) arc[start angle=0, end angle=90, radius=0.5] -- (0,0) -- (0.5, 0);
        \fill[fill=blue!35] (0.4, 0) arc[start angle=0, end angle=90, radius=0.4] -- (0,0) -- (0.4, 0);
        \fill[fill=blue!40] (0.3, 0) arc[start angle=0, end angle=90, radius=0.3] -- (0,0) -- (0.3, 0);
        \fill[fill=blue!45] (0.2, 0) arc[start angle=0, end angle=90, radius=0.2] -- (0,0) -- (0.2, 0);
        \fill[fill=blue!50] (0.1, 0) arc[start angle=0, end angle=90, radius=0.1] -- (0,0) -- (0.1, 0);

    \end{axis}
        
    \begin{axis}
    [   axis line style={latex-latex},
        axis y line=left,
        axis x line=left,
        xmode=linear,
        ymode=linear,
        xlabel = {$\bar{c}$},
        ylabel = {$\bar{e}$},
        xmin = 0, xmax = 1,
        ymin = 0, ymax = 1,
        xtick = {0,0.2,0.4,0.6,0.8,1.0},
        ytick = {0,0.2,0.4,0.6,0.8,1.0},
        minor x tick num=1,
        minor y tick num=1,
        legend cell align={left},
        legend style={at={(0.97, 0.97)},anchor=north east},
        clip mode=individual,
        x tick label style={/pgf/number format/.cd, fixed, fixed zerofill, precision=1, /tikz/.cd},
        y tick label style={/pgf/number format/.cd, fixed, fixed zerofill, precision=1, /tikz/.cd},
        label style={font=\large},
    ]
        \addplot[color=black, fill=white, style={thick}, only marks, mark=*, mark options={scale=1.2}]
        table[x = cost_mean, y = err_mean, col sep=comma]{./figs/data/euler_href_RP6_threshold.csv};
        
        \addplot[color=black, fill=red, style={thick}, only marks, mark=*, mark options={scale=1.2}]
        table[x = cost_mean, y = err_mean, col sep=comma]{./figs/data/euler_href_RP6_drl.csv};
        
    \end{axis}
\end{tikzpicture}}}
        \newline
        \subfloat[Case 12]{
        \adjustbox{width=0.32\linewidth,valign=b}{\begin{tikzpicture}[spy using outlines={rectangle, height=3cm,width=2.3cm, magnification=3, connect spies}]

    \begin{axis}
    [   axis lines = none,
        xmin = 0, xmax = 1,
        ymin = 0, ymax = 1
    ]
        \fill[fill=blue!5] (1.0, 0) arc[start angle=0, end angle=90, radius=1.0] -- (0,0) -- (1.0, 0);
        \fill[fill=blue!10] (0.9, 0) arc[start angle=0, end angle=90, radius=0.9] -- (0,0) -- (0.9, 0);
        \fill[fill=blue!15] (0.8, 0) arc[start angle=0, end angle=90, radius=0.8] -- (0,0) -- (0.8, 0);
        \fill[fill=blue!20] (0.7, 0) arc[start angle=0, end angle=90, radius=0.7] -- (0,0) -- (0.7, 0);
        \fill[fill=blue!25] (0.6, 0) arc[start angle=0, end angle=90, radius=0.6] -- (0,0) -- (0.6, 0);
        \fill[fill=blue!30] (0.5, 0) arc[start angle=0, end angle=90, radius=0.5] -- (0,0) -- (0.5, 0);
        \fill[fill=blue!35] (0.4, 0) arc[start angle=0, end angle=90, radius=0.4] -- (0,0) -- (0.4, 0);
        \fill[fill=blue!40] (0.3, 0) arc[start angle=0, end angle=90, radius=0.3] -- (0,0) -- (0.3, 0);
        \fill[fill=blue!45] (0.2, 0) arc[start angle=0, end angle=90, radius=0.2] -- (0,0) -- (0.2, 0);
        \fill[fill=blue!50] (0.1, 0) arc[start angle=0, end angle=90, radius=0.1] -- (0,0) -- (0.1, 0);

    \end{axis}
        
    \begin{axis}
    [   axis line style={latex-latex},
        axis y line=left,
        axis x line=left,
        xmode=linear,
        ymode=linear,
        xlabel = {$\bar{c}$},
        ylabel = {$\bar{e}$},
        xmin = 0, xmax = 1,
        ymin = 0, ymax = 1,
        xtick = {0,0.2,0.4,0.6,0.8,1.0},
        ytick = {0,0.2,0.4,0.6,0.8,1.0},
        minor x tick num=1,
        minor y tick num=1,
        legend cell align={left},
        legend style={at={(0.97, 0.97)},anchor=north east},
        clip mode=individual,
        x tick label style={/pgf/number format/.cd, fixed, fixed zerofill, precision=1, /tikz/.cd},
        y tick label style={/pgf/number format/.cd, fixed, fixed zerofill, precision=1, /tikz/.cd},
        label style={font=\large},
    ]
        \addplot[color=black, fill=white, style={thick}, only marks, mark=*, mark options={scale=1.2}]
        table[x = cost_mean, y = err_mean, col sep=comma]{./figs/data/euler_href_RP12_threshold.csv};
        
        \addplot[color=black, fill=red, style={thick}, only marks, mark=*, mark options={scale=1.2}]
        table[x = cost_mean, y = err_mean, col sep=comma]{./figs/data/euler_href_RP12_drl.csv};
        
    \end{axis}
\end{tikzpicture}}}
        \subfloat[Case 15]{
        \adjustbox{width=0.32\linewidth,valign=b}{\begin{tikzpicture}[spy using outlines={rectangle, height=3cm,width=2.3cm, magnification=3, connect spies}]

    \begin{axis}
    [   axis lines = none,
        xmin = 0, xmax = 1,
        ymin = 0, ymax = 1
    ]
        \fill[fill=blue!5] (1.0, 0) arc[start angle=0, end angle=90, radius=1.0] -- (0,0) -- (1.0, 0);
        \fill[fill=blue!10] (0.9, 0) arc[start angle=0, end angle=90, radius=0.9] -- (0,0) -- (0.9, 0);
        \fill[fill=blue!15] (0.8, 0) arc[start angle=0, end angle=90, radius=0.8] -- (0,0) -- (0.8, 0);
        \fill[fill=blue!20] (0.7, 0) arc[start angle=0, end angle=90, radius=0.7] -- (0,0) -- (0.7, 0);
        \fill[fill=blue!25] (0.6, 0) arc[start angle=0, end angle=90, radius=0.6] -- (0,0) -- (0.6, 0);
        \fill[fill=blue!30] (0.5, 0) arc[start angle=0, end angle=90, radius=0.5] -- (0,0) -- (0.5, 0);
        \fill[fill=blue!35] (0.4, 0) arc[start angle=0, end angle=90, radius=0.4] -- (0,0) -- (0.4, 0);
        \fill[fill=blue!40] (0.3, 0) arc[start angle=0, end angle=90, radius=0.3] -- (0,0) -- (0.3, 0);
        \fill[fill=blue!45] (0.2, 0) arc[start angle=0, end angle=90, radius=0.2] -- (0,0) -- (0.2, 0);
        \fill[fill=blue!50] (0.1, 0) arc[start angle=0, end angle=90, radius=0.1] -- (0,0) -- (0.1, 0);

    \end{axis}
        
    \begin{axis}
    [   axis line style={latex-latex},
        axis y line=left,
        axis x line=left,
        xmode=linear,
        ymode=linear,
        xlabel = {$\bar{c}$},
        ylabel = {$\bar{e}$},
        xmin = 0, xmax = 1,
        ymin = 0, ymax = 1,
        xtick = {0,0.2,0.4,0.6,0.8,1.0},
        ytick = {0,0.2,0.4,0.6,0.8,1.0},
        minor x tick num=1,
        minor y tick num=1,
        legend cell align={left},
        legend style={at={(0.97, 0.97)},anchor=north east},
        clip mode=individual,
        x tick label style={/pgf/number format/.cd, fixed, fixed zerofill, precision=1, /tikz/.cd},
        y tick label style={/pgf/number format/.cd, fixed, fixed zerofill, precision=1, /tikz/.cd},
        label style={font=\large},
    ]
        \addplot[color=black, fill=white, style={thick}, only marks, mark=*, mark options={scale=1.2}]
        table[x = cost_mean, y = err_mean, col sep=comma]{./figs/data/euler_href_RP15_threshold.csv};
        
        \addplot[color=black, fill=red, style={thick}, only marks, mark=*, mark options={scale=1.2}]
        table[x = cost_mean, y = err_mean, col sep=comma]{./figs/data/euler_href_RP15_drl.csv};
        
    \end{axis}
\end{tikzpicture}}}
        \subfloat[Case 17]{
        \adjustbox{width=0.32\linewidth,valign=b}{\begin{tikzpicture}[spy using outlines={rectangle, height=3cm,width=2.3cm, magnification=3, connect spies}]

    \begin{axis}
    [   axis lines = none,
        xmin = 0, xmax = 1,
        ymin = 0, ymax = 1
    ]
        \fill[fill=blue!5] (1.0, 0) arc[start angle=0, end angle=90, radius=1.0] -- (0,0) -- (1.0, 0);
        \fill[fill=blue!10] (0.9, 0) arc[start angle=0, end angle=90, radius=0.9] -- (0,0) -- (0.9, 0);
        \fill[fill=blue!15] (0.8, 0) arc[start angle=0, end angle=90, radius=0.8] -- (0,0) -- (0.8, 0);
        \fill[fill=blue!20] (0.7, 0) arc[start angle=0, end angle=90, radius=0.7] -- (0,0) -- (0.7, 0);
        \fill[fill=blue!25] (0.6, 0) arc[start angle=0, end angle=90, radius=0.6] -- (0,0) -- (0.6, 0);
        \fill[fill=blue!30] (0.5, 0) arc[start angle=0, end angle=90, radius=0.5] -- (0,0) -- (0.5, 0);
        \fill[fill=blue!35] (0.4, 0) arc[start angle=0, end angle=90, radius=0.4] -- (0,0) -- (0.4, 0);
        \fill[fill=blue!40] (0.3, 0) arc[start angle=0, end angle=90, radius=0.3] -- (0,0) -- (0.3, 0);
        \fill[fill=blue!45] (0.2, 0) arc[start angle=0, end angle=90, radius=0.2] -- (0,0) -- (0.2, 0);
        \fill[fill=blue!50] (0.1, 0) arc[start angle=0, end angle=90, radius=0.1] -- (0,0) -- (0.1, 0);

    \end{axis}
        
    \begin{axis}
    [   axis line style={latex-latex},
        axis y line=left,
        axis x line=left,
        xmode=linear,
        ymode=linear,
        xlabel = {$\bar{c}$},
        ylabel = {$\bar{e}$},
        xmin = 0, xmax = 1,
        ymin = 0, ymax = 1,
        xtick = {0,0.2,0.4,0.6,0.8,1.0},
        ytick = {0,0.2,0.4,0.6,0.8,1.0},
        minor x tick num=1,
        minor y tick num=1,
        legend cell align={left},
        legend style={at={(0.97, 0.97)},anchor=north east},
        clip mode=individual,
        x tick label style={/pgf/number format/.cd, fixed, fixed zerofill, precision=1, /tikz/.cd},
        y tick label style={/pgf/number format/.cd, fixed, fixed zerofill, precision=1, /tikz/.cd},
        label style={font=\large},
    ]
        \addplot[color=black, fill=white, style={thick}, only marks, mark=*, mark options={scale=1.2}]
        table[x = cost_mean, y = err_mean, col sep=comma]{./figs/data/euler_href_RP17_threshold.csv};
        
        \addplot[color=black, fill=red, style={thick}, only marks, mark=*, mark options={scale=1.2}]
        table[x = cost_mean, y = err_mean, col sep=comma]{./figs/data/euler_href_RP17_drl.csv};
        
    \end{axis}
\end{tikzpicture}}}
        \newline
        \caption{\label{fig:euler_href_pareto_rp}Pareto plot of normalized cost vs. error for $h$-refinement on the Euler equations with DynAMO and the threshold policy for the two-dimensional Riemann problem tests from \citet{Liska2003}. Contours of efficiency shown on background. }
    \end{figure}
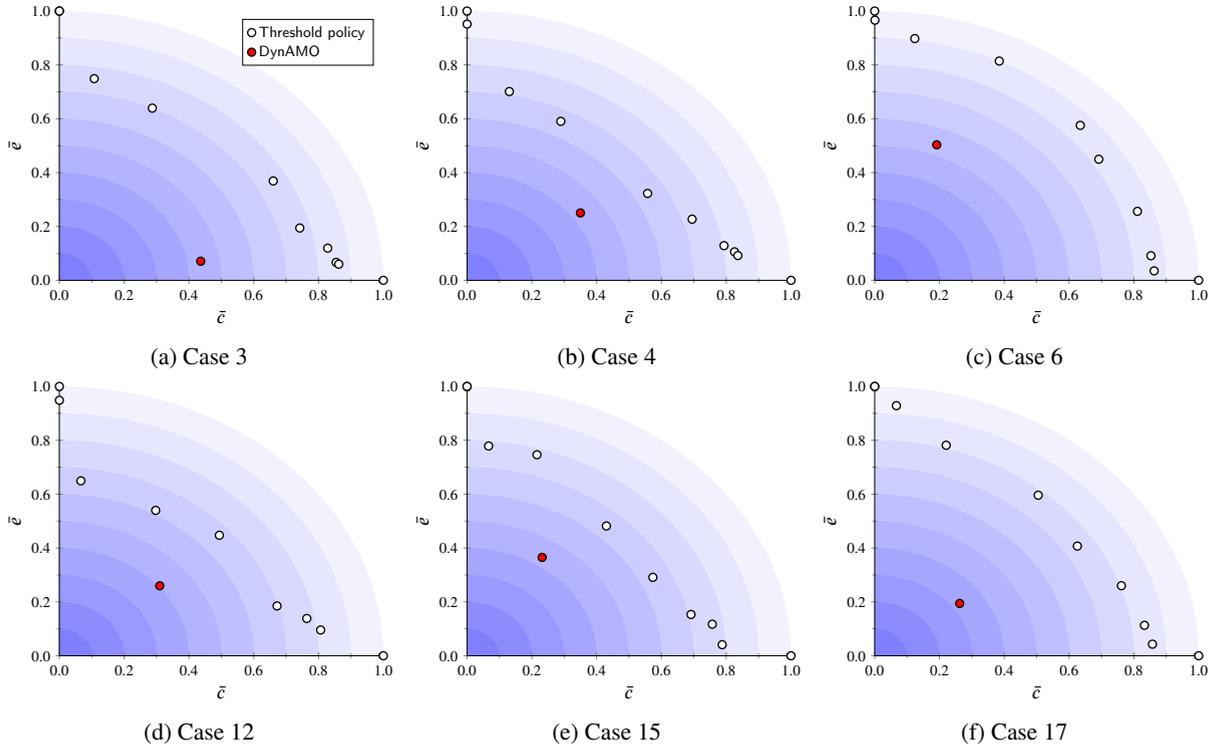

A comparison of the normalized cost vs. error between the threshold policy and the DynAMO approach for these six two-dimensional Riemann problem cases is shown through Pareto plots in \cref{fig:euler_href_pareto_rp}. Due to the previously mentioned inability of instantaneous error-based policies to appropriately refine the mesh at the initial state, the Pareto curve for the threshold policy stagnates at above a fixed error and does not show as well-behaved of a curve as with the smooth problems for $p$-refinement. As a result, the threshold policies did not achieve an efficiency of more than $0.4$ for any of the cases. Of the six cases, the threshold policy showed the worst performance for Case 6 and the best performance for Case 4 and Case 12. In comparison, the DynAMO policy showed significantly better performance, generally achieving an efficiency of approximately $0.6$ with minor variation across the cases. This increase in efficiency resulted in the policy achieving a similar error to the threshold policy with significantly lower computational cost, owing to its ability to preemptively refine mesh regions prior to features appearing. These results provide an initial quantitative evaluation of the proposed approach for providing more efficient adaptive mesh refinement strategies for complex problems of practical interest.

\begin{figure}[htbp!]
        \centering
        \subfloat[Remesh at $t = 0$]{
        \adjustbox{width=0.24\linewidth,valign=b}{\includegraphics{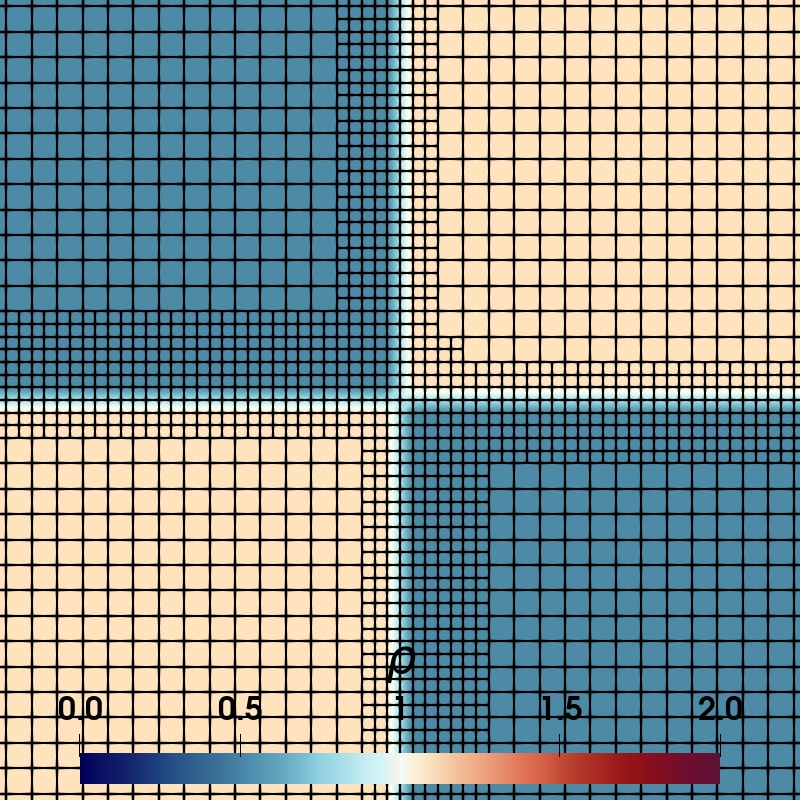}}}
        \subfloat[Solution at $t = T$]{
        \adjustbox{width=0.24\linewidth,valign=b}{\includegraphics{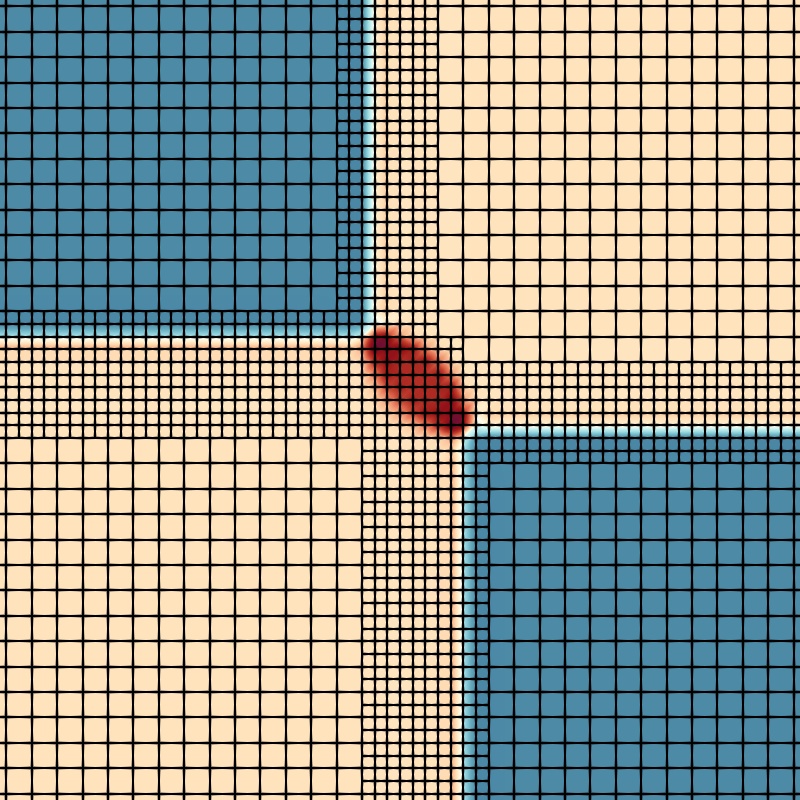}}}
        \subfloat[Remesh at $t = T$]{
        \adjustbox{width=0.24\linewidth,valign=b}{\includegraphics{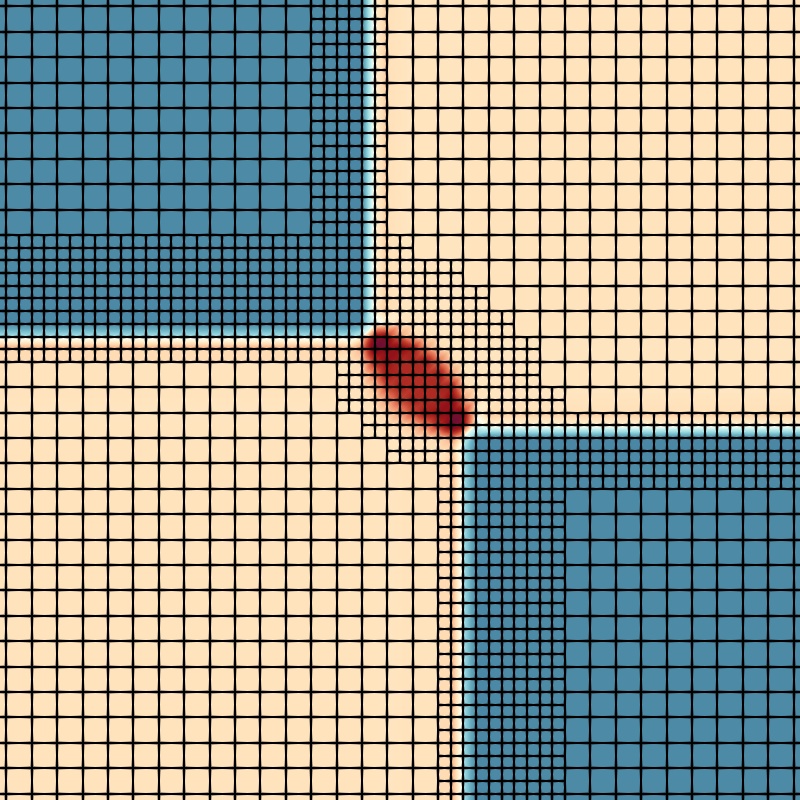}}}
        \subfloat[Solution at $t = 2T$]{
        \adjustbox{width=0.24\linewidth,valign=b}{\includegraphics{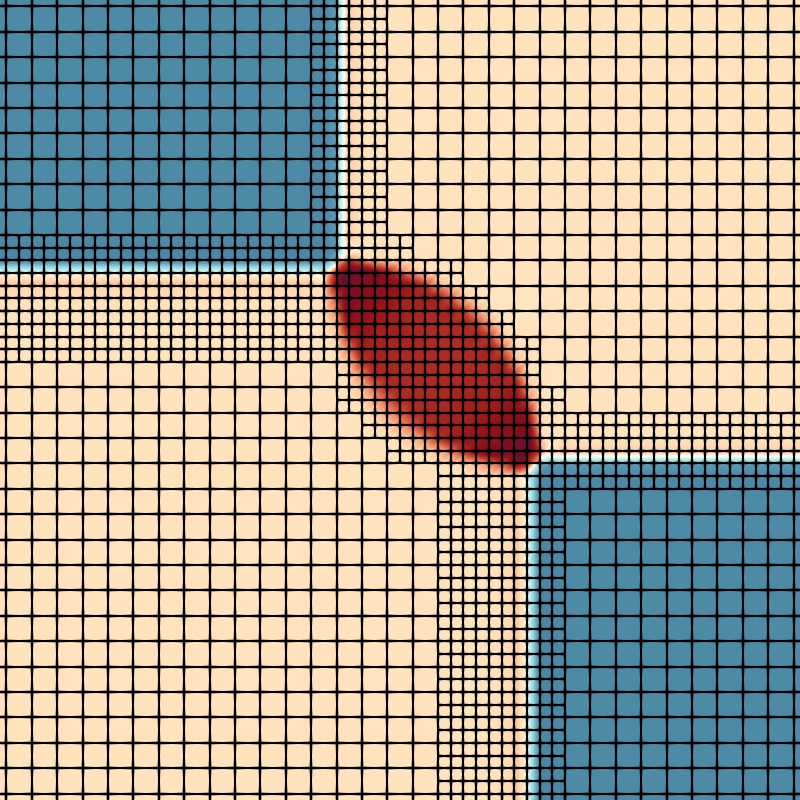}}}
        \newline
        \subfloat[Remesh at $t = 2T$]{
        \adjustbox{width=0.24\linewidth,valign=b}{\includegraphics{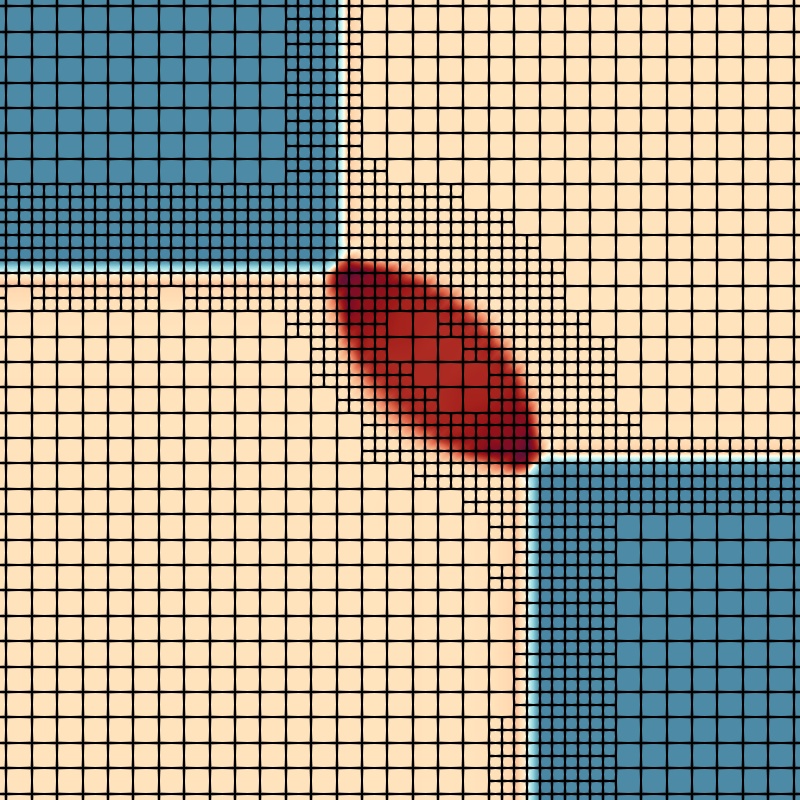}}}
        \subfloat[Solution at $t = 3T$]{
        \adjustbox{width=0.24\linewidth,valign=b}{\includegraphics{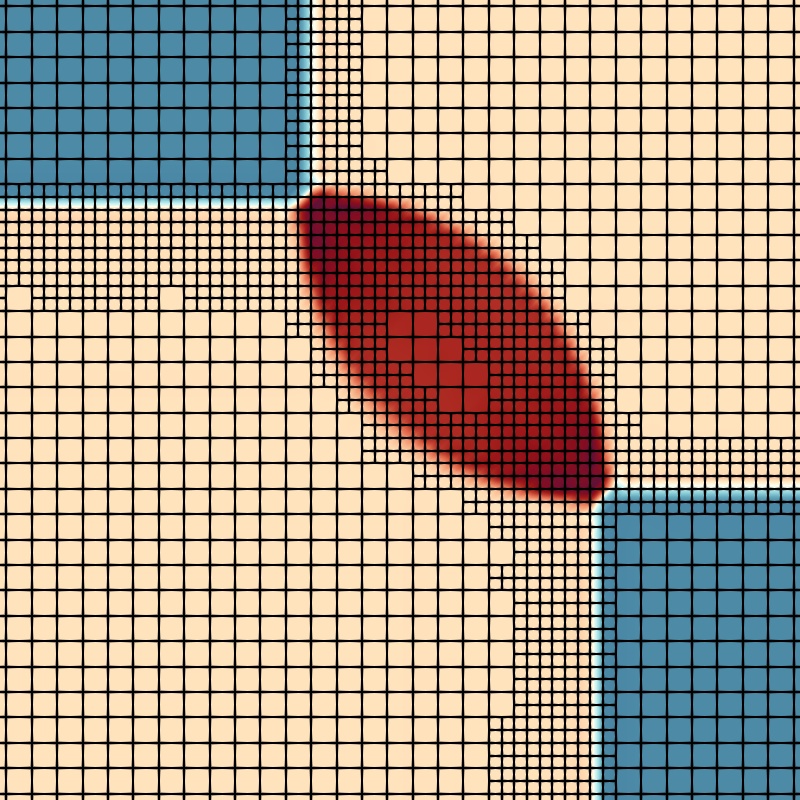}}}
        \subfloat[Remesh at $t = 3T$]{
        \adjustbox{width=0.24\linewidth,valign=b}{\includegraphics{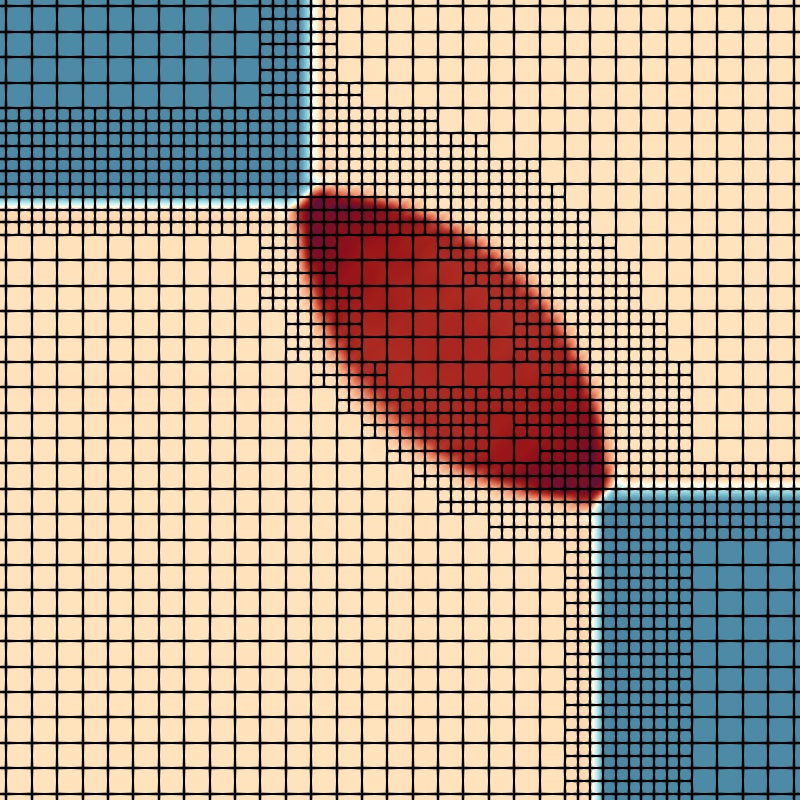}}}
        \subfloat[Solution at $t = 4T$]{
        \adjustbox{width=0.24\linewidth,valign=b}{\includegraphics{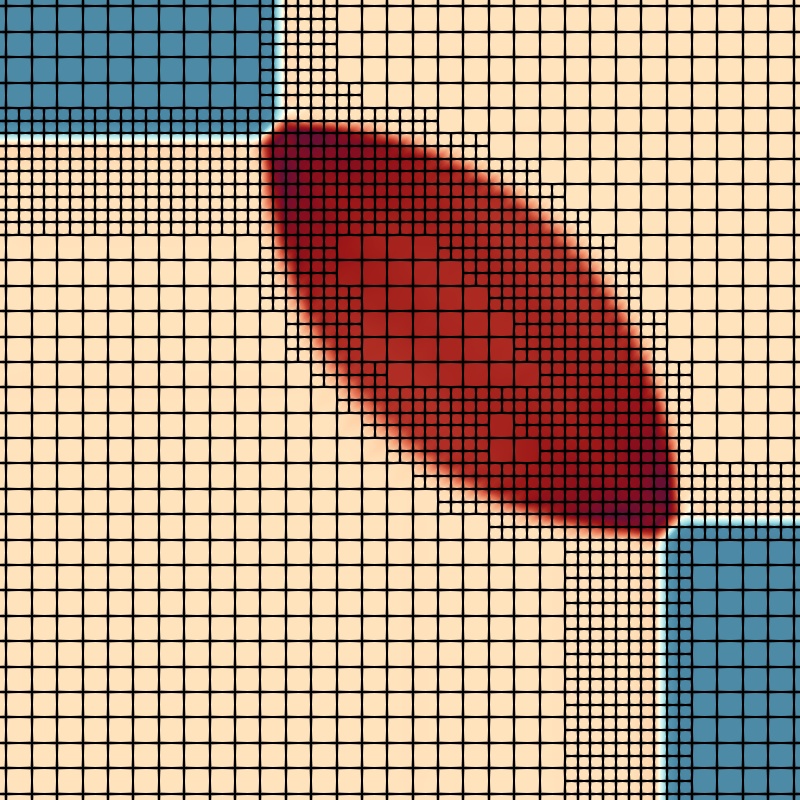}}}
        \newline
        \caption{\label{fig:rp4_amr_drl} Contours of density overlaid with $h$-adapted mesh at varying remesh intervals using DynAMO for the Case 4 two-dimensional Riemann problem from \citet{Liska2003}.}
    \end{figure}

\begin{figure}[htbp!]
        \centering
        \subfloat[Remesh at $t = 0$]{
        \adjustbox{width=0.24\linewidth,valign=b}{\includegraphics{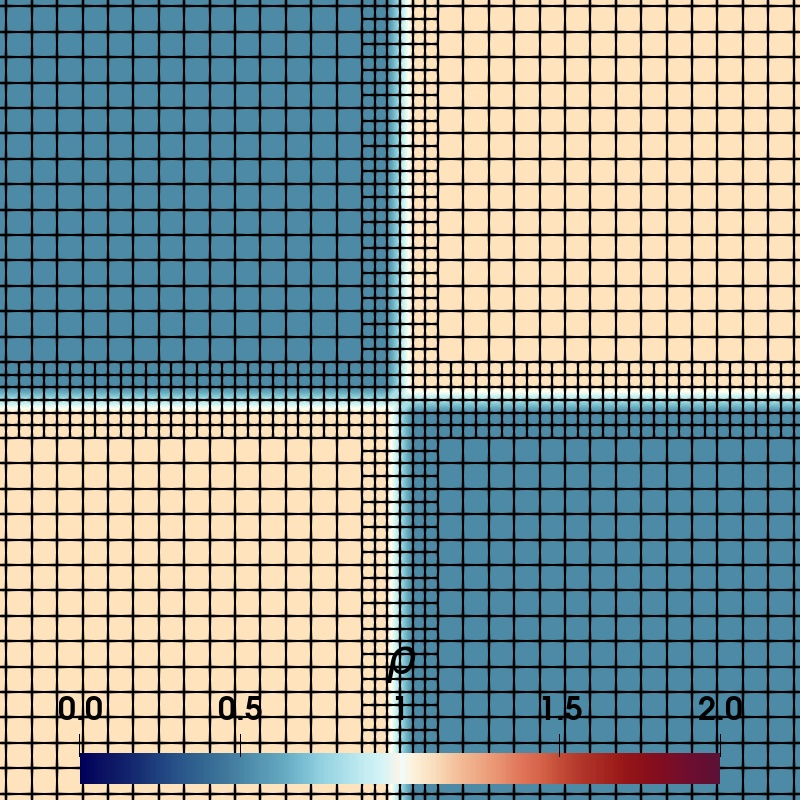}}}
        \subfloat[Solution at $t = T$]{
        \adjustbox{width=0.24\linewidth,valign=b}{\includegraphics{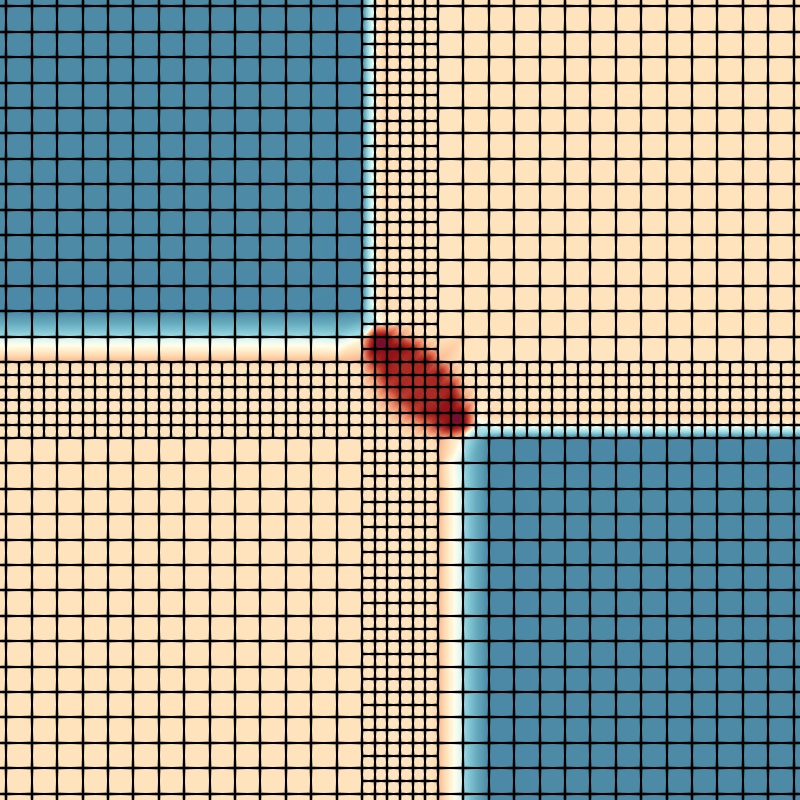}}}
        \subfloat[Remesh at $t = T$]{
        \adjustbox{width=0.24\linewidth,valign=b}{\includegraphics{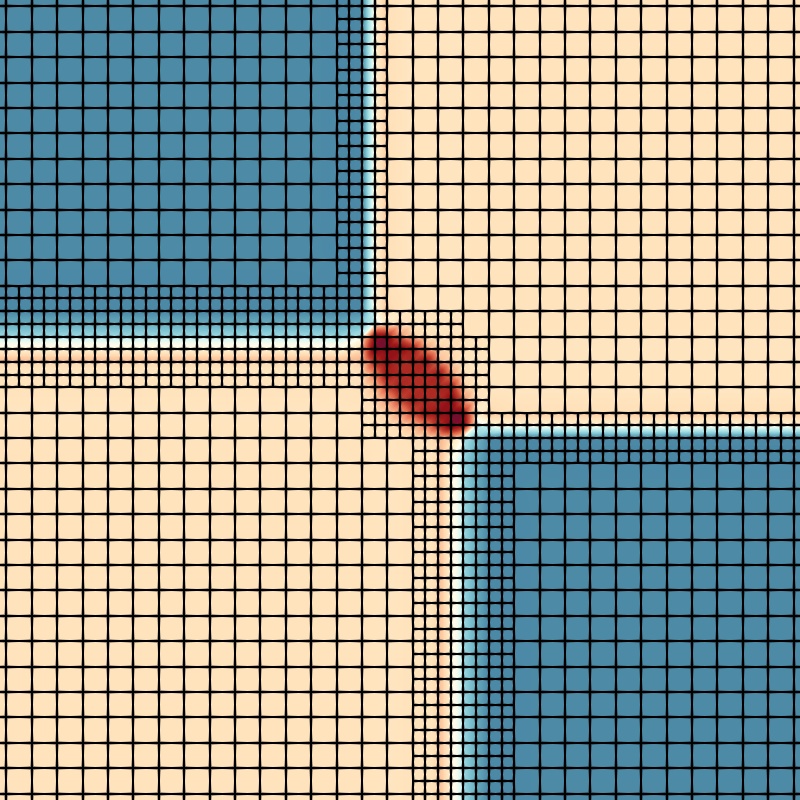}}}
        \subfloat[Solution at $t = 2T$]{
        \adjustbox{width=0.24\linewidth,valign=b}{\includegraphics{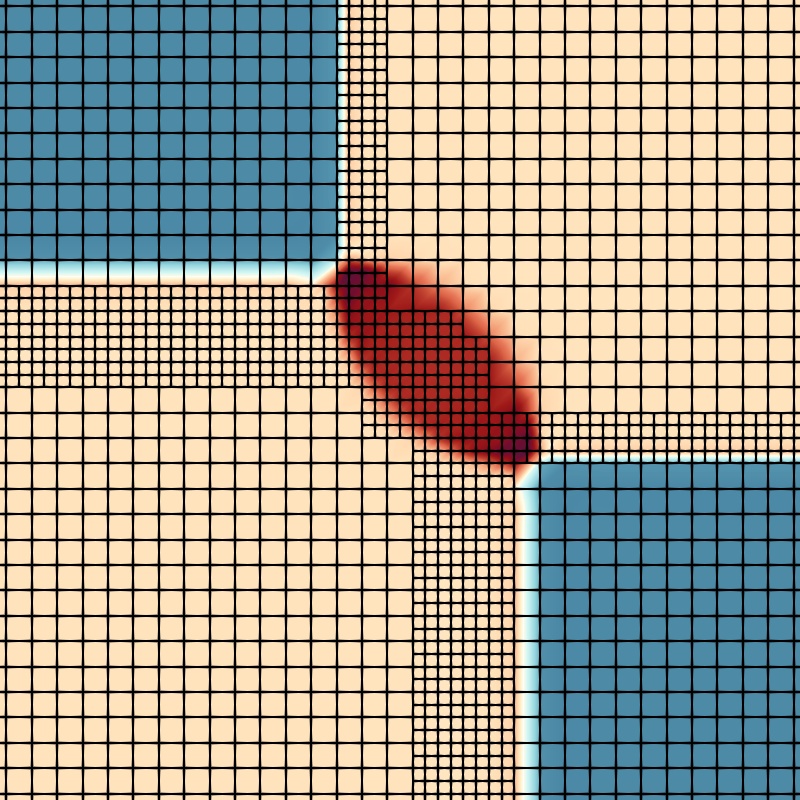}}}
        \newline
        \subfloat[Remesh at $t = 2T$]{
        \adjustbox{width=0.24\linewidth,valign=b}{\includegraphics{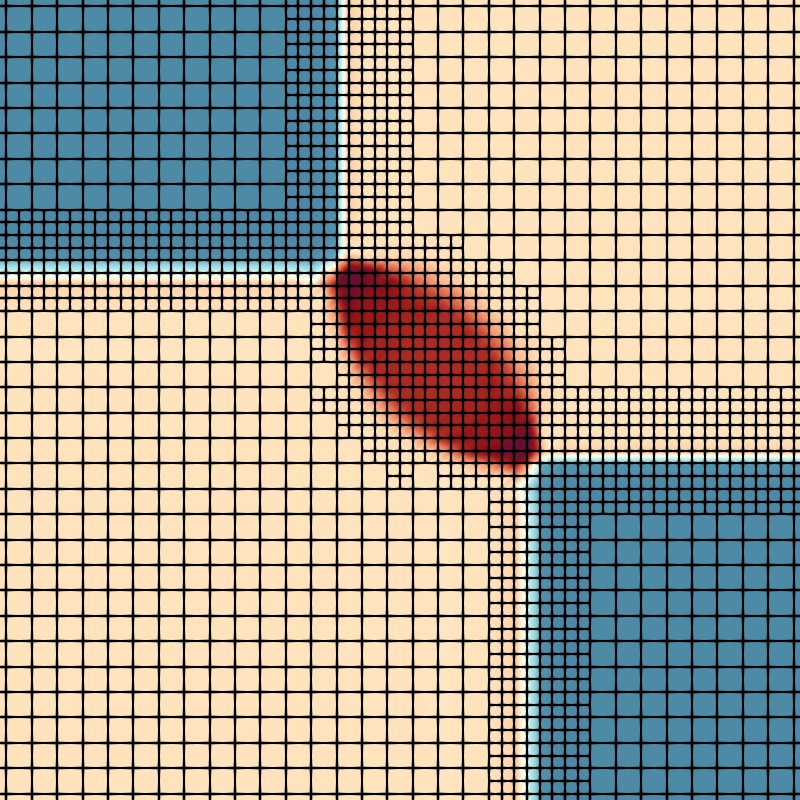}}}
        \subfloat[Solution at $t = 3T$]{
        \adjustbox{width=0.24\linewidth,valign=b}{\includegraphics{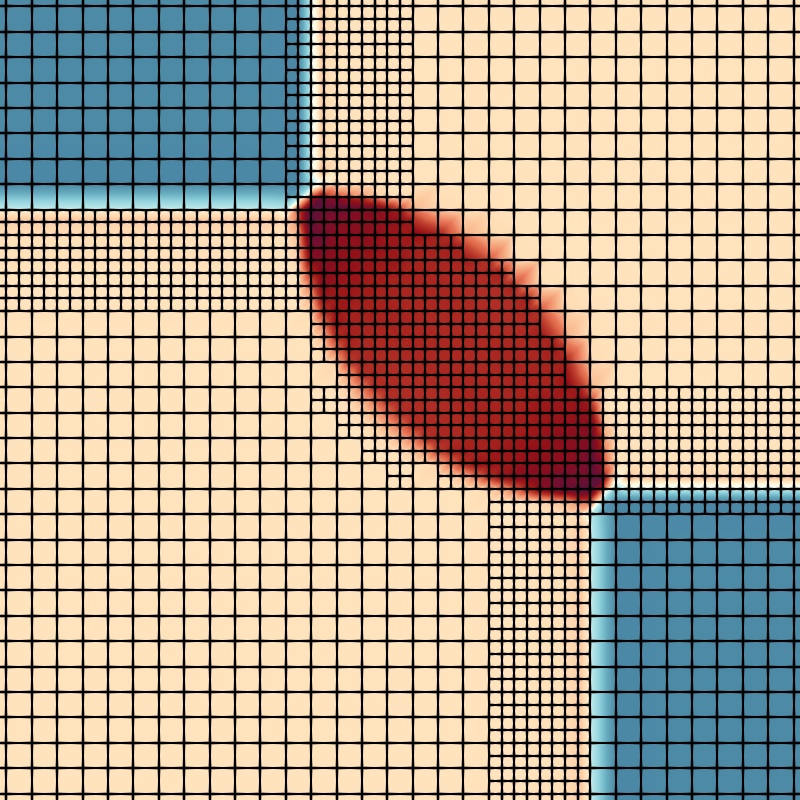}}}
        \subfloat[Remesh at $t = 3T$]{
        \adjustbox{width=0.24\linewidth,valign=b}{\includegraphics{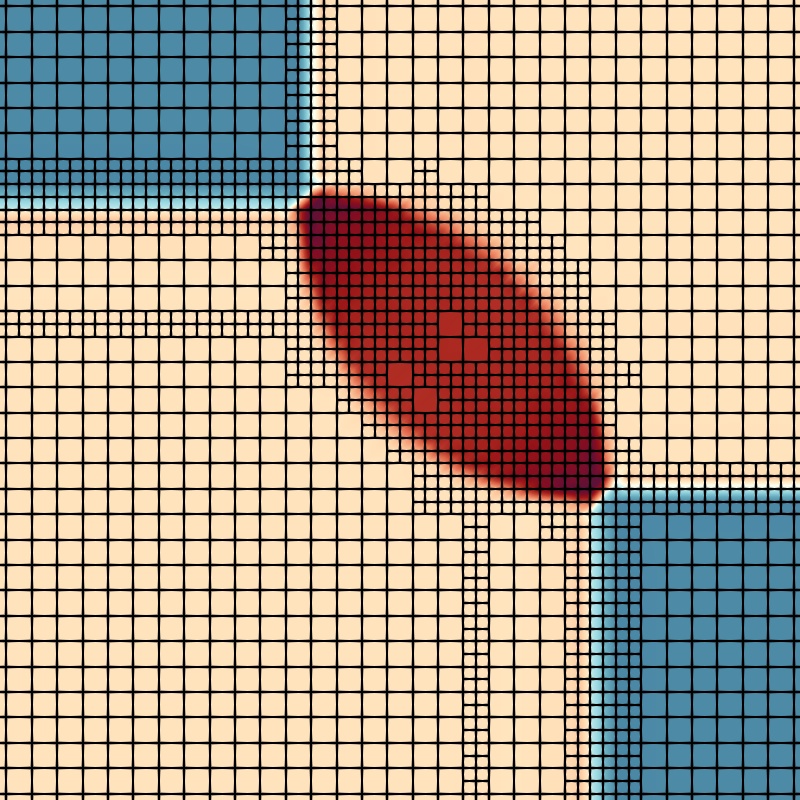}}}
        \subfloat[Solution at $t = 4T$]{
        \adjustbox{width=0.24\linewidth,valign=b}{\includegraphics{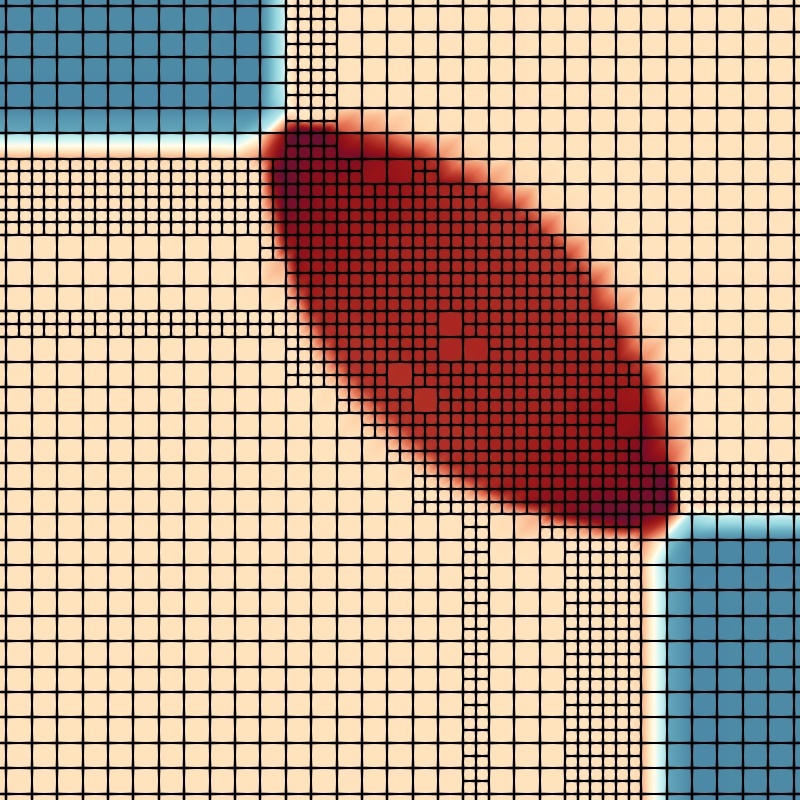}}}
        \newline
        \caption{\label{fig:rp4_amr_threshold} Contours of density overlaid with $h$-adapted mesh at varying remesh intervals using the threshold policy ($\theta = 10^{-2}$) for the Case 4 two-dimensional Riemann problem from \citet{Liska2003}.} 
    \end{figure}

The DynAMO policy and threshold policy was further qualitatively compared for these six two-dimensional Riemann problems by observing the mesh refinement decisions. A comparison was made for Case 4 and Case 12 as these are the cases for which the threshold policy showed the best performance. The threshold parameter value was set to $\theta = 10^{-2}$ which corresponds to the points of highest efficiency along the Pareto curve for both cases. The evolution of the density and $h$-adapted mesh for Case 4 as obtained by the DynAMO policy and threshold policy is shown in \cref{fig:rp4_amr_drl} and \cref{fig:rp4_amr_threshold}, respectively. Note that the density is shown instead of the total energy, the observed quantity, as this is the typical component used for comparison in the literature. It can be seen that at the initial time, the DynAMO policy refines several elements ahead of the discontinuities. This refinement region covered the entire expansion wave as well as the propagation of the discontinuities across the remesh interval. This behavior persisted at further times, showcasing the preemptive refinement capabilities of DynAMO. Furthermore, as the expansion wave region became sufficiently large, the DynAMO policy started to de-refine its interior which contains features with relatively low variation and, by extent, error. In comparison, the threshold policy could only refine the elements along the discontinuities and their immediate neighbors at the initial time. Therefore, by the next remesh time, the discontinuities propagated outside of the refinement region. This behavior became more egregious at the next remesh interval, where both the discontinuities as well as the expansion wave front moved outside of the refinement region, introducing a significant degree of discretization error. This error resulted in some irregular refinement patterns at later times which contributed to the cost of threshold policy without any benefits in the accuracy. 
    
    \begin{figure}[htbp!]
        \centering
        \subfloat[Remesh at $t = 0$]{
        \adjustbox{width=0.24\linewidth,valign=b}{\includegraphics{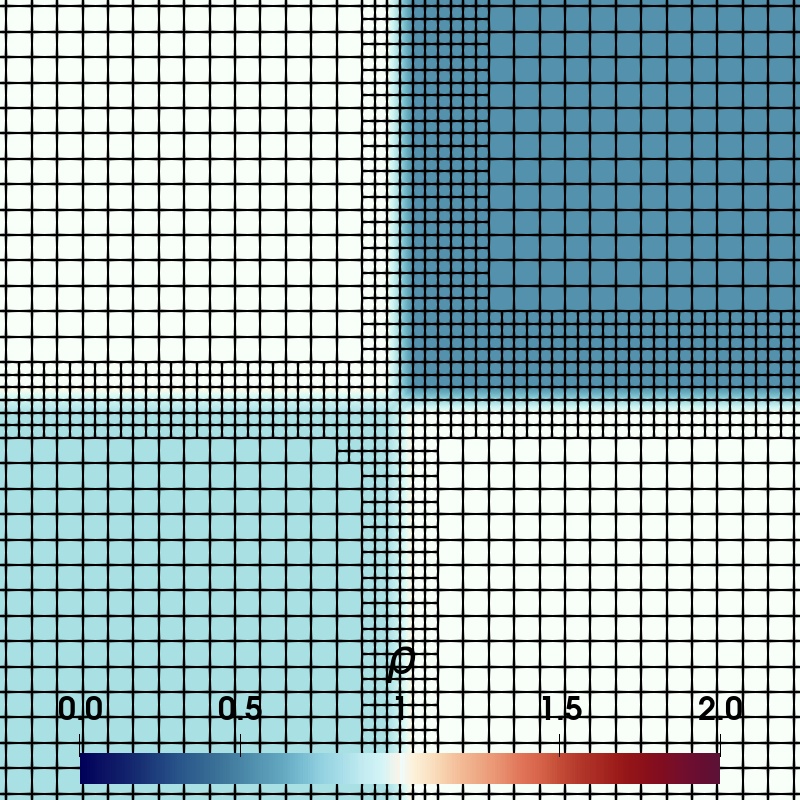}}}
        \subfloat[Solution at $t = T$]{
        \adjustbox{width=0.24\linewidth,valign=b}{\includegraphics{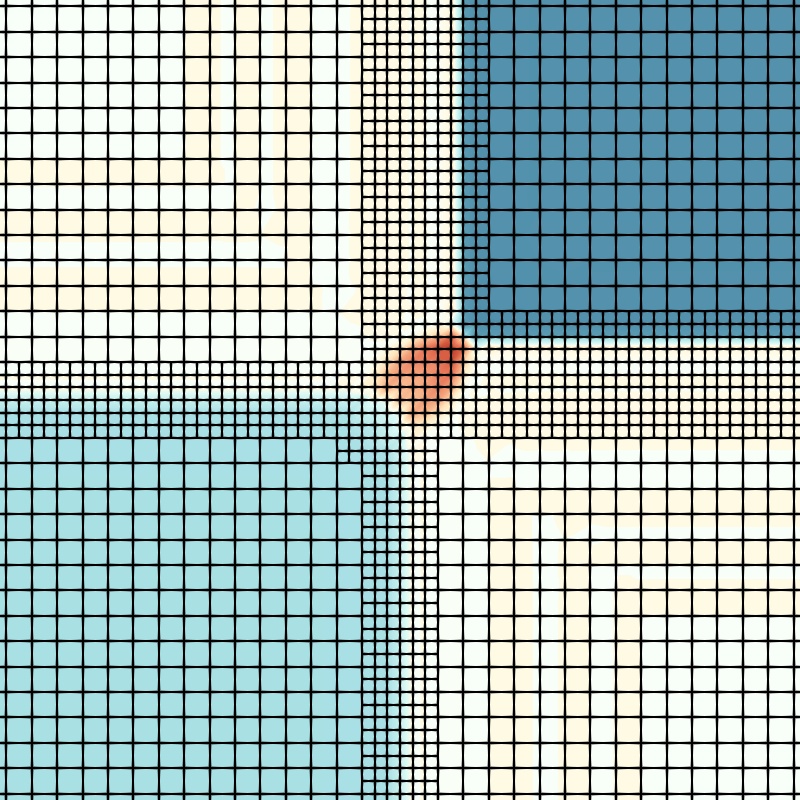}}}
        \subfloat[Remesh at $t = T$]{
        \adjustbox{width=0.24\linewidth,valign=b}{\includegraphics{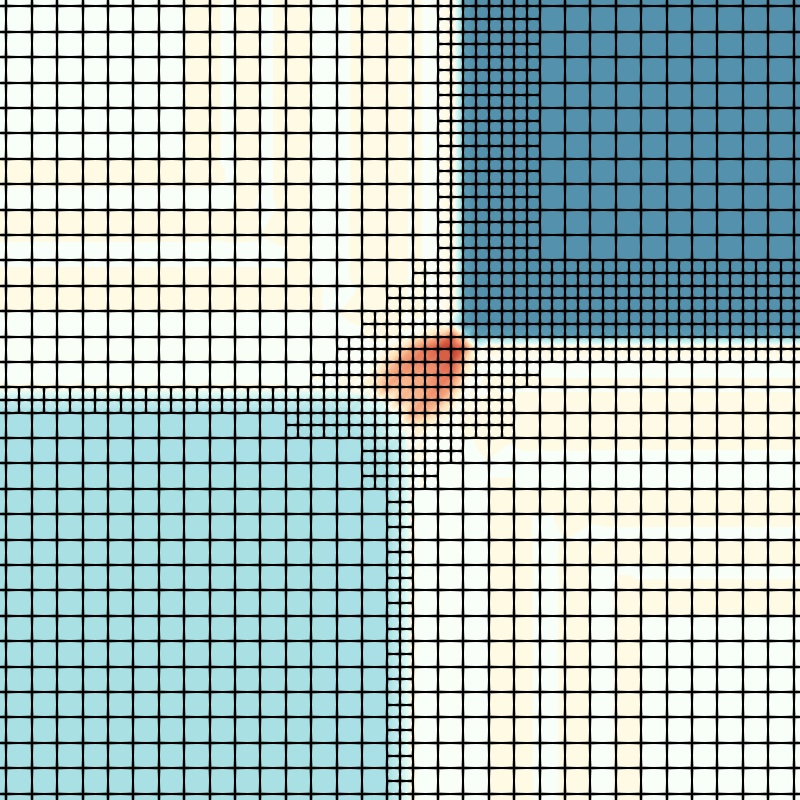}}}
        \subfloat[Solution at $t = 2T$]{
        \adjustbox{width=0.24\linewidth,valign=b}{\includegraphics{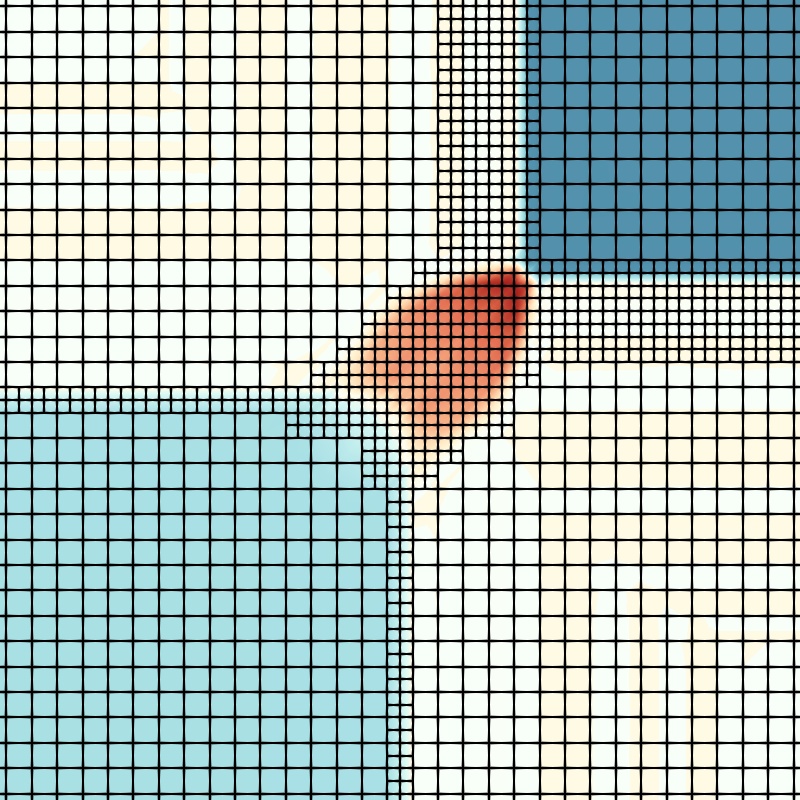}}}
        \newline
        \subfloat[Remesh at $t = 2T$]{
        \adjustbox{width=0.24\linewidth,valign=b}{\includegraphics{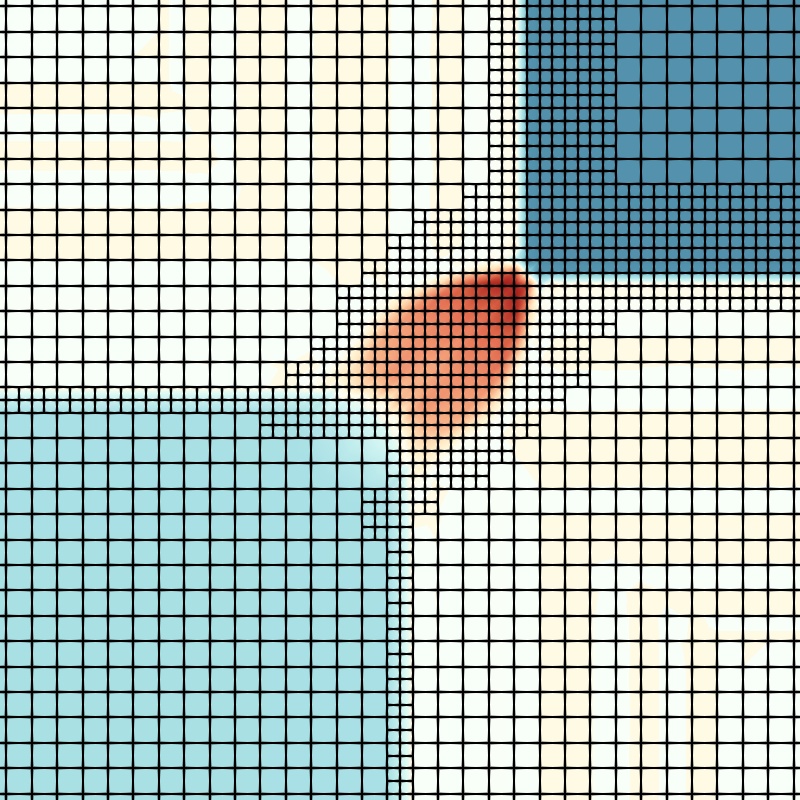}}}
        \subfloat[Solution at $t = 3T$]{
        \adjustbox{width=0.24\linewidth,valign=b}{\includegraphics{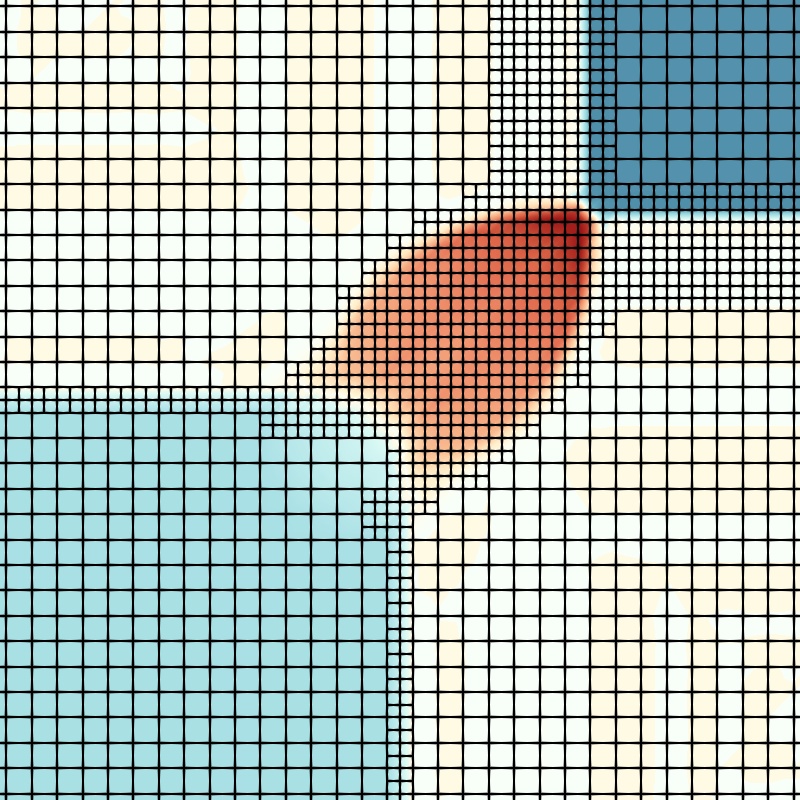}}}
        \subfloat[Remesh at $t = 3T$]{
        \adjustbox{width=0.24\linewidth,valign=b}{\includegraphics{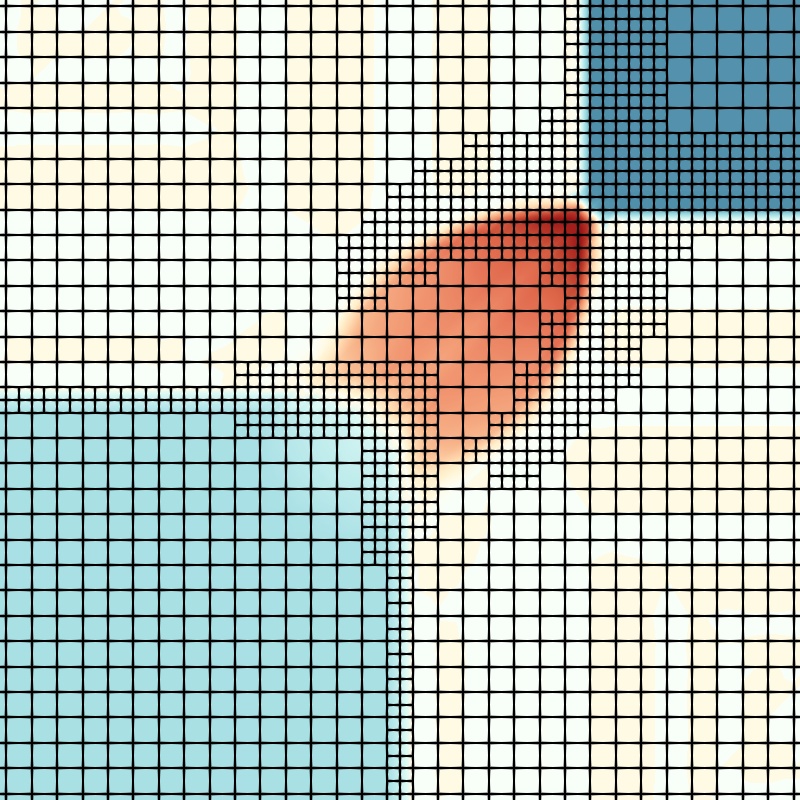}}}
        \subfloat[Solution at $t = 4T$]{
        \adjustbox{width=0.24\linewidth,valign=b}{\includegraphics{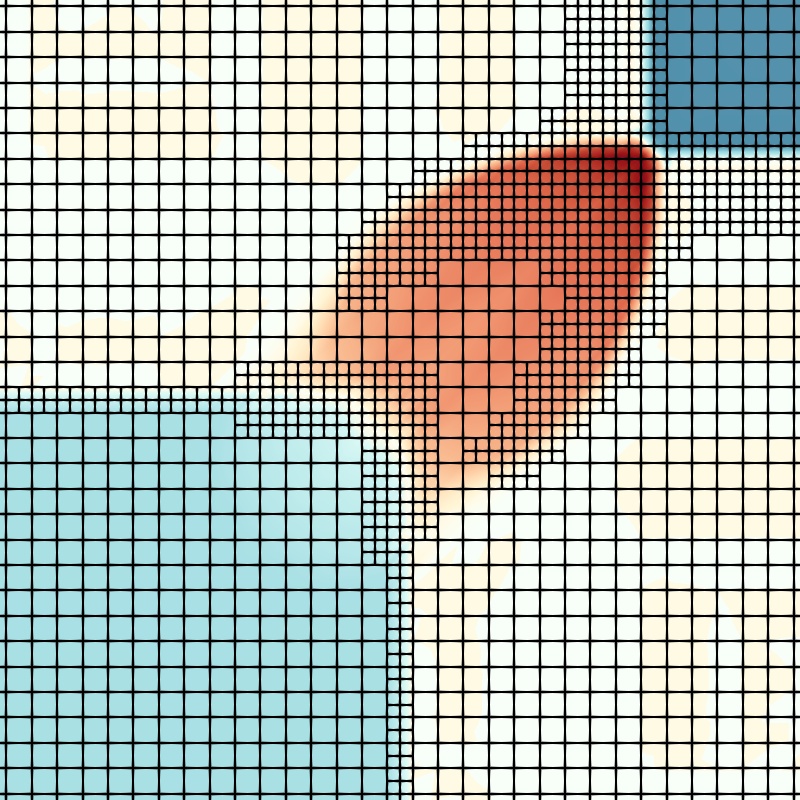}}}
        \newline
        \caption{\label{fig:rp12_amr_drl} Contours of density overlaid with $h$-adapted mesh at varying remesh intervals using DynAMO for the Case 12 two-dimensional Riemann problem from \citet{Liska2003}.}
    \end{figure}
    \begin{figure}[htbp!]
        \centering
        \subfloat[Remesh at $t = 0$]{
        \adjustbox{width=0.24\linewidth,valign=b}{\includegraphics{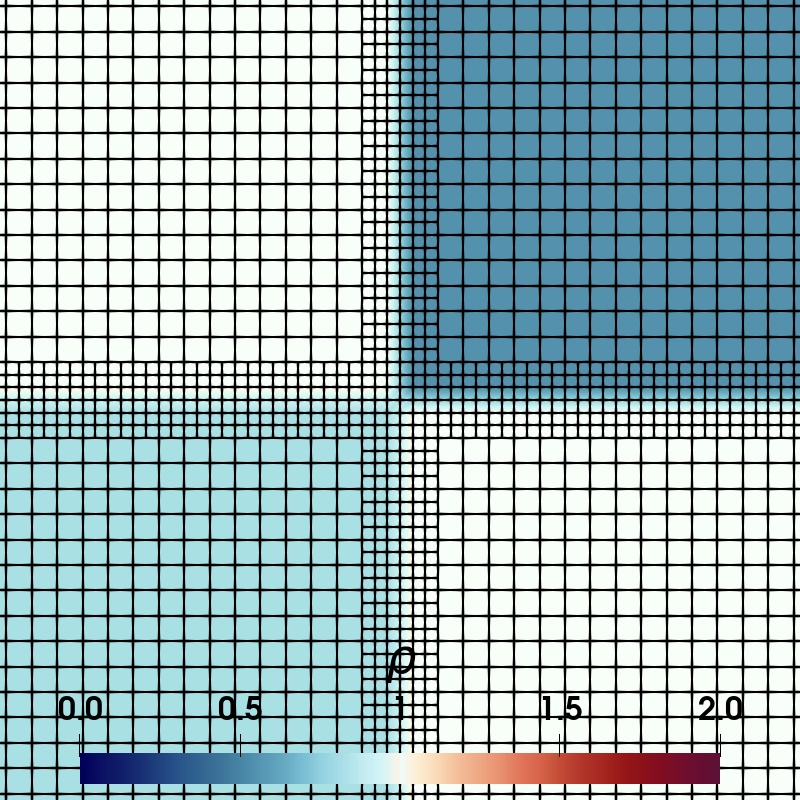}}}
        \subfloat[Solution at $t = T$]{
        \adjustbox{width=0.24\linewidth,valign=b}{\includegraphics{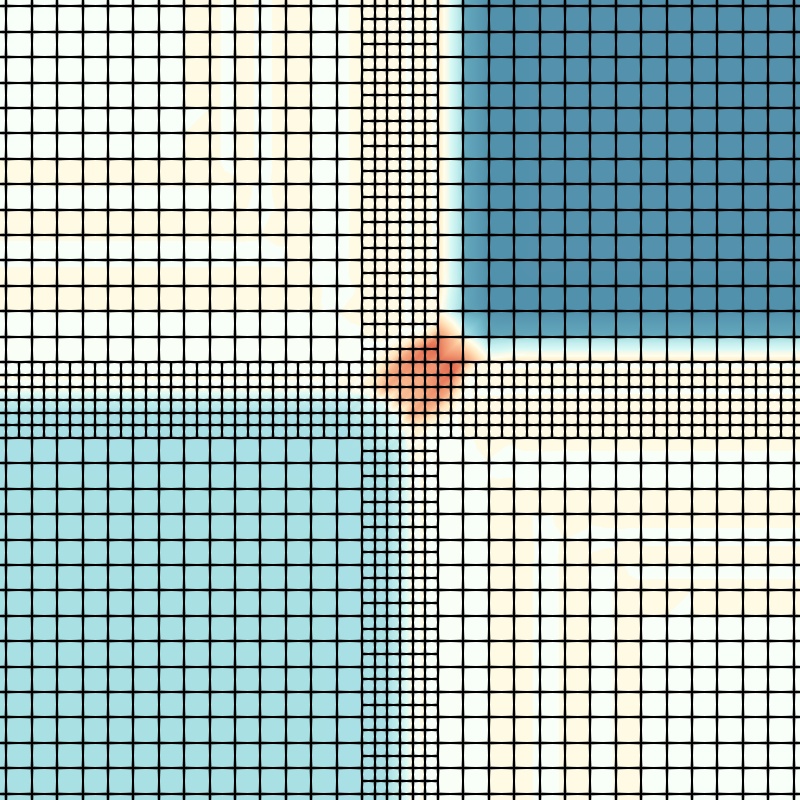}}}
        \subfloat[Remesh at $t = T$]{
        \adjustbox{width=0.24\linewidth,valign=b}{\includegraphics{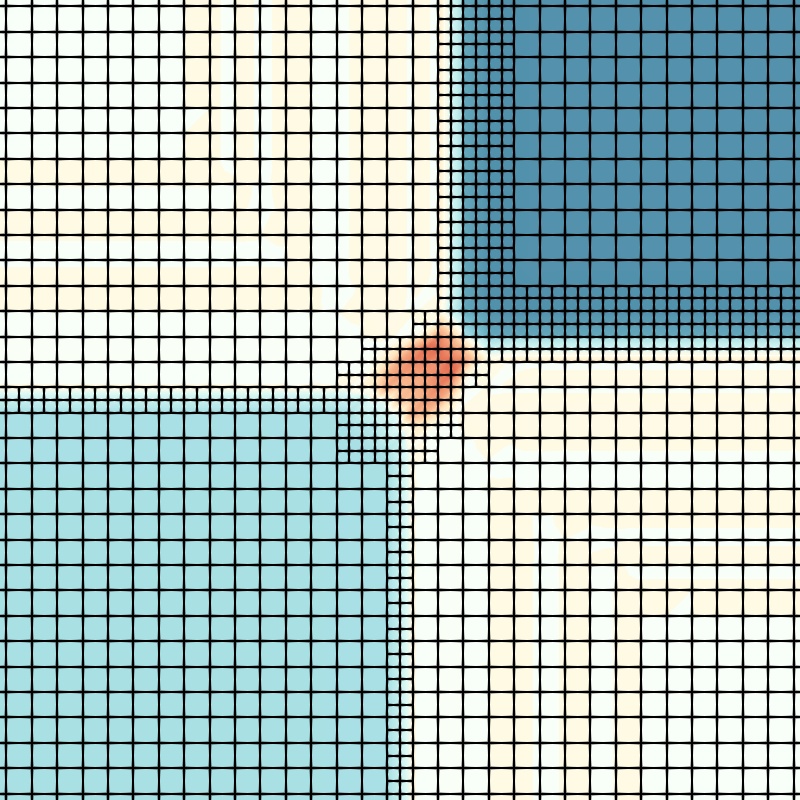}}}
        \subfloat[Solution at $t = 2T$]{
        \adjustbox{width=0.24\linewidth,valign=b}{\includegraphics{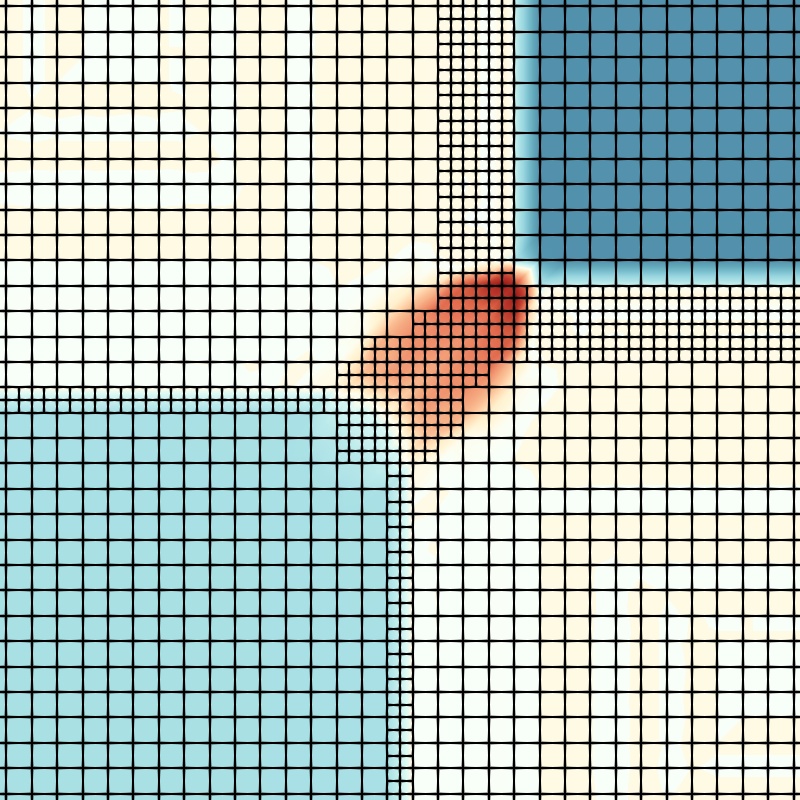}}}
        \newline
        \subfloat[Remesh at $t = 2T$]{
        \adjustbox{width=0.24\linewidth,valign=b}{\includegraphics{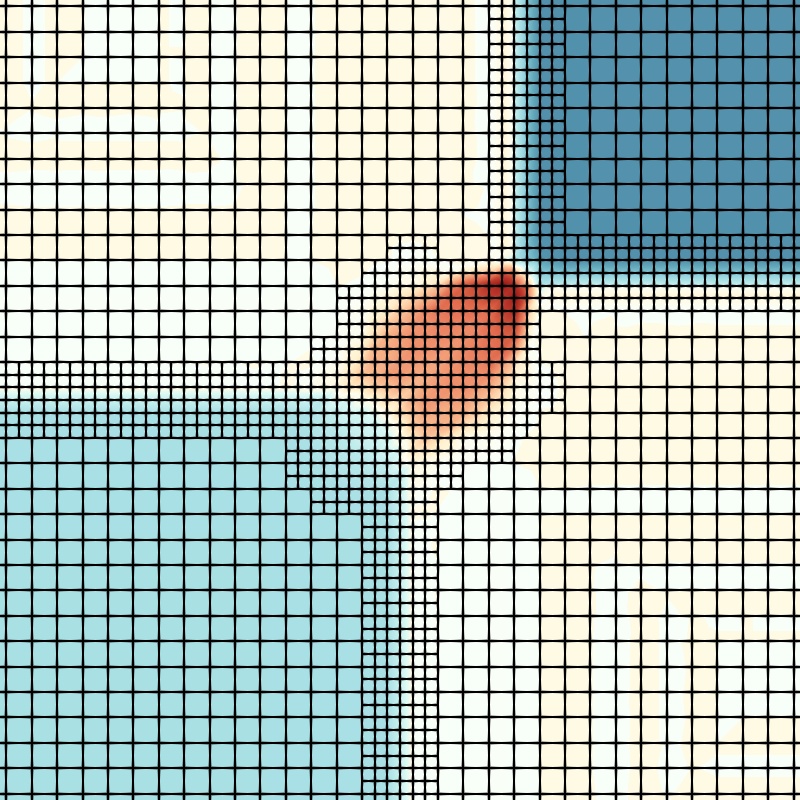}}}
        \subfloat[Solution at $t = 3T$]{
        \adjustbox{width=0.24\linewidth,valign=b}{\includegraphics{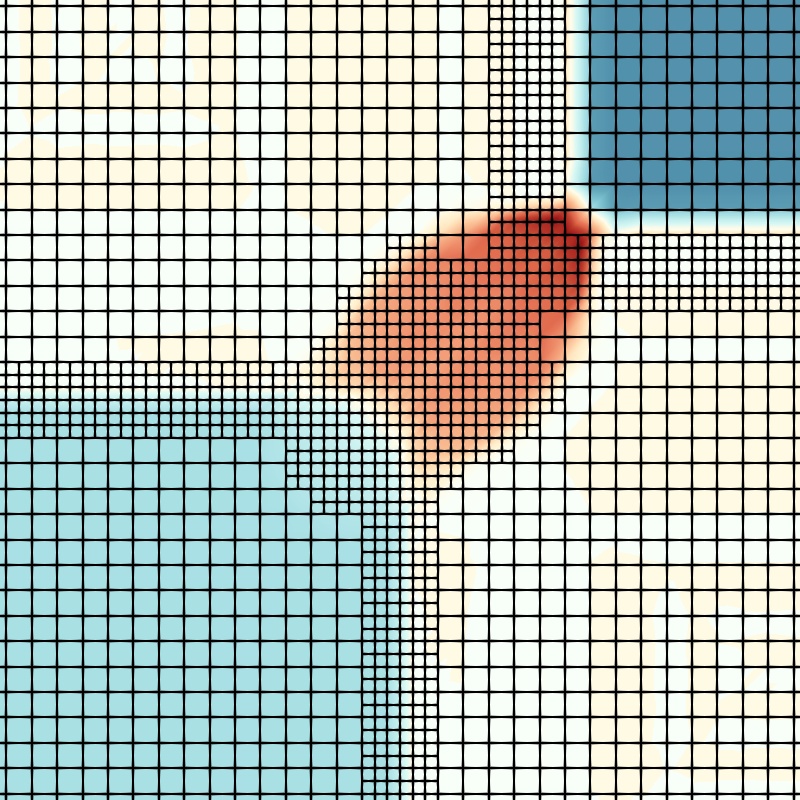}}}
        \subfloat[Remesh at $t = 3T$]{
        \adjustbox{width=0.24\linewidth,valign=b}{\includegraphics{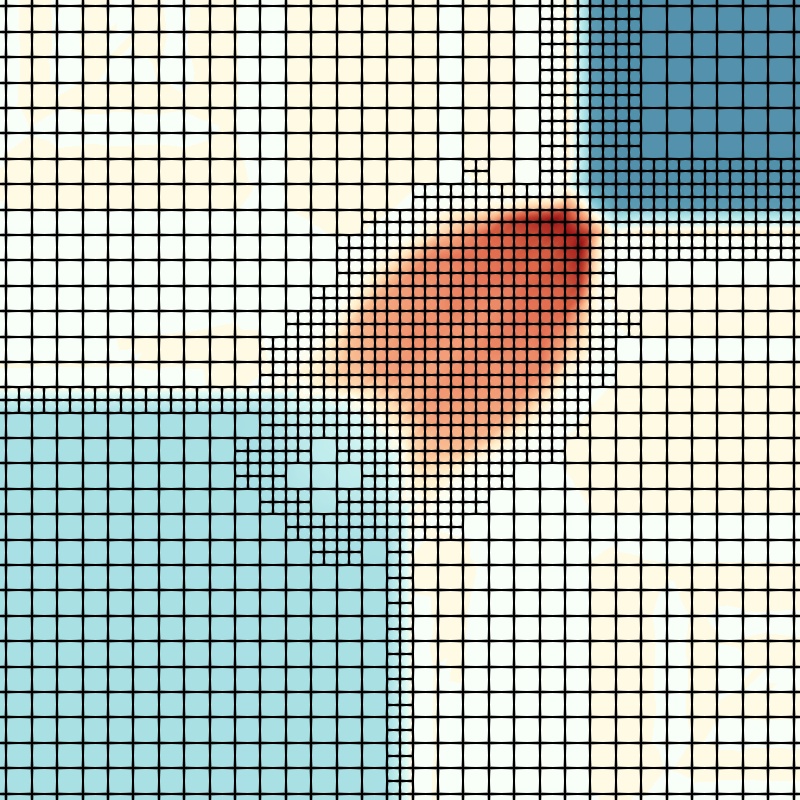}}}
        \subfloat[Solution at $t = 4T$]{
        \adjustbox{width=0.24\linewidth,valign=b}{\includegraphics{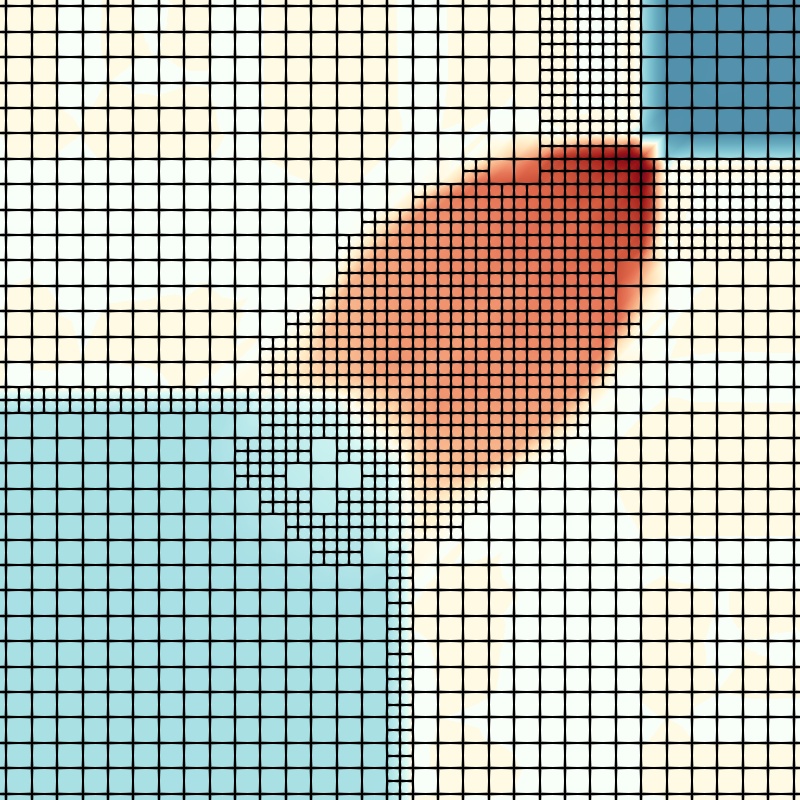}}}
        \newline
        \caption{\label{fig:rp12_amr_threshold} Contours of density overlaid with $h$-adapted mesh at varying remesh intervals using the threshold policy ($\theta = 10^{-2}$) for the Case 12 two-dimensional Riemann problem from \citet{Liska2003}.} 
    \end{figure}

The comparison between the DynAMO policy and threshold policy for Case 12 is also shown through the evolution of the density and $h$-adapted mesh in \cref{fig:rp12_amr_drl} and \cref{fig:rp12_amr_threshold}, respectively. Much like with Case 4, the DynAMO policy was able to preemptively refine the mesh at the initial time, such that the evolution of the discontinuities and expansion wave was contained within the refined region. However, it must be noted that since the contact discontinuities bordering the bottom-left quadrant are stationary, the DynAMO approach over-refined around these discontinuities. This behavior was later corrected, such that the refined region consisted of only a single element-wide region. The DynAMO policy did correctly predict the necessary refinement regions to cover the evolution of the discontinuities in the top-right quadrant. At later times, when the expansion wave was sufficiently developed, the DynAMO policy de-refined the interior region which could be sufficiently resolved on the coarse mesh level. For the threshold policy, an identical initial refinement pattern was observed with Case 12 as with Case 4 due to the similarity of the initial states. However, this refinement pattern also did not adequately cover the evolution of the discontinuities and expansion wave, such that a significant amount of discretization error was introduced around the discontinuity front. This behavior persisted throughout the simulation, with the threshold policy unable to adequately account for the propagation of the features.

\subsubsection{Generalization experiments}
The generalization capabilities of the DynAMO approach were further evaluated through \textit{out-of-distribution} experiments for problems on finer meshes, with longer remesh times, and with longer simulation times. The first generalization experiment was performed by increasing the base mesh resolution from $N = 32^2$ to $64^2$. A comparison between the mean efficiency of the DynAMO approach and the threshold policy on this finer base mesh is shown in \cref{tab:euler_href_ood} averaged over 100 runs. Both approaches showed an increase in mean efficiency, which is consist with previous generalization experiments, with the optimal threshold policy increasing from $0.251$ to $0.328$ and the DynAMO policy increasing from $0.486$ to $0.545$. This resulted in a relative efficiency increase of $66.2\%$ with the DynAMO policy, such that it significantly outperformed the threshold policy even on the finer mesh. 

    \begin{figure}[htbp!] 
        \centering
        \begin{tabular}{ccccc}
        \toprule
        Method & In-distribution &  Finer mesh  & Longer remesh time & Longer sim. time\\ 
        \midrule    
        Threshold ($\theta = 10^{-1}$) & 0.215 (0.083) & 0.280 (0.102) & 0.223 (0.082) &  \textbf{0.226 (0.096)} \\
        Threshold ($\theta = 10^{-2}$) & \textbf{0.251 (0.056)} & \textbf{0.328 (0.082)} &  \textbf{0.315 (0.052)} &  0.202 (0.065) \\
        Threshold ($\theta = 10^{-3}$) & 0.180 (0.031) & 0.245 (0.051) &  0.299 (0.044)  &  0.084 (0.022) \\
        Threshold ($\theta = 10^{-4}$) & 0.162 (0.027) & 0.212 (0.041) & 0.298 (0.042) &  0.075 (0.021) \\
        Threshold ($\theta = 10^{-5}$) & 0.157 (0.027) & 0.193 (0.035) & 0.298 (0.041) &  0.073 (0.021) \\
        Threshold ($\theta = 10^{-6}$) & 0.157 (0.027) & 0.180 (0.030) & 0.298 (0.041)  &  0.073 (0.021) \\
        \midrule
        DynAMO & \textbf{0.486 (0.073)} & \textbf{0.545 (0.069)} & \textbf{0.611 (0.048)} & \textbf{0.423 (0.113)} \\
        DynAMO/Optimal $\theta$ & \textcolor{green!70!black}{+93.6\%} & \textcolor{green!70!black}{+66.2\%} & \textcolor{green!70!black}{+93.9\%} & \textcolor{green!70!black}{+87.2\%} \\
        \bottomrule
        \end{tabular}
        \captionof{table}{\label{tab:euler_href_ood} Comparison of the mean efficiency for $h$-refinement on the Euler equations with DynAMO and the threshold policy for the two-dimensional Riemann problem over 100 \textit{out-of-distribution} runs using uniform random initial conditions with finer mesh resolution, longer remesh time, and longer simulation time. Standard deviation shown in parentheses. In-distribution results from \cref{tab:euler_href_indistribution} shown for comparison.}
    \end{figure}

The second generalization experiment was performed with respect to the remesh time $T$, where the remesh time was doubled from its training value. A comparison between the mean efficiency of the DynAMO approach and the threshold policy averaged over 100 runs for this longer remesh time is shown in \cref{tab:euler_href_ood}. The optimal threshold policy showed a marginal increase in efficiency in comparison to the \textit{in-distribution} results, from $0.251$ to $0.328$. The DynAMO policy also showed a similar increase, from $0.486$ to $0.545$. As such, the relative efficiency increase of the DynAMO policy over the optimal threshold policy stayed almost identical to the \textit{in-distribution} value, $93.9\%$, showcasing the ability of the approach to generalize to arbitrary remesh times. 

The third generalization experiment was performed with respect to the simulation time $t_f$, where the simulation time was doubled from its training value. \cref{tab:euler_href_ood} shows a comparison of the mean efficiency between the DynAMO approach and the threshold policy over 100 runs for this longer simulation time. For both approaches, the mean efficiency marginally decreased, likely as a result of the chaotic nature of the flow at later times for which the policies would tend towards more refinement. For the optimal threshold policy, this resulted in a mean efficiency decrease from $0.251$ to $0.226$, whereas for the DynAMO policy, this resulted in a mean efficiency decrease from $0.486$ to $0.423$. However, the DynAMO policy still showed a large relative efficiency increase over the threshold policy, $87.2\%$, which is on par with the \textit{in-distribution} results. These results indicate that the proposed approach can effectively generalize to longer simulation times and, by extent, flow physics which were not encountered during the training process.

    \begin{figure}[htbp!]
        \centering
        \subfloat[Remesh at $t = 0$]{
        \adjustbox{width=0.24\linewidth,valign=b}{\includegraphics{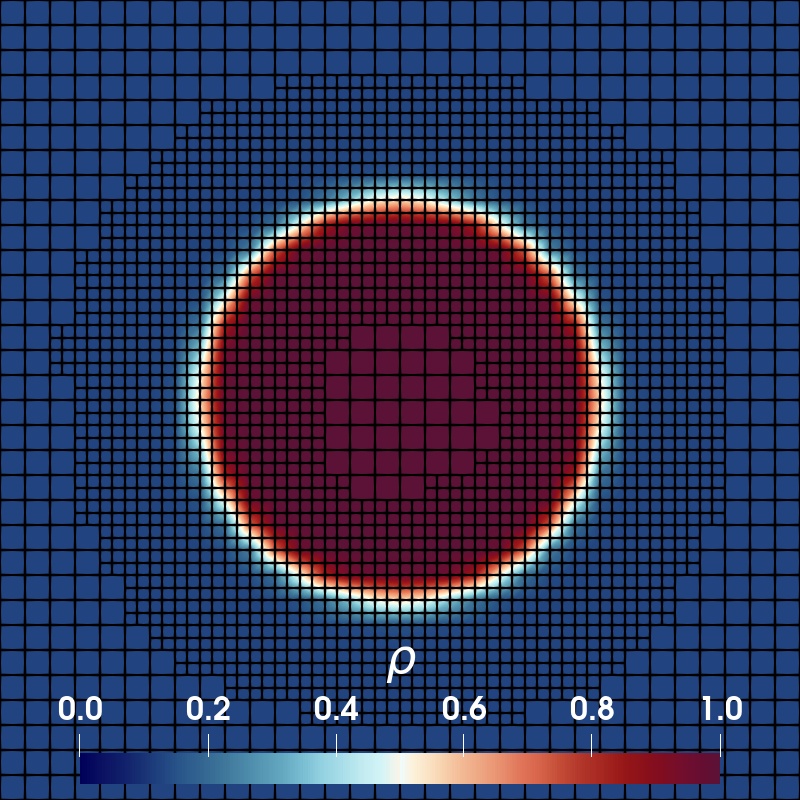}}}
        \subfloat[Solution at $t = T$]{
        \adjustbox{width=0.24\linewidth,valign=b}{\includegraphics{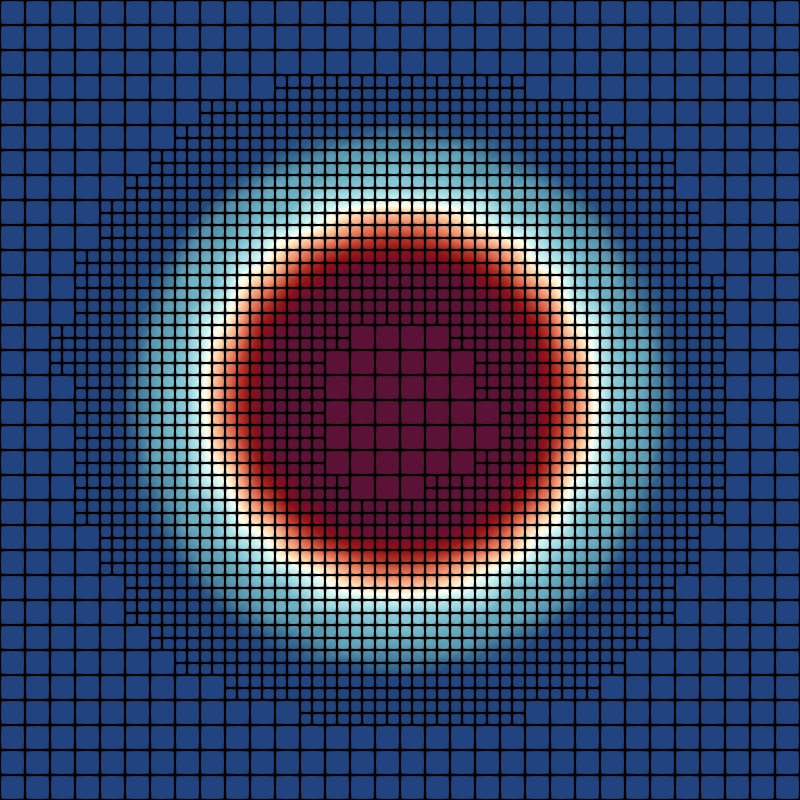}}}
        \subfloat[Remesh at $t = T$]{
        \adjustbox{width=0.24\linewidth,valign=b}{\includegraphics{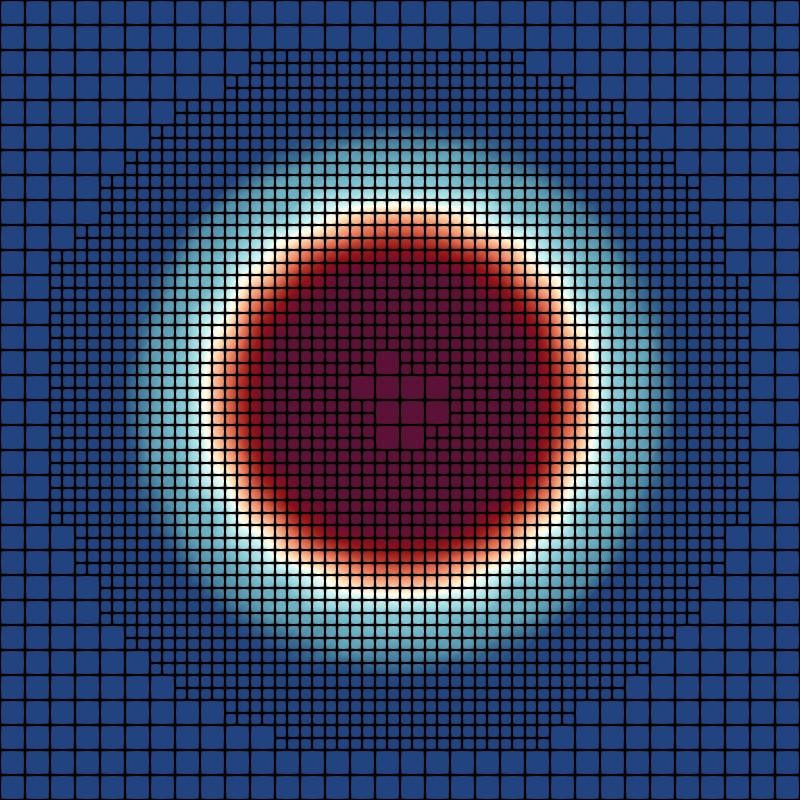}}}
        \subfloat[Solution at $t = 2T$]{
        \adjustbox{width=0.24\linewidth,valign=b}{\includegraphics{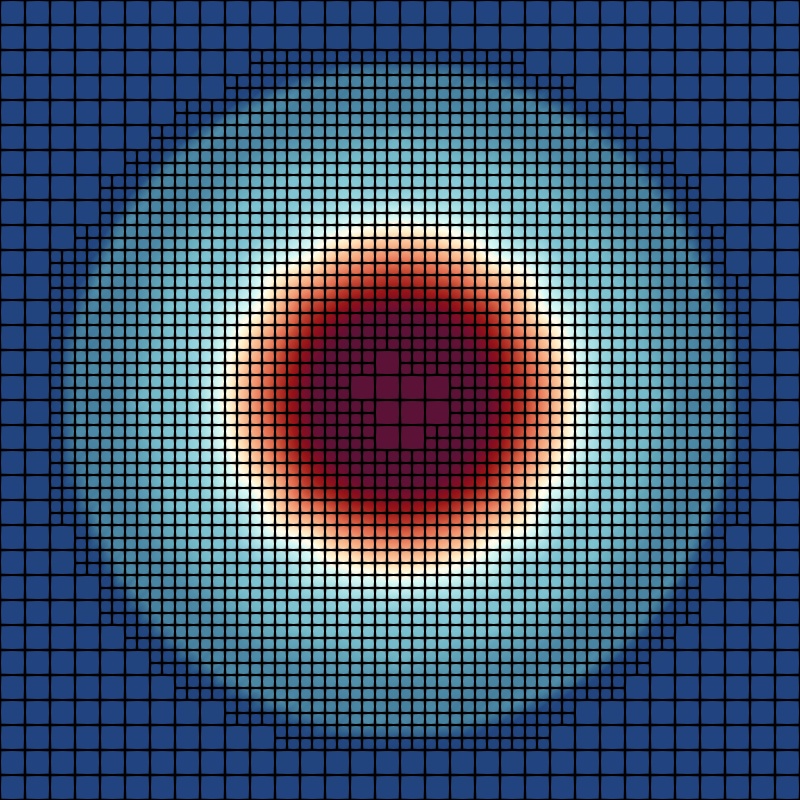}}}
        \newline
        \subfloat[Remesh at $t = 2T$]{
        \adjustbox{width=0.24\linewidth,valign=b}{\includegraphics{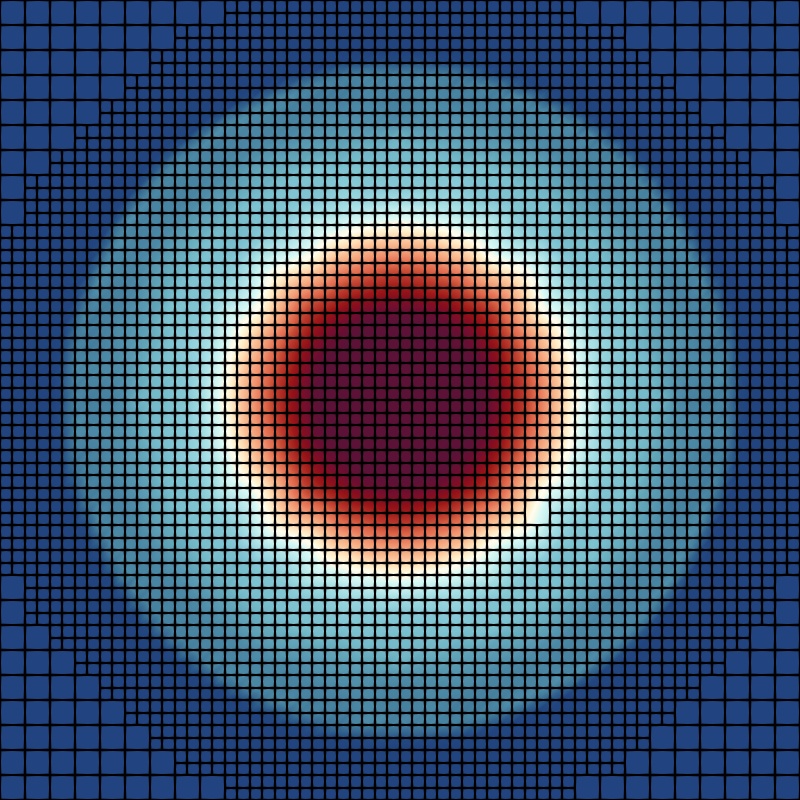}}}
        \subfloat[Solution at $t = 3T$]{
        \adjustbox{width=0.24\linewidth,valign=b}{\includegraphics{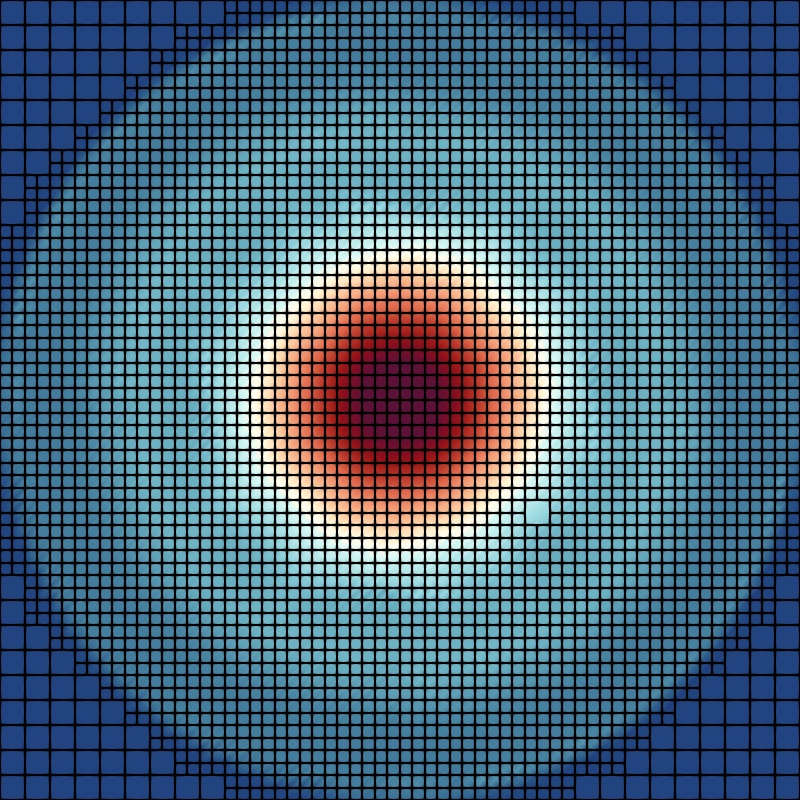}}}
        \subfloat[Remesh at $t = 3T$]{
        \adjustbox{width=0.24\linewidth,valign=b}{\includegraphics{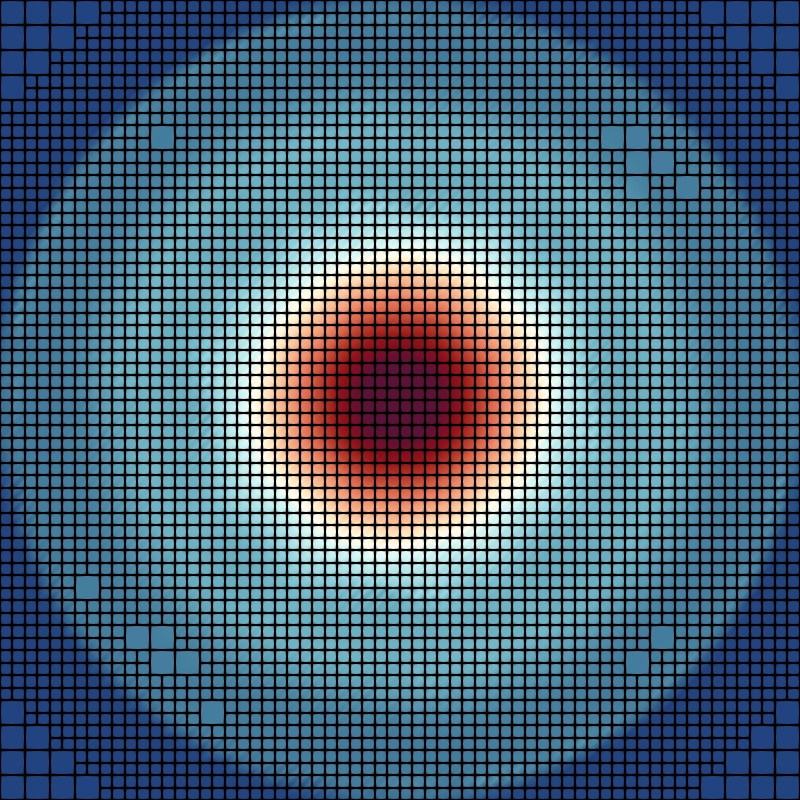}}}
        \subfloat[Solution at $t = 4T$]{
        \adjustbox{width=0.24\linewidth,valign=b}{\includegraphics{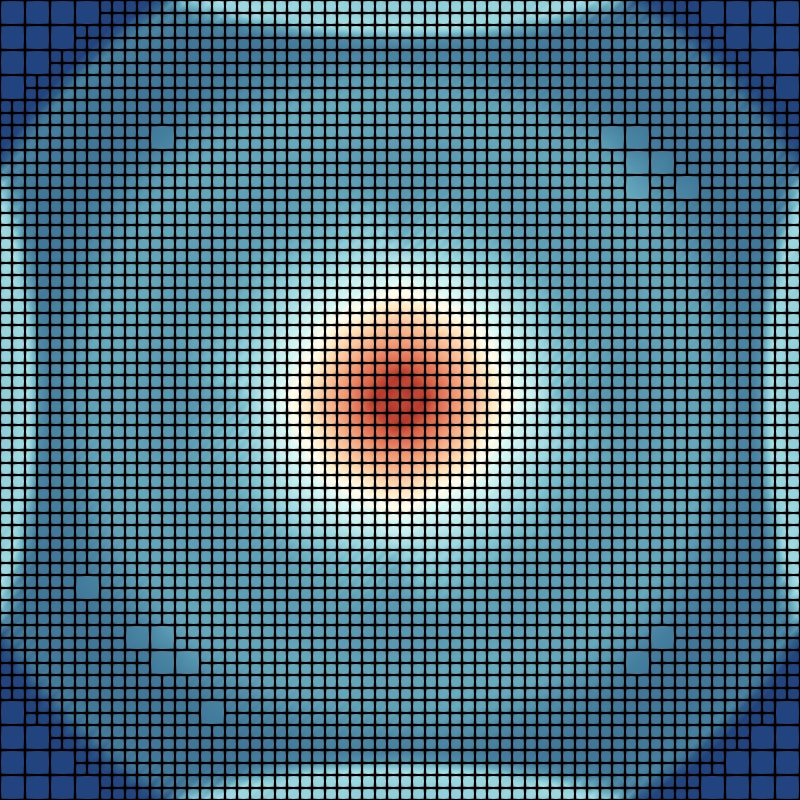}}}
        \newline
        \subfloat[Remesh at $t = 4T$]{
        \adjustbox{width=0.24\linewidth,valign=b}{\includegraphics{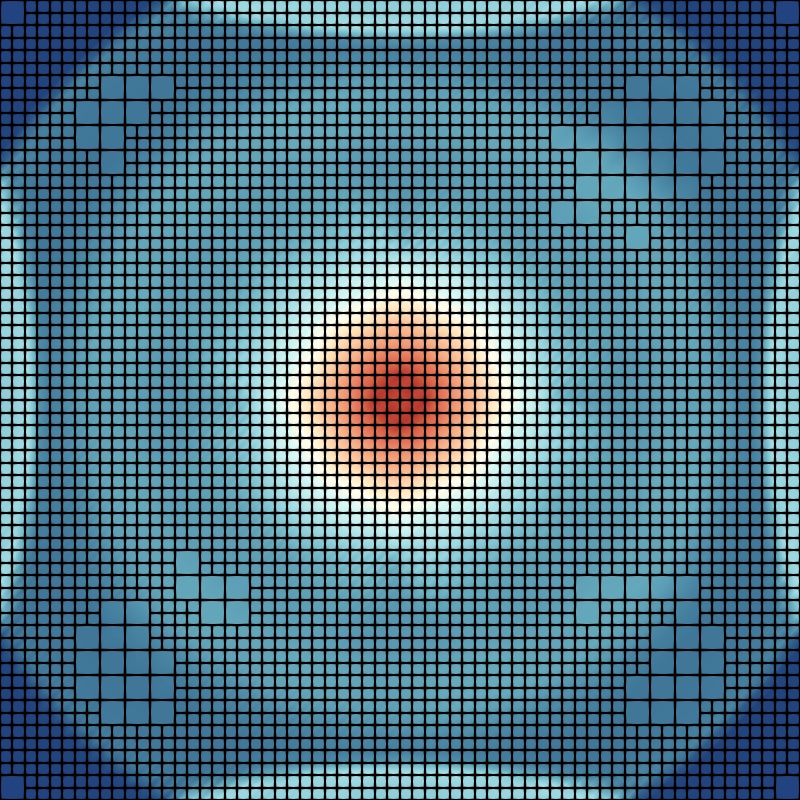}}}
        \subfloat[Solution at $t = 5T$]{
        \adjustbox{width=0.24\linewidth,valign=b}{\includegraphics{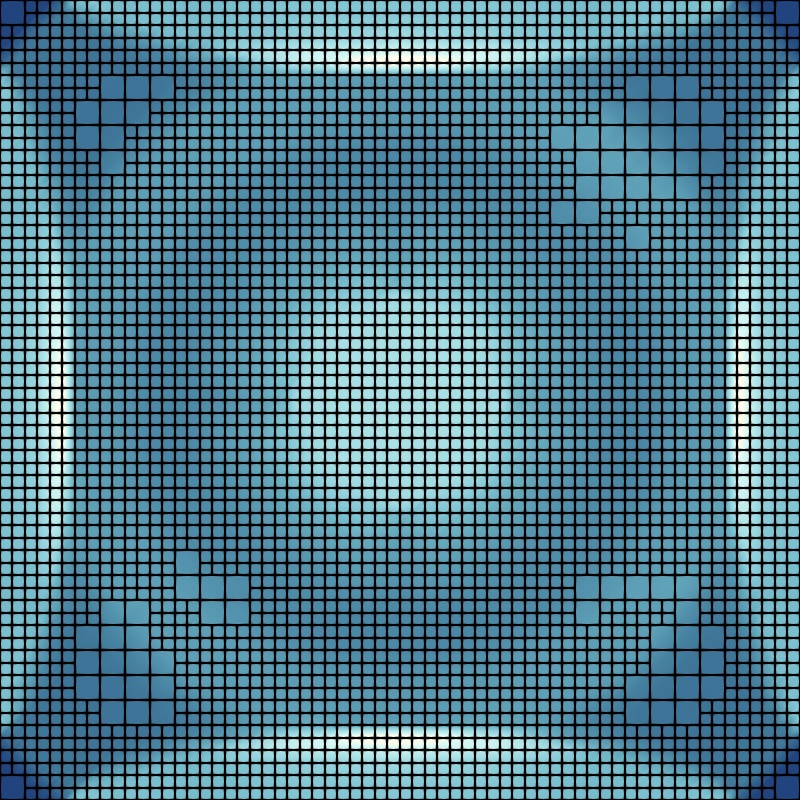}}}
        \subfloat[Remesh at $t = 5T$]{
        \adjustbox{width=0.24\linewidth,valign=b}{\includegraphics{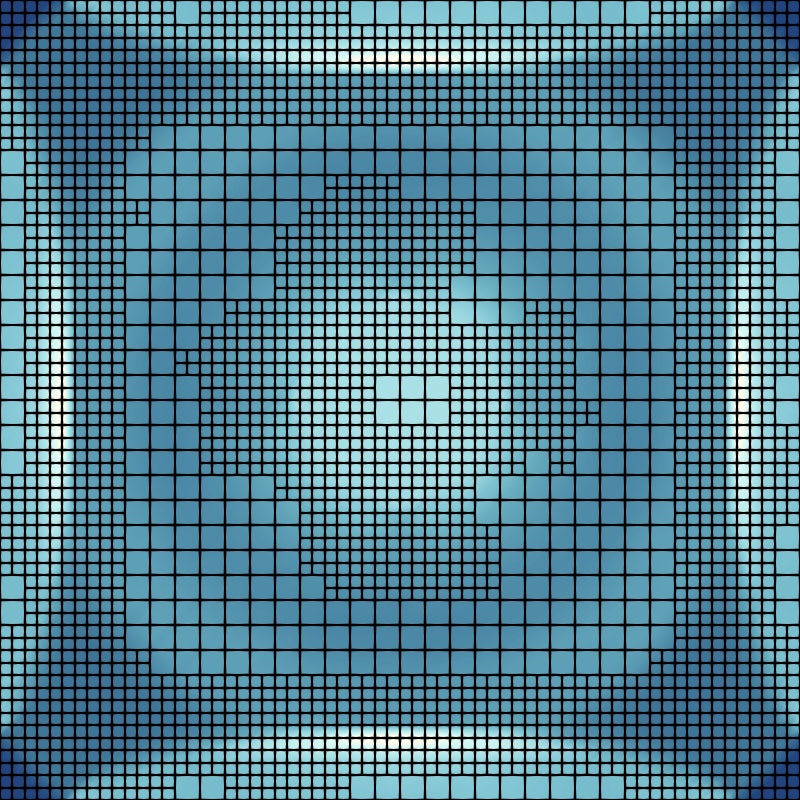}}}
        \subfloat[Solution at $t = 6T$]{
        \adjustbox{width=0.24\linewidth,valign=b}{\includegraphics{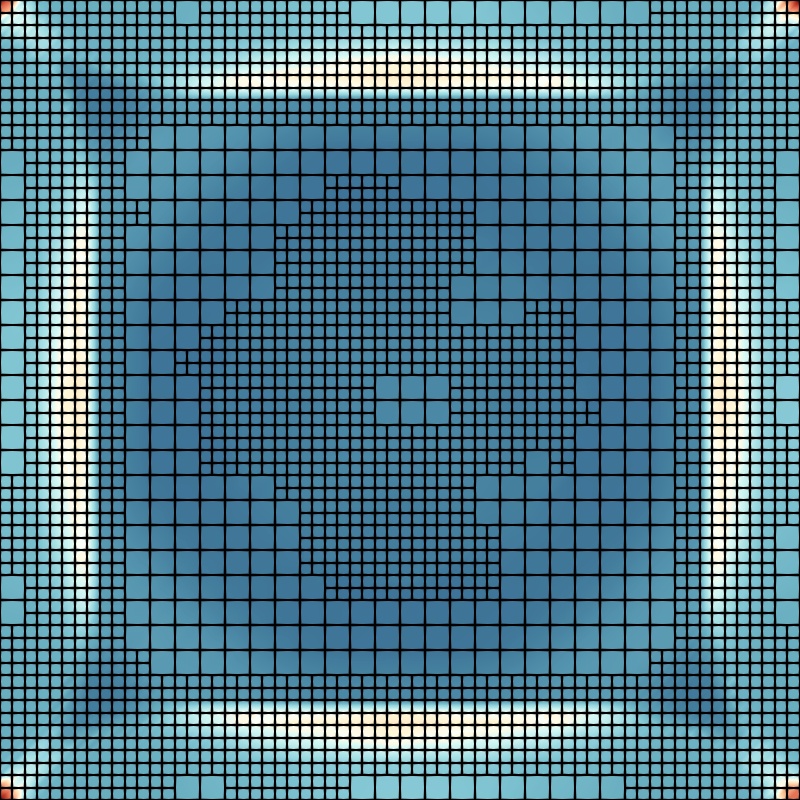}}}
        \newline
        \subfloat[Remesh at $t = 6T$]{
        \adjustbox{width=0.24\linewidth,valign=b}{\includegraphics{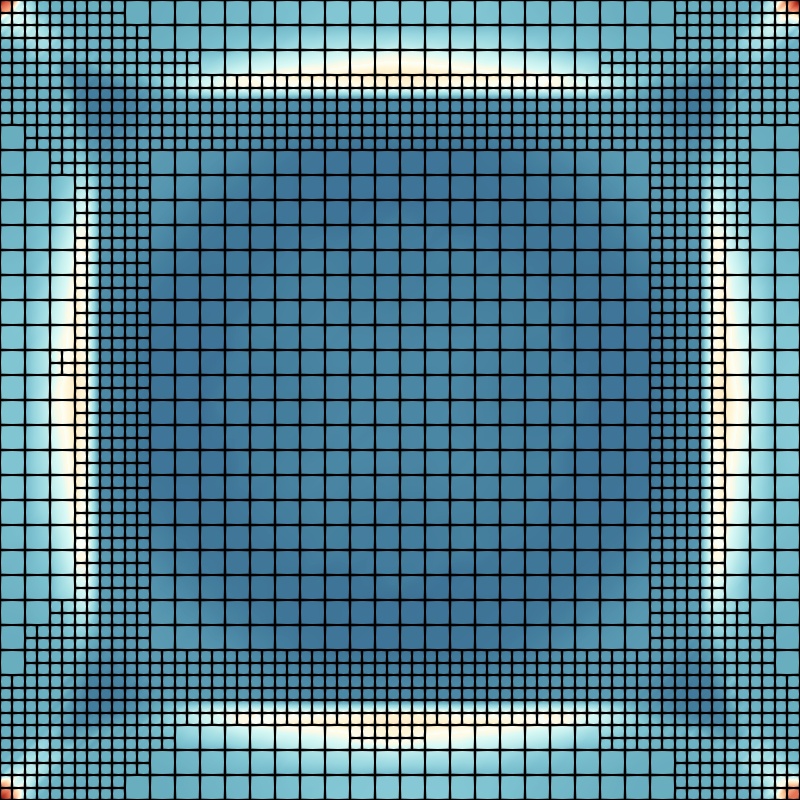}}}
        \subfloat[Solution at $t = 7T$]{
        \adjustbox{width=0.24\linewidth,valign=b}{\includegraphics{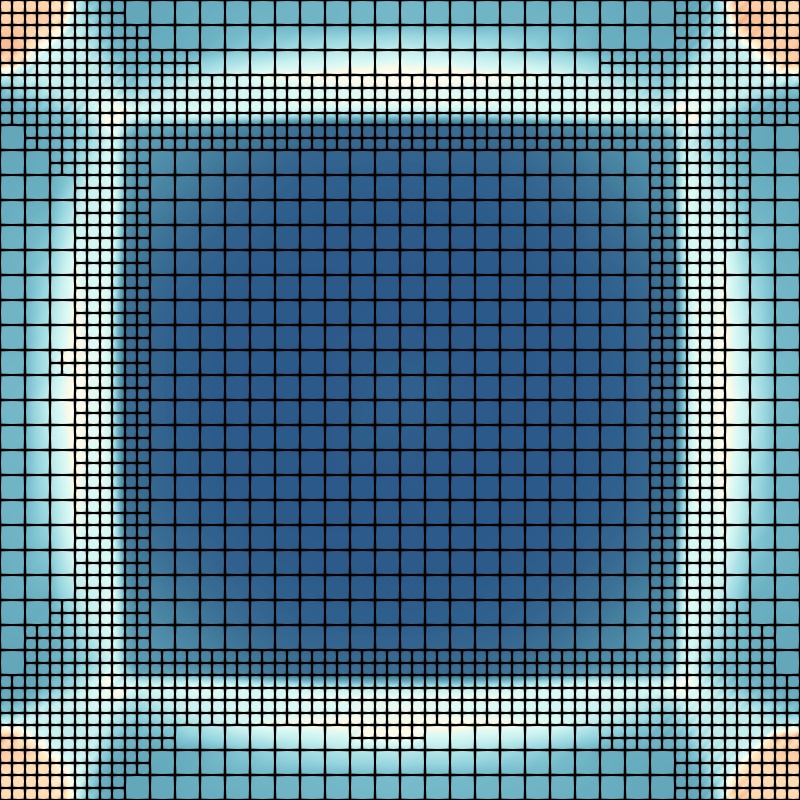}}}
        \subfloat[Remesh at $t = 7T$]{
        \adjustbox{width=0.24\linewidth,valign=b}{\includegraphics{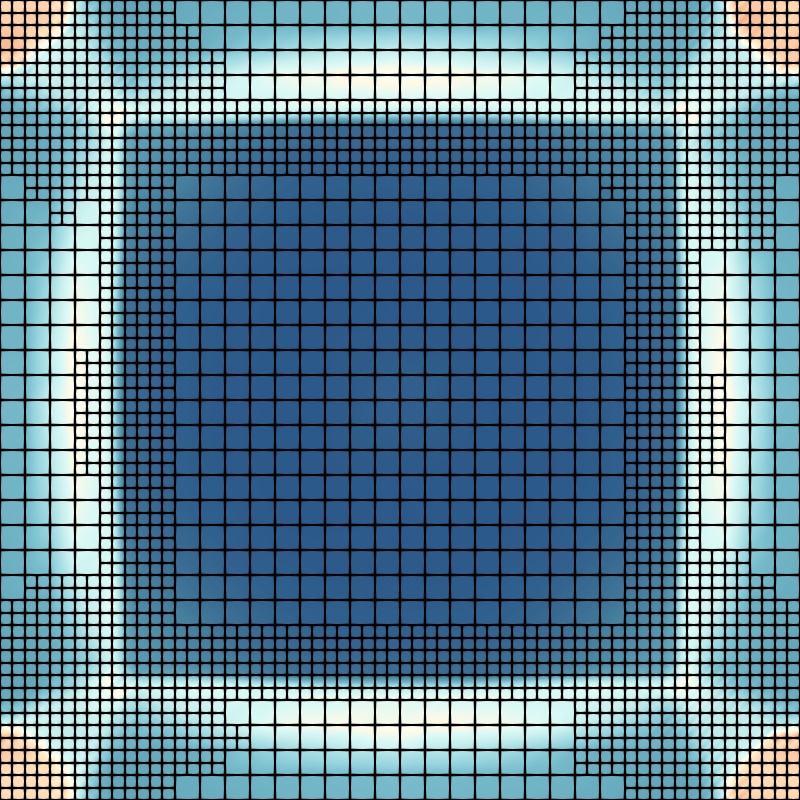}}}
        \subfloat[Solution at $t = 8T$]{
        \adjustbox{width=0.24\linewidth,valign=b}{\includegraphics{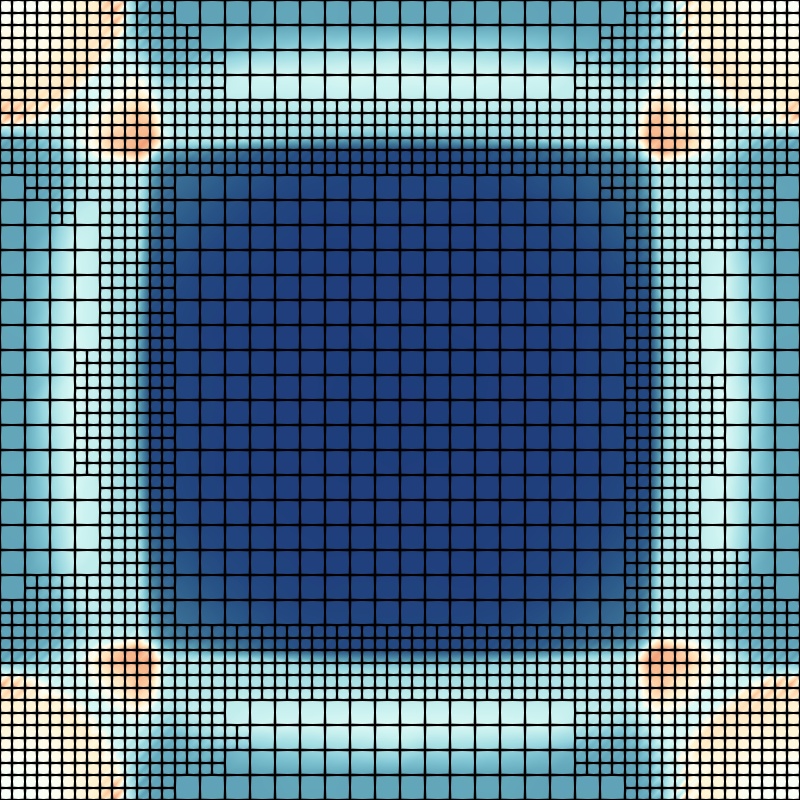}}}
        \newline
        \caption{\label{fig:sod_amr_drl} Contours of density overlaid with $h$-adapted mesh at varying remesh intervals using DynAMO for the two-dimensional Sod shock tube problem.} 
    \end{figure}
    
As a final validation of the DynAMO approach, we apply the DynAMO policy trained on two-dimensional Riemann problems to a particular \textit{out-of-distribution} problem: the two-dimensional radial Sod shock tube. This case is a more complex form of the quintessential example of the shock tube which exhibits the three main features of the Riemann problem, namely shock waves, contact discontinuities, and rarefaction waves. Furthermore, in its periodic form, it also exhibits shock-shock and shock-contact interactions at longer simulation times. The problem is solved on the domain $\Omega = [-0.5, 0.5]^2$, and the initial conditions are given as
\begin{equation}
    [\rho, u, v, P]^T = \begin{cases}
        [\rho_l, 0, 0, P_l]^T, \quad  \text{if } \sqrt{x^2 + y^2} \leq 0.25,\\
        [\rho_r, 0, 0, P_r]^T, \quad  \text{else,} 
    \end{cases}
\end{equation}
where $\rho_l = 1.0$, $P_l = 1.0$, $\rho_r = 0.125$, and $P_r = 0.1$. The base mesh resolution was set as $N=32^2$ and the remesh time was set as $T = 0.05$. The evolution of the density and $h$-adapted mesh as obtained by the DynAMO policy is shown in \cref{fig:sod_amr_drl} over 8 remesh intervals. At the initial time, the DynAMO policy preemptively refines the region both inside and outside of the initial discontinuity. This refinement region covers the propagation of the out-running shock front and the in-running expansion wave front, although with a slight overprediction in the refinement region necessary to cover the propagation distance of the shock. Furthermore, the interior region of the mesh was correctly left de-refined as the expansion wave front did not reach this region over the remesh interval, although with some minor asymmetry in the refinement patterns. Similar behavior was observed at the next few remesh intervals, with the DynAMO policy accurately predicting the preemptive refinement decisions necessary to account for the evolution of the shock front. As the shock front propagated towards the domain boundaries, certain spurious de-refinement decisions were chosen by the policy which were not consistent with the decisions at other similar locations in mesh. However, these decisions primarily acted to introduce some minor asymmetry in the refinement patterns without appreciably degrading the accuracy of the approach. After the shock fronts propagated through the periodic boundaries, the shock-shock and shock-contact interactions were correctly refined by the DynAMO policy. Furthermore, the interior region of the mesh, where the solution variation was low, was appropriately de-refined, such that the approach effectively balanced computational cost and accuracy. These results further showcase the ability of the proposed approach to be an effective method of performing anticipatory mesh refinement for complex nonlinear flows and generalizing to flow physics not encountered during training. 

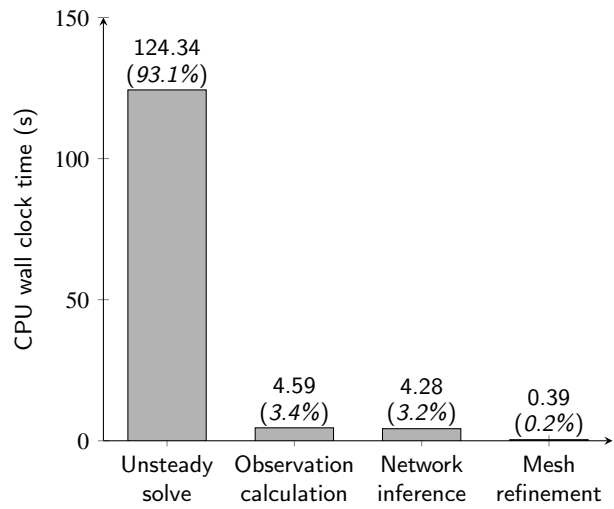
\begin{figure}[htbp!]
        \centering
        \adjustbox{width=0.5\linewidth,valign=b}{\begin{tikzpicture}[spy using outlines={rectangle, height=3cm,width=2.3cm, magnification=3, connect spies}]
\begin{axis} [ybar,
		axis line style={latex-latex},
	    axis x line=left,
        axis y line=left,
        bar width=30pt,
    	xmin=0.25, xmax=2.25,
    	ymin=0, ymax=150,
    	ylabel={CPU wall clock time (s)},
        xtick={0.5, 1, 1.5, 2},
        xticklabels={Unsteady solve, {Observation calculation}, Network inference, Mesh refinement},
        x tick label style={align=center, text width=1.75cm}, 
        clip mode=individual,
    	legend style={at={(0.03, 0.97)},anchor=north west},
    	legend cell align={left}]
     
\addplot[draw=black,fill=black!30]
    coordinates {
    	(0.5, 124.34) 
    	(1.0, 4.59) 
    	(1.5, 4.28) 
    	(2.0, 0.386) 
    };

\node at (0.5, 139.34) {124.34};
\node at (0.5, 129.34) {(\textit{93.1\%})};
\node at (1.0, 19.59) {4.59};
\node at (1.0, 9.59) {(\textit{3.4\%})};
\node at (1.5, 19.28) {4.28};
\node at (1.5, 9.28) {(\textit{3.2\%})};
\node at (2.0, 15.38) {0.39};
\node at (2.0, 5.38) {(\textit{0.2\%})};

\end{axis}

\end{tikzpicture}}
        \caption{\label{fig:cost} Computational cost distribution for the two-dimensional Sod shock tube problem in terms of unsteady solve, observation calculation, network inference, and mesh refinement. Labels represent absolute wall clock time (CPU-seconds) with relative cost shown in parentheses.}
    \end{figure}

The two-dimensional Sod shock tube problem was also used to evaluate the computational cost of the proposed approach. The cost, shown in absolute CPU wall clock time, was computed for the problem across four major categories: 1) advancing the unsteady solver between remesh intervals; 2) computing the observation (i.e., evaluating the error estimates, computing the observable quantities and metrics, etc.); 3) performing inference on the neural network; and 4) adapting the mesh and restructuring data. The absolute and relative computational costs of these four categories are shown in \cref{fig:cost}. It can be seen that the unsteady solver was the predominant factor in the overall computational cost, accounting for $93.1\%$ of the total compute time. In contrast, the cost of the proposed approach was only $6.6\%$ of the overall compute time, with the formation of the observation and network inference contributing approximately equally to the cost, and the cost of the mesh refinement process was essentially negligible. These results indicate that even with an RL framework that was not necessarily tailored for peak computational efficiency, the computational cost of the proposed approach did not significantly increase the overall cost of the solver for a complex system of equations.

\section{Conclusions}\label{sec:conclusion}
In this work, we presented DynAMO, a novel reinforcement learning paradigm for guiding anticipatory mesh refinement strategies for complex time-dependent PDEs. In comparison to standard adaptive mesh refinement approaches for time-dependent problems, which rely on instantaneous error indicators that typically cannot account for the spatio-temporal evolution of the error, the goal of the proposed approach is to predict future error propagation and preemptively refine the mesh to achieve superior accuracy and efficiency. To this end, we considered a multi-agent reinforcement learning approach to a decentralized partially observable Markov decision process model of mesh optimization, with mesh elements corresponding to independent agents that observe local surrounding information. To maximize the generalizability and robustness of the approach, we introduced novel observation and reward functions that utilized a user-adjustable error estimate normalization and a non-dimensionalized measure of information propagation. These proposed formulations enabled the extension of the approach to arbitrary nonlinear hyperbolic conservation laws, ensured that the approach was invariant to problem scale, mesh resolution, simulation time, and remeshing time interval, and allowed for the user to control error/cost targets at evaluation time. 

The DynAMO approach was applied to both $h$- and $p$-refinement for discontinuous Galerkin approximations of the linear advection and compressible Euler equations. We showed that due to its anticipatory refinement capability, the proposed approach could achieve significantly higher efficiency than conventional AMR approaches, such as threshold-based methods, yielding more accurate simulations at lower computational costs.
Furthermore, the use of DynAMO allowed for longer remesh intervals without sacrificing accuracy by producing meshes that were able to preemptively account for the underlying dynamics of the system. These benefits were observed over a wide variety of problems, ranging from simple linear transport to complex nonlinear shock interactions, including problems and numerical setups on which the policies were not trained. These results indicate that the proposed approach can effectively anticipate the spatiotemporal evolution of the error for complex nonlinear PDEs, generalize to unseen problems at evaluation time, and robustly extend to different meshes, simulation times, and remesh time intervals. 

Future developments for this approach will focus on mitigating the limitations present in its current form. In particular, this work showcases the method for one-level refinement on periodic meshes using a structured observation window. To extend this approach to higher refinement levels, problems with boundary conditions, and arbitrary unstructured meshes, the formulation can be modified to utilize a graph neural network as explored in the works of \citet{yang2023reinforcement} and \citet{foucart2023deep}. In fact, the proposed observation, which utilizes a non-dimensionalized measure of information propagation that is dependent on the mesh displacement vector, is highly amenable to graph-type implementations. Additionally, further work can be performed on exploring more complex mesh refinement techniques such as $hp$-refinement, for which optimal strategies are far from fully understood.

\section*{Acknowledgements}
\label{sec:ack}

This work was performed under the auspices of the U.S. Department of Energy by Lawrence Livermore National Laboratory under contract DE--AC52--07NA27344 and the LLNL-LDRD Program under Project tracking No.\ 21--SI--001. Release number LLNL--JRNL--855285.

\bibliographystyle{unsrtnat}
\bibliography{reference}

\clearpage
\begin{appendices}
\crefalias{section}{appendix}
\section{Euler flux Jacobian}\label{app:fluxjacobian}
Recall the solution and flux of the compressible Euler equations,
\begin{equation*}
    \mathbf{u} = \begin{bmatrix}
            \rho \\ \mathbf{m} \\ E
        \end{bmatrix} \quad  \mathrm{and} \quad \mathbf{F} = \begin{bmatrix}
            \mathbf{m}^T\\
            \mathbf{m}\otimes\mathbf{v} + P\mathbf{I}\\
        (E+P)\mathbf{v}^T
    \end{bmatrix}=\begin{bmatrix}
        \mathbf{m}^T\\
        \mathbf{m}\otimes\mathbf{m}/\rho+P\mathbf{I}\\
        \mathbf{m}^T(E+P)/\rho
    \end{bmatrix},
\end{equation*}
where $P=(\gamma-1)(E-\frac{1}{2}\|\mathbf{m}\|^2 / \rho)$.
The following relation is convenient to derive the flux Jacobian
\begin{equation*}
    \frac{\partial P}{\partial \mathbf{u}}=(\gamma-1)\begin{bmatrix}
        \frac{1}{2\rho^2}\|\mathbf{m}\|^2\\
        -\frac{1}{\rho}\mathbf{m}\\
        1
    \end{bmatrix}.
\end{equation*}
Let $\mathbf{F}_i$ be the $i$-th column of $\mathbf{F}(u)$. We can then calculate the flux Jacobian as
\begin{equation*}
\begin{aligned}
    \frac{\partial \mathbf{F}_i}{\partial \mathbf{u}}
    &=\begin{bmatrix}
        0&\mathbf{d}_i^T&0\\
        -\frac{m_i}{\rho^2}\mathbf{m}+\frac{1}{2\rho^2}\|\mathbf{m}\|^2\mathbf{d}_i&\frac{1}{\rho}(m_i\mathbf{I}+\mathbf{m}\otimes\mathbf{d}_i)-(\gamma-1)\frac{1}{\rho}\mathbf{d}_i\otimes\mathbf{m}&(\gamma-1)\mathbf{d}_i\\
        -\frac{m_i}{\rho^2}(E+P)+(\gamma-1)\frac{m_i}{2\rho^3}\|\mathbf{m}\|^2&\frac{1}{\rho}(E+P)\mathbf{d}_i^T-(\gamma-1)\frac{m_i}{\rho^2}\mathbf{m}^T&\gamma\frac{m_i}{\rho}
    \end{bmatrix},
\end{aligned}
\end{equation*}
where $\mathbf{d}_i$ is the standard basis in $\mathbb{R}^d$ with $(\mathbf{d}_i)_j=\delta_{ij}$ and $m_i$ is the $i$-th component of the momentum $\mathbf{m}$.

\section{Sample generalization experiments}\label{app:generalization}
\begin{figure}[htbp!]
        \centering
        \subfloat[Remesh at $t = 0$]{
        \adjustbox{width=0.24\linewidth,valign=b}{\includegraphics{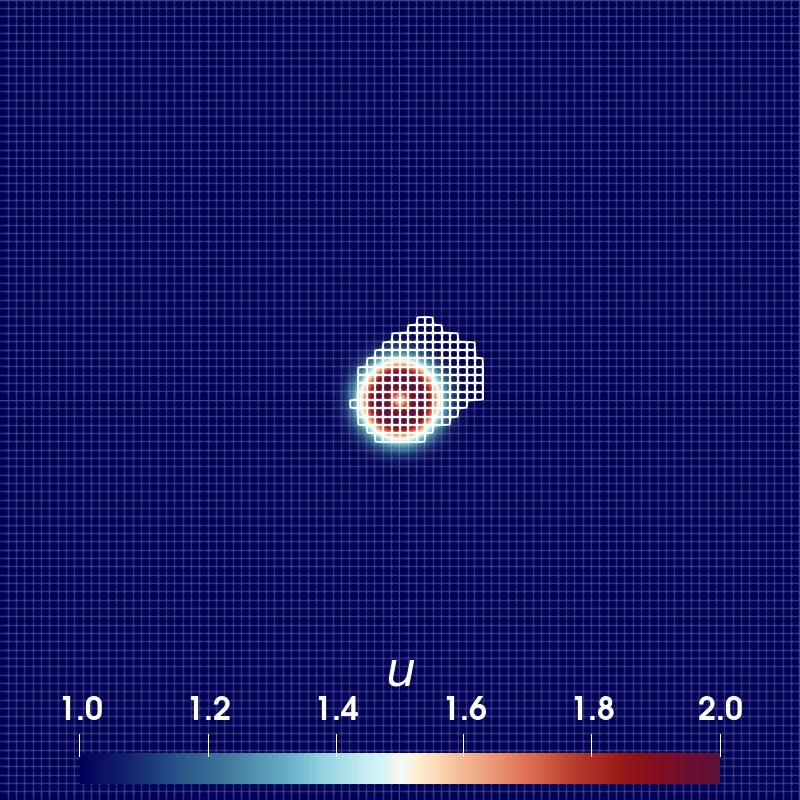}}}
        \subfloat[Solution at $t = T$]{
        \adjustbox{width=0.24\linewidth,valign=b}{\includegraphics{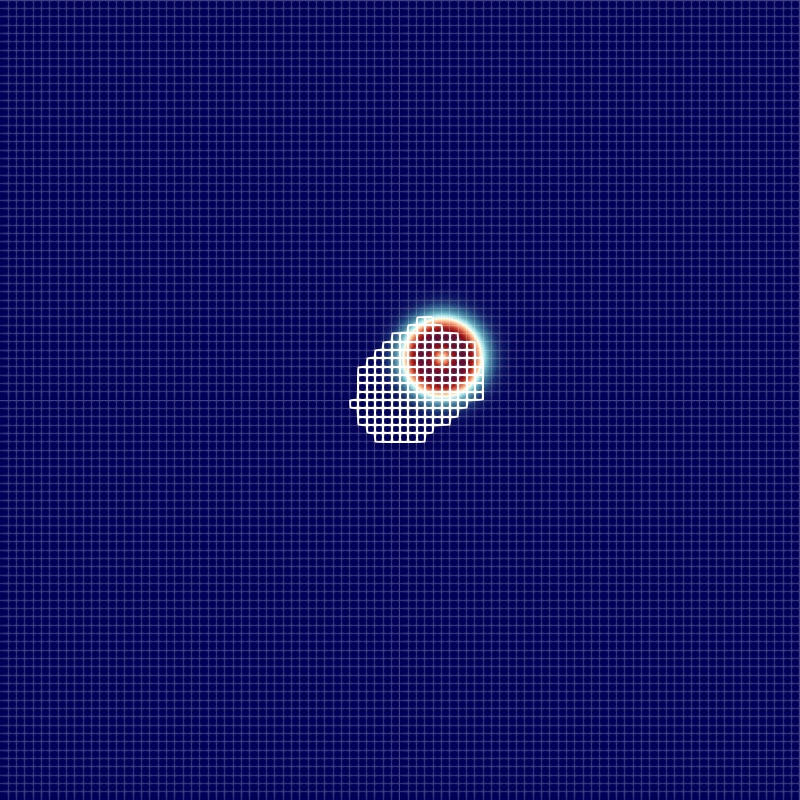}}}
        \subfloat[Remesh at $t = T$]{
        \adjustbox{width=0.24\linewidth,valign=b}{\includegraphics{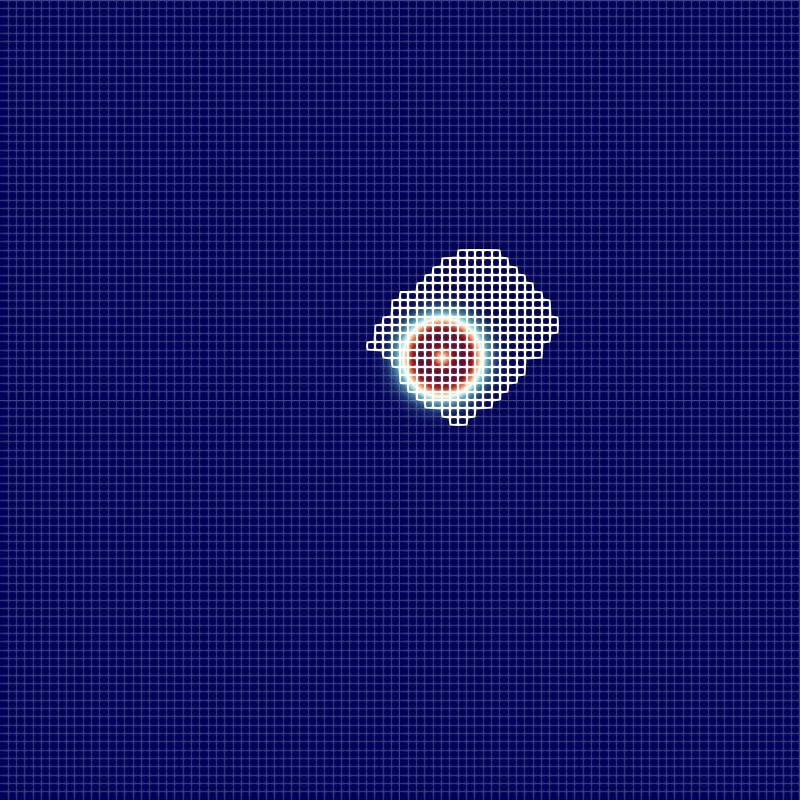}}}
        \subfloat[Solution at $t = 2T$]{
        \adjustbox{width=0.24\linewidth,valign=b}{\includegraphics{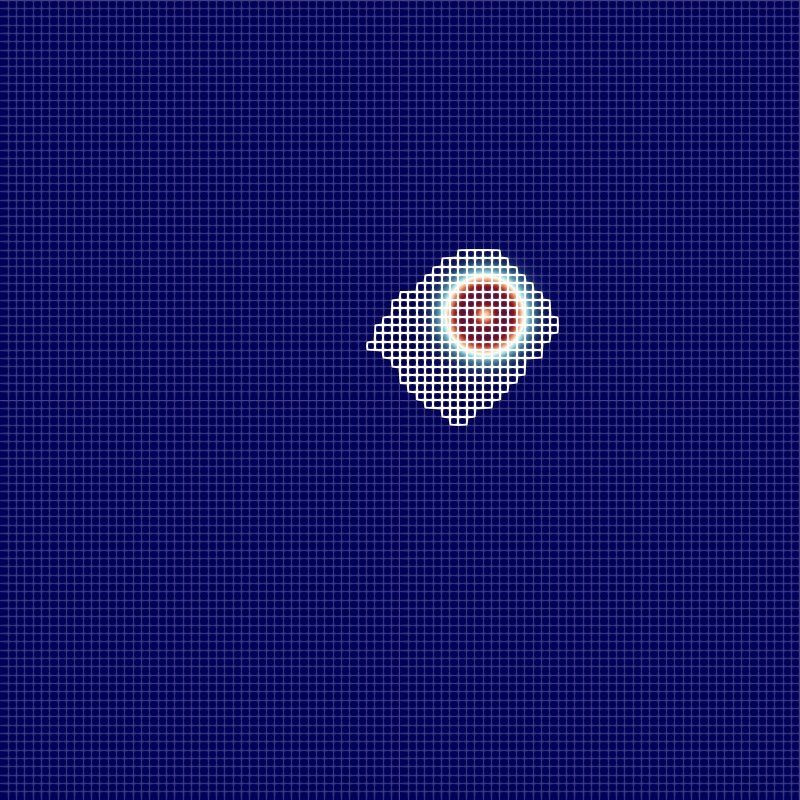}}}
        \newline
        \subfloat[Remesh at $t = 2T$]{
        \adjustbox{width=0.24\linewidth,valign=b}{\includegraphics{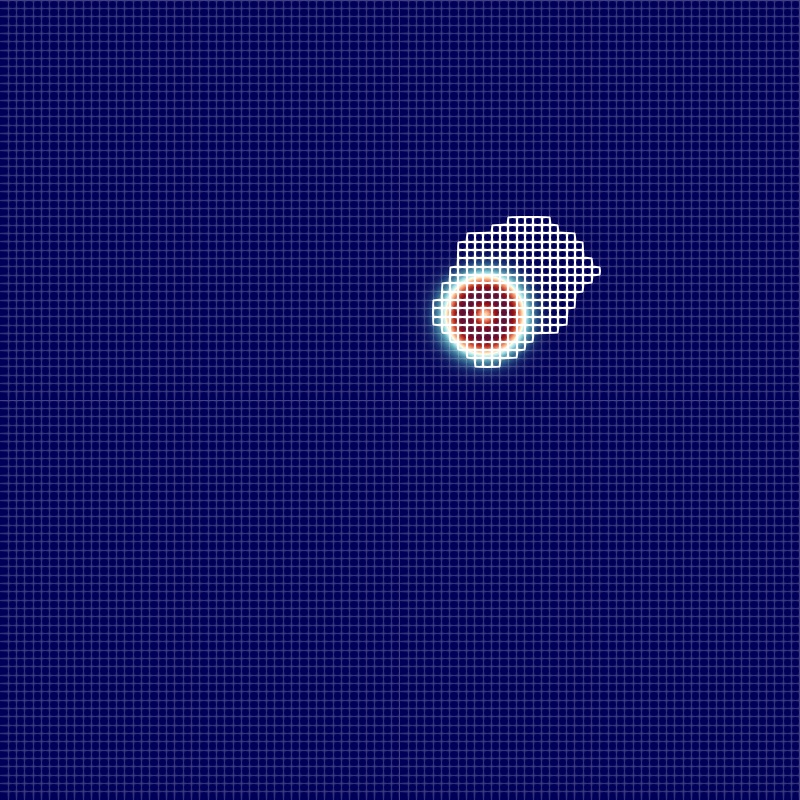}}}
        \subfloat[Solution at $t = 3T$]{
        \adjustbox{width=0.24\linewidth,valign=b}{\includegraphics{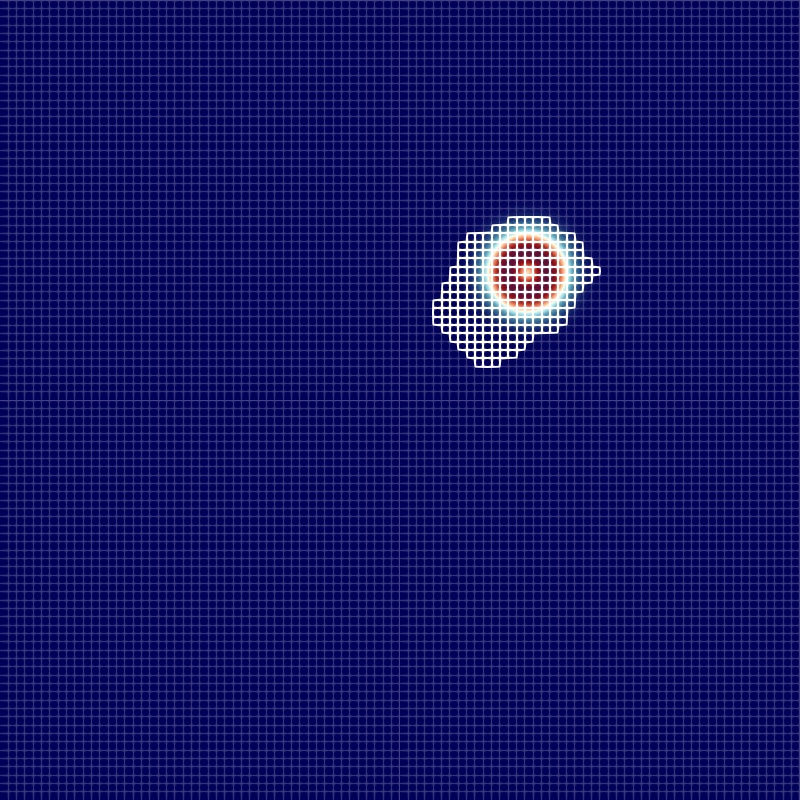}}}
        \subfloat[Remesh, $3T$]{
        \adjustbox{width=0.24\linewidth,valign=b}{\includegraphics{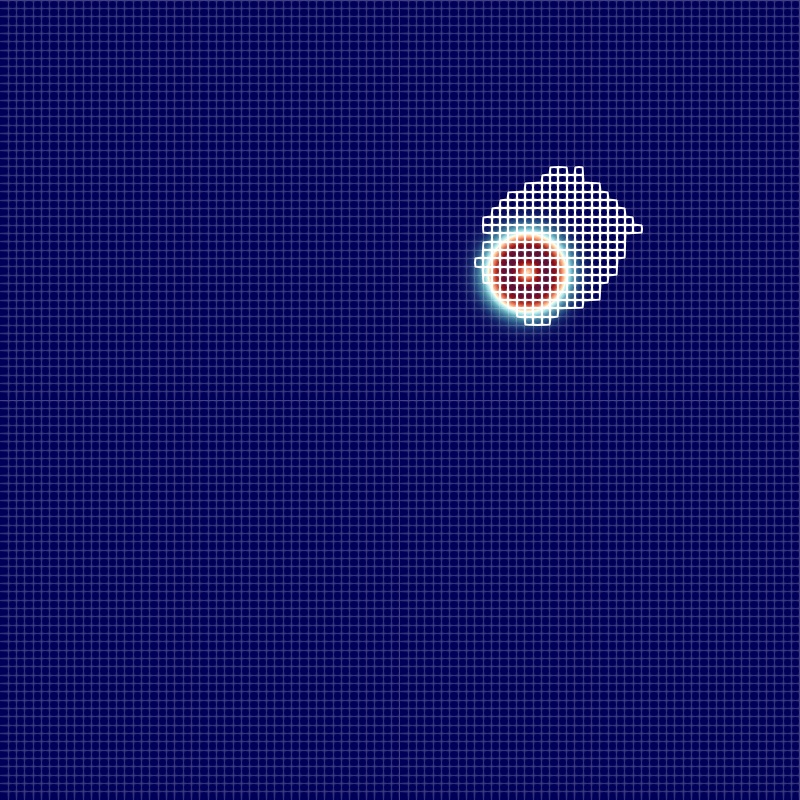}}}
        \subfloat[Solution, $4T$]{
        \adjustbox{width=0.24\linewidth,valign=b}{\includegraphics{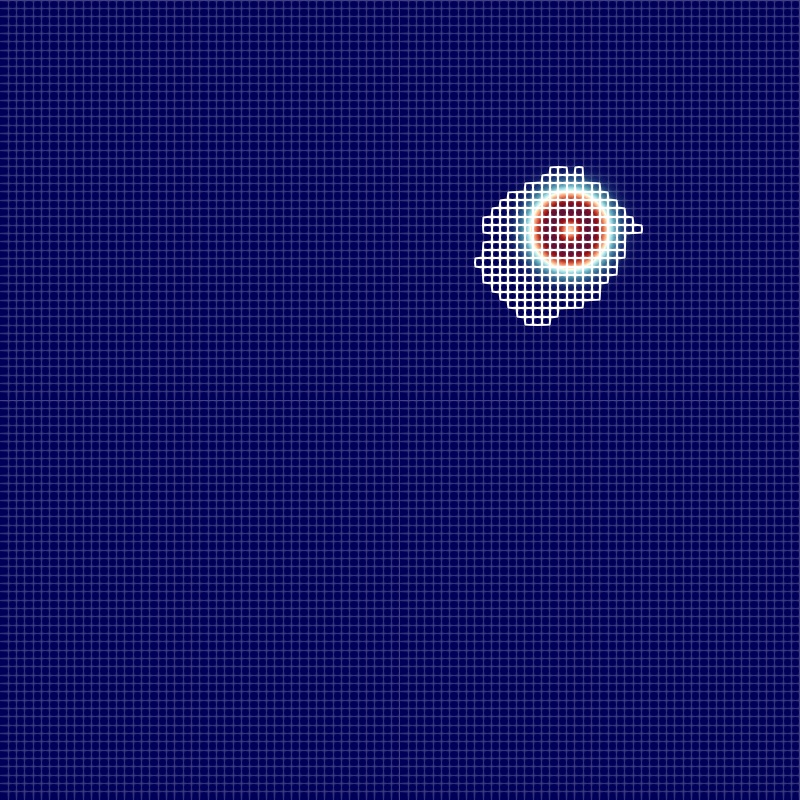}}}
        \newline
        \subfloat[Remesh, $4T$]{
        \adjustbox{width=0.24\linewidth,valign=b}{\includegraphics{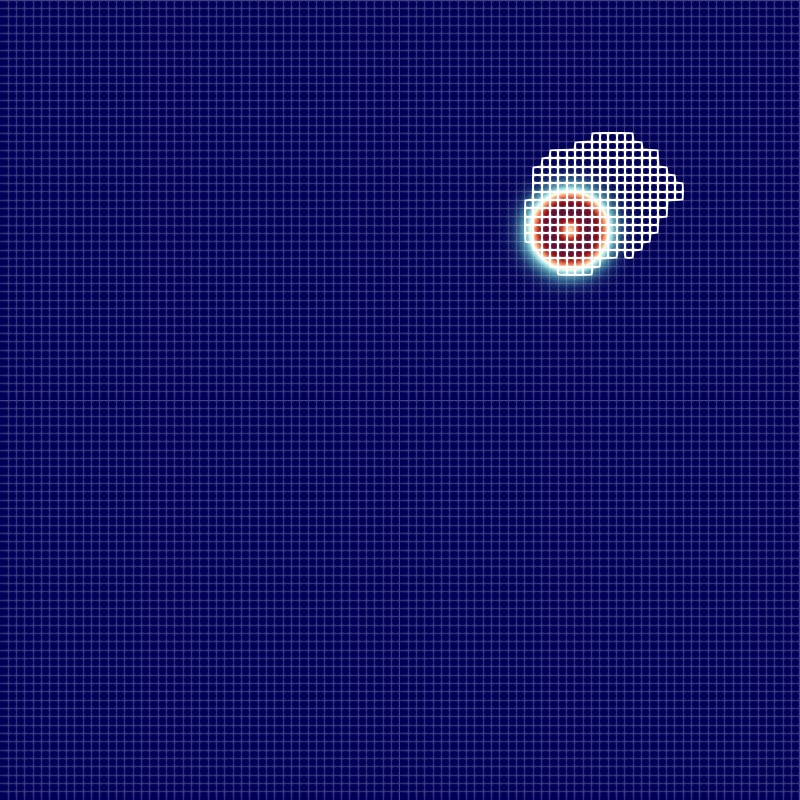}}}
        \subfloat[Solution, $5T$]{
        \adjustbox{width=0.24\linewidth,valign=b}{\includegraphics{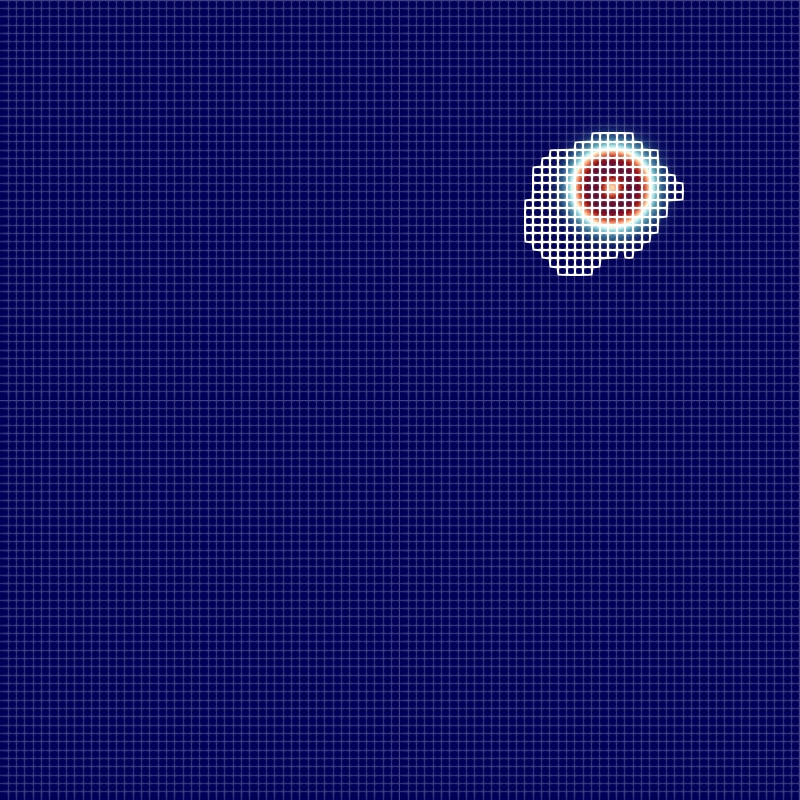}}}
        \subfloat[Remesh, $5T$]{
        \adjustbox{width=0.24\linewidth,valign=b}{\includegraphics{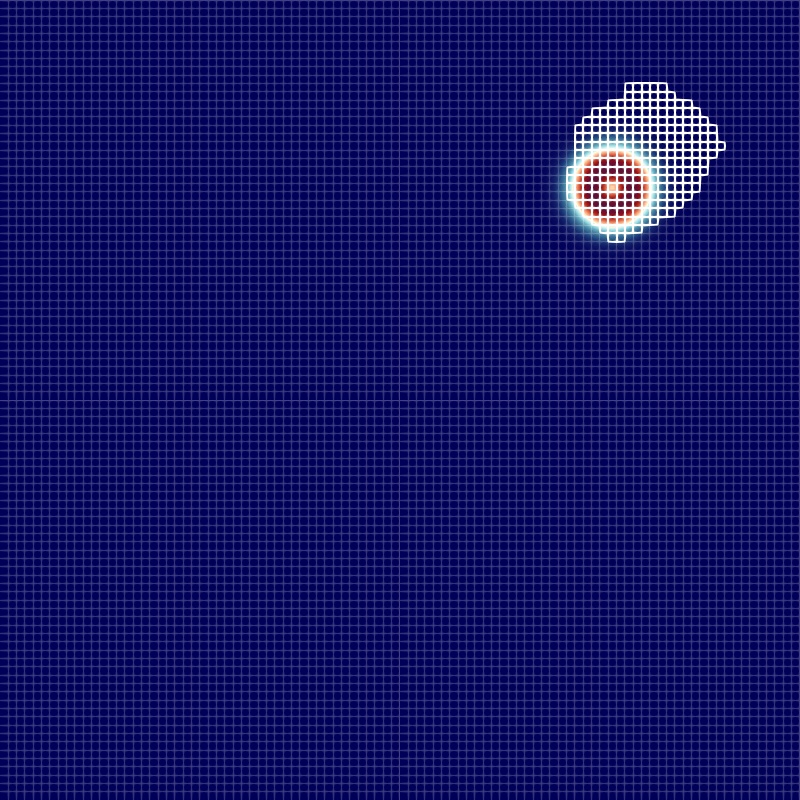}}}
        \subfloat[Solution, $6T$]{
        \adjustbox{width=0.24\linewidth,valign=b}{\includegraphics{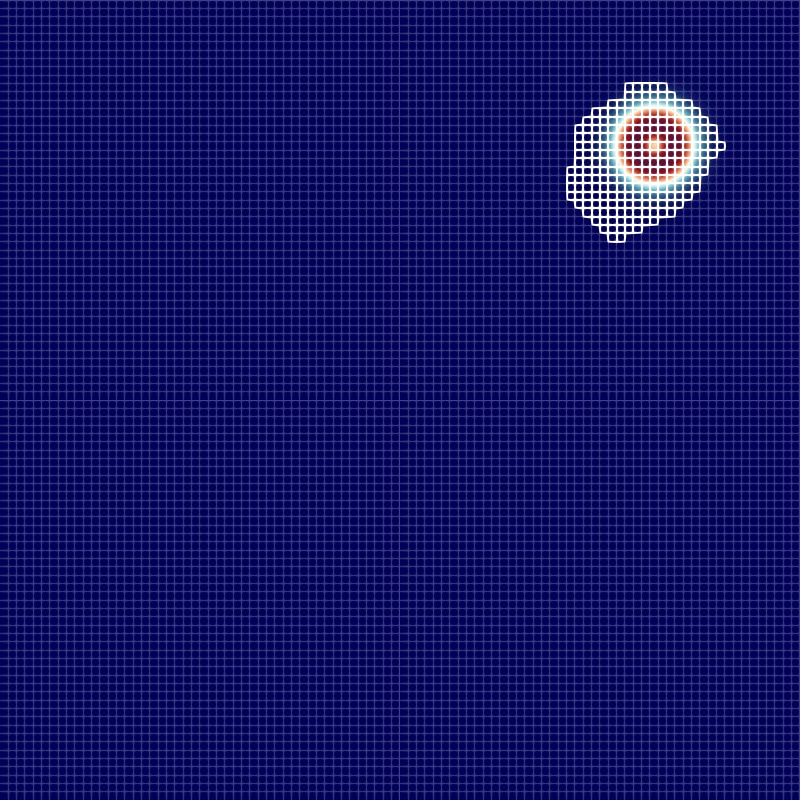}}}
        \newline
        \subfloat[Remesh, $6T$]{
        \adjustbox{width=0.24\linewidth,valign=b}{\includegraphics{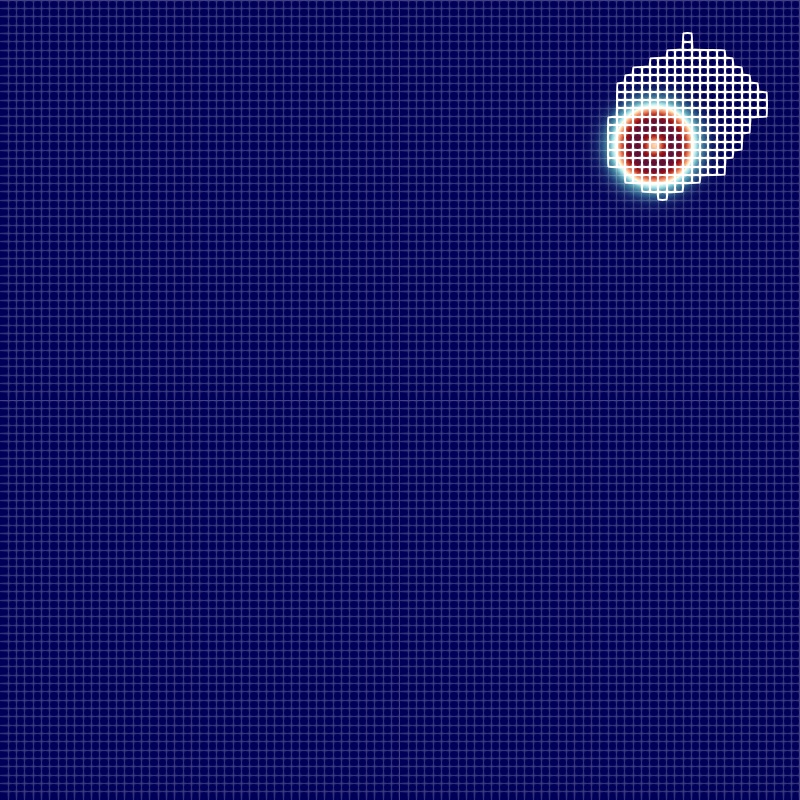}}}
        \subfloat[Solution, $7T$]{
        \adjustbox{width=0.24\linewidth,valign=b}{\includegraphics{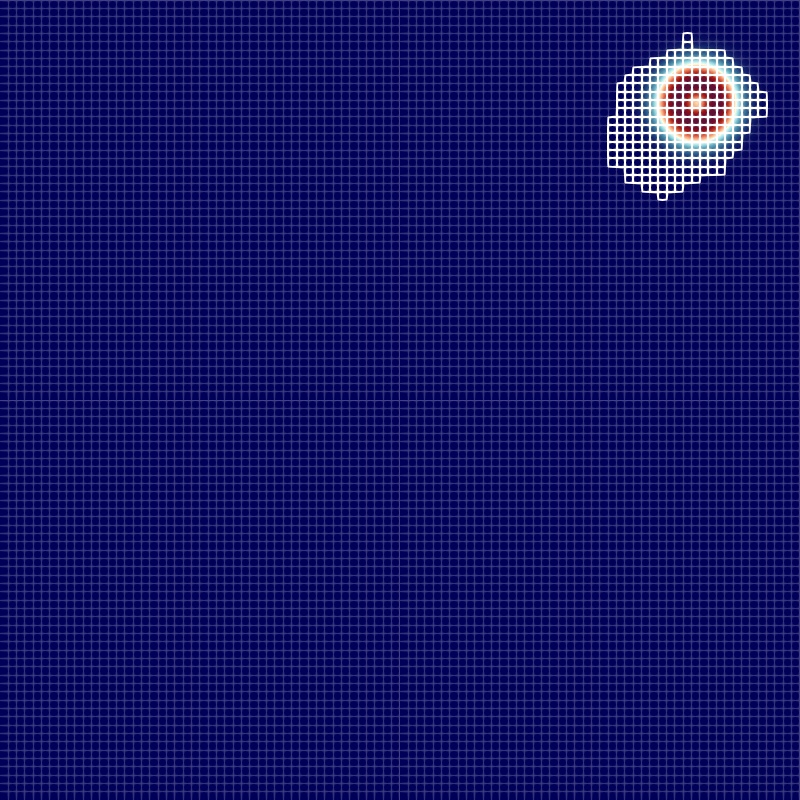}}}
        \subfloat[Remesh, $7T$]{
        \adjustbox{width=0.24\linewidth,valign=b}{\includegraphics{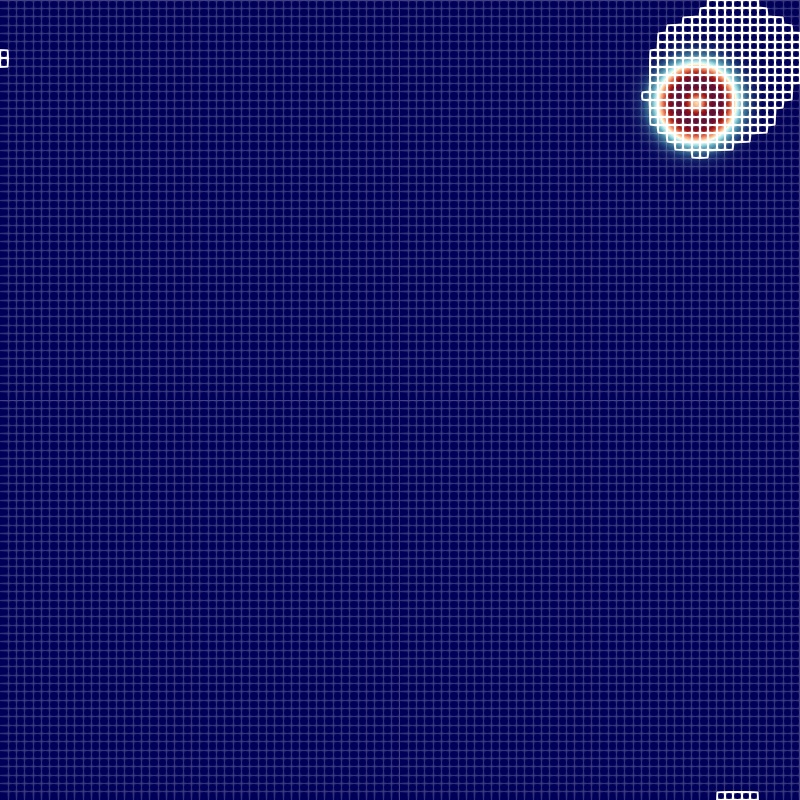}}}
        \subfloat[Solution, $8T$]{
        \adjustbox{width=0.24\linewidth,valign=b}{\includegraphics{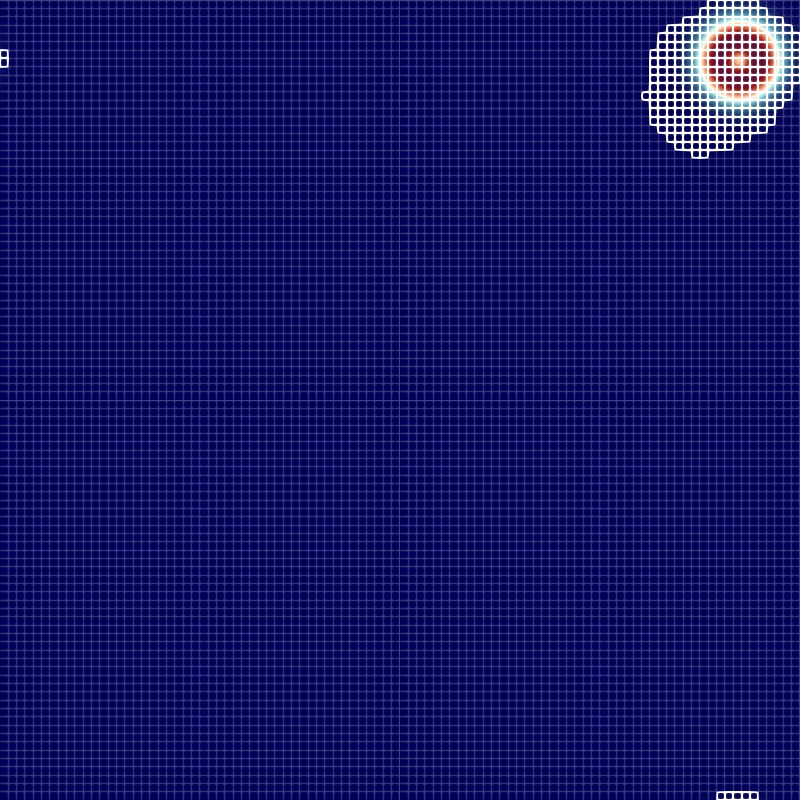}}}
        \newline
        \caption{\label{fig:app_pref} 
        Solution contours overlaid with $p$-adapted mesh at varying remesh intervals for a generalization experiment on the advection equation. Highlighted elements represent $p$-refinement.
        }
    \end{figure}

\begin{figure}[htbp!]
        \centering
        \subfloat[Remesh at $t = 0$]{
        \adjustbox{width=0.24\linewidth,valign=b}{\includegraphics{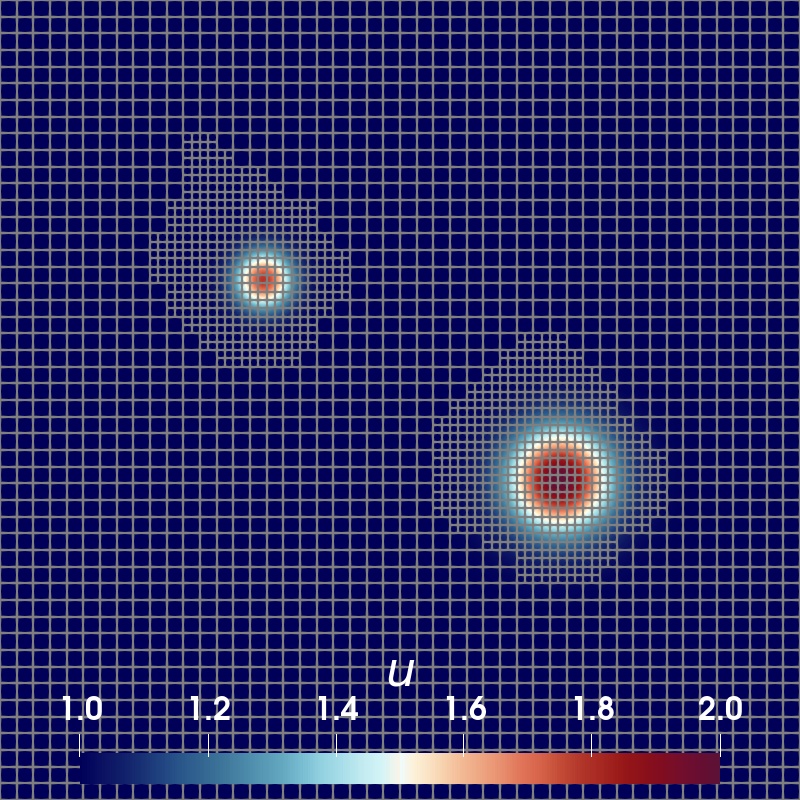}}}
        \subfloat[Solution at $t = T$]{
        \adjustbox{width=0.24\linewidth,valign=b}{\includegraphics{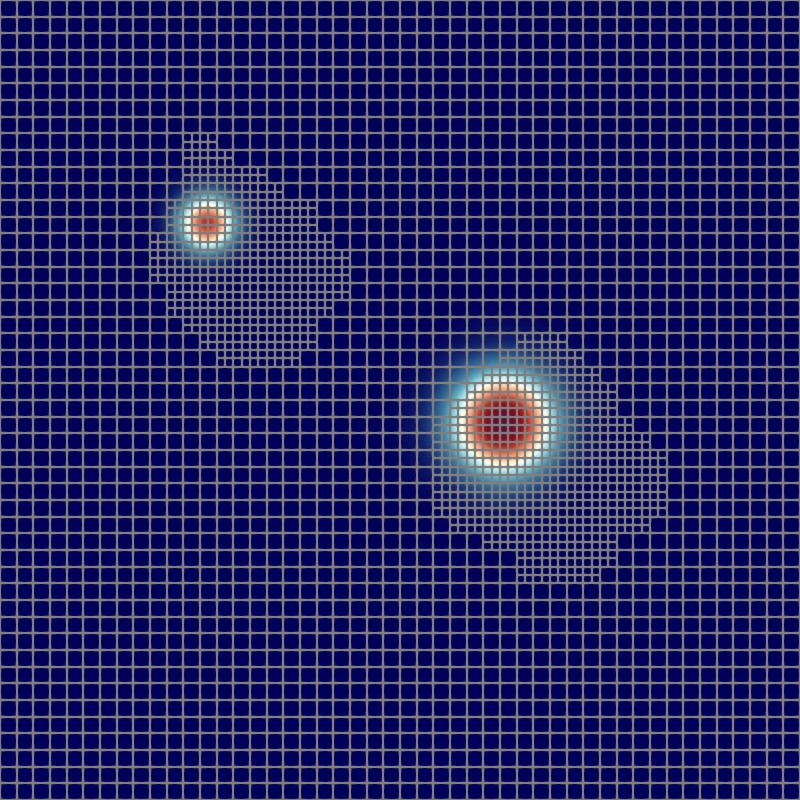}}}
        \subfloat[Remesh at $t = T$]{
        \adjustbox{width=0.24\linewidth,valign=b}{\includegraphics{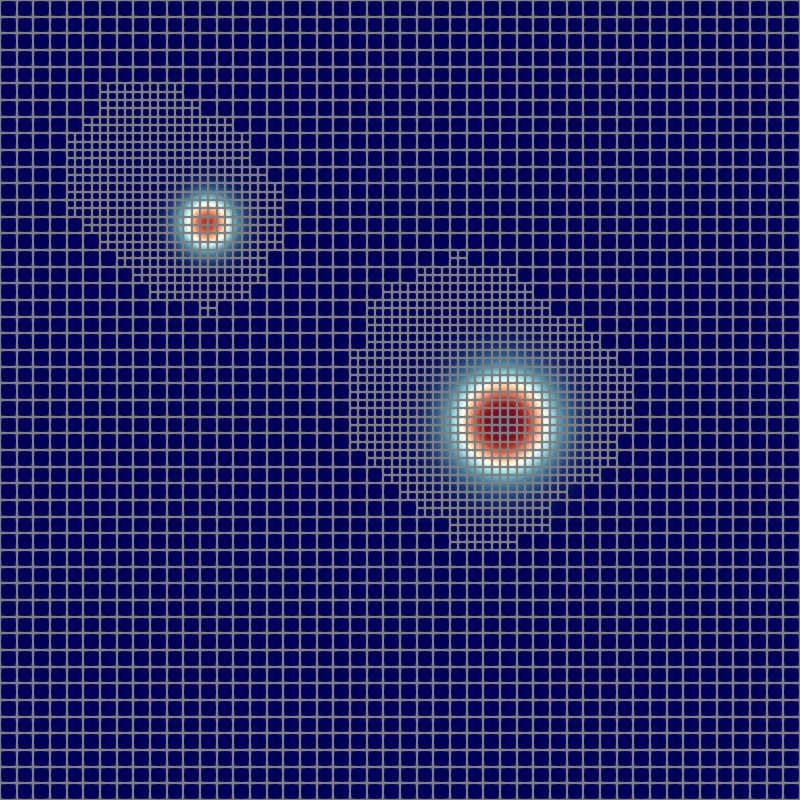}}}
        \subfloat[Solution at $t = 2T$]{
        \adjustbox{width=0.24\linewidth,valign=b}{\includegraphics{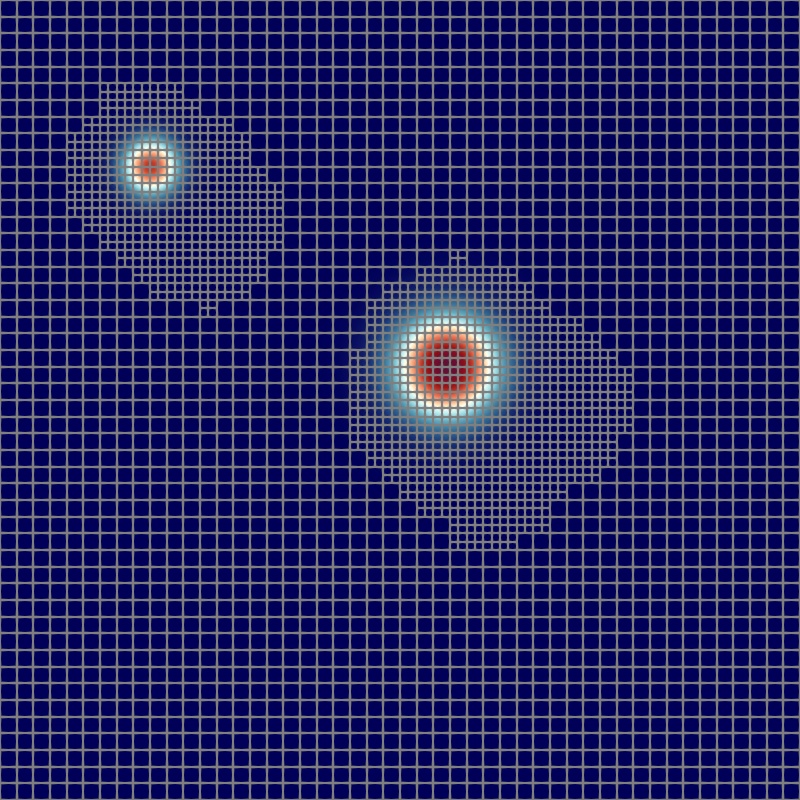}}}
        \newline
        \subfloat[Remesh at $t = 2T$]{
        \adjustbox{width=0.24\linewidth,valign=b}{\includegraphics{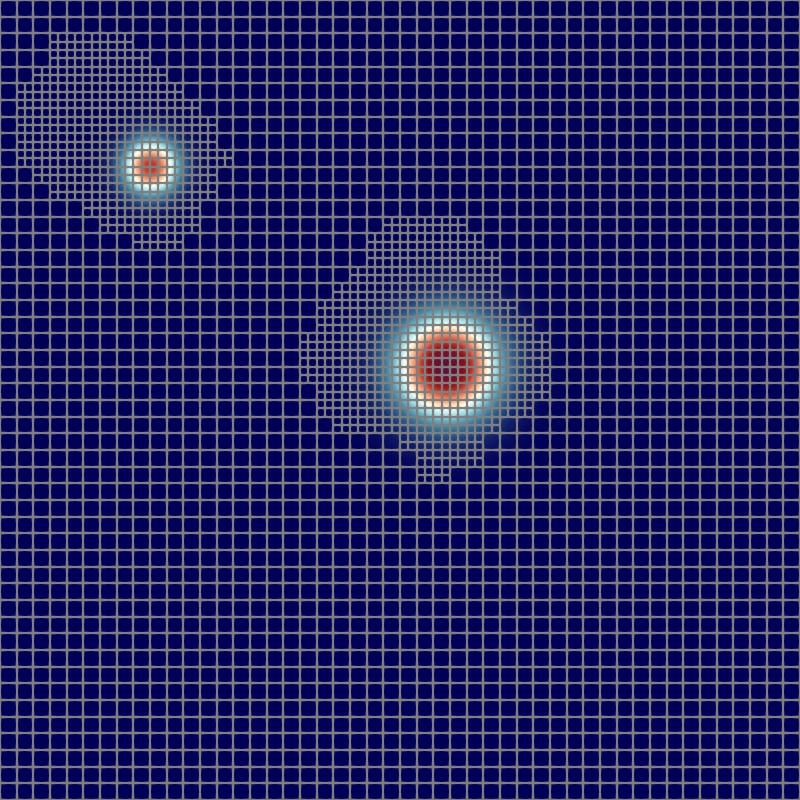}}}
        \subfloat[Solution at $t = 3T$]{
        \adjustbox{width=0.24\linewidth,valign=b}{\includegraphics{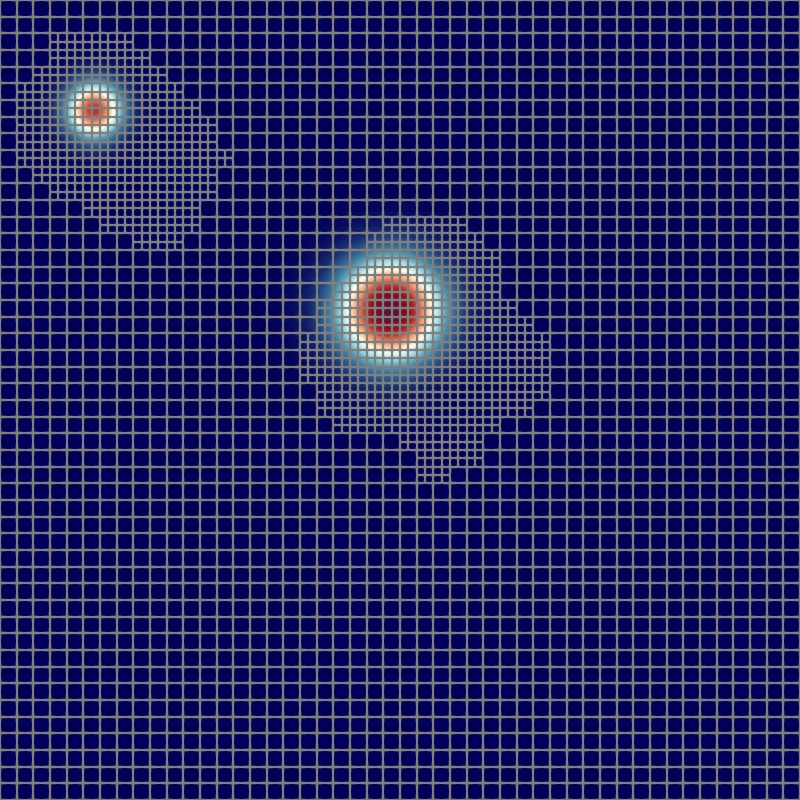}}}
        \subfloat[Remesh, $3T$]{
        \adjustbox{width=0.24\linewidth,valign=b}{\includegraphics{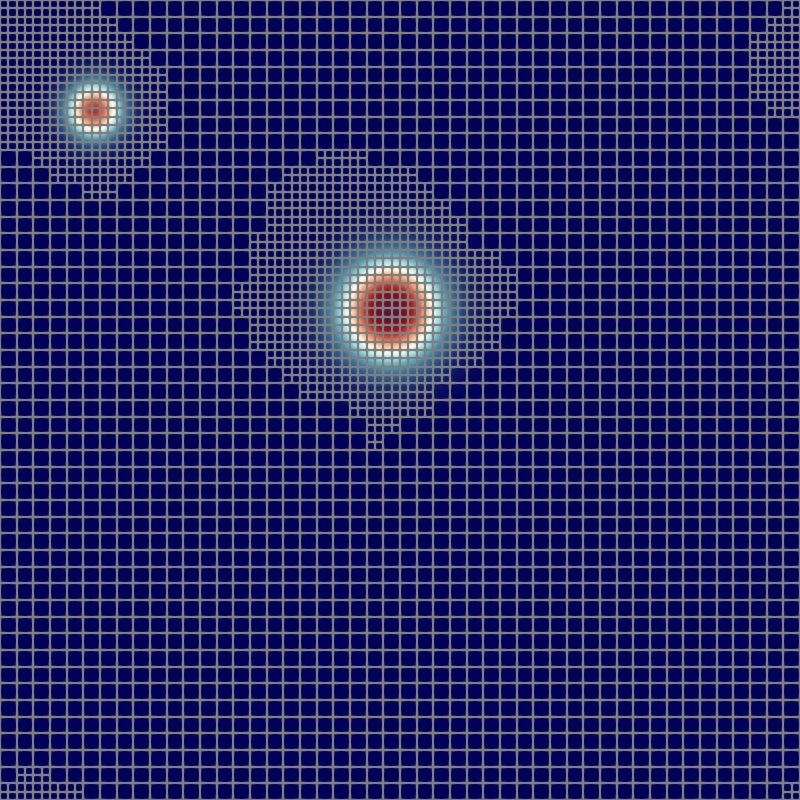}}}
        \subfloat[Solution, $4T$]{
        \adjustbox{width=0.24\linewidth,valign=b}{\includegraphics{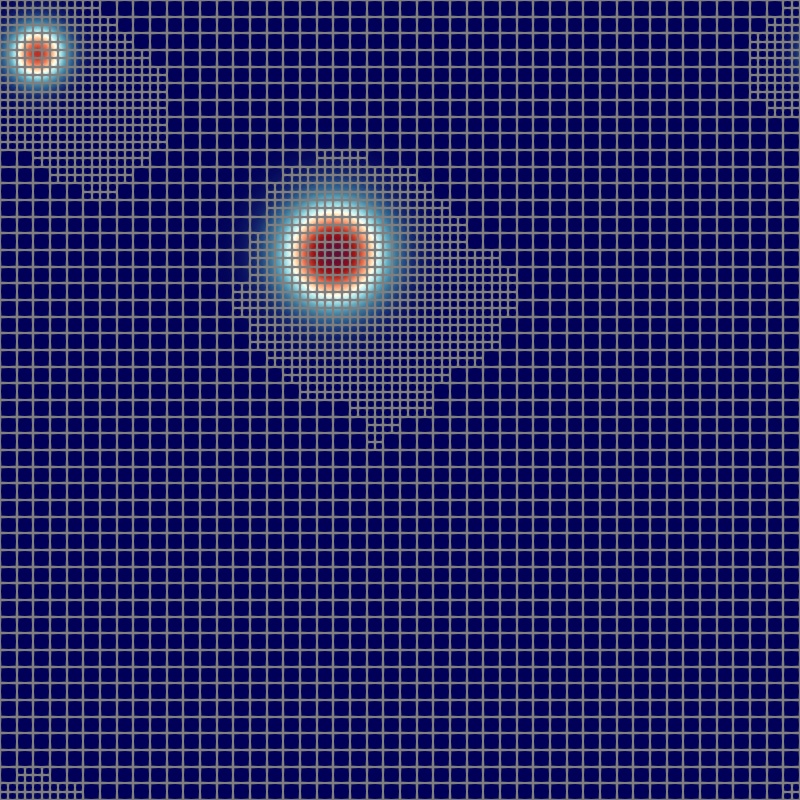}}}
        \newline
        \subfloat[Remesh, $4T$]{
        \adjustbox{width=0.24\linewidth,valign=b}{\includegraphics{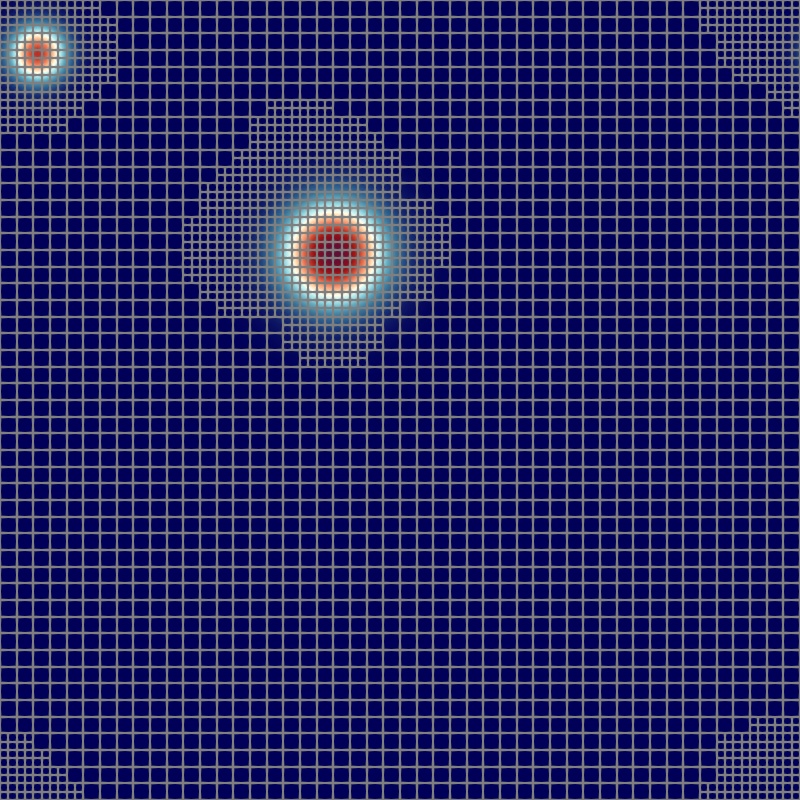}}}
        \subfloat[Solution, $5T$]{
        \adjustbox{width=0.24\linewidth,valign=b}{\includegraphics{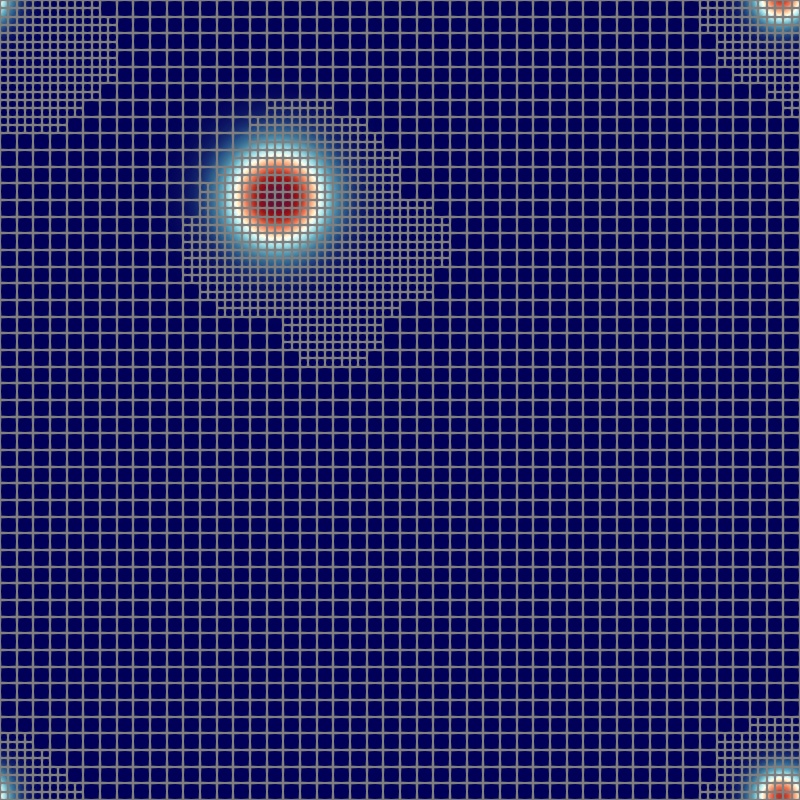}}}
        \subfloat[Remesh, $5T$]{
        \adjustbox{width=0.24\linewidth,valign=b}{\includegraphics{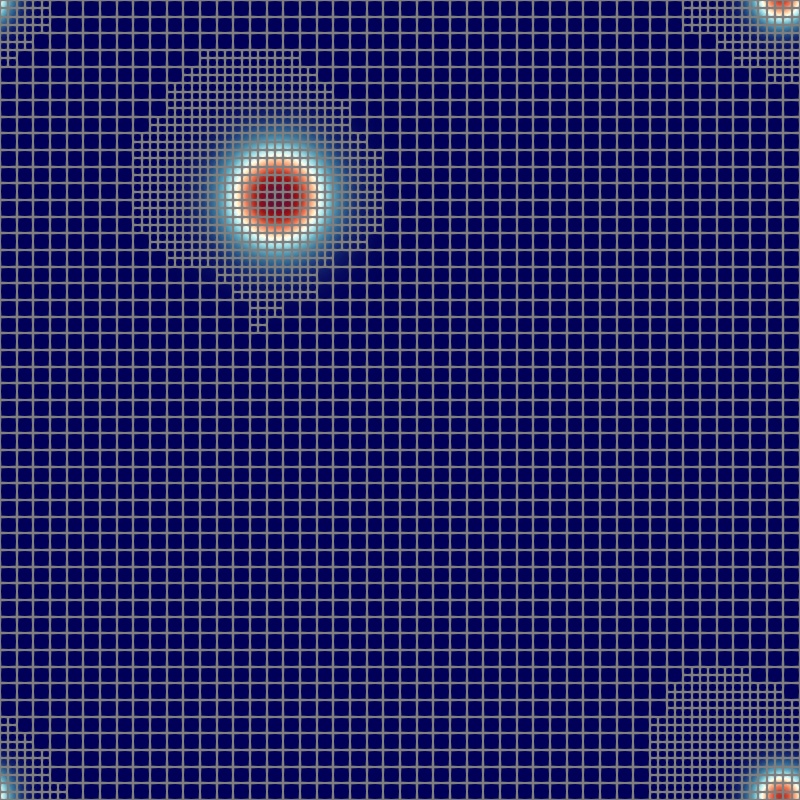}}}
        \subfloat[Solution, $6T$]{
        \adjustbox{width=0.24\linewidth,valign=b}{\includegraphics{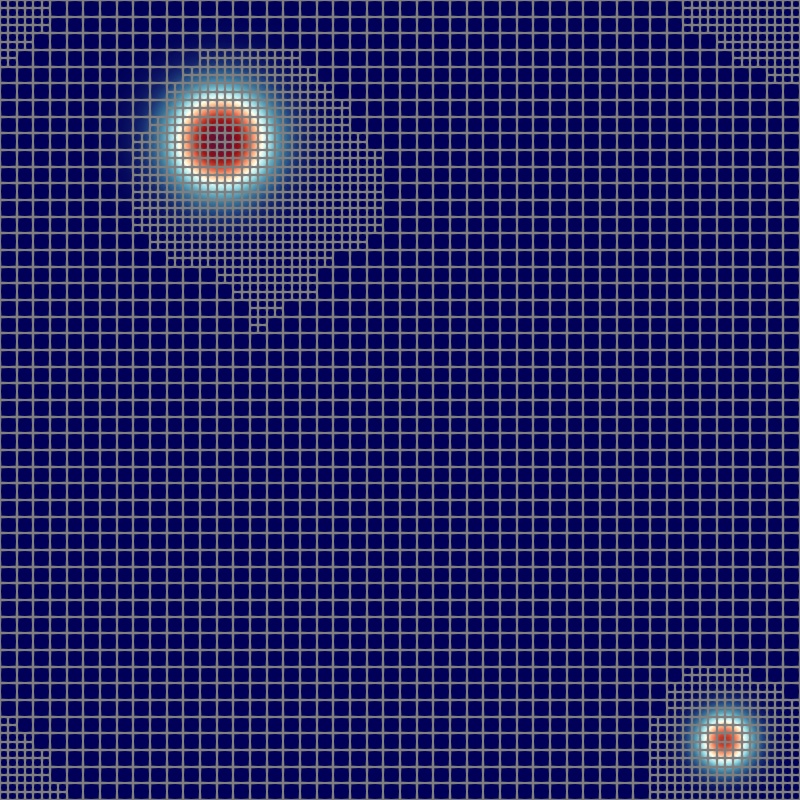}}}
        \newline
        \subfloat[Remesh, $6T$]{
        \adjustbox{width=0.24\linewidth,valign=b}{\includegraphics{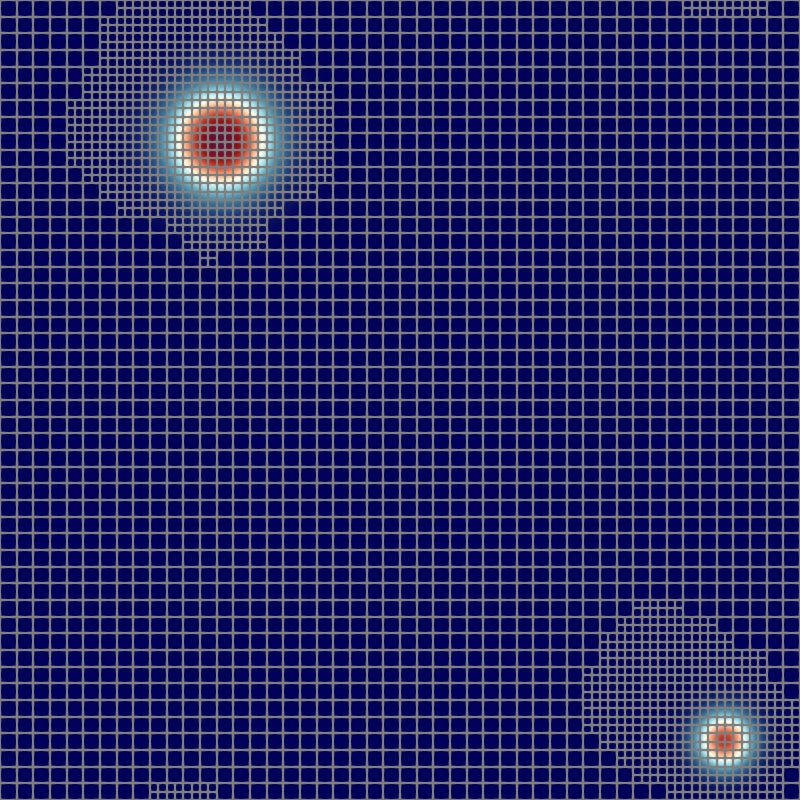}}}
        \subfloat[Solution, $7T$]{
        \adjustbox{width=0.24\linewidth,valign=b}{\includegraphics{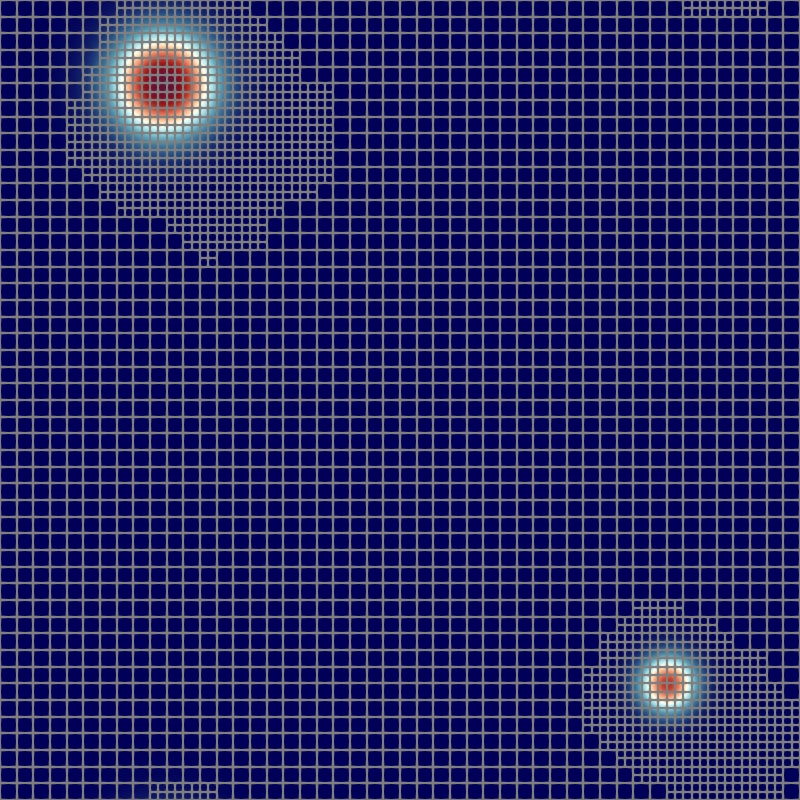}}}
        \subfloat[Remesh, $7T$]{
        \adjustbox{width=0.24\linewidth,valign=b}{\includegraphics{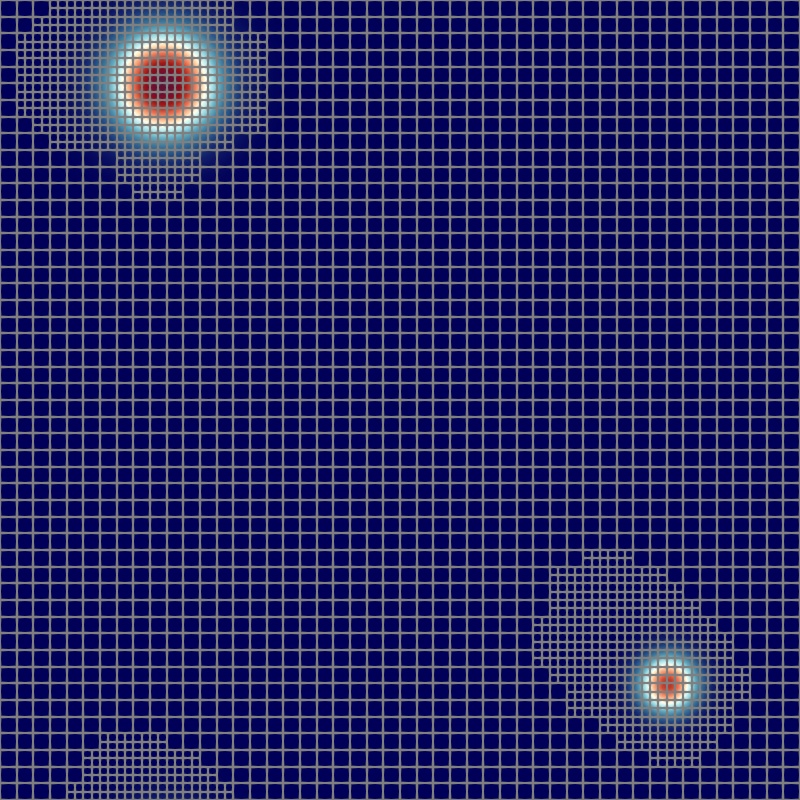}}}
        \subfloat[Solution, $8T$]{
        \adjustbox{width=0.24\linewidth,valign=b}{\includegraphics{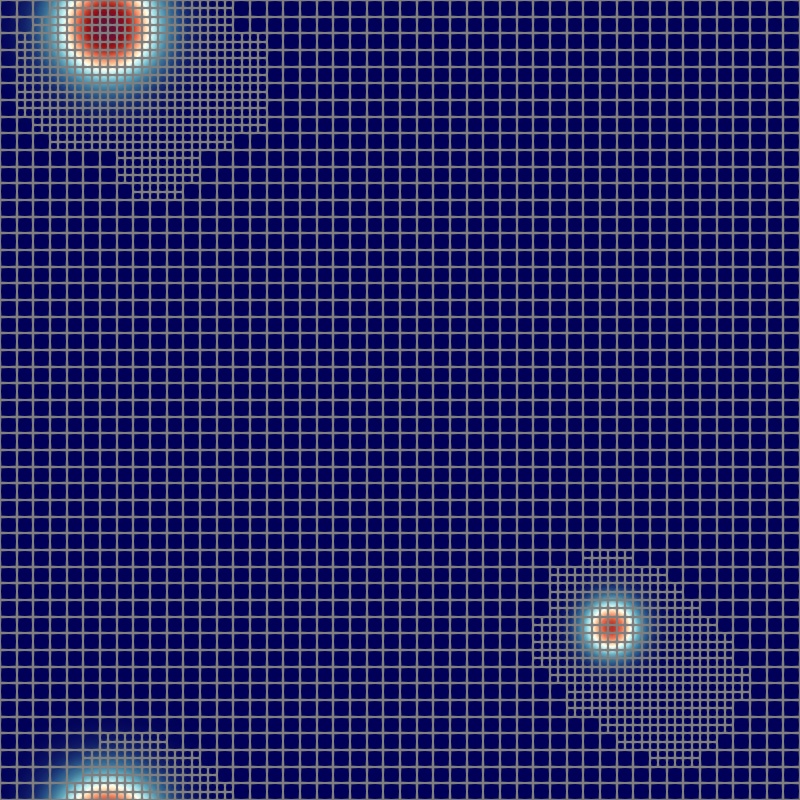}}}
        \newline
        \caption{\label{fig:app_href} 
        Solution contours overlaid with $h$-adapted mesh at varying remesh intervals for a generalization experiment on the advection equation.
        }
    \end{figure}

Some mesh refinement examples from the generalization experiments are presented in this section. \cref{fig:app_pref} showcases the adapted meshes for $p$-refinement on the advection equation. This setup utilizes a finer mesh (with $16$ times more elements) and a longer simulation time. \cref{fig:app_href} shows the adapted meshes for $h$-refinement on the advection equation, with a finer mesh (with $8$ times more elements), different solution profiles (multiple Gaussian bumps of varying width and height), and a longer simulation time.

\end{appendices}

\end{document}